\documentclass[12pt]{amsart}
\usepackage{amscd,amsxtra,amssymb, bbm}

\usepackage[dvips]{graphics}
\usepackage{epsfig,graphicx}

\input xy
\xyoption{all}
\xyoption{arc}
\SelectTips{cm}{10}

\newtheorem{theorem}{Theorem}[section]
\newtheorem{corollary}[theorem]{Corollary}
\newtheorem{lemma}[theorem]{Lemma}
\newtheorem{proposition}[theorem]{Proposition}

\theoremstyle{definition}
\newtheorem{definition}[theorem]{Definition}
\newtheorem{remark}[theorem]{Remark}
\newtheorem{algorithm}[theorem]{Algorithm}
\newtheorem{example}[theorem]{Example}
\newtheorem{conjecture}[theorem]{Conjecture}

\DeclareMathOperator{\Dbcoh}{\mathsf{D^{b}_{coh}}}
\DeclareMathOperator{\Perf}{\mathsf{Perf}}

\newcommand{\MPnode}{\mathsf{MP}_\mathrm{nd}}
\newcommand{\MPcusp}{\mathsf{MP}_\mathrm{cp}}
\DeclareMathOperator{\rk}{\mathrm{rk}}
\DeclareMathOperator{\cn}{\mathrm{cn}}
\DeclareMathOperator{\sn}{\mathrm{sn}}
\DeclareMathOperator{\dn}{\mathrm{dn}}
\DeclareMathOperator{\coker}{\mathrm{coker}}
\renewcommand{\ker}{\mathrm{ker}}
\newcommand{\im}{\mathrm{im}}
\newcommand{\conj}{\mathrm{cnj}}
\newcommand{\can}{\mathrm{can}}
\newcommand{\pr}{\mathrm{pr}}
\newcommand{\id}{\mathrm{id}}
\newcommand{\Tr}{\mathrm{Tr}}
\newcommand{\tr}{\mathrm{tr}}
\newcommand{\rank}{\mathrm{rank}}
\DeclareMathOperator{\cl}{\mathrm{cl}}
\DeclareMathOperator{\res}{\mathrm{res}}
\DeclareMathOperator{\ev}{\mathrm{ev}}

\DeclareMathOperator{\Coh}{\mathsf{Coh}}
\DeclareMathOperator{\BM}{\mathsf{BM}}

\DeclareMathOperator{\VB}{\mathsf{VB}}

\DeclareMathOperator{\Spl}{\mathsf{Spl}}

\DeclareMathOperator{\Tri}{\mathsf{Tri}}
\DeclareMathOperator{\AF}{\mathsf{AF}}
\DeclareMathOperator{\Pic}{\mathsf{Pic}}

\DeclareMathOperator{\Mf}{\mathsf{M}}

\DeclareMathOperator{\Hom}{\mathsf{Hom}}
\DeclareMathOperator{\RHom}{\mathsf{RHom}}
\DeclareMathOperator{\Mor}{\mathsf{Mor}}
\DeclareMathOperator{\Ext}{\mathsf{Ext}}
\DeclareMathOperator{\Lin}{\mathsf{Lin}}
\DeclareMathOperator{\GL}{\mathsf{GL}}
\DeclareMathOperator{\Aut}{\mathsf{Aut}}
\DeclareMathOperator{\End}{\mathsf{End}}
\DeclareMathOperator{\Mat}{\mathsf{Mat}}
\DeclareMathOperator{\Spec}{\mathsf{Spec}}

\DeclareMathOperator{\Ans}{\mathsf{Ans}}
\DeclareMathOperator{\Sets}{\mathsf{Sets}}
\DeclareMathOperator{\For}{\mathsf{For}}

\DeclareMathOperator{\Vect}{\mathsf{Vect}}

\newcommand{\bul}{\scriptstyle{\bullet}}
\newcommand{\kk}{\boldsymbol{k}}
\newcommand{\CC}{\mathbb{C}}
\newcommand{\ZZ}{\mathbb{Z}}
\newcommand{\fR}{\mathsf{R}}

\newcommand{\FF}{\mathbb{F}}
\newcommand{\GG}{\mathbb{G}}

\newcommand{\DD}{\mathbb{D}}

\newcommand{\PP}{\mathbb{P}}

\newcommand{\kA}{\mathcal{A}}
\newcommand{\kB}{\mathcal{B}}
\newcommand{\kE}{\mathcal{E}}
\newcommand{\kF}{\mathcal{F}}
\newcommand{\kG}{\mathcal{G}}

\newcommand{\kI}{\mathcal{I}}
\newcommand{\kO}{\mathcal{O}}
\newcommand{\kL}{\mathcal{L}}
\newcommand{\kP}{\mathcal{P}}
\newcommand{\kQ}{\mathcal{Q}}

\newcommand{\kM}{\mathcal{M}}
\newcommand{\kN}{\mathcal{N}}
\newcommand{\kV}{\mathcal{V}}
\newcommand{\kW}{\mathcal{W}}
\newcommand{\kT}{\mathcal{T}}
\newcommand{\kS}{\mathcal{S}}

\newcommand{\mm}{\mathsf{m}}
\newcommand{\nn}{\mathsf{n}}
\newcommand{\ii}{\mathsf{i}}
\newcommand{\jj}{\mathsf{j}}
\newcommand{\cc}{\mathsf{c}}
\newcommand{\pp}{\mathsf{p}}
\newcommand{\qq}{\mathsf{q}}
\newcommand{\rr}{\mathsf{r}}

\newcommand{\lar}{\longrightarrow}

\begin{document}

\title[Vector bundles and Yang--Baxter equations]{Vector bundles on
degenerations of elliptic curves
and  Yang--Baxter equations}

\author{Igor Burban}
\address{%
Mathematisches Institut,
Universit\"at Bonn,
Endenicher Allee , 53113 Bonn,
Germany
}
\email{burban@math.uni-bonn.de}

\author{Bernd Kreussler}
\address{%
Mary Immaculate College, South Circular Road, Limerick, Ireland
}
\email{bernd.kreussler@mic.ul.ie}

\begin{abstract}
In this paper we introduce the notion of a geometric associative $r$-matrix
attached  to a genus one  fibration with a section and irreducible fibres.
It allows us to study degenerations of solutions of the classical
Yang--Baxter equation using the
approach of Polishchuk. We also calculate  certain   solutions
of the classical, quantum and associative Yang--Baxter equations obtained from
moduli spaces of (semi-)stable vector bundles on Weierstra\ss{} cubic
curves.
\end{abstract}

\maketitle

\section{Introduction}

There are many indications (for example from homological mirror symmetry)
that the  formalism of derived categories provides a compact way to
formulate  and solve
complicated non-linear analytical problems. However, one would like to have
more concrete examples, in which one can follow the full path starting from
a categorical set-up and ending with an analytical output. In this article  we
study the interplay between the theory of the associative, classical and quantum
Yang--Baxter equations
and properties of vector bundles on projective curves of arithmetic genus one,
following the approach of Polishchuk \cite{Polishchuk1}.

Let $\mathfrak{g}$ be the Lie algebra $\mathfrak{sl}_n(\mathbb{C})$ and
$A = U(\mathfrak{g})$ its universal enveloping algebra.
The classical Yang--Baxter equation (CYBE)  is
$$
[r^{12}(x), r^{13}(x+y)] + [r^{13}(x+y), r^{23}(y)] + [r^{12}(x), r^{23}(y)] = 0,
$$
where $r(z)$ is the germ of a meromorphic function of one  variable
in a neighbourhood of $0$ taking values in
$\mathfrak{g} \otimes \mathfrak{g}$.
The upper indices in this equation   indicate various embeddings of
$\mathfrak{g}\otimes \mathfrak{g}$ into $A \otimes A \otimes A$. For example,
the function $r^{13}$ is defined as
$$
r^{13}: \mathbb{C} \stackrel{r}\lar \mathfrak{g}\otimes \mathfrak{g}
\stackrel{\tau_{13}}\lar A \otimes A  \otimes A,
$$
where $\tau_{13}(x\otimes y) = x \otimes 1 \otimes y$. Two other
maps $r^{12}$ and $r^{23}$ have a similar meaning.

In the physical literature, solutions of (CYBE)  are frequently
called $r$--\emph{matrices}.
They   play an important role in  mathematical physics, representation theory,
integrable  systems and statistical mechanics.

By a famous  result of Belavin and Drinfeld \cite{BelavinDrinfeld}, there
exist exactly three types of non-degenerate solutions of the classical
Yang--Baxter equation:  elliptic (two-periodic),
trigonometric (one-periodic) and rational.
This trichotomy corresponds to three models in statistical mechanics:
XYZ (elliptic), XXZ (trigonometric) and XXX (rational), see \cite{Baxter}.

Belavin and Drinfeld have also obtained a complete classification of elliptic
and trigonometric solutions, see
\cite[Proposition 5.1 and Theorem 6.1]{BelavinDrinfeld}.
A certain classification of rational solutions was given by Stolin
\cite[Theorem 1.1]{Stolin}.

\medskip

This article is devoted to a study of degenerations of elliptic $r$-matrices
into trigonometric and then into rational ones. We hope that this sort of
questions will be interesting from the point of view of applications in
mathematical physics.
In order to attack this problem we use a construction
of  Polishchuk \cite{Polishchuk1}.
After  certain modifications of his original presentation,
the core of this  method can be described as follows.

\medskip

Let $E$ be  a Weierstra\ss{} cubic curve, $\breve{E} \subset E$ the open subset of smooth points, $M = M_E^{(n,d)}$ the moduli space
of stable bundles of rank $n$ and degree $d$, assumed to be coprime.
Let $\kP = \kP(n,d) \in \VB(E \times M)$ be a universal family of the
moduli functor $\underline{\Mf}_E^{(n,d)}$. For a point $v \in M$ we denote by
$\kV = \kP|_{E \times v}$ the corresponding vector bundle on $E$.   Consider
the following data:
\begin{itemize}
\item two distinct  points $v_1, v_2 \in M$ in the moduli space;
\item two distinct   points  $y_1, y_2 \in \breve{E}$  such that
      $\kV_1(y_2) \not\cong \kV_2(y_1)$.
\end{itemize}

\noindent
Using Serre Duality, the triple Massey product
$$
\Hom_E(\kV_1, \mathbb{C}_{y_1}) \otimes \Ext_E^1(\mathbb{C}_{y_1}, \kV_2)
\otimes \Hom_E(\kV_2, \mathbb{C}_{y_2}) \lar
\Hom_E(\kV_1, \mathbb{C}_{y_2}),
$$
induces a linear map
$$
r^{\kV_1, \kV_2}_{y_1, y_2}:
\Hom_E(\kV_1, \mathbb{C}_{y_1}) \otimes
\Hom_E(\kV_2, \mathbb{C}_{y_2}) \lar
\Hom_E(\kV_2, \mathbb{C}_{y_1}) \otimes
\Hom_E(\kV_1, \mathbb{C}_{y_2})
$$
and satisfying  the so-called \emph{associative Yang--Baxter equation} (AYBE)
$$
\left(r^{\kV_3, \kV_2}_{y_1, y_2}\right)^{12}
\left(r^{\kV_1, \kV_3}_{y_1, y_3}\right)^{13} -
\left(r^{\kV_1, \kV_3}_{y_2, y_3}\right)^{23}
\left(r^{\kV_1, \kV_2}_{y_1, y_2}\right)^{12} +
\left(r^{\kV_1, \kV_2}_{y_1, y_3}\right)^{13}
\left(r^{\kV_2, \kV_3}_{y_2, y_3}\right)^{23} = 0
$$
viewed as a map
$$
\Hom_E(\kV_1, \mathbb{C}_{y_1}) \otimes \Hom_E(\kV_2, \mathbb{C}_{y_2})
\otimes \Hom_E(\kV_3, \mathbb{C}_{y_3}) \lar
$$
$$
\lar \Hom_E(\kV_2, \mathbb{C}_{y_1}) \otimes \Hom_E(\kV_3, \mathbb{C}_{y_2})
\otimes \Hom_E(\kV_1,
\mathbb{C}_{y_3}).
$$
This map  can be rewritten as the germ of a tensor-valued
meromorphic  function in  four variables, defined in a neighbourhood of a
smooth point $o$ of the moduli space $M \times M \times E \times E$ (the
choice of $o$ will be explained in Corollary \ref{C:MainCorol})
$$r(\kV_1, \kV_2; y_1, y_2):
(\mathbb{C}^2 \times \mathbb{C}^2, 0) \cong
\bigl((M \times M) \times (E \times E), o\bigr) \lar
 \Mat_{n \times n}(\mathbb{C}) \otimes \Mat_{n \times n}(\mathbb{C}).$$

Since the complex manifold $M_E^{(n,d)}$ is a homogeneous space  over the
algebraic group $J= \Pic^0(E)$, it turns out that
$$
r(v_1, v_2; y_1, y_2) \sim r(v_1 - v_2; y_1 , y_2) = r(v;  y_1, y_2),
$$
with respect to  a certain equivalence relation on the set of solutions. We
show that this equivalence relation corresponds to a change of a
trivialization of the universal family $\kP$.

Let  $e$ be the neutral element of $J$.
It was shown by Polishchuk \cite[Lemma 1.2]{Polishchuk1}
that the function of two variables
$$
\bar{r}(y_1, y_2) = \lim_{v\to e}(\pr \otimes \pr) r(v; y_1, y_2)
\in \mathfrak{sl}_n(\mathbb{C}) \otimes \mathfrak{sl}_n(\mathbb{C})
$$
is a non-degenerate unitary solution of the classical Yang--Baxter equation.
Moreover, under certain restrictions (which are always fulfilled at least
for elliptic curves and Kodaira cycles of projective lines), for
any fixed value $g \ne e$ from a small neighbourhood $U_e \subseteq J$
of   $e$, the   tensor-valued function
$$r: \bigl(\{g\}\times E \times E, e\bigr)  \lar   \Mat_{n \times n}(\mathbb{C}) \otimes
\Mat_{n \times n}(\mathbb{C})
$$
satisfies the  \emph{quantum} Yang--Baxter equation, see
\cite[Theorem 1.4]{Polishchuk2}. Hence, this approach gives an explicit method
to quantize some known solutions of the classical Yang--Baxter equation.

Moreover, as was pointed out by Kirillov \cite{Kirillov},
a solution $r(v;y_1, y_2)$ of the associative Yang--Baxter equation defines an interesting
family of pairwise commuting first order differential operators,  generalizing Dunkl operators
studied by Buchstaber, Felder and Veselov \cite{BFV}, see Proposition \ref{P:Kirillov}.

The aim of our article is to study a \emph{relative} version of
Polishchuk's construction.  Although  most of the results can be generalized
to the case of arbitrary reduced projective curves of arithmetic genus one
having trivial dualizing sheaf, in this
article we shall concentrate mainly on the case of irreducible curves.

Let $E$ be a Weierstra\ss{} cubic curve, i.e.~a plane projective curve
  given by the equation
$zy^2 = 4 x^3  - g_2 x z^2  - g_3 z^3$.
It is singular if any only if $\Delta := g_2^3 -  27 g_3^2 = 0$.
$$
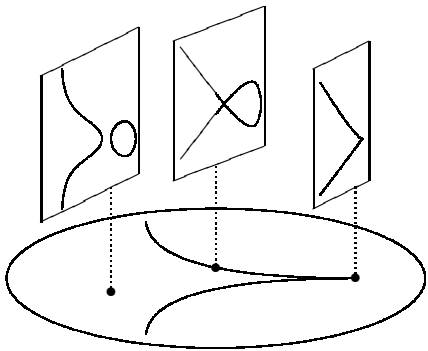
{\xy 0;/r0.22pc/:
\POS(40,0);
{(25,0)\ellipse(30,10){-}},%{\ellipse(,.75){}}
{(12,20)\ellipse(1.5,2.5)_,=:a(180){-}},
{(12,20)\ellipse(2,2.5)^,=:a(180){-}},
\POS(3,30);
@(,
\POS(3,25)@+,  \POS(8,22)@+, \POS(9,20)@+,
\POS(8,18)@+, \POS(3,15)@+, \POS(3,10)@+,
**\crvs{-}
,@i @);
{\ar@{-}(0,8);(0,29)},
{\ar@{-}(0,8);(14,15)},
{\ar@{-}(0,29);(14,36)},
{\ar@{-}(14,15);(14,36)},
%nodal curve
\POS(20,33);@-
{%(28,28)@+,
(20,30)@+,  (30,20)@+, (32,25)@+
,(30,30)@+, (20,20)@+, \POS(20,17)@-,
,**\qspline{}};  %@)
{\ar@{-}(19,14);(19,34)},
{\ar@{-}(19,14);(32,19)},
{\ar@{-}(19,34);(32,39)},
{\ar@{-}(32,19);(32,39)},
%cusp curve
\POS(40,28); \POS(46,20);
**\crv{(40,28);
(30,25); (42,22); (45,20) };
\POS(40,12); \POS(46,20);
**\crv{(40,18);
(45,5);   (42,8); (45,20); (46,20)};
%\POS(46.5,20)*{\bul}; \POS(48,20);
{\ar@{-}(39,10);(39,30)},
{\ar@{-}(39,10);(47,14)},
{\ar@{-}(39,30);(47,34)},
{\ar@{-}(47,14);(47,34)},
%%
%%cusp base
\POS(15,8); \POS(45,0)*{\bul} ;
**\crv{(16,0)};
\POS(15,-8); \POS(45,0);
**\crv{(16,0)};
\POS(15,8);
{\ar@{.}(45,13); (45,0) }
\POS(25,1.5)*{\bul};
{\ar@{.}(25,2); (25,16) }
\POS(10,-2)*{\bul};
{\ar@{.}(10,-2); (10,13) }
\endxy}
$$
Unless $g_2 = g_3 = 0$, the singularity is a node, whereas for $g_2 = g_3 = 0$
it is a cusp.

\medskip

A connection between the theory  of vector bundles on cubic curves and exactly
solvable models
of mathematical physics was observed a long time ago, see for example
\cite[Chapter 13]{McKean} and
\cite{Mulase} for a link  with
KdV equation, \cite{Cherednik} for applications to integrable systems and
\cite{BenZviNevins} for an interplay with
Calogero-Moser systems. In particular, the correspondence
\begin{center}
\begin{tabular}{l||l}
elliptic & elliptic \\
\hline
trigonometric & nodal \\
\hline
rational  & cuspidal
\end{tabular}
\end{center}
was discovered at  the very beginning
of the algebraic theory of completely integrable systems.

\medskip

In this article we follow another strategy. Instead of looking at each curve
of arithmetic genus one individually,  we consider the relative case, so that
all solutions will be considered as specializations of one \emph{universal}
solution.  Our main result can be stated  as follows.

\medskip
\noindent
\emph{Let $E \to T$ be a genus one fibration with a section having reduced and
irreducible fibres, $M = M_{E/T}^{(n,d)}$ the moduli space of relatively stable
vector bundles of rank $n$ and degree $d$. We construct a meromorphic function
$$
r: (M \times_T \times M\times_T E \times_T E, o) \lar
\Mat_{n \times n}(\mathbb{C}) \otimes \Mat_{n\times n}(\mathbb{C})
$$
in a neighbourhood of a smooth point $o$ of
$M \times_T \times M\times_T E \times_T E$,
which satisfies the associative Yang--Baxter equation for each fixed value $t
\in T$ and $(v_1,v_2,y_1,y_2) \in \bigl((M_{E_t} \times M_{E_t}) \times (E_t
\times E_t), o\bigr)$.
Moreover, $r_t(v_1, v_2, y_1, y_2)$ depends analytically on $t$, is compatible
with base change of the given family $E \to T$  and the corresponding solution
of the classical Yang--Baxter equation $\bar{r}_t(y)$ is
\begin{itemize}
\item elliptic if $E_t$ is smooth;
\item trigonometric if $E_t$ is nodal;
\item   rational  if $E_t$ is cuspidal.
\end{itemize}
}

\medskip
We also  carry  out explicit calculations  for vector bundles of rank two and
degree one on irreducible Weierstra\ss{} cubic curves.
In the case of an elliptic curve $E = E_\tau$ the corresponding solution is

\medskip
\begin{align*}
   r_{\mathrm{ell}}(v; y)  =
\frac{\theta'_1(0|\tau)}{\theta_1(y|\tau)}
&\left[
\frac{\theta_1(y + v|\tau)}{\theta_1(v|\tau)}
\mathbbm{1} \otimes \mathbbm{1} +
\frac{\theta_2(y + v|\tau)}{\theta_2(v|\tau)} h  \otimes h+\right.
\\
&\left.+
\frac{\theta_3(y + v|\tau)}{\theta_3(v|\tau)} \sigma  \otimes \sigma
+
\frac{\theta_4(y + v|\tau)}{\theta_4(v|\tau)} \gamma
\otimes \gamma\right],
\end{align*}
where $\mathbbm{1} =  e_{11} + e_{22}$, $h = e_{11} - e_{22},
\sigma = i(e_{21} - e_{12})$ and  $\gamma = e_{21} + e_{12}$.

\medskip
\noindent
In the case of a nodal cubic curve we get
\begin{align*}
  r_{\mathrm{trg}}(v; y) &=
\frac{\sin(y+v)}{\sin(y)\sin(v)}(e_{11}\otimes e_{11} +
e_{22} \otimes e_{22}) +  \frac{1}{\sin(v)}(e_{11} \otimes e_{22} + e_{22} \otimes
e_{11}) +
\\
&+ \frac{1}{\sin(y)}(e_{12}\otimes e_{21} + e_{21}\otimes e_{12}) + \sin(y+v)
e_{21} \otimes e_{21}
\end{align*}
and in  the case of a cuspidal cubic curve,  the associative $r$--matrix is
\begin{align*}
  r_{\mathrm{rat}}(v; y_1, y_2)  &=
\frac{1}{v} \mathbbm{1} \otimes \mathbbm{1} +
\frac{2}{y_2 - y_1}(e_{11} \otimes e_{11} + e_{22}\otimes e_{22} +
e_{12} \otimes e_{21} + e_{21}\otimes e_{12}) +
\\
&+ (v - y_1)  e_{21} \otimes h +
(v + y_2)  h  \otimes e_{21} -
v(v - y_1)(v + y_2)  e_{21}\otimes e_{21}.
\end{align*}
Our results imply that up to a gauge transformation the trigonometric and
rational solutions $r_{\mathrm{trg}}(v; y)$ and $r_{\mathrm{rat}}(v; y_1,y_2)$
are degenerations of $r_{\mathrm{ell}}(v; y)$, which seems to be difficult to
show by a direct computation.

Moreover, for a generic $v$ the tensors
$r_{\mathrm{ell}}(v; y)$, $r_{\mathrm{trg}}(v; y)$ and
$r_{\mathrm{rat}}(v; y_1, y_2)$ satisfy the quantum Yang--Baxter equation and
are quantizations   of the following classical $r$--matrices:

\medskip
\noindent
$\bullet$  Elliptic solution found and studied by Baxter, Belavin and Sklyanin:
$$
\bar{r}_{\mathrm{ell}}(y) =
\frac{\mathrm{cn}(y)}{\mathrm{sn}(y)} h \otimes h +
\frac{1+ \mathrm{dn}(y)}{\mathrm{sn}(y)}
(e_{12}\otimes e_{21} + e_{21} \otimes e_{12})  +
\frac{1 -\mathrm{dn}(y)}{\mathrm{sn}(y)}
(e_{12} \otimes e_{12} + e_{21} \otimes e_{21}).
$$

\medskip
\noindent
$\bullet$  Trigonometric solution of Cherednik:
$$
\bar{r}_{\mathrm{trg}}(y) = \frac{1}{2} \cot(y)h\otimes h +
\frac{1}{\sin(y)}(e_{12}\otimes e_{21} + e_{21} \otimes e_{12}) +
\sin(y) e_{21} \otimes e_{21}.
$$

\medskip
\noindent
$\bullet$ Rational solution
$$
\bar{r}_{\mathrm{rat}}(y) =
\frac{1}{y}\Bigl(\frac{1}{2} h \otimes h + e_{12} \otimes e_{21} +
e_{21} \otimes e_{12}\Bigr) +
y(e_{21} \otimes h + h \otimes e_{21}) - y^3 e_{21} \otimes e_{21},
$$
which is gauge equivalent to a solution found by Stolin \cite{Stolin}.

This paper is organized as follows. In Section \ref{S:YBE} we
collect  some results about the associative  Yang--Baxter
equation and its relations with Dunkl operators as well as with the
 classical and quantum Yang--Baxter equations.
Section \ref{S:PolConstr} gives a short introduction into a
construction of Polishchuk which provides a method to obtain solutions
of Yang--Baxter equations from triple Massey products in a derived category.

In order to be able to calculate solutions explicitly, this construction has
to be translated into another language, involving residue maps. In Section
\ref{S:Massey} we explain the corresponding result of Polishchuk whereby we
provide some details which are only implicit in \cite{Polishchuk1}. The
understanding of these details is crucial for the study of the relative case,
which is carried out in Sections \ref{S:RelativeConstr} and
\ref{S:GeomR}. Theorem \ref{C:MainCorol} is the main result of this article.
It states that for any genus one  fibration $E \rightarrow T$ satisfying certain
restrictions and any pair  of coprime integers $0 < d < r$  one can attach
a family of solutions  of the associative Yang-Baxter equation $r^\xi(v_1, v_2; y_1, y_2)$
depending \emph{analytically} on the parameter of the basis and functorial with
respect to the base change.  This solutions actually  depend
on the choice of a  trivialization $\xi$ of the universal family $\kP(n,d)$ of stable vector bundles
of rank $n$ and degree $d$.  However, in Proposition \ref{P:gauge} we show that
a choice of another trivialization $\zeta$ leads to a \emph{gauge equivalent} solution
$r^\zeta(v_1, v_2; y_1, y_2)$

In Section \ref{S:GeomYB} we prove  that in the case of a Weierstra\ss{} cubic curve
$E$  there exists  a trivialization  $\xi$  of the universal family $\kP(n,d)$  such that
the corresponding solution $r^\xi(v_1, v_2; y_1, y_2)$ is invariant
under simultaneous shifts $v_1 \mapsto v_1 + v, v_2 \mapsto v_2 + v$.
In other words, the solution $r^\xi(v_1, v_2; y_1, y_2)$, also called \emph{geometric associative $r$-matrix},  depends on the difference
$v_2 -v_1$ of the first pair of spectral parameters. Hence, the obtained
solution $r^\xi(v; y_1, y_2)$ also satisfies the quantum Yang--Baxter equation and defines an interesting
quantum integrable system. The key point  of the  proof is to
show  equivariance  of triple Massey products with respect to the action of the
Jacobian $J$ on the moduli space $M_E^{(n,d)}$.

Since it is indispensable  for carrying out explicit calculations of $r$--matrices, in the following
sections  we elaborate foundations of the theory of vector bundles on genus one curves.
In  Section \ref{S:elliptrm} we recall some classical results about
holomorphic  vector
bundles on a smooth elliptic curve.
Using the methods described before, we explicitly compute the  solution of the
associative Yang--Baxter equation and the classical $r$--matrix corresponding
to a  universal family of stable vector bundles of rank two and degree one.
These  solutions were obtained  by Polishchuk in \cite[Section 2]{Polishchuk1}
using homological mirror symmetry and formulas for higher products in the
Fukaya category of an elliptic curve. Our direct computation, however, is
independent of homological mirror symmetry. We are lead directly to express the
resulting associative $r$--matrix in terms of Jacobi's theta-functions and the
corresponding classical $r$--matrix in terms of the elliptic functions
$\mathrm{sn}(z)$, $\mathrm{cn}(z)$ and $\mathrm{dn}(z)$.

Sections \ref{S:triples} and \ref{S:singrm} are devoted to similar
calculations for nodal and cuspidal Weierstra\ss{} curves. Our computations
are based on the description of vector bundles on singular genus one curves  in terms of
so-called matrix problems, which was given by Drozd and Greuel
\cite{DrozdGreuel} and Burban \cite{Thesis}.  We show that their description
of canonical forms of matrix problems corresponds precisely to a very explicit
presentation  of universal families of stable vector bundles.  We explicitly
compute geometric $r$--matrices  coming from  universal families of stable vector bundles
of rank two and degree one on a nodal and cuspidal cubic curves and the
$r$--matrix coming from the universal family of semi--stable vector bundles of
rank two and degree zero on a nodal cubic curve.

This article is concluded  with a brief summary of analytical results in Section \ref{S:summary}.

\medskip
\noindent
\emph{Notation}.
Throughout this paper we work in the category of analytic spaces over
the field of complex numbers  $\mathbb{C}$, see \cite{PetRem}.
However, most of the results remain valid in the category of algebraic
varieties over an algebraically closed field $\kk$ of characteristic zero.
If $V,W$ are two complex vector spaces, $\Lin(V,W)$ denotes the vector space
of complex linear maps from $V$ to $W$. For  an additive category
$\mathsf{C}$, a pair of objects $X, Y \in \mathsf{C}$
and a pair of isomorphisms $X \stackrel{f}\lar X'$ and $Y \stackrel{g}\lar Y'$
we denote $\conj(f, g)$ the morphism of abelian groups
$\Hom_{\mathsf{C}}(X, Y) \lar \Hom_{\mathsf{C}}(X', Y')$ mapping
a morphism $h$ to the composition $g \circ h \circ f^{-1}$.

If $X$ is a complex projective variety, we denote by $\Coh(X)$ the category of
coherent $\mathcal{O}_{X}$-modules and by $\VB(X)$ its full subcategory of
locally free sheaves (holomorphic vector bundles). The torsion sheaf of length one,
supported at a closed point $y\in X$, is always denoted by $\CC_{y}$.
By $\Dbcoh(X)$ we denote the full subcategory of the derived category of the
abelian category of all $\mathcal{O}_{X}$-modules whose objects are those
complexes which have bounded and coherent cohomology.
The notation $\Perf(X)$ is used  for the full subcategory of $\Dbcoh(X)$ whose
objects are isomorphic to bounded complexes of locally free sheaves.
For a morphism of reduced complex spaces $E \stackrel{p}\lar T$ we denote
by $\breve{E}$ the regular locus of $p$.

A Weierstra\ss{} curve is a plane cubic curve given in homogeneous coordinates
by an equation $zy^2 = 4 x^3 -  g_2 xz^2 - g_3 z^3$, where $g_1, g_2 \in \CC$.
Such a curve is always irreducible. It is a smooth elliptic curve if and only
if $\Delta(g_2, g_3) = g_2^3 - 27 g_3^2  \ne 0$.

\medskip
\noindent
\emph{Acknowledgement}. The first-named author would like to thank A.~Kirillov,
A.~Stolin  and D.~van Straten for fruitful discussions and to T.~Henrich
for pointing out on numerous misprints in the previous version of this paper.
  The main work  on this article was carried out during the authors stay at Max-Planck Institut
f\"ur Mathematik in Bonn, at the Mathematical
Research Institute  Oberwolfach within  the  ``Research in Pairs'' programme
and during the visits of the second-named
author at the Johannes-Gutenberg University of Mainz and Friedrich-Wilhelms
University of Bonn supported
by Research Seed Funding at Mary Immaculate College. The first-named author
 was also supported by  the DFG grants Bu 1866/1-1 and Bu 1866/1-2.
\clearpage
\tableofcontents

\clearpage
\section{Yang--Baxter equations}\label{S:YBE}
In this section we are going to recall some standard results about Yang--Baxter
equations.
Let $\mathfrak{g}$ be a simple complex Lie algebra (throughout this paper
$\mathfrak{g} = \mathfrak{sl}_n(\mathbb{C})$),
$\langle \quad,\quad\rangle: \mathfrak{g} \times \mathfrak{g} \to \CC$
the Killing form.  The classical Yang--Baxter equation
is
\begin{equation}\label{E:CYBE2}
[r^{12}(y_1,y_2), r^{23}(y_2,y_3)] + [r^{12}(y_1,y_2), r^{13}(y_1,y_3)] +
[r^{13}(y_1,y_3), r^{23}(y_2,y_3)] = 0,
\end{equation}
where $r(x,y)$ is the germ of a meromorphic function of two complex variables
in a neighbourhood of $0$,  taking values in $\mathfrak{g} \otimes \mathfrak{g}$.
A solution of (\ref{E:CYBE2}) is called \emph{unitary}  if
$$
r^{12}(y_1,y_2) = - r^{21}(y_2,y_1)
$$
and \emph{non-degenerate} if $r(y_1, y_2) \in \mathfrak{g} \otimes
\mathfrak{g} \cong \mathfrak{g}^* \otimes
\mathfrak{g} \cong \End(\mathfrak{g})$ is invertible for generic $(y_1, y_2)$.
On the set of  solutions of (\ref{E:CYBE2})  there exists  a natural action of
the group of holomorphic function germs
$\phi: (\mathbb{C},0) \lar {\mathrm{Aut}}(\mathfrak{g})$ given by the rule
\begin{equation}\label{E:equiv}
r(y_1, y_2) \mapsto
\bigl(\phi(y_1) \otimes \phi(y_2)\bigr) r(y_1,y_2).
\end{equation}

\begin{proposition}[see \cite{BelavinDrinfeld2}]
Modulo the equivalence relation (\ref{E:equiv}) any non-degenerate unitary
solution of the equation (\ref{E:CYBE2})  is equivalent to  a solution $r(u,v)
= r(u-v)$ depending on the difference (or  the quotient) of spectral
parameters only.
\end{proposition}

\noindent
This means that equation (\ref{E:CYBE2}) is essentially equivalent to the
equation
\begin{equation}\label{E:CYBE1}
[r^{12}(x), r^{13}(x+y)] + [r^{13}(x+y), r^{23}(y)] + [r^{12}(x), r^{23}(y)] = 0.
\end{equation}
Although  the classical Yang--Baxter equation with one spectral parameter is
better adapted  for applications in mathematical physics, it seems that from
a geometric point of view equation (\ref{E:CYBE2}) is more natural.

Let $m = \dim(\mathfrak{g})$, $e_1, e_2, \dots,e_m$ be a basis
of $\mathfrak{g}$  and $e^1, e^2, \dots,e^m$ be the dual basis of $\mathfrak{g}$
with respect to the Killing form $\langle \quad,\quad\rangle$. Then
$\Omega = \sum_{i = 1}^m e^i \otimes e_i \in \mathfrak{g} \otimes
\mathfrak{g}$
is independent of the choice of a basis and is called the
\emph{Casimir element}.

\begin{theorem}[see Proposition 2.1 and Proposition 4.1 in
  \cite{BelavinDrinfeld}]
Let $r(y)$ be a non-constant non-degenerate solution of (\ref{E:CYBE1}). Then the tensor
$r(y)$
\begin{itemize}
\item has a simple pole at $0$ and
$
\res_{y=0}\bigl(r(y)\bigr) = \alpha \Omega \in \mathfrak{g} \otimes \mathfrak{g}, \alpha \in \mathbb{C}^*;
$
\item is automatically unitary, i.e. $r^{12}(y) = -r^{21}(-y)$.
\end{itemize}
\end{theorem}

As it was already mentioned in the introduction, there is the following
classification of non-degenerate solutions of (CYBE) due to Belavin and
Drinfeld.

\begin{theorem}[see Proposition 4.5 and Proposition 4.7 in
  \cite{BelavinDrinfeld}]
There are three types of non-degenerate solutions  of the classical
Yang--Baxter equation (\ref{E:CYBE1}):  elliptic, trigonometric and rational.
\end{theorem}

\noindent
Let us now consider some examples. Fix the following basis
$$
h =
\left(
\begin{array}{cc}
1 & 0 \\
0 & -1
\end{array}
\right),  \quad
e_{12} =
\left(
\begin{array}{cc}
0 & 1 \\
0 & 0
\end{array}
\right),  \quad
e_{21} =
\left(
\begin{array}{cc}
0 & 0 \\
1 & 0
\end{array}
\right)
$$
of the Lie algebra  $\mathfrak{g} = \mathfrak{sl}_2(\mathbb{C})$. Note
that $\Omega = \dfrac{1}{2} h\otimes h + e_{12} \otimes e_{21} + e_{21} \otimes
e_{12}$  is the Casimir element of $\mathfrak{sl}_2(\mathbb{C})$.

\medskip
\noindent
$\bullet$  Historically,  the first solution ever found was the rational
solution of Yang
$$
r_{\mathrm{rat}}(y) =
\frac{1}{y}\left(\frac{1}{2} h\otimes h + e_{12}\otimes e_{21} + e_{21} \otimes
e_{12}\right).
$$

\noindent
$\bullet$
A few years later,  Baxter discovered the  trigonometric solution
$$
r_{\mathrm{trg}}(y) = \frac{1}{2} \cot(y)h\otimes h +
\frac{1}{\sin(y)}(e_{12}\otimes e_{21} + e_{21} \otimes e_{12}).
$$

\noindent
$\bullet$   The following solution of elliptic type was found and studied by Baxter, Belavin and Sklyanin:
$$
r_{\mathrm{ell}}(y) =   \frac{\mathrm{cn}(y)}{\mathrm{sn}(y)} h \otimes h +
\frac{1+ \mathrm{dn}(y)}{\mathrm{sn}(y)}(e_{12}\otimes e_{21} +
e_{21} \otimes e_{12})  +
\frac{1 -  \mathrm{dn}(y)}
{\mathrm{sn}(y)}(e_{12} \otimes e_{12} + e_{21} \otimes e_{21}),
$$
where $\mathrm{cn}(y)$, $\mathrm{sn}(y)$ and $\mathrm{dn}(y)$ are doubly
periodic meromorphic functions on $\mathbb{C}$ with periods  $2$ and
$2\tau$. These functions also satisfy identities of the form
$f(y+1) = \varepsilon f(y)$ and $f(y + \tau) = \varepsilon f(y)$ with
$\varepsilon=\pm1$.

\medskip
\noindent
At first glance, all these solutions seem to be  completely
different. However, it is easy to see that
$$
\lim_{t  \to \infty} \frac{1}{t}r_{\mathrm{trg}}\left(\frac{y}{t}\right) =
r_{\mathrm{rat}}(y),
$$
hence the solution of Yang is a degeneration of Baxter's solution.
Moreover, there exist degenerations
 $\mathrm{dn}(y) \to 1$, $\mathrm{cn}(y) \to \cos(y)$ and
$\mathrm{sn}(y) \to \sin(y)$, when  the imaginary period $\tau$  tends to
infinity, see for example \cite[Section 2.6]{Lawden}.
Hence, both solutions of Baxter and Yang are degenerations
of the elliptic  solution. However, as we shall see later,
the theory of degenerations of $r$--matrices
 is more complicated as it might look like at first sight.

\medskip
\medskip
In this article we deal with a new type of Yang--Baxter equation, called
\emph{associative Yang--Baxter equation} (AYBE). It
appeared for the first time in a paper
of Fomin and Kirillov  \cite{FominKirillov}. Later, it was
studied by Aguiar \cite{Aguiar} in the framework of the theory
of infinitesimal Hopf algebras. The following version of the associative
Yang--Baxter equation with \emph{spectral parameters} is due to Polishchuk
\cite{Polishchuk1}. A special case of this equation was also considered by
Odesski and Sokolov \cite{OdeSok}.

\begin{definition}
An associative $r$-matrix is the germ of a meromorphic function in
 four variables
$$
r: \bigl(\mathbb{C}^4_{(v_1, v_2; y_1, y_2)}, 0\bigr) \lar
\Mat_{n \times n}(\mathbb{C}) \otimes \Mat_{n \times n}(\mathbb{C})
$$
holomorphic on
$\Bigl(\mathbb{C}^4\setminus V\bigl((y_1 -y_2)(v_1 -v_2)\bigr), 0\Bigr)$
and
satisfying the equation
\clearpage
\begin{equation}\label{E:AYBE1}
r(v_1, v_2; y_1, y_2)^{12} r(v_1, v_3; y_2, y_3)^{23} =
r(v_1, v_3; y_1, y_3)^{13} r(v_3, v_2; y_1, y_2)^{12} +
\end{equation}
$$
+ r(v_2, v_3; y_2, y_3)^{23} r(v_1, v_2; y_1, y_3)^{13}.
$$
Such a matrix  is called unitary if
\begin{equation}\label{E:AYBEunitary}
r(v_1, v_2; y_1, y_2)^{12} = - r(v_2, v_1; y_2, y_1)^{21}.
\end{equation}
\end{definition}

\noindent
On the set of solutions of the equation (\ref{E:AYBE1})
there exists a natural
equivalence relation.

\begin{definition}[see section 1.2 in \cite{Polishchuk1}]\label{D:gaugeAYBE}
Let $\phi: (\mathbb{C}^2, 0) \to \GL_n(\mathbb{C})$ be the germ of a
holomorphic function and $r(v_1, v_2; y_1, y_2)$ be a solution of (AYBE) then
\begin{equation}\label{E:gaugeAYBE}
r'(v_1, v_2; y_1, y_2) =
\bigl(\phi(v_1; y_1) \otimes \phi(v_2; y_2)\bigr) r(v_1, v_2; y_1, y_2)
\bigl(\phi(v_2; y_1)^{-1} \otimes \phi(v_1; y_2)^{-1}\bigr)
\end{equation}
is again  a  solution of (\ref{E:AYBE1}). Two such  tensors $r$ and $r'$ are
called gauge equivalent. Note that if the matrix $r$ is unitary
then $r'$ is unitary, too.
\end{definition}

\vspace{2mm}

\begin{example}\label{Ex:specialgage}
Let $r(v_1, v_2; y_1, y_2) \in \Mat_{n \times n}(\CC) \otimes \Mat_{n \times n}(\CC)$ be a solution
of (\ref{E:AYBE1}),
$c \in \CC$ and  $\phi = \exp(c vy) \cdot \mathbbm{1}: (\CC^2, 0) \lar \GL_n(\CC)$ be a gauge transformation. Then
$
\exp\bigl(c(v_2 - v_1)(y_2 - y_1)\bigr)r(v_1, v_2; y_1, y_2)$
 is a gauge equivalent solution of (AYBE).
\end{example}

\vspace{2mm}

\begin{lemma}\label{L:useful}
Let $r(v_1, v_2; y_1, y_2)$ be a unitary solution of the
associative Yang--Baxter equation (\ref{E:AYBE1}). Then $r$ also
satisfies the ``dual'' equation
\begin{equation}\label{E:AYBEdual}
r(v_2, v_3; y_2, y_3)^{23} r(v_1, v_3; y_1, y_2)^{12} =
r(v_1, v_2; y_1, y_2)^{12} r(v_2, v_3; y_1, y_3)^{13} +
\end{equation}
$$
+ r(v_1, v_3; y_1, y_3)^{13} r(v_2, v_1; y_2, y_3)^{23}.
$$
\end{lemma}

\begin{proof}
Let $\tau$ be the linear automorphism of  $\Mat_{n \times n}(\CC) \otimes \Mat_{n \times n}(\CC)$
defined by $\tau(a \otimes b) = b \otimes a$.  Applying the
operator $\tau \otimes \mathbbm{1}$ to the equation (\ref{E:AYBE1}),  we obtain:
$$
r(v_1, v_2; y_1, y_2)^{21} r(v_1, v_3; y_2, y_3)^{13} =
r(v_1, v_3; y_1, y_3)^{23} r(v_3, v_2; y_1, y_2)^{21} +
$$
$$
+
r(v_2, v_3; y_2, y_3)^{13} r(v_1, v_2; y_1, y_3)^{23}.
$$
Using the unitarity condition (\ref{E:AYBEunitary}) we get:
$$
-r(v_2, v_1; y_2, y_1)^{12} r(v_1, v_3; y_2, y_3)^{13} =
- r(v_1, v_3; y_1, y_3)^{23} r(v_2, v_3; y_2, y_1)^{12} +
$$
$$
+ r(v_2, v_3; y_2, y_3)^{13} r(v_1, v_2; y_1, y_3)^{23}.
$$
After  the change of variables $v_1 \leftrightarrow v_2,
v_3 \leftrightarrow v_3$ and
$y_1 \leftrightarrow y_2, y_3 \leftrightarrow y_3$,
we obtain the equation (\ref{E:AYBEdual}).
\end{proof}

\medskip
\noindent
Assume a unitary solution $r(v_1, v_2; y_1, y_2)$ of the associative
Yang--Baxter equation (\ref{E:AYBE1}) depends on the difference $v= v_1 - v_2$
of the first pair of spectral parameters only. For the sake of simplicity, we
shall use the notation:
$r(v_1, v_2; y_1, y_2) = r(v_1 - v_2; y_1,y_2) = r(v;y_1,y_2).$
  Then the equation
(\ref{E:AYBE1}) can be rewritten as
\clearpage

\begin{equation}\label{E:AYBE2}
r(u; y_1, y_2)^{12} r(u+v; y_2, y_3)^{23}  =
r(u+v;y_1, y_3)^{13} r(-v;y_1, y_2)^{12}
\end{equation}
$$
+ r(v; y_2, y_3)^{23} r(u; y_1, y_3)^{13}.
$$
\begin{remark}
It will be shown in Theorem  \ref{T:diffinvbpar} that any solution $r$ of the
associative Yang--Baxter equation (\ref{E:AYBE1}) obtained from a universal
family of stable vector bundles on an irreducible genus one curve, is  gauge
equivalent to a solution $r'$ depending on the difference $v_1 - v_2$ only.
\end{remark}

\medskip
\noindent
Let $A$ be the algebra of germs of meromorphic functions
$f: \bigl(\mathbb{C}^4_{(v_1, v_2; w_1, w_2)}, 0\bigr) \lar \CC$
holomorphic on
$\Bigl(\mathbb{C}^4\setminus V\bigl((v_1 -v_2)(w_1 -w_2)\bigr), 0\Bigr)$.
A solution of the equation (\ref{E:AYBE2}) defines an element
$$
r \in \Mat_{n \times n}(A) \otimes_A \Mat_{n \times n}(A) \cong
A \otimes_\CC \Mat_{n \times n}(\CC) \otimes_\CC   \Mat_{n \times n}(\CC).
$$
In a similar way, for any integer  $m \ge 3$ denote by $B$  the algebra of
germs of meromorphic functions
$f: \bigl(\mathbb{C}^{2m}_{(x_1, \dots, x_m; y_1, \dots, y_m)},0\bigr) \lar \CC$
holomorphic on $\bigl(\mathbb{C}^{2m}\setminus  D, 0\bigr)$, where
$D$ is the divisor
$$
D = V\left(\prod\limits_{i \ne j}(x_i-x_j)(y_i -
y_j)\right).
$$
Next, for any pair of indices $1 \le i \ne j \le m$ we have
\begin{itemize}
\item
  a ring homomorphism
  $\psi^{ij}: A \lar B$ which sends a function  $f(v_1, v_2; w_1, w_2)$ to
  $f(x_i, x_j; y_i, y_j)$;
\item a ring homomorphism $k^{ij}: B \lar B$ defined as
  $$
  f(\dots, x_i, \dots, x_j, \dots; \underline{y}) \mapsto
  f(\dots, x_j, \dots, x_i, \dots; \underline{y});
  $$
\item a morphism
  $\varrho_{ij}: \Mat_{n \times n}(\CC)^{\otimes 2} \lar \Mat_{n \times n}(\CC)^{\otimes m}$
  mapping a simple tensor
  $a \otimes b$ to $1 \otimes \dots 1 \otimes a \otimes 1
  \otimes \dots 1 \otimes b \otimes 1 \otimes \dots \otimes 1$,
  where $a$ and $b$ belong to  the $i$-th and $j$-th components respectively.
\end{itemize}
In these notations, consider
$$
\Psi^{ij} := \psi_{ij} \otimes \varrho_{ij} \, : \,\,
A \otimes_\CC \Mat_{n \times n}(\CC)^{\otimes 2} \lar B \otimes_\CC
\Mat_{n \times n}(\CC)^{\otimes m}.
$$
For  example $\Psi^{13}\bigl(f(v_1, v_2; w_1, w_2) \otimes a \otimes b\bigr)
 = f(x_1, x_3; y_1, y_3) \otimes a \otimes 1 \otimes b$.
Next, we set
  $r^{ij} := \Psi^{ij}(r) \in B \otimes_\CC
\Mat_{n \times n}(\CC)^{\otimes m}$ and
$$
K^{ij} = k^{ij} \otimes \mathbbm{1} \, : \,\,  B \otimes_\CC \Mat_{n \times n}(\CC)^{\otimes m} \lar
B \otimes_\CC \Mat_{n\times n}(\CC)^{\otimes m}.
$$
Consider the linear operator
$$
\tilde{r}^{ij} = r^{ij} \circ K^{ij} \, : \,\, B \otimes_\CC \Mat_{n \times n}(\CC)^{\otimes m} \lar
B \otimes_\CC \Mat_{n\times n}(\CC)^{\otimes m},
$$
which is  the composition of $K^{ij}$ and the multiplication with the element
$r^{ij}$. For any $1 \le i \le m$  consider the differential operator
$$
\partial_i = \frac{\partial}{\partial x_i} \otimes \mathbbm{1} \,: \, \,
B \otimes_\CC \Mat_{n \times n}(\CC)^{\otimes m} \lar
B \otimes_\CC \Mat_{n \times n}(\CC)^{\otimes m}.
$$
Next, for any $\kappa \in \CC$ let
$$
\theta_i := \kappa \partial_i + \sum\limits_{j \ne i} \tilde{r}^{ij}
\,: \, \,
B \otimes_\CC \Mat_{n \times n}(\CC)^{\otimes m} \lar
B \otimes_\CC \Mat_{n \times n}(\CC)^{\otimes m}
$$
be the \emph{Dunkl operator} of level $\kappa$. The following result
was explained to the first-named author by Anatoly Kirillov, see also \cite{Kirillov}.

\begin{proposition}\label{P:Kirillov}
Let $r(v; y_1,y_2) \in \Mat_{n \times n}(A) \otimes_A \Mat_{n \times n}(A)$
be a unitary solution of the equation (\ref{E:AYBE2}),  $\kappa \in \CC$
be a scalar  and $\theta_i$ be the Dunkl operator of level $\kappa$   defined above.
Then for all $
1 \le i, j \le m$  we have:  $
\bigl[\theta_i, \theta_j\bigr] = 0$.
\end{proposition}

\begin{proof} First note that we have:
$$
\Bigl(\frac{\partial}{\partial x_i} +
\frac{\partial}{\partial x_j}\Bigr) r(x_i - x_j; y_i, y_j) = 0,
$$
which implies the equality
$
\bigl[\partial_i + \partial_j, \tilde{r}^{ij}\bigr] = 0.
$
Next, the Yang--Baxter relations
 (\ref{E:AYBE1}) and (\ref{E:AYBEdual}) yield  that for any
triple of  mutually different  indices $1 \le i <  j <  k \le m$ we have:
$$
\tilde{r}^{ij} \tilde{r}^{jk} = \tilde{r}^{jk} \tilde{r}^{ik} +
\tilde{r}^{ik} \tilde{r}^{ij} \quad \mbox{\textrm{and}} \quad
\tilde{r}^{jk} \tilde{r}^{ij} = \tilde{r}^{ik} \tilde{r}^{jk} +
\tilde{r}^{ij} \tilde{r}^{ik}.
$$
From the  unitarity of $r$  it follows that $\tilde{r}^{ij} = - \tilde{r}^{ji}$ for
all $1 \le i \ne j \le m$.
Finally, the following two equalities are obvious:
$$
\bigl[\tilde{r}^{ij}, \tilde{r}^{kl}\bigr] = 0, \quad
\bigl[\partial_i, \tilde{r}^{kl}\bigr] = 0
$$
where $1 \le i,j,k,l \le m$ are mutually distinct. Combining these
equalities together, we  obtain the claim.
\end{proof}

\begin{remark} The above proposition means that
to any unitary solution of the associative Yang--Baxter equation
$(\ref{E:AYBE2})$ one can attach
a very interesting second order differential operator
$$
H = \theta_1^2 + \theta_2^2 + \dots + \theta_m^2:  \quad
B \otimes_\CC
\Mat_{n\times n}(\CC)^{\otimes m} \lar B \otimes_\CC
\Mat_{n \times n}(\CC)^{\otimes m}.
$$
These operators are ``matrix versions'' of  the Hamiltonians  considered
in the work of Buchstaber, Felder and Veselov \cite{BFV}.
\end{remark}

\medskip
\noindent
Another motivation to study solutions of the equation (\ref{E:AYBE2})
is provided  by their  close connection with the theory of the
classical Yang--Baxter equation.

\begin{lemma}[see Lemma 1.2 in \cite{Polishchuk1}]\label{L:AYBEandCYBE}
Let $r(v; y_1, y_2)$ be a unitary solution of the associative Yang--Baxter
equation (\ref{E:AYBE2}) and
$\pr: \Mat_{n\times n}(\mathbb{C}) \to \mathfrak{sl}_n(\mathbb{C})$
be the projection along the scalar matrices, i.e.\/
$\pr(A) = A - \frac{\displaystyle \tr(A)}{\displaystyle n}  \cdot \mathbbm{1}$.
Assume that $(\pr \otimes \pr)r(v; y_1, y_2)$ has a limit as
$v \to 0$. Then
$$
\bar{r}(y_1, y_2)  := \lim_{v\to 0}(\pr \otimes \pr) r(v; y_1, y_2)
$$
is a unitary solution of  the classical Yang--Baxter equation (\ref{E:CYBE2}).
\end{lemma}

\begin{proof}
 First note that   (\ref{E:AYBEdual}) implies  the equality
$$
r(v; y_2, y_3)^{23} r(u+v; y_1, y_2)^{12}  =
r(u;y_1, y_2)^{12} r(v;y_1, y_3)^{13} +
$$
$$
+ r(u+v; y_1, y_3)^{13} r(-u; y_2, y_3)^{23}.
$$
Using the change of variables $u \mapsto -v$ and $v \mapsto u + v$,
 we obtain  the relation
$$
r(u + v; y_2, y_3)^{23} r(u; y_1, y_2)^{12}  =
r(-v; y_1, y_2)^{12} r(u+v; y_1, y_3)^{13} +
$$
$$
+
r(u; y_1, y_3)^{13} r(v; y_2, y_3)^{23}.
$$
Subtracting this equation from (\ref{E:AYBE2}) we  get
\begin{equation}\label{E:derivingCYBE}
[r(-v; y_1, y_2)^{12}, r(u+v; y_1, y_3)^{13}] + [r(u; y_1, y_2)^{12},
r(u+v; y_2, y_3)^{23}] +
\end{equation}
$$
+ [r(u; y_1, y_3)^{13}, r(v; y_2, y_3)^{23}] = 0.
$$
By definition,  the  function $r(v; y_1, y_2)$ is meromorphic, hence we can
write its Laurent expansion:
$
r(v; y_1, y_2) = \sum_{\alpha \in \mathbb{Z}}
r_\alpha(y_1, y_2) v^\alpha,
$
where $r_\alpha(y_1, y_2)$ are meromorphic and $r_\alpha = 0 $ for $\alpha
 \ll 0$.
Since we have assumed that  $(\pr \otimes \pr)r(v; y_1, y_2)$  is regular
with respect to $v$  in a neighbourhood of $v = 0$, we have:
$(\pr \otimes \pr) r_\alpha(y_1, y_2) = 0$ for all $\alpha \le -1$.
This implies, if $\alpha\le-1$, that
$$
r_\alpha(y_1, y_2) = s_\alpha(y_1, y_2) \otimes \mathbbm{1} +
\mathbbm{1}  \otimes t_\alpha(y_1, y_2)
$$
for some matrix-valued functions $s_\alpha(y_1, y_2)$ and $t_\alpha(y_1,
y_2)$.
Hence,
$$
(\pr \otimes \pr \otimes \pr)\bigl[r_\alpha^{ij}, r_\beta^{lk}\bigr] = 0
$$
for arbitrary permutations $(ij) \ne (lk)$ and all indices
$\alpha \le -1$ and
$\beta \in \mathbb{Z}$.
The claim of Lemma \ref{L:AYBEandCYBE} follows by applying $\pr \otimes \pr \otimes \pr$  to the
equation (\ref{E:derivingCYBE})
and taking the limit $u, v \to 0$.
\end{proof}

\medskip
\noindent
The statement of
Lemma \ref{L:AYBEandCYBE} leads to the following question. Let
$r= r(v;y_1,y_2)$ be a unitary solution of the associative Yang--Baxter equation (\ref{E:AYBE2})
satisfying the conditions of Lemma \ref{L:AYBEandCYBE} and
$s=s(v; y_1, y_2)$ be an equivalent solution in the sense of
Definition \ref{D:gaugeAYBE}.
Are the corresponding solutions $\bar{r}(y_1, y_2)$ and
$\bar{s}(y_1, y_2)$ of the classical Yang--Baxter equation also
gauge equivalent?

The answer on this question is affirmative, if
one imposes an additional restriction on the function $r$.
Namely, we assume that the Laurent expansion of  $r$  has the form:
\begin{equation}\label{E:AnsatzAYBE}
r(v; y_1, y_2) =
\frac{\mathbbm{1} \otimes \mathbbm{1}}{v} + r_{0}(y_1, y_2) + v r_1(y_1, y_2) +
v^2 r_2(y_1, y_2) +\dots
\end{equation}
Then the following proposition  is true.

\begin{proposition}
Let $r: (\CC^3_{(v; y_1, y_2)}, 0) \to \Mat_{n \times n}(\CC) \otimes \Mat_{n
  \times n}(\CC)$ be a unitary solution of the associative Yang--Baxter
equation (\ref{E:AYBE2}) having  a Laurent expansion of the form
(\ref{E:AnsatzAYBE}) and $\bar{r}_0(y_1, y_2)$ be the corresponding solution
of the classical Yang--Baxter equation.
If $\phi: (\CC^2_{(v;y)}, 0) \to \GL_n(\CC)$ is the germ of a holomorphic
function such that
\begin{equation}\label{E:somegauge}
s(v_1, v_2; y_1, y_2) :=
\bigl(\phi(v_1; y_1) \otimes \phi(v_2; y_2)\bigr) r(v; y_1, y_2)
\bigl(\phi(v_2; y_1)^{-1} \otimes \phi(v_1; y_2)^{-1}\bigr)
\end{equation}
is again a function of $v = v_2 - v_1$, then we have
$$
s(v;y_1, y_2) = \frac{1}{v}\mathbbm{1} \otimes \mathbbm{1} + s_0(y_1, y_2) +
v s_1(y_1, y_2) +
v^2 s_2(y_1, y_2) + \dots
$$
and moreover, $\bar{r}_0(y_1, y_2)$ and $\bar{s}_0(y_1, y_2)$ are equivalent
in the sense of the relation  (\ref{E:equiv}).
\end{proposition}

\begin{proof} We denote $v = v_1 - v_2$ and $h = v_2$. Then
$v_1 = v +h$ and using the Taylor expansion of $\phi$ with respect to $v$, we
may rewrite (\ref{E:somegauge}) in the form
\begin{multline*}
  \left(\phi(h;y_1) + v \phi'(h;y_1) + \frac{v^2}{2} \phi''(h; y_1) +
    \dots\right) \otimes \phi(h; y_2)\cdot   \\
  \cdot \left(\frac{\mathbbm{1} \otimes \mathbbm{1}}{v} +
    r_0(y_1, y_2) + v r_1(y_1, y_2) + \dots\right)
\end{multline*}
$$
= \left(\sum\limits_{i \in \mathbb{Z}} s_i(y_1, y_2) v^i\right) \cdot
\left(\phi(h; y_1) \otimes \left(\phi(h; y_2) + v\phi'(h; y_2) + \frac{v^2}{2}
\phi''(h; y_2) + \dots\right)\right),
$$
where $\phi'(v;y)$ and $\phi''(v;y)$ are the partial derivatives of
$\varphi(v;y)$ of the first and the second order  with respect to $v$.
This equality implies that $s_i(y_1, y_2) = 0$ for  $i \le -2$ and
$s_{-1}(y_1, y_2) = \mathbbm{1} \otimes \mathbbm{1}$. Moreover, we have:
$$
  \bigl(\phi(h; y_1) \otimes \phi(h; y_2)\bigr) \cdot r_0(y_1, y_2) \cdot
  \bigl(\phi(h; y_1)^{-1} \otimes \phi(h; y_2)^{-1}\bigr) =
$$
$$
=
s_0(y_1, y_2) + \mathbbm{1}\otimes \phi'(h; y_2) \phi(h; y_2)^{-1}
- \phi'(h; y_1) \phi(h; y_1)^{-1} \otimes \mathbbm{1}.
$$
Let $\bar{\phi}(y) := \phi(0; y)$. Applying the operator
$\pr \otimes \pr$ to the last equality and putting $h = 0$ we obtain:
$$
\bar{s}_0(y_1, y_2) = \bigl(\bar{\phi}(y_1) \otimes \bar{\phi}(y_2)\bigr) \cdot
\bar{r}_0(y_1, y_2)\cdot  \bigl(\bar{\phi}(y_1)^{-1} \otimes \bar{\phi}(y_2)^{-1}\bigr).
$$
This  implies the claim.
\end{proof}

\begin{remark}
It was proven by Polishchuk \cite[Theorem 2]{Polishchuk1} that
all solutions of the associative Yang--Baxter equation
(\ref{E:AYBE2}) arising from a universal family of \emph{stable}
vector bundles on an (irreducible) projective cubic curve
satisfy the conditions of Lemma \ref{L:AYBEandCYBE}. However, we do not know
a  conceptual
explanation  of the fact that all these  solutions have a Laurent expansion of the form
(\ref{E:AnsatzAYBE}),  although in turns out to be so in all the examples known so far.
Later,  we shall see that the equation (\ref{E:AYBE2}) has many  solutions $r(v; y_1, y_2)$
with  higher order poles with respect to $v$.  Some of them can be obtained
by the same geometric method, applied to certain families  of
\emph{semi-stable} vector  bundles, see Subsection \ref{SS:solfromSS}.
However, they do not project
to a solution of the classical Yang--Baxter equation.
\end{remark}

\medskip
\noindent
Finally,  assume that a solution of (\ref{E:AYBE1})  has the form
  $r(v_1, v_2; y_1, y_2) = r(v_2 - v_1; y_2 - y_1)$. Then
 the associative Yang--Baxter equation can be rewritten as
\begin{equation}\label{E:AYBE3}
r(u; x)^{12} r(u+v; y)^{23}  = r(u + v; x+y)^{13} r(-v; x)^{12} + r(v; y)^{23} r(u; x+y)^{13}.
\end{equation}
This is the form of the associative Yang--Baxter equation introduced and studied by Polishchuk
in \cite{Polishchuk1} and  \cite{Polishchuk2}.

\begin{definition}
A  solution $r(y)$ of the classical Yang--Baxter equation (\ref{E:CYBE1})
has an infinitesimal symmetry, if there exists an element $a \in \mathfrak{g}$
such that
$$
\bigl[r(y), a \otimes \mathbbm{1} + \mathbbm{1} \otimes a\bigr] = 0.
$$
\end{definition}

\noindent
For example, let $r(y) = r_{\mathrm{rat}}(y) = \frac{\displaystyle 1}{\displaystyle y}\Omega$
be Yang's solution for $\mathfrak{sl}_2(\mathbb{C})$,
then $$[r(y), a \otimes \mathbbm{1} + \mathbbm{1} \otimes a] = 0$$
for any $a \in \mathfrak{g}$.

An important reason to study unitary solutions of  the  equation
\eqref{E:AYBE3} satisfying \eqref{E:AnsatzAYBE} is explained by the following
theorem.

\begin{theorem}[see Theorem 1.4 of  \cite{Polishchuk2} and Theorem 6
of \cite{Polishchuk1}]\label{T:QYBEandCYBE}
Let $r(v;y)$ be a non-degenerate unitary solution
of the associative Yang--Baxter equation (\ref{E:AYBE3}) with Laurent expansion of the form
(\ref{E:AnsatzAYBE}). Then we have:

\noindent
$\bullet$  The function
$
\bar{r}_0(y) := (\pr \otimes \pr)\bigl(r_0(y)\bigr)
$
is a non-degenerate unitary solution of  the classical Yang--Baxter equation
(\ref{E:CYBE1}).

\noindent
$\bullet$
 If $\bar{r}_0(y)$ is either periodic (elliptic or trigonometric),
or without infinitesimal symmetries, then for a fixed $v = v_0$ the tensor  $r(v_0;y)$
satisfies the  \emph{quantum
Yang--Baxter} equation:
\begin{equation}\label{E:QYBE}
r(v_0; x)^{12} r(v_0; x+y)^{13} r(v_0; y)^{23} =
r(v_0; y)^{23} r(v_0; x+y)^{13} r(v_0; x)^{12}
\end{equation}

\noindent
$\bullet$
 If $\bar{r}_0(y)$ does not have infinitesimal symmetries and if
  $s(v;y)$ is another  solution of (\ref{E:AYBE3}) of the form
  (\ref{E:AnsatzAYBE}) and such that $\bar{r}_0(y) = \bigl(\pr\otimes\pr\bigr)(s_0(y))$,
  then there exist $\alpha_1 \in \mathbb{C}^*$ and $\alpha_2 \in \CC$ such
  that $s(v;y) = \alpha_1 \exp(\alpha_2 vy) r(v;y)$. In other words, under
  these conditions, a solution of the associative Yang--Baxter equation  $r(v;y)$
is uniquely determined by the corresponding solution of the classical Yang--Baxter equation
$\bar{r}_0(y)$  up  a gauge transformation described in Example
\ref{Ex:specialgage} and  rescaling.
\end{theorem}

\begin{remark} Theorem \ref{T:QYBEandCYBE} gives an explicit recipe to lift a non-degenerate
solution of the classical Yang--Baxter equation to a
 solution of the quantum
Yang--Baxter equation.  Existence of such quantization is known due to a result of
Etingof and Kazhdan \cite{EtingofKazhdan}.
Moreover,
it was proven by Polishchuk in \cite{Polishchuk1} that any elliptic solution
of the classical Yang--Baxter equation (\ref{E:CYBE1}) can be lifted to a
solution of (\ref{E:AYBE3}) having a Laurent expansion of the form
(\ref{E:AnsatzAYBE}). However, Schedler showed in \cite{Schedler} that there exist
trigonometric solutions of (CYBE), which  \emph{can not}  be lifted to a
solution of the associative Yang--Baxter equation (\ref{E:AYBE3}) of the form (\ref{E:AnsatzAYBE}).
\end{remark}

\section{Polishchuk's construction}\label{S:PolConstr}

Let $X$ be a connected projective Gorenstein variety of dimension $n$
over a field $\kk$ and $\Perf(X)$ be the triangulated category of perfect complexes, i.e.
a full subcategory of the derived category $D(X)= \Dbcoh(X)$ consisting of complexes quasi-isomorphic to
 bounded complexes of locally free $\kO_X$--modules.

We denote by  $\omega_X$
 the   dualizing sheaf on $X$.  This means (see for example
\cite[Section III.7]{Hartshorne}) that we have an isomorphism
$t: H^n(\omega_X) \to \kk$, also called  a \emph{trace map},
 such that for any coherent sheaf
$\kF \in \Coh(X)$ the pairing
$$
H^n(\kF) \times \Hom_X(\kF, \omega_X) \lar H^n(\omega_X) \stackrel{t}\to
\kk
$$
is non-degenerate.

\begin{remark}
Such a  map $t$ is defined only up to a non-zero constant. However, it will be
explained later that in the case of reduced projective Gorenstein curves with
trivial canonical sheaf there exists a ``canonical''  choice for $t$, see
subsection \ref{subsec:geometricMassey}, directly after Theorem \ref{T:Kunz}.
\end{remark}

\noindent
Let  $\mathbb{S}: \Dbcoh(X) \lar \Dbcoh(X)$ be the functor
given by the rule $\mathbb{S}(\mathbb{\kF})
= \kF \stackrel{\mathbb{L}}\otimes \omega_X[n]$.
For  a perfect complex $\kF$,
let $\tr_{\kF}: \Hom_{D(X)}\bigl(\kF, \mathbb{S}(\kF)\bigr) \to \kk$ be the
morphism
$$\Hom_{D(X)}(\kF, \kF \stackrel{\mathbb{L}}\otimes \omega_X[n])
\stackrel{\cong}\lar
\Hom_{D(X)}(\kO, \kF^\vee \stackrel{\mathbb{L}}\otimes \kF
\stackrel{\mathbb{L}}\otimes \omega_X[n])
\lar H^n(\omega_X) \stackrel{t}\lar \kk,$$
 where the first arrow is a canonical isomorphism and the second
is induced by the canonical evaluation morphism
$\kF^\vee \stackrel{\mathbb{L}}\otimes \kF \lar \kO$.

The following theorem seems to be  well-known, however we were not able to
find its proof in the literature and therefore sketch  it here.

\begin{theorem}\label{T:SerreDuality}
Let $p: X \to Y = \Spec(\kk)$ be a connected projective Gorenstein variety of dimension $n$
over a field $\kk$. Then for any $\kF \in \Perf(X)$ and $\kG \in \Dbcoh(X)$ the pairing
$$
\langle \,\, ,\,\rangle_{\kF, \, \kG}: \Hom_{D(X)}(\kF, \kG) \times
\Hom_{D(X)}\bigl(\kG, \mathbb{S}(\kF)\bigr)
\lar \kk,
$$
given by the formula $\langle f, g\rangle_{\kF, \,\kG} := \tr_{\kF}(gf)$,  is
non-degenerate.\footnote{We thank Amnon Neeman for helping us at this place.}
Moreover, if the complex $\kG$ is also perfect, we have: $\tr_{\kF}(gf) =
\tr_{\kG}\bigl(\mathbb{S}(f) g\bigr)$.
\end{theorem}

\begin{proof}
 Recall that by the Grothendieck duality the functor
$\mathbb{R}p_*: \Dbcoh(X) \to \Dbcoh(Y)$ has a right adjoint
$p^!: \Dbcoh(Y) \to \Dbcoh(X)$, see \cite{RD, Neeman, Conrad}.
Moreover, $p^!(\kk) \cong \omega_X[n]$ and the adjunction  morphism
$\tilde{t}: \mathbb{R}p_* p^!(\kk) \to \kk$ coincides up to a non-zero
constant with the morphism
$$
\mathbb{R}p_* p^!(\kk) \stackrel{\can}\lar
 H^0\bigl(\mathbb{R}p_* p^!(\kk)\bigr) \cong H^n(\omega_X)
\stackrel{t}\lar \kk,
$$
see the proof of \cite[Theorem 5.12]{Conrad}).
Next, we have the following canonical isomorphisms:
$$
\Hom_{D(X)}\bigl(\kG, \kF \stackrel{\mathbb{L}}\otimes p^!(\kk)\bigr)
\cong \Hom_{D(X)}\bigl(\kF^\vee \stackrel{\mathbb{L}}\otimes \kG,  p^!(\kk)\bigr)
\cong
$$
$$
\cong
\Hom_{D(X)}\bigl({\mathcal R}{\mathcal Hom}_X(\kF, \kG),  p^!(\kk)\bigr)
\cong
\Hom_{D(Y)}\bigl(\mathbb{R}p_*({\mathcal R}{\mathcal Hom}_X(\kF, \kG)),  \kk\bigr)
\cong
$$
$$
\cong
\Hom_Y\Bigl(H^0\bigl(\mathbb{R}p_*({\mathcal R}{\mathcal Hom}_X(\kF, \kG))\bigr),  \kk\Bigr)
\cong \Hom_{D(X)}(\kF, \kG)^*.
$$
One can check that the image of a  morphism
$g \in \Hom_{D(X)}\bigl(\kG, \mathbb{S}(\kF)\bigr)$
under this chain of isomorphisms
is the functional $\tr_\kF(g \circ -): \Hom_{D(X)}\bigl(\kF, \kG\bigr) \to
\kk$ up to a constant not depending on $f$.
This implies that the bilinear form
$\langle \,\, ,\,\rangle_{\kF, \,  \kG}$ is non-degenerate. In particular,
$\mathbb{S}$ is a Serre functor of the category  $\Perf(X)$ and the equality
 $\tr_{\kF}(gf) = \tr_{\kG}\bigl(\mathbb{S}(f) g\bigr)$ is automatically true,
see  \cite[Lemma I.1.1]{ReitenvandenBergh}.
\end{proof}

\begin{corollary}
Let $X$ be a connected projective Gorenstein variety over $\kk$  of dimension
$n$ such that its dualizing sheaf is trivial. Let $\omega \in  H^0(\omega_X)$
be a non-zero global section of $\omega_X$. Then for any pair of perfect
complexes $\kF, \kG$ on $X$ the pairing
$\langle \,\,,\,\rangle_{\kF, \, \kG}^{\omega}:
\Hom_{D(X)}(\kF, \kG) \times
\Hom_{D(X)}(\kG, \kF[n])  \lar
\Hom_{D(X)}(\kF, \kF[n]) \xrightarrow{\omega_{*}}
\Hom_{D(X)}\bigl(\kF, \mathbb{S}(\kF)\bigr) \stackrel{\tr_\kF}\lar
\kk
$
is non-degenerate.
\end{corollary}

Let $X$ be a reduced Gorenstein curve over the field $\mathbb{C}$.
 By \cite[Chapter VIII]{AltmanKleiman}  or by
\cite[Appendix B]{Conrad} the dualising sheaf  $\omega_X$
is isomorphic to the sheaf
of  \emph{regular} or \emph{Rosenlicht's} differential 1-forms
$\Omega_X = \Omega^{1,R}_X$.
If $X$ is smooth, then $\Omega_X$ coincides with the sheaf of holomorphic
1-forms.  For $X$ singular the definition is as follows.

\begin{definition}\label{D:Rosenlicht}
Let $X$ be a reduced projective Gorenstein curve,  $n: \widetilde{X} \to X$
its normalization. Denote by   $\Omega^M_X$ and $\Omega^M_{\widetilde{X}}$
the sheaves  of meromorphic differential  1-forms on $X$ and $\widetilde{X}$
respectively.  Observe that
$\Omega^M_X = n_*(\Omega^M_{\widetilde{X}}).
$
 Then $\Omega_X$ is defined to be the subsheaf of $\Omega^M_X$ such that for
any open subset $U \subseteq X$ one has
$$
\Omega_X(U) = \left\{\omega   \in \Omega^{M}_{\widetilde{X}}\bigl(n^{-1}(U)\bigr)\left|
\forall p \in U, \forall f \in \kO_{X}(U):\;
\sum\limits_{i= 1}^t \res_{p_i}\bigl((f\circ n)\omega\bigr) = 0
\right.\right\},
$$
where $\{p_1,p_2,\dots,  p_t\} = n^{-1}(p)$.
\end{definition}

A reduced projective curve $E$  whose dualizing sheaf  $\omega_E$ is
isomorphic to the  structure sheaf is Gorenstein and has arithmetic genus
one. For example, reduced plane cubics, Kodaira cycles and generic
configurations of $n+1$ lines in $\PP^n$ passing through a given point are of
this type.
In  what follows, for such a curve $E$, we identify $\omega_{E}$ wiht
$\Omega_{E}$ and fix a global section $\omega \in H^0(\Omega_E)$ giving an
isomorphism $\omega: \kO \to \Omega_E$ and a  trace map
$t^{\omega}: H^1(\kO) \xrightarrow{\omega_\ast} H^1(\Omega_E)
\xrightarrow{t} \mathbb{C}$.

A characteristic property of reduced projective curves
with trivial dualizing sheaf is a very special form of the Serre duality.
By Theorem
\ref{T:SerreDuality} we have the following result.

\begin{proposition}\label{P:SD}
  Let $E$ be a reduced projective curve with trivial dualizing sheaf and
  $\kE, \kF \in \Perf(E)$. Then the map
  $$
  \langle \,\, , \,  \rangle_{\kE, \,  \kF}^{\omega}:
  \Hom(\kE, \kF) \otimes \Hom(\kF, \kE[1]) \xrightarrow{\circ}
  \Hom(\kE, \kE[1]) \stackrel{\Tr_\kE}\lar H^1(\kO) \xrightarrow{t^{\omega}} \CC
  $$
  where
  $
  \Tr_\kE: \Hom(\kE, \kE[1]) \xrightarrow{\cong}
  \Hom(\kO, \kE^\vee \stackrel{\mathbb{L}}\otimes \kE[1]) \lar
  \Hom(\kO, \kO[1]) = H^1(\kO),
  $
  is a non-degenerate pairing. This pairing coincides with the composition
  $$
  \langle \,\, , \,  \rangle_{\kF, \, \kE}^{\omega}:
  \Hom(\kE, \kF) \otimes \Hom(\kF, \kE[1]) \xrightarrow{\cong}
  \Hom(\kF, \kE[1]) \otimes \Hom(\kE[1], \kF[1]) \xrightarrow{\circ}
  $$
  $$
  \lar \Hom(\kF, \kF[1])\xrightarrow{\Tr_\kF} H^1(\kO)\xrightarrow{t^{\omega}}\CC.
  $$
\end{proposition}

\begin{remark}
  The choice of non-degenerate pairings $\langle \,\, , \,  \rangle_{\kE, \, \kF}$
  is actually not unique, see the proof of Proposition I.2.3 in
  \cite{ReitenvandenBergh}.
  In particular, $\langle \,\, , \,  \rangle_{\kE, \, \kF}$ depends on the
  choice of a global section of the dualizing sheaf $\omega_E$.
  If $\omega_{\ast}:\Hom(\kF, \kE[1]) \lar \Hom(\kF, \kE\otimes\omega_{E}[1])$
  denotes the isomorphism induced by
  $\omega:\kO_{E}\xrightarrow{\sim}\omega_{E}$, we obtain
  $\langle \,\underline{\phantom{A}}\, ,
  \,\underline{\phantom{A}}\,  \rangle_{\kE, \,  \kF}^{\omega} =
  \langle \,\underline{\phantom{A}}\, ,
  \omega_{\ast}(\underline{\phantom{A}})  \rangle_{\kE, \,  \kF}$.
\end{remark}

In \cite{PolishchukMP} Polishchuk showed how
a construction of Merkulov \cite{Merkulov}, which fairly explicitly
provides an $A_{\infty}$-structure, applied to the category of
Hermitian vector bundles on a complex manifold with a hermitian metric, gives
an $A_{\infty}$-structure on the bounded derived category of this manifold.
He shows that this construction provides a
cyclic $A_{\infty}$-structure \cite[Theorem 1.1]{PolishchukMP}, which is
crucial for the proof of the Yang-Baxter equations in the situation considered
here. Therefore, we formulate his result explicitly in the following
proposition.

\begin{proposition}
  If $E$ is a smooth elliptic curve, there exits an   $A_\infty$-structure on
  the category $\Dbcoh(E)$, which is cyclic with respect to the pairing
  described in Proposition \ref{P:SD}. In particular, this means
  $$
  \bigl\langle m_3(f_1, g_1, f_2), g_2\bigr\rangle =
  - \bigl\langle f_1,  m_3(g_1, f_2, g_2) \bigr\rangle =
  - \bigl\langle m_3(f_2, g_2, f_1), g_1\bigr\rangle
  $$
  for any objects $\kE_1, \kE_2$ and $\kF_1, \kF_2 \in \Dbcoh(E)$,  and
  any morphisms  $$f_i \in \Hom(\kE_i, \kF_i) \text{ and } g_i \in
  \Hom(\kF_i, \kE_{3-i}[1]), i = 1,2.
  $$
\end{proposition}

\noindent
Because the proof of Polishchuk uses harmonic forms and Hodge theory, it
heavily depends on the smoothness assumption for $E$.
  If $E$ is singular, the same result can be derived with some effort using
  more recent methods from non-commutative symplectic geometry
  \cite[Theorem 10.2.2]{KoSo}.

\medskip
Now we recall the main construction of \cite{Polishchuk1}.
Take a reduced  projective curve $E$ with trivial dualising sheaf and  fix the
following data:

\begin{itemize}
\item Two  vector bundles $\kV_1$ and $\kV_2$ on $E$ having  the same rank $n$ and such that
  $\Hom_E(\kV_1, \kV_2) = 0 = \Ext^1_E(\kV_1, \kV_2)$.
\item Two distinct  smooth points $y_1, y_2 \in \breve{E}$
  lying on the same irreducible component of $E$ and such that
      $\Hom_E\bigl(\kV_1(y_2), \kV_2(y_1)\bigr) = 0 =
      \Ext^1_E\bigl(\kV_1(y_2), \kV_2(y_1)\bigr)$.
\end{itemize}

\medskip
\begin{remark}
This ``orthogonality'' assumption on vector bundles $\kV_1$ and $\kV_2$  might
seem to be quite artificial. The natural example of such data is the
following.
Let  $(n,d) \in \mathbb{N} \times \mathbb{Z}$ be a pair of coprime integers,
$M_{E}^{(n,d)}$ the moduli space of stable vector bundles of rank $n$ and degree
$d$ on a Weierstra\ss{} curve $E$.  Let $\kP(n,d)$ be a universal family
on $E \times M_E^{(n,d)}$. For a point $v_{i} \in M_E^{(n,d)}$ denote by
$\kV_{i}$ the corresponding stable vector bundle $\kP(n,d)|_{E\times v_{i}}$
on the curve $E$.
Then  for any two distinct points $v_1, v_2 \in M_E^{(n,d)}$ we have
$\Hom_E(\kV_1, \kV_2) = 0 = \Ext^1_E(\kV_1, \kV_2)$.
\end{remark}

Actually, one can also consider a more general situation.
Namely, for any pair
$(n,d) \in \mathbb{N} \times \mathbb{Z}$, not necessarily coprime, one can
take indecomposable semi-stable vector bundles of rank $n$ and degree $d$
having \emph{locally free}  Jordan-H\"older factors. The orthogonality
condition between non-isomorphic bundles of this type follows from the
following lemma.

\begin{lemma}[see \cite{BK3}]
Let $(n,d) \in \mathbb{N} \times \mathbb{Z}$, $m = \gcd(n,d)$ and $n =  mn', d = md'$.
Let $\kV$ be an indecomposable semi-stable vector bundle of rank
$n$ and degree $d$ on a Weierstra\ss{} curve $E$ with locally free
Jordan-H\"older factors.
Then all these factors  are isomorphic to a single stable vector bundle
$\kV' \in M_E^{(n', d')}$. Moreover,  we have
$
\kV \cong \kV' \otimes \kA_m,
$
where $\kA_m$ is the
 indecomposable vector bundle of rank $m$ and degree $0$ defined
recursively by the non-split extension sequences
$$
 0 \lar \kA_m \lar \kA_{m+1} \lar \kO \lar 0 \quad m \ge 1,
$$
where $\kA_1 = \kO$.
\end{lemma}

\medskip
\noindent
Let us return to Polishchuk's construction.
Since $\Hom_E(\kV_1, \kV_2) = 0 = \Ext_E^1(\kV_1, \kV_2)$, we have a linear map
$$
m_3: \Hom_E(\kV_1, \CC_{y_1}) \otimes \Ext_E^1(\CC_{y_1}, \kV_2)
\otimes \Hom_E(\kV_2, \CC_{y_2}) \lar
\Hom_E(\kV_1, \CC_{y_2})
$$
called the \emph{triple Massey product}
 and defined as follows. Let $a \in \Ext_E^1(\CC_{y_1}, \kV_2)$,
$g \in \Hom_E(\kV_1, \CC_{y_1})$, $h \in \Hom_E(\kV_2, \CC_{y_2})$ and let
$$
0 \lar \kV_2 \stackrel{\alpha}\lar \kA \stackrel{\beta}\lar \CC_{y_1} \lar 0
$$
be an exact sequence representing  the element $a$.
The vanishing of $\Hom_E(\kV_1, \kV_2)$ and $\Ext_E^1(\kV_1, \kV_2)$ implies that
we can uniquely lift the morphisms $g$ and $h$ to morphisms
 $\tilde{g}: \kV_1 \lar \kA$ and $\tilde{h}: \kA \lar \CC_{y_2}$. So, we
obtained a diagram
$$
\xymatrix{
         &              &           & \kV_1 \ar[d]^g  \ar[dl]_{\tilde g}    &  \\
a: 0 \ar[r] & \kV_2 \ar[r]^\alpha \ar[d]_h & \kA \ar[r]^\beta \ar[dl]^{\tilde h}
& \CC_{y_1} \ar[r] &  0 \\
         & \CC_{y_2}    &           &                  & \\
}
$$
and the triple Massey product is defined as  $m_3(g\otimes a \otimes h) =
\tilde{h}\tilde{g}$.
Note that $a$ determines
an extension only up to an automorphism of the middle term,
but the action of $\Aut(\kA)$ leads to the same answer  for
$m_3(g \otimes a \otimes h)  =: m^{\kV_1, \kV_2}_{y_1, y_2}(g \otimes a \otimes h).$

\medskip

Now one can  use a sequence of canonical isomorphisms in order
to rewrite $m^{\kV_1, \kV_2}_{y_1, y_2}$
in another form:
$$
\Lin\bigl(\Hom_E(\kV_1, \CC_{y_1}) \otimes \Ext_E^1(\CC_{y_1}, \kV_2)
\otimes \Hom_E(\kV_2, \CC_{y_2}),
\Hom_E(\kV_1, \CC_{y_2})\bigr)
\cong
$$
$$
\Lin\bigl(\Hom_E(\kV_1, \CC_{y_1}) \otimes \Hom_E(\kV_2, \CC_{y_2}),
\Ext_E^1(\CC_{y_1}, \kV_2)^* \otimes
\Hom_E(\kV_1, \CC_{y_2})\bigr)
\cong
$$
$$
\Lin\bigl(\Hom_E(\kV_1, \CC_{y_1}) \otimes \Hom_E(\kV_2, \CC_{y_2}),
\Hom_E(\kV_2, \CC_{y_1}) \otimes
\Hom_E(\kV_1, \CC_{y_2})\bigr),
$$
where we use  the Serre duality formula
$\Ext_E^1(\CC_{y_1}, \kV_2)^* \cong \Hom_E(\kV_2, \CC_{y_1})$
given by the bilinear form $\langle \,\,,\, \rangle_{\kV_2, \, \CC_{y_1}}^{\omega}$ from
Proposition \ref{P:SD}. Let $\widetilde{m}^{\kV_1, \kV_2}_{y_1, y_2}$ be the
image of $m^{\kV_1, \kV_2}_{y_1, y_2}$ under this chain of isomorphisms.

\begin{theorem}[see Theorem 1 in \cite{Polishchuk1}]
If $E$ is smooth, then
$\widetilde{m}^{\kV_1, \kV_2}_{y_1, y_2}$ satisfies
the following ``triangle equation'' (associative Yang--Baxter equation)
\begin{equation}\label{E:yb1}
(\widetilde{m}^{\kV_3, \kV_2}_{y_1, y_2})^{12} (\widetilde{m}^{\kV_1, \kV_3}_{y_1, y_3})^{13} -
(\widetilde{m}^{\kV_1, \kV_3}_{y_2, y_3})^{23} (\widetilde{m}^{\kV_1, \kV_2}_{y_1, y_2})^{12} +
(\widetilde{m}^{\kV_1, \kV_2}_{y_1, y_3})^{13} (\widetilde{m}^{\kV_2, \kV_3}_{y_2, y_3})^{23} = 0.
\end{equation}
The left-hand side of this equation is a linear map
$$
\Hom_E(\kV_1, \CC_{y_1}) \otimes \Hom_E(\kV_2, \CC_{y_2})  \otimes
\Hom_E(\kV_3, \CC_{y_3}) \lar
$$
$$
\lar \Hom_E(\kV_2, \CC_{y_1}) \otimes \Hom_E(\kV_3, \CC_{y_2})  \otimes
\Hom_E(\kV_1, \CC_{y_3}).
$$
Moreover, the tensor $\widetilde{m}^{\kV_1, \kV_2}_{y_1, y_2}$ is
non-degenerate and skew-symmetric:
$$
\tau(\widetilde{m}^{\kV_1, \kV_2}_{y_1, y_2}) = - \widetilde{m}^{\kV_2, \kV_1}_{y_2, y_1},
$$
where $\tau$ is the isomorphism $$\Hom_E(\kV_1, \CC_{y_1}) \otimes
\Hom_E(\kV_2, \CC_{y_2}) \lar
\Hom_E(\kV_2, \CC_{y_2}) \otimes \Hom_E(\kV_1, \CC_{y_1})$$
 given by
$\tau(f\otimes g) = g \otimes f$.
\end{theorem}

\noindent
\emph{Idea of the proof}.
This equality is a consequence of the $A_{\infty}$\,-\,constraint
$$
m_3 \circ (m_3 \otimes \mathbbm{1} \otimes \mathbbm{1} +
\mathbbm{1} \otimes m_3 \otimes \mathbbm{1} + \mathbbm{1} \otimes
\mathbbm{1} \otimes m_3) = 0,
$$
and  skew-symmetry of $\widetilde{m}^{\kV_1, \kV_2}_{y_1, y_2}$ follows from the
cyclicity of the $A_\infty$\,-\,structure.
\qed

\medskip
Note that for a vector bundle $\kV$ and a smooth point $y \in E$
we have canonical isomorphisms
$$
\Hom_E(\kV, \CC_y) \cong \Hom_{E}(\kV \otimes \CC_y, \CC_y)  =
\Hom_{\CC}(\kV|_y, \CC) = \kV|_y^{*}.
$$
In these  terms,  $\widetilde{m}^{\kV_1, \kV_2}_{y_1, y_2}$ is a linear map
$$
\widetilde{m}^{\kV_1, \kV_2}_{y_1, y_2}:
\kV_1|_{y_1}^{*} \otimes \kV_2|_{y_2}^{*}
\lar
\kV_2|_{y_1}^{*} \otimes \kV_1|_{y_2}^{*}.
$$
Now we use the canonical isomorphism
$$
\alpha:
\Hom_\CC(\kV_2|_{y_1}, \kV_1|_{y_1}) \otimes \Hom_\CC(\kV_1|_{y_2}, \kV_2|_{y_2})
\lar
\Hom_\CC(\kV_1|_{y_1}^{*} \otimes \kV_2|_{y_2}^{*},
        \kV_2|_{y_1}^{*} \otimes \kV_1|_{y_2}^{*})
$$
mapping a simple tensor $f_1 \otimes f_2$ to $f_1^t \otimes f_2^t$.
Then the tensor
$$
r^{\kV_1, \kV_2}_{y_1, y_2} := \alpha^{-1}(\widetilde{m}^{\kV_1, \kV_2}_{y_1, y_2}) \in
\Hom_\CC(\kV_2|_{y_1}, \kV_1|_{y_1}) \otimes
\Hom_\CC(\kV_1|_{y_2}, \kV_2|_{y_2})
$$
satisfies the equation
\begin{equation}\label{E:yb2}
(r^{\kV_1, \kV_3}_{y_1, y_3})^{13} (r^{\kV_3, \kV_2}_{y_1, y_2})^{12} -
(r^{\kV_1, \kV_2}_{y_1, y_2})^{12} (r^{\kV_1, \kV_3}_{y_2, y_3})^{23} +
(r^{\kV_2, \kV_3}_{y_2, y_3})^{23} (r^{\kV_1, \kV_2}_{y_1, y_3})^{13}  = 0
\end{equation}
and the unitarity condition
\begin{equation}\label{E:ybunit}
\tau(r^{\kV_1, \kV_2}_{y_1, y_2}) = - r^{\kV_2, \kV_1}_{y_2, y_1}.
\end{equation}

\begin{remark}\label{R:formsofAYBE}
Since the functorial isomorphism of vector spaces
 $\Hom_\CC(U, V) \to \Hom_\CC(V^*, U^*)$ is contravariant, the
 tensors  $r^{\kV_1, \kV_2}_{y_1, y_2}$ and
$\widetilde{m}^{\kV_1, \kV_2}_{y_1, y_2}$ appear in inverse order
in Equations (\ref{E:yb1}) and (\ref{E:yb2}).
\end{remark}

\medskip

\noindent
Note that the bilinear map
$$
\tr: \Hom_\CC(U, V) \times \Hom_\CC(V, U) \lar \mathbb{C}, \quad (f, g)
\mapsto \tr(f\circ g)
$$
is non-degenerate and induces an isomorphism
$
\Hom_\CC(U, V)^* \cong \Hom_\CC(V, U).
$
Using this,  we get a chain of canonical isomorphisms
$$
\Hom_\CC(\kV_2|_{y_1}, \kV_1|_{y_1}) \otimes \Hom_\CC(\kV_1|_{y_2}, \kV_2|_{y_2}) \cong
\Hom_\CC(\kV_1|_{y_1}, \kV_2|_{y_1})^* \otimes \Hom_\CC(\kV_1|_{y_2}, \kV_2|_{y_2}) \cong
$$
$$
\cong
\Lin\bigl(\Hom_\CC(\kV_1|_{y_1}, \kV_2|_{y_1}),
          \Hom_\CC(\kV_1|_{y_2}, \kV_2|_{y_2})\bigr).
$$
We let  $\tilde{r}^{\kV_1, \kV_2}_{y_1, y_2} \in
\Lin\bigl(\Hom_\CC(\kV_1|_{y_1}, \kV_2|_{y_1}),
 \Hom_\CC(\kV_1|_{y_2}, \kV_2|_{y_2})\bigr)$ be the image of $r^{\kV_1, \kV_2}_{y_1, y_2}$.

\medskip
\begin{remark}
Note that the triple Massey product $m_3$ is canonical, however the
tensor $r^{\kV_1, \kV_2}_{y_1, y_2}$ and the linear map $\tilde{r}^{\kV_1, \kV_2}_{y_1, y_2}$
\emph{depend} on the choice
of a global section $\omega \in H^0(\Omega_E)$.  Indeed, passing
from $m_3$ to $r^{\kV_1, \kV_2}_{y_1, y_2}$ and  $\tilde{r}^{\kV_1, \kV_2}_{y_1, y_2}$
we use the bilinear form $\langle \,\,,\, \rangle_{\kV_2, \, \CC_{y_1}}^{\omega}$, which
depends on the choice of $\omega$.
\end{remark}

\vspace{0.2cm}
\noindent
Our next aim is to answer the following questions:
\begin{itemize}
\item[\textbf{Q1}] What is a geometrical interpretation of the equivalence
  relation given in Definition \ref{D:gaugeAYBE}?
\item[\textbf{Q2}] How can we view $\kV_1$ and $\kV_2$ as variables?
\item[\textbf{Q3}] Is there  a practical way to compute $r^{\kV_1, \kV_2}_{y_1, y_2}$?
\end{itemize}

\section{Geometric description of triple Massey products}\label{S:Massey}

Let $E$ be a reduced projective curve with trivial dualizing sheaf and  $\breve{E}$ be the open subset
of the smooth points of $E$.
As in the previous section, we fix the following data:
\begin{itemize}
\item two vector bundles $\kV_1$ and $\kV_2$ on $E$ of rank $n$ such that
      $\Hom_E(\kV_1, \kV_2) = 0$ and $\Ext_E^1(\kV_1, \kV_2) = 0$.
\item two distinct   points $y_1, y_2 \in \breve{E}$ lying
      on the same irreducible component of $E$ such that
      $\Hom_E\bigl(\kV_1(y_2), \kV_2(y_1)\bigr)$ $ =
      \Ext_E^1\bigl(\kV_1(y_2), \kV_2(y_1)\bigr) = 0$.
\item  A non-zero global regular differential one-form $\omega \in
  H^0(\Omega_E)$.
\end{itemize}

\noindent
The main goal of this section is to get
an alternative description of the linear map
$$
\tilde{r}^{\kV_1, \kV_2}_{y_1, y_2}(\omega):
\Hom_\CC(\kV_1|_{y_1}, \kV_2|_{y_1}) \lar
 \Hom_\CC(\kV_1|_{y_2}, \kV_2|_{y_2})
$$
introduced in Section \ref{S:PolConstr}. To do this, we  first elaborate the
theory of residues and evaluation morphisms on reduced complex Gorenstein
curves.

\subsection{Residue map for vector bundles}\label{SS:resmap}
In this subsection,  our set-up is rather general. We fix the following data:

\begin{itemize}
\item a reduced Gorenstein analytic curve $X$ (not necessarily compact);
\item a subset $D\subset\breve{X}$, where $\breve{X}$ denotes the
  open subset of smooth points of $X$, such that $D\subset X$ is locally
  finite\footnote{%
  This means that for any $p\in X$, there exists an open neighbourhood of $p$
  in $X$ which contains only a finite number of points from $D$. According to
  \cite[Ch.~1.1]{GH}, $D$ defines a divisor and a corresponding line bundle
  $\kO_{X}(D)$ on $X$.};
\item a point $x \in D$; the set $D'=D\setminus\{x\}$ may be empty;
\item a pair of vector bundles $\kV$ and $\kW$ on $X$;
\item a germ of a differential 1-form  $\omega \in \Omega_x$, not vanishing at
  $x$ (i.e.~$\omega\not\in\Omega(-x)_{x}$), where $\Omega = \Omega_X$
  is the sheaf of regular differential 1-forms  on $X$.
\end{itemize}
The only application of this general set-up with $D'\ne\emptyset$ occurs in
Section \ref{S:elliptrm}, where $Y=\CC\xrightarrow{\pi}\CC/\Lambda=X$
(cf.~Prop.~\ref{P:functofres}, Prop.~\ref{L:functofev} and
Thm.~\ref{C:Masseyprod}) is the universal cover of an elliptic curve and
$D=y+\Lambda$ is an infinite $\Lambda$-invariant subset of $\CC$. Therefore,
we resist the temptation to include $D$ or $D'$ into the decorations of the
residue maps below.
In all other applictions of Theorem \ref{C:Masseyprod}, $\pi:Y\to X$ is the
normalisation, hence restricts to an isomorphism $\pi:\breve{Y}\to \breve{X}$
and $D$ consists of a single point.

Consider the canonical short exact sequence
\begin{equation}\label{E:residue}
0 \to \Omega_X \to \Omega_X(x) \stackrel{\underline{\res}_x}\lar \CC_x \to 0.
\end{equation}

\begin{proposition}\label{P:propofres}
The following diagram is commutative:
$$
\xymatrix
{
{\mathcal Hom}_X\bigl(\kV \otimes \Omega, \kW \otimes \Omega(D)\bigr)
\ar[rr]^-{(1 \otimes \underline{\res}_x)_*} & &
{\mathcal Hom}_X\bigl(\kV \otimes \Omega,\kW\otimes\kO(D')\otimes\CC_x\bigr)
\\
{\mathcal Hom}_X\bigl(\kV, \kW \otimes \kO(D)\bigr)
\ar[u]^{\,-\, \otimes \Omega} & &
{\mathcal Hom}_X\bigl(\kV \otimes \Omega,\kW\otimes\CC_x\bigr)
\ar[u]_{\left(r_{D'}\right)_{\ast}}
\\
{\mathcal Hom}_X(\kV, \kW) \otimes \kO(D)
\ar[rr]^-{\underline{\res}'_x} \ar[u]^{\can}  & &
{\mathcal Hom}_X\bigl(\kV \otimes \Omega \otimes \CC_x,\kW\otimes\CC_x\bigr).
\ar[u]_\can
}
$$
The map $\left(r_{D'}\right)_{\ast}$ is induced by the inclusion
$r_{D'}:\kO\to\kO(D')$ of subsheaves of the sheaf of meromorphic functions on
$X$ and $\underline{\res}'_x$ is determined by the following morphism of
presheaves.
Let $x\in U \subseteq X$ be an open subset, $s \in
\Hom_{\kO(U)}\bigl(\kV(U), \kW(U)\bigr)$, $f \in \Gamma\bigl(U, \kO(D)\bigr),
v \in \kV(U), \delta \in \Omega(U)$. Then
$$
\underline{\res}'_x(s \otimes f)\bigl[v \otimes \delta\bigr] =
\res_x(f \delta) [s(v)],
$$
where  $[v \otimes \delta] := v \otimes \delta \otimes 1$. If $x\not\in U$,
$\underline{\res}'_x|_{U}$ is the zero map.
The upper horizontal morphism is induced by the short exact sequence
\eqref{E:residue} and all other morphisms are standard canonical
isomorphisms.
\end{proposition}

\begin{proof}
  After observing that $r_{D'}$ is the identity map over any open set $U$ which
  does not intersect $D'$, the proof  is an easy  diagram chase.
\end{proof}

\noindent
Note that a  germ  $\omega \in \Omega_x$, which is not in $\Omega(-x)_{x}$,
i.e.~not equal to zero in $\Omega\otimes\CC_{x}$, induces an isomorphism of
sheaves
$$
{\mathcal Hom}_X\bigl(\kV\otimes\Omega\otimes\CC_x,\kW\otimes\CC_x\bigr)
\stackrel{\omega^{*}}\lar
{\mathcal Hom}_X\bigl(\kV\otimes\CC_x,\kW\otimes\CC_x\bigr).
$$

\begin{definition}
Consider the morphism of sheaves of $\kO_X$--modules
$$
\underline{\res}_x^{\kV, \kW}(\omega):
{\mathcal Hom}_X\bigl(\kV, \kW \otimes \kO(D)\bigr) \lar
{\mathcal Hom}_X\bigl(\kV  \otimes \CC_x, \kW \otimes \CC_x\bigr)
$$
defined as the composition of morphisms
$$
{\mathcal Hom}_X\bigl(\kV, \kW \otimes \kO(D)\bigr)
\xrightarrow{-\otimes \Omega}
{\mathcal Hom}_X\bigl(\kV \otimes \Omega, \kW \otimes \Omega(D)\bigr)
$$
$$
\xrightarrow{(1 \otimes \underline{\res}_x)_{*}}
{\mathcal Hom}_X\bigl(\kV \otimes \Omega,\kW\otimes\kO(D')\otimes\CC_x\bigr)
$$
$$
\xrightarrow{\left(\left(r_{D'}\right)_{\ast}\circ\can\right)^{-1}}
{\mathcal Hom}_X\bigl(\kV \otimes \Omega \otimes \CC_x, \kW \otimes \CC_x\bigr)
\stackrel{\omega^{*}}\lar
{\mathcal Hom}_X\bigl(\kV  \otimes \CC_x, \kW \otimes \CC_x\bigr),
$$
coinciding with the composition
$
{\mathcal Hom}_X\bigl(\kV, \kW \otimes \kO(D)\bigr)
\xrightarrow{\can}
{\mathcal Hom}_X(\kV, \kW) \otimes \kO(D)$
$
\xrightarrow{\underline{\res}'_x}
{\mathcal Hom}_X\bigl(\kV \otimes \Omega \otimes \CC_x, \kW \otimes \CC_x\bigr)
\stackrel{\omega^{*}}\lar
{\mathcal Hom}_X\bigl(\kV  \otimes \CC_x, \kW \otimes \CC_x\bigr).
$
\end{definition}

\begin{remark}\label{rem:thirdres}
  The fact that the composition
  $$\mathcal{H}om(\kV,\kW(D)) \xrightarrow{-\otimes \Omega}
  \mathcal{H}om(\kV\otimes \Omega, \kW\otimes \Omega(D))
  \xrightarrow{\omega^{*}} \mathcal{H}om(\kV, \kW\otimes \Omega(D))$$
  coincides with the map $\omega_{*}$ implies that there is a third
  description for $\underline{\res}_x^{\kV, \kW}(\omega)$ as:
  $\mathcal{H}om(\kV,\kW(D)) \xrightarrow{\omega_{*}}
  \mathcal{H}om(\kV, \kW\otimes \Omega(D))
  \xrightarrow{(1 \otimes \underline{\res}_x)_{*}}
  \mathcal{H}om(\kV, \kW\otimes \CC_{x})
  \xrightarrow{\left(\left(r_{D'}\right)_{\ast}\circ\can\right)^{-1}}
  \mathcal{H}om(\kV\otimes \CC_{x}, \kW\otimes \CC_{x})$.
  This will be used in the proof of Proposition \ref{P:comdiag}.
\end{remark}

\begin{definition}\label{def:res-absolut}
In the above notation we define:
$$
\res^{\kV, \kW}_x(\omega) := H^0\bigl(\underline{\res}_x^{\kV, \kW}(\omega)\bigr):
\Hom_X\bigl(\kV, \kW(D)\bigr) \lar  \Hom_X(\kV \otimes \CC_x, \kW \otimes \CC_x).
$$
\end{definition}

\noindent
Proposition \ref{P:propofres} yields an explicit
 formula for  the morphisms
$\res_x^{\kV, \kW}(\omega)$ in the following simple situation.

\begin{lemma}\label{L:resexplit}
Let $U\subseteq \CC$ be open, $O = H^0(U, \kO_\CC)$, $D\subset U$ a
locally finite subset, $O(D) = H^0\bigl(U, \kO_\CC(D)\bigr)$ and $x\in D$.
Let $\kV = \kO_U^n, \kW = \kO_U^m$ and   $\omega = f(z)dz$ be
a  meromorphic  one-form on $U$, holomorphic at $x$ with $f(x)\ne0$.
Then the morphism $\res_x$, which is defined as the  composition of canonical
morphisms
$$
\xymatrix
{
\Hom_U\bigl(\kV, \kW(D)\bigr) \ar[rrr]^-{\res_x^{\kV, \kW}(\omega)} & & &
\Hom_U(\kV \otimes \CC_x, \kW \otimes \CC_x) \ar[d]^\can \\
\Mat_{m \times n}\bigl(O(D)\bigr) \ar[rrr]^-{\res_x} \ar[u]^\can & & &
\Mat_{m \times n}(\CC)
}
$$
is given by the formula
$
\res_x(F) = \res_x(F \cdot \omega).
$
\end{lemma}

\begin{proof}
Let $F(z) = \bigl(f_{ij}(z)\bigr) \in \Mat_{m \times n}\bigl(O(D)\bigr)$ be a matrix
whose entries are meromorphic functions on $U$ having at most simple poles
along $D$.
Then we can write $F(z) = \frac{\displaystyle G(z)}{\displaystyle z-x}$, where
the entries of the matrix $G(z)$ are meromorphic functions which are
holomorphic at $x$. Using the definition of $\res_x$ in terms
of $\underline{\res}'_x$ we obtain
$$
\res_x(F)(a) = \left(\res_x \frac{\omega}{z-x}\right) G(x)a =
\res_x(F\omega)a
$$
for all $a \in \CC^n$.
\end{proof}

\noindent
In general, the morphism $\res^{\kV, \kW}_x(\omega)$ is neither surjective nor
injective.
However, there is an important  special case where it is an isomorphism.

\begin{proposition}\label{P:isores}
Let $E$ be a reduced projective curve with trivial dualizing sheaf. Let $\kV$
and $\kW$ be a pair of vector bundles on $E$ such that
$\Hom_E(\kV, \kW) = 0$ and $\Hom_E(\kW, \kV) = 0$
and let $x\in \breve{E}$, $\omega \in \Omega_{x}$ be as above.
Then, the morphism
$$
\res_x^{\kV, \kW}(\omega):
\Hom_E\bigl(\kV, \kW(x)\bigr) \lar \Hom_E(\kV\otimes\CC_x,\kW\otimes\CC_x),
$$
defined with $D=\{x\}$, is an isomorphism.
\end{proposition}

\begin{proof}
First note that by Serre  duality we have:
$
\Ext^1_E(\kV, \kW) \cong \Hom_E(\kW, \kV)^* = 0.
$
Hence, then short exact sequence
$
0 \to \kW \otimes \Omega \to  \kW \otimes \Omega(x) \to \kW \otimes \CC_x \to 0
$
induces an isomorphism
$
H^{0}\bigl((1\otimes\underline{\res}_{x})_{*}\bigr):
\Hom_E\bigl(\kV \otimes \Omega, \kW \otimes \Omega(x)\bigr) \stackrel{\cong}\lar
\Hom_E(\kV \otimes \Omega, \kW \otimes \CC_x).
$
As a result, the morphism
$
\Hom_E\bigl(\kV, \kW(x)\bigr) \xrightarrow{\res_x^{\kV, \kW}(\omega)}
\Hom_E(\kV \otimes \CC_x, \kW \otimes \CC_x)
$
is an isomorphism, too.
\end{proof}

\begin{remark}
  The vanishing conditions of Proposition \ref{P:isores} are satisfied if $E$
  is an irreducible projective Weierstra\ss{} curve and $\kV$ and $\kW$ are
  two stable vector bundles of the same rank and degree and such that $\kV
  \not\cong \kW$.
\end{remark}

\noindent
The next goal is to show that the morphism
$\res^{\kV, \kW}_x(\omega)$ has nice functorial properties.

\begin{proposition}\label{P:functofres}
  Let $Y \stackrel{\pi}\lar X$ be a morphism of reduced Gorenstein  curves,
  $y \in \breve{Y}$ and $x = \pi(y) \in \breve{X}$ be such that
  $f$ is unramified over $x$.  Let $D=\pi^{-1}(x)$ and $\kV, \kW$ be a pair of
  vector bundles on $X$, $\widetilde\kV = \pi^*\kV$ and
  $\widetilde\kW = \pi^*\kW$.
  Finally, let $\omega \in \Omega_{X, x}$ be the germ of  a  regular
  differential one-form, not vanishing at $x$, and $\tilde\omega =
  \pi^*(\omega) \in \Omega_{Y, y}$ be the corresponding germ on $Y$.
  Then the following diagram is commutative:
  $$
  \xymatrix
  {
    \Hom_X\bigl(\kV, \kW(x)\bigr)
    \ar[rrr]^-{\res_x^{\kV, \kW}(\omega)} \ar[d]_{\pi^*}
    & & &
    \Hom_X\bigl(\kV \otimes \CC_x, \kW \otimes \CC_x\bigr) \ar[d]^{\pi^*}
    \\
    \Hom_Y\bigl(\widetilde\kV, \widetilde\kW(D)\bigr)
    \ar[rrr]^-{\res_y^{\widetilde\kV, \widetilde\kW}(\tilde\omega)}
    & & &
    \Hom_Y\bigl(\widetilde\kV \otimes \CC_y, \widetilde\kW \otimes \CC_y\bigr).
  }
  $$
\end{proposition}

\begin{proof}
Let $i: U \hookrightarrow X$ be an open embedding containing the point $x$,
$\kV' = \kV|_U$ and  $\kW' = \kW|_U$. Then the diagram
$$
\xymatrix
{
\Hom_X\bigl(\kV, \kW(x)\bigr)  \ar[rrr]^-{\res_x^{\kV, \kW}(\omega)} \ar[d]_{i^*}
& & &
\Hom_X\bigl(\kV \otimes \CC_x, \kW \otimes \CC_x\bigr) \ar[d]^{i^*}\\
\Hom_U\bigl(\kV', \kW'(x)\bigr)  \ar[rrr]^-{\res_x^{\kV', \kW'}(\omega)}
& & &
\Hom_U\bigl(\kV' \otimes \CC_x, \kW' \otimes \CC_x\bigr).
}
$$
is commutative: this is a consequence of the ``local'' definition of the
morphism
$\res^{\kV, \kW}_x(\omega)$ as $H^0\bigl(\underline{\res}_x^{\kV, \kW}(\omega)\bigr)$.

In order to pass to the general case, recall that any unramified
morphism of smooth Riemann surfaces is locally biholomorphic.
The assumptions imply $\pi^{\ast}\kO(x)=\kO(D)$.
Since both points $x \in X$ and $y \in Y$ are smooth, there
exist open neighbourhoods $x \in U \stackrel{i}\hookrightarrow X$ and
$y \in V \stackrel{j}\hookrightarrow Y$ such that $\pi: V \to U$ is an
isomorphism. In particular, $V\cap D=\{y\}$ and we have a diagram
$$
\xymatrix
{
\Hom_X\bigl(\kV, \kW(x)\bigr)  \ar[rrr]^-{\res_x^{\kV, \kW}(\omega)} \ar[d]_{i^*} \ar@/_5pc/[ddd]_{\pi^*}
& & &
\Hom_X\bigl(\kV \otimes \CC_x, \kW \otimes \CC_x\bigr) \ar[d]^{i^*} \ar@/^7pc/[ddd]^{\pi^*}\\
\Hom_U\bigl(\kV', \kW'(x)\bigr) \ar[d]_{\pi^*} \ar[rrr]^-{\res_x^{\kV', \kW'}(\omega)}
& & &
\Hom_U\bigl(\kV' \otimes \CC_x, \kW' \otimes \CC_x\bigr) \ar[d]^{\pi^*} \\
\Hom_V\bigl(\widetilde\kV', \widetilde\kW'(y)\bigr)  \ar[rrr]^-{\res_y^{\widetilde\kV',
\widetilde\kW'}(\tilde\omega)}
& & &
\Hom_V\bigl(\widetilde\kV' \otimes \CC_y, \widetilde\kW' \otimes \CC_y\bigr) \\
\Hom_Y\bigl(\widetilde\kV, \widetilde\kW(D)\bigr)  \ar[rrr]^-{\res_y^{\widetilde\kV,
\widetilde\kW}(\tilde\omega)}  \ar[u]^{j^*}
& & &
\Hom_Y\bigl(\widetilde\kV \otimes \CC_y, \widetilde\kW \otimes \CC_y\bigr) \ar[u]_{j^*},
}
$$
in which  all three cental  squares and both  exterior squares are commutative.
To conclude the claim, it remains to note that all vertical
morphisms on the right hand side are isomorphisms.
\end{proof}

Note that the constructed morphism $\res^{\kV, \kW}_x(\omega)$ is
functorial in $\kV$ and $\kW$. In what follows, we shall use the following
bifunctoriality.

\begin{definition}
  Let $\kV_1, \kV_2$ and $\kW_1, \kW_2$ be coherent sheaves on $X$,
  $f: \kV_1 \to \kV_2$ and $g: \kW_1 \to \kW_2$ be isomorphisms in
  $\Coh(X)$. Then we have an morphism
  $$
  \conj(f, g): \Hom_X(\kV_1, \kW_1) \lar  \Hom_X(\kV_2, \kW_2)
  $$
  given by the rule $\Hom_X(\kV_1, \kW_1) \ni h  \mapsto \hspace{0.1cm}  g
  \circ h \circ f^{-1} \in \Hom_X(\kV_2, \kW_2)$.
\end{definition}

\noindent
Having this notation, the proof of the following lemma is straightforward.

\begin{lemma}\label{L:conjofres}
  Let $X$ be a reduced Gorenstein curve, $D\subset\breve{X}$, $x \in D$  and
  $\omega \in \Omega_x$ as before.
  Let $f: \kV_1 \to \kV_2$ and $g: \kW_1 \to \kW_2$ be isomorphisms of vector
  bundles on $X$.
  Then the following diagram, in which $g(D)$ denotes $g\otimes 1$, is
  commutative:
  $$
  \xymatrix
  {
    \Hom_X\bigl(\kV_1, \kW_1(D)\bigr)
    \ar[rrr]^-{\res_x^{\kV_1, \kW_1}(\omega)}
    \ar[d]_{\conj\left(f, \, g(D)\right)}
    & & &
    \Hom_X\bigl(\kV_1 \otimes \CC_x, \kW_1 \otimes \CC_x\bigr)
    \ar[d]^{\conj(\bar{f}, \, \bar{g})}
    \\
    \Hom_X\bigl(\kV_2, \kW_2(D)\bigr)
    \ar[rrr]^-{\res_x^{\kV_2, \kW_2}(\omega)}
    & & &
    \Hom_X\bigl(\kV_2 \otimes \CC_x, \kW_2 \otimes \CC_x\bigr).
  }
  $$
\end{lemma}

\medskip

\subsection{Evaluation map  for vector bundles}\label{SS:evmap}
As in the previous subsection, we fix the following notation:

\begin{itemize}
\item a reduced Gorenstein analytic curve $X$ (not necessarily compact);
\item a subset $D\subset\breve{X}$, locally finite in $X$, and a smooth point
  $y \in \breve{X}$, $y\not\in D$;
\item a pair of vector bundles $\kV$ and $\kW$ on $X$.
\end{itemize}

\noindent
Consider the short exact sequence
\begin{equation}\label{E:evaluat}
0 \to \kO(-y) \to \kO \stackrel{\underline{\ev}_y}\lar \CC_y \to 0.
\end{equation}
It induces a short exact sequence of coherent sheaves
$$
0 \to \kW(D)\otimes\kO(-y) \to \kW(D)
\xrightarrow{1 \otimes \underline{\ev}_y}
\kW(D) \otimes \CC_y \to 0
$$
and a morphism of sheaves $\underline{\ev}^{\kV, \kW(D)}_y$ making the
following diagram commutative
$$
\xymatrix
{
  {\mathcal Hom}_X\bigl(\kV, \kW(D)\bigr)
  \ar[d]_{\underline{\ev}^{\kV, \kW(D)}_y}
  \ar[rr]^-{(1 \otimes \underline{\ev}_y)_*} & &
  {\mathcal Hom}_X\bigl(\kV, \kW(D) \otimes \CC_y\bigr)
  \\
  {\mathcal Hom}_X\bigl(\kV \otimes \CC_y, \kW \otimes \CC_y\bigr)
  \ar[rr]^-{(r_D)_*} & &
  {\mathcal Hom}_X\bigl(\kV \otimes \CC_y, \kW(D) \otimes \CC_y\bigr)
  \ar[u]_\cong
}
$$
where the lower horizontal morphism $(r_D)_*$ is induced by the canonical
inclusion $r_D: \kO \to \kO(D)$ (here we view both sheaves as subsheaves of
the sheaf of meromorphic functions on $X$).
Note that $(r_{D})_{*}$ is an isomorphism because $y\not\in D$ by assumption.

\begin{definition}\label{def:ev-absolut} In the notation as above, we set:
  $$
  \ev^{\kV, \kW(D)}_y := H^0\bigl(\underline{\ev}^{\kV, \kW(D)}_y\bigr):
  \Hom_X\bigl(\kV, \kW(D)\bigr) \lar \Hom_X(\kV\otimes \CC_y, \kW\otimes \CC_y).
  $$
\end{definition}

\noindent
Similar to the case of the residue map, the following statements are true.

\begin{proposition}\label{L:functofev}
  Let $\pi: Y \to X$ be a morphism of reduced Gorenstein curves,
  $x_{1}\in \breve{X}$ a point over which $\pi$ is unramified,
  $D_{1}=\pi^{-1}(x_{1})$ and $y_2 \in \breve{Y}$ a smooth point, such that
  $x_{1}\ne x_2 = \pi(y_2)\in \breve{X}$.
  For a pair of vector bundles $\kV$ and $\kW$ on $X$ we denote
  $\widetilde\kV = \pi^*\kV$ and $\widetilde\kW = \pi^*\kW$. Then the
  following diagram is commutative:
  $$
  \xymatrix
  {
    \Hom_X\bigl(\kV, \kW(x_1)\bigr)
    \ar[rrr]^-{\ev_{x_2}^{\kV, \kW(x_1)}} \ar[d]_{\pi^*}
    & & &
    \Hom_X\bigl(\kV \otimes \CC_{x_2}, \kW \otimes \CC_{x_2}\bigr)
    \ar[d]^{\pi^*}
    \\
    \Hom_Y\bigl(\widetilde\kV, \widetilde\kW(D_1)\bigr)
    \ar[rrr]^-{\ev_{y_2}^{\widetilde\kV, \widetilde\kW(D_1)}}
    & & &
    \Hom_Y\bigl(\widetilde\kV\otimes\CC_{y_2},\widetilde\kW\otimes\CC_{y_2}\bigr).
  }
  $$
  Moreover, the constructed morphism of vector spaces
  $\ev_{y_2}^{\widetilde\kV, \widetilde\kW(D_1)}$ is natural in
  $\widetilde\kV$ and $\widetilde\kW$. In particular, if
  $f: \widetilde\kV_1 \to \widetilde\kV_2$ and
  $g: \widetilde\kW_1 \to \widetilde\kW_2$ are isomorphisms of vector
  bundles on $Y$ then the following diagram is commutative:
  $$
  \xymatrix
  {
    \Hom_Y\bigl(\widetilde\kV_1, \widetilde\kW_1(D_1)\bigr)
    \ar[rrr]^-{\ev_{y_2}^{\widetilde\kV_1, \widetilde\kW_1(D_1)}}
    \ar[d]_{\conj\left(f, \,  g(D_{1})\right)}
    & & &
    \Hom_Y\bigl(\widetilde\kV_1\otimes\CC_{y_2},
    \widetilde\kW_1\otimes\CC_{y_2}\bigr)
    \ar[d]^{\conj(\bar{f}, \, \bar{g})}
    \\
    \Hom_Y\bigl(\widetilde\kV_2, \widetilde\kW_2(D_1)\bigr)
    \ar[rrr]^-{\ev_{y_2}^{\widetilde\kV_2, \widetilde\kW_2(D_1)}}
    & & &
    \Hom_Y\bigl(\widetilde\kV_2 \otimes \CC_{y_2},
    \widetilde\kW_2 \otimes \CC_{y_2}\bigr).
  }
  $$
\end{proposition}

\noindent
The proof of the following formula for the evaluation morphism $\ev_y$ is
straightforward.

\begin{lemma}\label{L:evexplit}
  Let $U\subseteq \CC$ be open, $D\subset U$
  locally finite and $y\in U\setminus D$.
  Define $O(D) = H^0\bigl(U, \kO_\CC(D)\bigr)$ and let
  $\kV = \kO_U^n, \kW = \kO_U^m$.  The formula $\ev_y(F) = F(y)$ defines
  a morphism
  $\ev_y: \Mat_{n \times n}\bigl(O(D)\bigr) \lar \Mat_{n\times n}(\CC)$
  which fits into the following commutative diagram
  $$
  \xymatrix
  {
    \Hom_U\bigl(\kV, \kW(D)\bigr)
    \ar[rrr]^-{\ev_y^{\kV, \kW(D)}} & & &
    \Hom_U(\kV \otimes \CC_y,\kW \otimes \CC_y) \ar[d]^\can
    \\
    \Mat_{m \times n}\bigl(O(D)\bigr)
    \ar[rrr]^-{\ev_y} \ar[u]^\can & & &
    \Mat_{m \times n}(\CC).
  }
  $$
\end{lemma}

\noindent
In general, the evaluation morphism $\ev_{x_2}^{\kV, \kW(x_1)}$ is neither
injective nor surjective. However, there is one particular case, when it is an
isomorphism.

\begin{proposition}\label{P:isoev}
  Let $E$ be a reduced projective curve with trivial dualizing sheaf.
  Let $\kV, \kW$ be vector bundles on $E$ and $x, y \in \breve{E}$ be such
  that $x \ne y$ and
  \[
  \Hom_E\bigl(\kV, \kW(x-y)\bigr) = 0 = \Ext^1_E\bigl(\kV, \kW(x-y)\bigr).
  \]
  Then the morphism
  $
  \ev^{\kV, \kW(x)}_y: \Hom_E\bigl(\kV, \kW(x)\bigr) \lar
  \Hom_E(\kV \otimes \CC_y, \kW\otimes \CC_y)
  $
  is an isomorphism of vector spaces.
\end{proposition}

\begin{proof}
Applying  the functor $\Hom_E(\kV, \,-\,)$ to the short exact sequence
$$
0 \to \kW(x-y)  \to \kW(x) \xrightarrow{1 \otimes \underline{\ev}_y}
 \kW(x) \otimes \CC_y \to 0
$$
we see that  the morphism defined as the composition
$$
\Hom_E\bigl(\kV, \kW(x)\bigr) \xrightarrow{(1 \otimes \underline{\ev}_y)_*}
 \Hom_E\bigl(\kV \otimes \CC_y, \kW(x) \otimes \CC_y\bigr)
\xrightarrow{(r_x)^{-1}_*} \Hom_E\bigl(\kV \otimes \CC_y, \kW \otimes \CC_y\bigr)
$$
coincides with $\ev^{\kV, \kW(x)}_y$ and is an isomorphism of vector spaces.
\end{proof}

\begin{remark}
The vanishing conditions  of Proposition \ref{P:isoev} are satisfied if $E$ is an irreducible
Weierstra\ss{} projective curve and $\kV$ and $\kW$ are two stable vector bundles
of the same rank and degree and such that $\kV \not\cong \kW(x-y)$.
\end{remark}

\medskip

\subsection{Geometric description of triple Massey products on genus one curves}
\label{subsec:geometricMassey}
In this subsection,  $E$ is  a reduced projective Gorenstein curve with
trivial dualizing sheaf. In particular, the sheaf $\Omega_E$ of regular
differential 1-forms on $E$ is trivial.
For any smooth point $x \in E$  consider the coboundary map  $\delta_x:
H^0(\CC_x) \to H^1(\Omega_E)$  of the short exact sequence (\ref{E:residue}).
This is an isomorphism.
Define $w_x := \delta_x(1_x) \in H^1(\Omega_E)$. By a result of Kunz, which is
also true without the assumption $\Omega_E \cong \kO_{E}$, we get:

\begin{theorem}[see  Satz 4.1 in \cite{Kunz}]\label{T:Kunz}
The element $w_x$ does not depend on $x$.
\end{theorem}

\noindent
Let $w = w_x$ and $t: H^1(\Omega_E) \to \CC$ be the isomorphism which maps
$w$ to $1$.  We fix a global regular differential form
$\omega: \kO_E  \to \Omega_E$. For any two perfect complexes
$\kE, \kF \in \Perf(E)$, $\omega$ induces a non-degenerate pairing (see
Proposition \ref{P:SD})
$$
\langle \,\, , \,  \rangle_{\kE, \, \kF}^{\omega}: \Hom_{D(E)}(\kE, \kF)
\otimes \Hom_{D(E)}(\kF, \kE[1]) \lar
 \mathbb{C}.
$$
Recall that, when passing from the triple Massey product
$m^{\kV_1, \kV_2}_{y_1, y_2}$ to the tensor  ${\widetilde m}^{\kV_1,  \kV_2}_{y_1, y_2}$,
we have already used these  bilinear forms.

\medskip
\noindent
The alternative description of $\tilde{r}^{\kV_1, \kV_2}_{y_1, y_2}$ involves
the two isomorphisms $\res^{\kV_1, \kV_2}_{y_1}(\omega)$ and
$\ev^{\kV_1, \kV_2(y_1)}_{y_2}$ constructed
in Subsections \ref{SS:resmap} and \ref{SS:evmap} respectively.
The following theorem (see also \cite[Theorem 4]{Polishchuk1})
is the key statement to explicitly compute the tensor
$\tilde{r}^{\kV_1, \kV_2}_{y_1, y_2}$ describing triple Massey products.

\begin{theorem}\label{T:GeomMassey} Let $E$ be a reduced projective
curve with trivial dualizing sheaf, $\kV_1, \kV_2 \in \VB(E)$,
$y_1, y_2 \in \breve{E}$ and $\omega \in H^0(\Omega_E)$ be as at the beginning
of Section \ref{S:Massey}. By $\omega_{y_{1}}$ we denote the germ of $\omega$
at $y_{1}\in E$. Then the diagram
$$
\xymatrix@C-C@R+3mm
{
& \Hom_E\bigl(\kV_1, \kV_2(y_1)\bigr)
\ar[ld]_{\res^{\kV_1, \kV_2}_{y_1}(\omega_{y_{1}})\phantom{x}}
\ar[rd]^{\ev^{\kV_1, \kV_2(y_{1})}_{y_2}} &
\\
\Hom_{\CC}(\kV_1|_{y_1}, \kV_2|_{y_1})
\ar[rr]^{\tilde{r}^{\kV_1, \kV_2}_{y_1, y_2}(\omega)} & &
\Hom_{\CC}(\kV_1|_{y_2}, \kV_2|_{y_2})
}
$$
is commutative.
\end{theorem}

\noindent
An important message from this theorem is: \emph{only} if $\omega_{y_{1}}$ is
the germ of a \emph{global} holomorphic 1-form $\omega\in H^{0}(\Omega_{E})$,
we can guarantee that
$\tilde{r}=\ev_{y_{2}} \circ \left(\res_{y_{1}}(\omega_{y_{1}})\right)^{-1}$.

Since this result plays a crucial role in our approach to degeneration
problems, we decided to give a detailed proof of this statement, stressing
those  points which are implicit in \cite{Polishchuk1}.  As a preparation,
several technical lemmas have to be proven.

\begin{lemma}
Let $E$ be a reduced projective curve with trivial dualizing sheaf and
$x \in \breve{E}$. Then we have  an isomorphism of functors
$\VB(E) \lar \Vect_{\CC}$:
$$
T_x:  \Ext^1_E(\CC_x, \,-\,) \lar \Hom_E(\CC_x, \,- \otimes \CC_x)
$$
\end{lemma}

\begin{proof}  Let $\kV$ be a vector bundle on $E$ of rank $n$.
The short exact sequence \eqref{E:residue}
and the isomorphism $\omega:\kO_{C}\xrightarrow{\sim}\Omega_{E}$
yield the short exact sequence
$
0 \rightarrow \kV \rightarrow\kV(x) \rightarrow \kV \otimes \CC_x \rightarrow 0,
$
which  induces the long exact sequence
$$
0 \rightarrow \Hom_E(\CC_x, \kV\otimes \CC_x) \stackrel{\delta_x}\lar
\Ext^1_E(\CC_x, \kV) \rightarrow \Ext^1_E(\CC_x, \kV(x)) \rightarrow
\Ext^1_E(\CC_x, \kV\otimes \CC_x) \rightarrow 0.
$$
Because $\Ext^1_E\bigl(\CC_x, \kV(x)\bigr) \cong H^0\bigl({\mathcal
  Ext}^1(\CC_x, \kV(x))\bigr)$ and $\Ext^1_E(\CC_x, \kV\otimes \CC_x)$ are
both of dimension $n=\rank(\kV)$,
we conclude that $\delta_x$ is an isomorphism. Moreover, this map is
functorial and we can put $T_x = \delta_x^{-1}$.
\end{proof}

\begin{remark}\label{R:connmap}
Due to the construction of the functor $T_x$  we have a commutative diagram
$$
\xymatrix
{
0 \ar[r] & \kV\ar[r]\ar@{=}[d]&\kA\ar[r]\ar[d] & \CC_x\ar[r]\ar[d]^{T_x(a)} & 0 \\
0 \ar[r] & \kV\ar[r] &\kV(x)\ar[r]^{\underline{\res}^\kV_x} &\kV\otimes\CC_x \ar[r]& 0,
}
$$
where the upper exact sequence corresponds to the element
$a \in \Ext^1_E(\CC_x, \kV)$.
\end{remark}

\medskip
In order to justify our calculations in Sections \ref{S:elliptrm} and
\ref{S:singrm}   we need to establish an explicit link between the
``categorical trace map'' of Proposition \ref{P:SD} and the usual trace from
linear algebra.

Let $X$ be a reduced projective Gorenstein curve,
$x \in \breve{X}$
a smooth point,  $\kV$ a vector bundle on $X$. From the exact sequence
(\ref{E:residue})
we get a commutative diagram
$$
\xymatrix
{\Ext^1_X(\kV, \kV \otimes \Omega_X) \ar[rr]^{\Tr^{1}_\kV} & &
H^1(\Omega_X) \ar[r]^t &  \CC \\
\Hom_X(\kV, \kV \otimes \CC_x) \ar[rr]^{\Tr^{0}_\kV} \ar[u]^{\delta_x}  & &
H^0(\CC_x) \ar[u]_{\delta_x} \ar[ur]_\can &  \\
\Hom_\CC(\kV|_x, \kV|_x) \ar[rr]^{\tr} \ar[u] & &
\CC \ar[u]_\can \ar@/_1pc/[uur]_{=}& \\
}
$$
Here  $t$ is the trace map from  Theorem \ref{T:Kunz} and $\tr$ is the ordinary
trace of an endomorphism of the vector space $\kV|_x$. The morphism
$\Tr^{0}_\kV$ is the composition $\Hom_X(\kV, \kV \otimes \CC_x) \rightarrow
H^0\bigl(\kV^\vee \otimes \kV \otimes \CC_x\bigr) \rightarrow H^0(\CC_x)$ and
 $\Tr^{1}_\kV$ is defined in a similar way, so that $t\circ \Tr^{1}_\kV =
 \tr_{\kV}$, see Theorem \ref{T:SerreDuality}.
The commutativity of this diagram gives us the following result.

\begin{lemma}\label{C:Traces}
For an element $f \in \Hom_X(\kV, \kV \otimes \CC_x)$ we have:
$$
t\Bigl(\Tr^{1}_{\kV}\bigl(\delta_x(f)\bigr)\Bigr) = \tr(f_x),
$$
which is the required  link between the categorical trace and the usual trace
for vector spaces.
\end{lemma}

\begin{lemma}\label{L:SvsT}
Let $E$ be a reduced projective curve with trivial dualizing sheaf,
$x \in E$ a smooth point,  $\kV \in \VB(E)$ a vector bundle  and
$ S: \Ext^1_E(\CC_x, \kV) \lar \Hom_E(\kV, \CC_x)^*$
 the isomorphism
induced by the bilinear form $\langle \,\, ,\,\rangle_{\kV,
  \,\CC_x}^{\omega}$, defined in Proposition \ref{P:SD}.
Then the following diagram is commutative:
$$
\xymatrix
{
\Ext^1_E(\CC_x, \kV) \ar[rrr]^S \ar[d]_{T_x}  & &  & \Hom_E(\kV, \CC_x)^* \ar[d] \\
\Hom_E(\CC_x, \kV\otimes\CC_x) \ar[rrr]^{\tr} &&&\Hom_E(\kV\otimes\CC_x,\CC_x)^{*},
}
$$
where  $\tr$ is induced by the canonical isomorphism of vector spaces
$\Hom_\CC(U, V)^* \cong \Hom_\CC(V, U)$, which is given by the usual trace of
endomorphisms.
\end{lemma}

\begin{proof} Let $\xi \in \Hom_E(\kV, \CC_x)$ and $a\in \Ext^1_E(\CC_x, \kV)$.
Then $T_x(a): \CC_x \to \kV \otimes \CC_x$ and
$\xi_x: \kV \otimes \CC_x \to \CC_x$ satisfy:
\begin{align*}
  \tr\bigl(\xi_x \circ T_x(a)\bigr) &=
\tr\bigl(T_x(a) \circ \xi_x\bigr) =
t\bigl(\Tr^{1}_\kV\bigl(\delta_x(T_x(a) \circ \xi)\bigr)\bigr) =
t\bigl(\Tr^{1}_\kV\bigl(\omega_{*}(a \circ \xi)\bigr)\bigr)\\
&= t\bigl(\omega_{*}\bigl(\Tr_\kV(a \circ \xi)\bigr)\bigr) = S(a)(\xi),
\end{align*}
where the second equality holds by Lemma  \ref{C:Traces} and the others come
from straightforward commutative diagrams.
\end{proof}

\medskip
Now, after proving these preliminary statements we are ready to prove Theorem
\ref{T:GeomMassey}.
Let $(\kV_1, \kV_2, y_1, y_2, \omega)$ be the data fixed at  the beginning of
Section \ref{S:Massey}.
Recall that we have to compare the triple Massey product
$$
m^{\kV_1, \kV_2}_{y_1, y_2}: \Hom_E(\kV_1, \CC_{y_1})
\otimes \Ext_E^1(\CC_{y_1}, \kV_2) \otimes \Hom_E(\kV_2, \CC_{y_2}) \lar
\Hom_E(\kV_1, \CC_{y_2})
$$
with the map
$$
\tilde{\tilde{r}}^{\kV_1, \kV_2}_{y_1, y_2} := \ev_{y_2}^{\kV_1, \kV_2(y_1)}
\circ \bigl(\res_{y_1}^{\kV_1, \kV_2}(\omega_{y_{1}})\bigr)^{-1}:
\Hom_\CC(\kV_1|_{y_1}, \kV_2|_{y_1}) \lar
\Hom_\CC(\kV_1|_{y_2}, \kV_2|_{y_2}).
$$

\begin{proposition}\label{P:comdiag}
If $g \in\! \Hom_E(\kV_1, \CC_{y_1})$, $h \in\! \Hom_E(\kV_2, \CC_{y_2})$ and
$a \in\! \Ext_E^1(\CC_{y_1}, \kV_2)$,  then
$$
h_{y_2} \circ \tilde{\tilde{r}}^{\kV_1, \kV_2}_{y_1, y_2}\bigl(T_{y_1}(a) g_{y_1}\bigr)
= \bigl(m^{\kV_1, \kV_2}_{y_1, y_2}(g \otimes a \otimes h)\bigr)_{y_2}.
$$
\end{proposition}

\begin{proof} Let us first explain our notation. We have a composition map
$$
\kV_1|_{y_1} \xrightarrow{g_{y_1}} \CC
\xrightarrow{T_{y_1}(a)} \kV_2|_{y_1},
$$
hence we may consider
$$
\kV_1|_{y_2}
\xrightarrow{\tilde{\tilde{r}}^{\kV_1, \kV_2}_{y_1, y_2}\bigl(T_{y_1}(a)
g_{y_1}\bigr)}
\kV_2|_{y_2}.
$$
Let
$
0 \to \kV_2 \stackrel{\alpha}\lar \kA \stackrel{\beta}\lar \CC_{y_1} \to 0
$
be an exact sequence representing $a \in \Ext^1_E(\CC_{y_1}, \kV_2)$.
Then we have a commutative diagram
$$
\xymatrix@R-=20pt
{
&    0   \ar[d]          &              &   0 \ar[d] \\
&     \kV_2   \ar[d]_\alpha \ar@{=}[rr]  & &   \kV_2 \ar[d]^r \\
\kV_1 \ar[r]^{\tilde g} \ar[rd]_g \ar[d] & \kA \ar[d]_\beta \ar[rr]^\varepsilon & &
\kV_2(y_1)  \ar[d]^{\underline{\res}^{\kV_2}_{y_1}}\\
\kV_1\otimes \CC_{y_1} \ar[r]_-{g_{y_1}}  & \CC_{y_1} \ar[rr]_-{T_{y_1}(a)} \ar[d] & &
\kV_2 \otimes \CC_{y_1}\ar[d]  \\
&    0                 &          &   0
}$$
where $\tilde{g}$ is the unique lift of $g$ and the two columns on the right
form a transposed version of the diagram from  Remark \ref{R:connmap}. Since
$$
\res^{\kV_1, \kV_2}_{y_1}(\omega_{y_{1}}):
\Hom_E\bigl(\kV_1, \kV_2(y_1)\bigr) \lar \Hom_E(\kV_1 \otimes
\CC_{y_1},  \kV_2 \otimes \CC_{y_1})
$$
is an isomorphism, by definition we have (see Remark \ref{rem:thirdres})
$$
\bigl(\res^{\kV_1, \kV_2}_{y_1}(\omega_{y_{1}})\bigr)^{-1}
\bigl(T_{y_1}(a)g_{y_1}\bigr)
= \varepsilon \tilde{g}.
$$
Moreover, tensoring the whole diagram with $\CC_{y_2}$ we obtain
a new commutative diagram
$$
\xymatrix
{
   &  \kV_2 \otimes \CC_{y_2} \ar[rr]^-=
\ar[d]_{\alpha_{y_2}}& & \kV_2 \otimes \CC_{y_2} \ar[d]^{r_{y_2}}\\
\kV_1 \otimes \CC_{y_2} \ar[r]^{\tilde{g}_{y_2}} & \kA \otimes \CC_{y_2}
\ar[rr]^-{\varepsilon_{y_2}} & & \kV_2(y_1) \otimes \CC_{y_2},
}
$$
from which the identity
$
\ev^{\kV_1, \kV_2(y_1)}_{y_2}(\varepsilon \tilde{g}) = \alpha_{y_2}^{-1} \tilde{g}_{y_2}
$
follows.
By the definition of Massey products we have a commutative diagram
$$
\xymatrix
{
0 \ar[r] & \kV_2 \ar[rr]^\alpha \ar[d]_h & & \kA \ar[r] \ar[lld]_{\tilde{h}}
& \CC_{y_1} \ar[r] & 0 \\
& \CC_{y_2}  & & \kV_1 \ar[ll]^-{m_3(g \otimes a \otimes h)}
\ar[u]^{\tilde g}
}
$$
which finally implies
$$
\bigl(m^{\kV_1, \kV_2}_{y_1, y_2}(g\otimes a \otimes h)\bigr)_{y_2} =
h_{y_2} \circ \alpha_{y_2}^{-1} \circ \tilde{g}_{y_2} =
h_{y_2} \circ \tilde{\tilde{r}}^{\kV_1, \kV_2}_{y_1, y_2}(T_{y_1}(a) g_{y_1}).
$$
\qed

\medskip
\medskip
Now we are ready to  finish  the proof of Theorem \ref{T:GeomMassey}.
Our goal is to  keep track of the linear map
$m^{\kV_1, \kV_2}_{y_1, y_2}$ under a long chain of canonical  isomorphisms.
Let us do it step by step.
Each linear map
 $$m \in \Lin\bigl(\Hom_E(\kV_1, \CC_{y_1}) \otimes
\Ext^1_E(\CC_{y_1}, \kV_2) \otimes
\Hom_E(\kV_2, \CC_{y_2}), \Hom_E(\kV_1, \CC_{y_2})\bigr)$$
corresponds to an element
$$
n
\in \Lin\bigl(\Hom_E(\kV_1, \CC_{y_1}) \otimes
\Hom_E(\kV_2, \CC_{y_1})^*,
\Lin\bigl(\Hom_E(\kV_2, \CC_{y_2}), \Hom_E(\kV_1, \CC_{y_2})\bigr)\bigr)
$$
which is related to $m$ by the formula
$$
n\bigl(g \otimes S(a)\bigr)(h) = m(g \otimes a \otimes h),
$$
where $S: \Ext^1_E(\CC_{y_1}, \kV_2) \to \Hom_E(\kV_2, \CC_{y_1})^*$ is
given by the bilinear form $\langle \,\, ,\,\,\rangle_{\kV_2, \CC_{y_1}}^{\omega}$ from
Proposition \ref{P:SD}.
By  Lemma \ref{L:SvsT},  the element $S(a) \in \Hom_E(\kV_2, \CC_{y_1})^*$
is mapped to  $T_{y_{1}}(a) \in \Hom_\CC(\CC, \kV_2|_{y_1})$ under the chain of
isomorphisms
$$
\Hom_E(\kV_2, \CC_{y_1})^* \lar
\Hom_E(\kV_2 \otimes \CC_{y_1}, \CC_{y_1})^* \lar
\Hom_\CC(\kV_2|_{y_1}, \CC)^* \lar
\Hom_\CC(\CC, \kV_2|_{y_1}).
$$
This  implies that the linear map
$n$ corresponds to
$$
l \in \Lin\bigl(\Hom_\CC(\kV_1|_{y_1}, \CC) \otimes
\Hom_\CC(\CC, \kV_2|_{y_1}), \Lin\bigl(\Hom_\CC(\kV_2|_{y_2}, \CC),
\Hom_\CC(\kV_1|_{y_2}, \CC)\bigr)\bigr)
$$
given by $l(g_{y_1} \otimes T_{y_{1}}(a))(h_{y_2}) =
m(g \otimes a \otimes h)_{y_2}$.
But since $$\Hom_\CC(\kV_1|_{y_1}, \CC) \otimes
\Hom_\CC(\CC, \kV_2|_{y_1}) \stackrel{\circ}\lar \Hom_\CC(\kV_1|_{y_1}, \kV_2|_{y_1})$$
is an isomorphism and
$$
\Lin\bigl(\Hom_\CC(\kV_2|_{y_2}, \CC),
\Hom_\CC(\kV_1|_{y_2}, \CC)\bigr) \cong
\Hom_\CC(\kV_1|_{y_2}, \kV_2|_{y_2}),
$$
we obtain a linear map
$$
k \in \Lin\bigl(\Hom_\CC(\kV_1|_{y_1}, \kV_2|_{y_1}),
\Hom_\CC(\kV_1|_{y_2}, \kV_2|_{y_2})\bigr)
$$
such that for any elements $g$, $a$ and $h$ the following diagram commutes
$$
\xymatrix
{
\kV_1|_{y_2} \ar[rr]^{k\bigl(T_{y_{1}}(a)g_{y_1}\bigr)}
\ar[rrdd]_{m(g \otimes a \otimes h)_{y_2}}& &  \kV_2|_{y_2} \ar[dd]^{h_{y_2}}\\
& & \\
 & &  \CC. \\
}
$$
Using Proposition \ref{P:comdiag}, we obtain with $m = m^{\kV_1, \kV_2}_{y_1, y_2}$
the identity
$$k =
\tilde{\tilde{r}}^{\kV_1, \kV_2}_{y_1, y_2} = \ev_{y_2}\circ \res_{y_1}^{-1},
$$
A tedious diagram chase shows that $k$ is equal to $\tilde{r}^{\kV_1,
  \kV_2}_{y_1, y_2}$, which  was defined directly after Remark
\ref{R:formsofAYBE}. This completes the proof.
\end{proof}

\noindent
The following Theorem explains how  triple Massey products on a genus one curve
can be computed in a practical way.

\begin{theorem}\label{C:Masseyprod}
Let $E$ be a reduced projective curve with trivial dualizing sheaf,
$x_1, x_2 \in \breve{E}$ be a pair of distinct smooth points lying on the same
irreducible component of $E$. Let $\kV_1$ and $\kV_2$ be a pair of vector
bundles on $E$ satisfying both vanishing conditions from the beginning of
Section \ref{S:Massey}. Let $\pi: Y \to E$ be the  normalization morphism if
$E$ is singular or the universal covering $\CC \to E = \CC/\langle 1, \tau
\rangle$ if $E$ is smooth.
Take a point $y_2$ on $Y$ such that $\pi(y_2) = x_2$, let
$y_{1}\in D_{1}=\pi^{-1}(x_{1})$  and denote $\widetilde\kV_i = \pi^*\kV_i$
for $i=1,2$.
Let $\omega \in H^0(\Omega_E)$ be a global regular differential form on $E$
and $\tilde\omega$ be its (possibly meromorphic) lift on $Y$.
Then the diagram
$$
  \xymatrix@C-C
  {
    \Hom_E\bigl(\kV_1 \otimes \CC_{x_1}, \kV_2 \otimes \CC_{x_1}\bigr)
    \ar[dd]_{\tilde{r}^{\kV_1, \kV_2}_{x_1,x_2}(\omega)}
    \ar[rr]^{\pi^*} & &
    \Hom_Y\bigl(\widetilde\kV_1 \otimes \CC_{y_1},
    \widetilde\kV_2 \otimes \CC_{y_1}\bigr)
    \\
    & \Hom_E\bigl(\kV_1, \kV_2(x_1)\bigr)
    \ar[lu]^{\res_{x_1}^{\kV_1, \kV_2}(\omega_{x_{1}})}
    \ar[r]^{\pi^*}
    \ar[ld]_{\ev_{x_2}^{\kV_1, \kV_2(x_1)}} &
    \Hom_Y\bigl(\widetilde\kV_1, \widetilde{\kV}_2(D_1)\bigr)
    \ar[u]_{\res_{y_1}^{\widetilde\kV_1,\widetilde\kV_2}(\tilde\omega_{y_{1}})}
    \ar[d]^{\ev_{y_2}^{\widetilde\kV_1, \widetilde{\kV}_2(D_1)}}
    \\
    \Hom_E\bigl(\kV_1 \otimes \CC_{x_2}, \kV_2\otimes \CC_{x_2}\bigr)
    \ar[rr]^{\pi^*} & &
    \Hom_Y\bigl(\widetilde\kV_1 \otimes \CC_{y_2},
    \widetilde\kV_2 \otimes \CC_{y_2}\bigr)
  }
$$
is commutative. This shows in particular that the computation of triple Massey
products on elliptic curves  (resp.\/ on singular genus one curves)
can be expressed by computations on the universal
covering (resp.\/ on the normalization).
\end{theorem}

\begin{proof}
The left triangle is commutative by Theorem \ref{T:GeomMassey}.
The two squares are commutative by Propositions \ref{P:functofres} and
\ref{L:functofev}.
\end{proof}

\section{A relative construction  of geometric triple Massey products}
\label{S:RelativeConstr}

Our next goal is to extend the definition of the morphism
$r^{\kV_1, \kV_2}_{y_1,  y_2}(\omega)$, constructed in the previous section,
to genus one fibrations. We achieve this by generalizing the
construction of Theorem \ref{T:GeomMassey} to the relative case.
Throughout this section we work either in the category of locally Noetherian
algebraic schemes over an algebraically closed field $\kk$ of characteristic
zero  or in the category of complex analytic spaces.

\subsection{The relative residue map}
Let $p: X \lar S$ be a \emph{smooth} map of complex analytic  spaces or of
algebraic schemes. Assume $p$ has a section $\imath: S \lar X$,  let $D$ be the
image of $\imath$ equipped with the  ringed space structure induced from $S$.
Recall that the sheaf of relative differentials $\Omega^1_{X/S}$ is defined via
the exact sequence
$$
p^*\Omega^1_S \lar \Omega^1_X \lar \Omega^1_{X/S} \lar 0,
$$
see \cite[Chapter 7]{AltmanKleiman},
\cite[Section II.8 and Section III.10]{Hartshorne} and \cite{PetRem} for
definitions and basic properties of smooth morphisms
and K\"ahler differential forms. In particular, for any closed point $s \in S$
we have: $\Omega^1_{X/S}|_{X_s} \cong \Omega^1_{X_s}$ and $\Omega^1_{X/S}$
is locally free.

Assume additionally that $p$ has relative dimension one and
$X$ itself is smooth.
Our aim is to define  a canonical
epimorphism of $\kO_X$--modules
$$
\underline{\res}_D: \Omega^1_{X/S}(D) \lar \kO_D,
$$
later called \emph{the residue map}. We shall explain our  construction in the
case of algebraic schemes, whereas its generalisation on the case of complex
analytic spaces is straightforward.

\medskip
Let $x \in D \subset X$ be a closed point, then we can find affine
neighbourhoods $U = \Spec(B)$ of $x \in X$ and  $V = \Spec(A)$ of $f(x) \in S$
such that the map $p|_U: U \to V$ is induced by a ring homomorphism
$p^*: A \to B$:
$$
\xymatrix
{ U = \Spec(B) \ar@{^{(}->}[r] \ar[d]_{p|_U} & X \ar[d]^{p}\\
  V = \Spec(A) \ar@{^{(}->}[r]    & S.
}
$$
Then the sheaf $\Omega^1_{X/S}|_U$ is isomorphic to the sheafification of the
$B$--module
of K\"ahler differentials $\Omega_{B/A}$.

Let $\imath^*: B \to A$ be the ring homomorphism corresponding to the section
$\imath$ and $I = \ker(\imath^*)$. Then the map
$C:= B/I \stackrel{\imath^*}\lar  A$ is an isomorphism and
$I$ is the ideal,  locally defining  the subscheme $D$.
By Krull's Hauptidealsatz,  since $U$ is smooth and $V(I)\subset U$ has
codimension one, shrinking the open sets $U$ and $V$ if necessary, we can
achieve that $I$ is generated by a single element $a  \in B$.
From the exact sequence
$$
I/I^2 \stackrel{\delta}\lar \Omega_{B/A} \otimes_B C \lar \Omega_{C/A} \lar 0
$$
where $\delta([b]) = d(b) \otimes 1$ and the fact that $\Omega_{C/A} = 0$ it
follows that the $C$--module $\Omega_{B/A} \otimes_B C$ is generated by a
single element, namely $d(a) \otimes 1$. The smoothness of $p$ implies that
$\Omega_{B/A} \otimes_B C$ is a free $C$--module with this generator.

\begin{definition}\label{D:ressmooth}
Let $p: X \lar S$ be a smooth map of relative dimension one,
$\imath: S \lar X$ a section of $p$ and $D = \imath(S)$.
We define the sheaf homomorphism
$$
\underline{\res}_D:  \Omega^1_{X/S}(D) \lar \kO_D
$$ to be the composition of the canonical map
$
\Omega^1_{X/S}(D) \lar \Omega^1_{X/S}|_D \otimes \kO(D)|_D
$
and the morphism
$\Omega^1_{X/S}|_D \otimes \kO(D)|_D \to  \kO_D
$
locally defined as follows.
In the above notation let
$M = \left\{\left.\dfrac{u}{a}\;\right|\; u \in B\right\} =
\Gamma\bigl(U,  \kO_X(D)\bigr) \subset Q(B)$, where
$Q(B)$ is the field of fractions of the integral domain $B$. The map
$$
\rho:  (\Omega_{B/A} \otimes_B C) \otimes (M \otimes_B C) \lar C
$$
is given by the formula
$
\bigl(d(a)  \otimes 1\bigr) \otimes \bigl(\frac{u}{a} \otimes 1\bigr)
\mapsto u \otimes 1 = \overline{u} := u \mod I.
$
\end{definition}

It is easy to see that the morphism $\rho$ is $C$-linear, surjective
and  does not depend on the choice of a  generator  of the ideal $I$.

\begin{remark}\label{rem:compareres}
  If $S$ is a point and $X$ a smooth complex curve, the residue map in
  Definition \ref{D:ressmooth} coincides with the classical residue map, which
  was used in sequence \eqref{E:residue} at the beginning of subsection
  \ref{SS:resmap}.
  To see this, we let $D=\{x\}$ and $U$ a neighbourhood of $x$ in $X$ with a
  coordinate $z$ centred at $x$. Then, in the notation of Definition
  \ref{D:ressmooth}, $a=z$ and $\frac{f(z)}{z}dz \in \Omega^{1}_{X}(x)$ is
  first sent to $\bigl(dz  \otimes 1\bigr) \otimes \bigl(\frac{f(z)}{z}
  \otimes 1\bigr)$ and then to $f\mod I_{x} = f(0)$, which is equal to the
  ordinary residue of $\frac{f(z)}{z}dz$.
\end{remark}

\medskip

\begin{proposition}\label{P:resbasechange}
Let $p: X \lar  S$ be a smooth map as above, $\imath: S \lar  X$ a section of
$p$  and $g: S' \lar  S$ any morphism. Let $X' = X \times_S S'$ and
 $\imath': S' \lar
X'$ be the section defined by  the universal property of pull-backs:
$$
\xymatrix
{
S' \ar@/_/[ddr]_{=} \ar@/^/[drr]^{\imath g} \ar@{.>}[dr]|-{\imath'} &  &  \\
   & X' \ar[r]^f \ar[d]^{p'} & X \ar[d]_p \\
   & S' \ar[r]^g        & S \ar@/_10pt/[u]_\imath
}
$$
and $D' = \imath'(S')$. Then the following diagram is commutative:
$$
\xymatrix
{
f^*\bigl(\Omega^1_{X/S}(D)\bigr)
\ar[rr]^-{f^*(\underline{\res}_D)}
\ar[d]_{\cong} & & f^*(\kO_D)
\ar[d]^\cong \\
\Omega^1_{X'/S'}(D')
\ar[rr]^-{\underline{\res}_{D'}}  &  &
\kO_{D'}
}
$$
where the vertical arrows  are canonical isomorphisms.
\end{proposition}

\noindent
\emph{Proof}. The problem is local, so we can assume, without loss of
generality, $X = \Spec(B), X' = \Spec(B'), S = \Spec(A)$ and
$S' = \Spec(A')$. Then, we have a  Cartesian diagram of rings and ring
homomorphisms
$$
\xymatrix
{
B' \ar@/_10pt/[d]_{\imath'^*} & & \ar[ll]_{f^*}  B \ar@/^10pt/[d]^{\imath^*} \\
A' \ar[u]_{p'^*} & & A \ar[ll]^{g^*} \ar[u]^{p^*}
}
$$
where $B' = B \otimes_A A'$,
$p'^{*}(a') = 1\otimes a'$ and $f^*(b) = b \otimes 1$.
Denote  $C := B/\ker(\imath^*)$ and $C' := B'/\ker(\imath'^{*})$
then we have an isomorphism
of $C'$--modules   $C \otimes_B B' \lar  C'$.

Let $d: B \to \Omega_{B/A}$ and $d': B' \to \Omega_{B'/A'}$ be the
universal derivations  from the definition of K\"ahler differentials.
By the  universal property we obtain a uniquely determined
$B$--module homomorphism $\Omega_{B/A} \lar  \Omega_{B'/A'}$ and an induced
$B'$--module isomorphism
$\tilde{f}^*:\Omega_{B/A} \otimes_B B' \lar  \Omega_{B'/A'}$
making the following  diagram
$$
\xymatrix
{
B \ar[d]_{f^*}  \ar[r]^-{d}  & \Omega_{B/A} \ar@{.>}[d] \ar[r] &
\Omega_{B/A} \otimes_B B' \ar@{:>}[ld]^{\tilde{f}^*} \\
B' \ar[r]^-{d'} & \Omega_{B'/A'}
}
$$
commutative, in particular
$\tilde{f}^*\bigl(d(b) \otimes 1\bigr) = d'\bigl(f^*(b)\bigr)$.
Moreover, we have a canonical homomorphism
$f^*_M:  M \otimes_B B' \lar M'$, given by
$f^*_M\Bigl(\dfrac{u}{a}\otimes 1\Bigr) = \dfrac{f^*(u)}{f^*(a)}$, where
$$
M = \left\{\left.\frac{u}{a}\right| u \in B\right\} \subset Q(B)
\, \, \text{ and }\, \,
M' = \left\{\left.\frac{v}{f^*(a)}\right| v \in B'\right\} \subset Q(B').
$$
We know that the $C$--module
$\Omega_{B/A} \otimes C$ is generated by the single element $d(a) \otimes
1$. Hence the commutativity of the diagram
$$
\xymatrix
{
\bigl(\Omega_{B/A} \otimes_B B'\bigr) \otimes \bigl(M \otimes_B B'\bigr)
\ar[d]_{\tilde{f}^* \otimes f_M^*} \ar[rrrr]^-{f^*(\underline{\res}_D)} & & & &
C \otimes_B B' \ar[d] \\
\Omega_{B'/A'} \otimes M'  \ar[rrrr]^-{\underline{\res}_{D'}} & & & & C'
}
$$
can be checked on elements of the following form:
$$
\xymatrix
{
\bigl(d(a) \otimes 1\bigr)
\otimes \bigl(\frac{\displaystyle u}{\displaystyle a}
\otimes 1\bigr) \ar@{|->}[rrr]^-{f^*(\underline{\res}_D)}
\ar@{|->}[d]_{\tilde{f}^* \otimes f^*_M}  & & &
\bar{u} \otimes 1 \ar@{|->}[d] \\
d'\bigl(f^*(a)\bigr)
\otimes \dfrac{f^*(u)}{f^*(a)}
\ar@{|->}[rrr]^-{\underline{\res}_{D'}} & & &
\overline{f^*(u)}.
}
$$
and the proposition is proven.
\qed

\subsection{On the sheaf of relative differential forms  of a Gorenstein fibration}

\medskip
\noindent
Let $p: X \lar S$ be a \emph{proper} and \emph{flat}  morphism of relative
dimension \emph{one}, either in the category of complex analytic spaces or of
algebraic schemes over an algebraically closed field $\kk$ of characteristic zero. Assume
additionally that for all closed points $s \in S$ the fibres $X_s$ are reduced
and we have an embedding
$$
\xymatrix
{
X  \ar[rd]_{p} \ar@{^{(}->}[rr] &  &  Y \ar[ld]^{q} \\
 &  S &
}
$$
where $q: Y \lar S$ is a \emph{proper} and \emph{smooth}
  morphism of the relative dimension \emph{two}.

\begin{remark}
Since for any $s \in S$ the surface $Y_s$ is smooth and $X_s \subset Y_s$
has codimension one,
the curve $X_s$ has hypersurface singularities and is in particular Gorenstein.
\end{remark}

\noindent
Recall that for any morphism $q$ we have an exact sequence
$$
q^* \Omega^1_S \lar \Omega^1_Y \lar \Omega^1_{Y/S} \lar 0
$$
and that $\Omega^1_{Y/S}$ is a locally free $\kO_Y$--module of rank two,
because $q$ is a smooth morphism.

\begin{definition}
The relative dualizing sheaf is defined by the formula
$$
\omega_{X/S} :=
\left(\left. {\textstyle\bigwedge^2}\; \Omega^1_{Y/S} \otimes
\kO_Y(X)\right)\right|_{X}.
$$
\end{definition}

\begin{proposition}[see Chapter II in \cite{BarthPetersVen}]
For any $s \in S$ the sheaf $\omega_{X/S}|_{X_s}$ is the dualising sheaf
of the projective curve  $X_s$.
\end{proposition}

\begin{remark}
It can be shown that up to the  pull-back of a line bundle on $S$
this definition of $\omega_{X/S}$ does not depend on the embedding $X
\hookrightarrow Y$.
\end{remark}

\medskip
\noindent
Let $\breve{X}$ be the regular locus of $p$, then
$\jmath: \breve{X} \lar   X$ is an open embedding and the morphism
$\breve{X} \stackrel{p}\lar S$  is flat but in general not proper.
Our aim is to define an injective map of $\kO_X$--modules
$\underline{\cl}_S: \omega_{X/S} \lar \jmath_*(\Omega^1_{\breve{X}/S})$.

For a closed point
$x \in \breve{X}$ let  $U \subset Y$ be an open neighbourhood of  $x$ and
$S_0$  an  open neighbourhood of $f(x)$ in $S$.   Choose
local coordinates $(u, v, s)$ on $U$ such that we have a commutative diagram
$$
\xymatrix
{ \bigl(\CC^2 \times S_0\bigr) = U \ar[d]_{\pr} \ar[rr]  & & Y \ar[d]^q \\
  S_0 \ar[rr] & &  S
}
$$
where $\pr(u,v,s) = s$ and $du|_{\breve{X}} \not\equiv 0$, $dv|_{\breve{X}}
\not\equiv 0$. Assume that the closed subset $X \cap U$ is given in $U$  by
an equation $f(u,v,s) = 0$. Then
$$
\left(\frac{\partial  f}{\partial  u} du + \left.
\frac{\partial  f}{\partial  v} dv\right)\right|_{\breve{X}} = 0.
$$
where the left-hand side of this  equality is viewed as a local section
of $\Omega^1_{\breve{X}/S}$.

\medskip
\noindent
Consider  the composition map
$
\ell: \breve{X} \stackrel{\jmath}\lar    X \lar   Y.
$
\begin{definition}[see Section II.1 in \cite{BarthPetersVen}]
The Poincar\'e residue map is the morphism of $\kO_Y$--modules
$$
\underline{\res}^P: \wedge^2\Omega^1_{Y/S}(X) \lar \ell_*\Omega^1_{\breve{X}/S}
$$
locally defined as follows.
Let $U \subseteq Y$ be an open neighbourhood of $x \in \breve{X}$ as above and
$V := U \cap \breve{X}$,  then the map
$$
\res^P: \Gamma\bigl(U, \, \wedge^2\Omega^1_{Y/S}(X)\bigr) \lar
\Gamma\bigl(V, \, \Omega^1_{\breve{X}/S}\bigr) =
\Gamma\bigl(U, \, \ell_{*}\Omega^1_{\breve{X}/S}\bigr)
$$
is given by the formula
$$
\frac{h \, du \wedge dv}{f} \mapsto
\begin{cases}
  \phantom{-}\left. \displaystyle\frac{h du}{\partial_v f}\right|_{V} & \text{ if } \;
 \displaystyle\frac{\partial f}{\partial v}(u,v,s) \ne 0, \\[3mm]
\left.  -\displaystyle\frac{ h dv}{\partial_u f}\right|_{V} & \text{ if }\;
 \displaystyle\frac{\partial f}{\partial u}(u,v,s) \ne 0.
\end{cases}
$$
\end{definition}

\begin{remark}
Since for any point $s \in S$  the fibre $\breve{X}_s$ is a  smooth curve,
the set $V(f, \partial_u f, \partial_v f) \subseteq \breve{X}_s$ is empty and
the map $\res^P$ is well-defined. Moreover,  $\underline{\res}^P$ is independent of the
choice of a local equation $f \in \kO_Y(U)$ for $X \subset Y$ and also of
the choice of local coordinates $(u, v, s)$ on $Y$,   see for example
\cite[Section II.1]{BarthPetersVen}.
\end{remark}

\medskip
\noindent
From what was said above it follows:

\begin{corollary}\label{C:classmap}
The commutative diagram of $\kO_Y$--modules
$$
\xymatrix
{
0 \ar[r] &  \wedge^2 \Omega^1_{Y/S} \ar@{->}[rd]_0  \ar[r] &
\wedge^2 \Omega^1_{Y/S}(X)
\ar[r] \ar[d]^{\underline{\res}^{P}} & \left.
\wedge^2 \Omega^1_{Y/S}(X)\right|_{X}  \ar[r] \ar@{.>}[ld]  & 0 \\
         &                       & \ell_{*} \Omega^1_{\breve{X}/S} &  &
}
$$
induces an injective morphism of $\kO_X$--modules
$$
\underline{\cl}_S: \omega_{X/S} = \left. \wedge^2 \Omega^1_{Y/S}(X)\right|_{X}   \lar \jmath_* \Omega^1_{\breve{X}/S}.
$$
\end{corollary}

\medskip

\begin{remark}
In what follows  the  morphism $\underline{\cl}_S$ will be called the \emph{class map}.
For a Gorenstein projective variety $X$ of dimension $n$ let $\kM_X$ denote
the sheaf of meromorphic functions on $X$.
Ang\'eniol and Lejeune-Jalabert construct a morphism $\Omega^n_{X} \lar
\omega_X$ which induces  an isomorphism
$
\Omega^n_{X} \otimes \kM_X  \stackrel{\cong}\lar
\omega_X \otimes \kM_X,$
also called ``class map'', see \cite{AngeniolLejeune}.
The relationship between this class map and the class map from  Corollary \ref{C:classmap}
will be discussed elsewhere.
\end{remark}

\noindent
The following proposition can be shown on the lines of \cite[Section
II.1]{BarthPetersVen}.

\begin{proposition}\label{P:classchange}
Let $p: X \lar S$ be a Gorenstein fibration of relative dimension one
satisfying the conditions from the beginning of this subsection.
If $g: S' \lar S$ is any  base change, we obtain the Cartesian diagram
$$
\xymatrix
{
 X' \ar[r]^f \ar[d]_{p'} & X \ar[d]^p \\
 S'  \ar[r]^g  & S.
}
$$
Then, the following diagram is commutative
$$
\xymatrix
{
  f^*\left(\left.\wedge^2 \Omega^1_{Y/S}(X)\right|_X\right) \ar[rr]^-\cong
  \ar[d]_{f^*(\underline{\cl}_S)}
  & & \left. \wedge^2 \Omega^1_{Y'/S'}(X')\right|_{X'} \ar[d]^{\underline{\cl}_{S'}} \\
  f^*(\jmath_*\Omega^1_{\breve{X}/S}) \ar[rr]^-\cong& &\jmath'_*\Omega^1_{\breve{X'}/S'}
}
$$
where the upper  horizontal isomorphism is canonical and the lower one is
induced by the  base-change property.
\end{proposition}

\noindent
The reason to introduce the map $\underline{\cl}_S$ is explained by the following
proposition.

\begin{proposition}[see Proposition 6.2 in \cite{BarthPetersVen}]
\label{prop:classmap}
Let $p: X \lar S$ be as in Proposition \ref{P:classchange},  $t \in S$ a
closed point and
$
\underline{\cl}_t: \omega_{X_t} \lar \jmath_{t*}\Omega^1_{{\breve X}_t}
$
the class map constructed in Corollary \ref{C:classmap}. Then we have:
\begin{enumerate}
\item If the fibre $X_t$ is smooth, then the image of $\underline{\cl}_t$ is the sheaf
$\Omega^1_{X_t}$ of holomorphic differential one-forms on $X_t$.
\item In the case $X_t$ is singular, the image of $\underline{\cl}_t$
is the sheaf of Rosenlicht's differential forms, see Definition
\ref{D:Rosenlicht}.
In particular, $\im(\underline{\cl}_t)$  is a subsheaf of
the sheaf of \emph{meromorphic} differential one-forms on $X_t$
 regular at smooth points of $X_t$.
\end{enumerate}
\end{proposition}

\medskip
The following definition is
central for our  construction of associative geometric
$r$--matrices.
Let $p: X \lar S$ be flat and proper morphism  such that

\vspace{0.2mm}

\noindent
$\bullet$ All fibres $X_t$, $t \in S$ are reduced
projective Gorenstein curves.

\noindent
$\bullet$ There exists an embedding
$$
\xymatrix
{
X  \ar[rd]_{p} \ar@{^{(}->}[rr] &  &  Y \ar[ld]^{q} \\
 &  S &
}
$$
where $q:Y \lar S$ is a proper and smooth morphism of relative dimension two.

\begin{definition}\label{D:ResMap}
Let $\jmath: \breve{X} \lar X$ be the inclusion of the  smooth locus of $p$,
$\imath: S \lar \breve{X}$ a section of $p$ and $D = \imath(S)$.
Then the residue map
$$
\underline{\res}_D: \omega_{X/S}(D) \stackrel{\underline{\cl}_S}\lar
\jmath_*\bigl(\Omega^1_{\breve{X}/S}(D)\bigr) \lar \kO_D
$$
is defined as the composition of the class map $\underline{\cl}_S$ from
Corollary \ref{C:classmap} and the residue map for smooth morphisms of
relative dimension one from Definition \ref{D:ressmooth}.
\end{definition}

\begin{remark}\label{rem:comparecl}
  If $S$ is a point, $X$ a complex curve and $D=\{x\}$, the residue map in
  Definition \ref{D:ResMap} fits as the top horizontal arrow into the
  following commutative diagram
  \[
  \xymatrix@C+10mm{
    \omega_{X}(x) \ar[r]^-{\underline{\res}_x}\ar[d]^{\underline{\cl}_{x}} &
    \CC_{x} \ar@{=}[d]\\
    \Omega_{X}(x) \ar[r]^-{\underline{\res}_x} & \CC_{x}
  }
  \]
  in which the lower horizontal map is the classical residue map used in the
  sequence \eqref{E:residue} at the beginning of subsection \ref{SS:resmap}.
  This follows from Remark \ref{rem:compareres} and Proposition
  \ref{prop:classmap}.
\end{remark}

\noindent
Propositions \ref{P:resbasechange} and \ref{P:classchange} imply the following
corollary.

\begin{proposition}\label{P:Resbasechange}
Let $p: X \lar  S$ and $\imath: S \lar  X$ be as in Definition \ref{D:ResMap}
and $g: S' \lar  S$ be any base change. Denote $X' = X \times_S S'$,
$f: X' \lar X$,  $\imath': S' \lar X'$ the pull-back of $\imath$ and $D' =
\imath'(S')$.
Then the following diagram is commutative
$$
\xymatrix
{
f^*\bigl(\omega_{X/S}(D)\bigr)
\ar[rr]^-{f^*(\underline{\res}_{D})}
\ar[d]_{\cong} & & f^*(\kO_D)
\ar[d]^\cong \\
\omega_{X'/S'}(D')
\ar[rr]^-{\underline{\res}_{D'}}  &  &
\kO_{D'}
}
$$
where the vertical maps are canonical isomorphisms.
\end{proposition}

\subsection{Geometric triple Massey products}

Let $E \stackrel{p}\lar S$ be a genus one fibration embedded into
a smooth fibration of surfaces, i.e we have a commutative diagram
$$\xymatrix
{E \ar@{^{(}->}[rr] \ar[rd]_p & & Y \ar[ld]^q \\
                   & S &
}
$$ where $p$ is a proper and flat map,  for all
$t \in S$ the fibre $E_t$ is  a reduced projective curve
with trivial dualizing sheaf and $q$ is a smooth and proper map
of relative dimension two.

Let $\breve{E}$ be the regular locus of $p$. Assume $S$ is chosen sufficiently
small, so that  $\omega_{E/S} \cong \kO_E$. Fix the following data:
\begin{itemize}
\item A nowhere vanishing  global section $\omega \in H^0(\omega_{E/S})$.
\item Two holomorphic vector bundles $\kV$ and $\kW$ on the
total space $E$ having the same rank and
such that for all $t\in S$ we have:
$$
\Hom_{E_t}\bigl(\kV_{t}, \kW_{t}\bigr) = 0 =
\Ext^1_{E_t}\bigl(\kV_{t}, \kW_{t}\bigr).
$$
Here and in the sequel we denote $\kF_t = \kF|_{E_t}$ for any vector
bundle $\kF$ on $E$.
\item Two sections $h_1, h_2: S \lar \breve{E}$ of $p$
such that for all $t \in S$ we have:
$h_1(t) \ne h_2(t)$ and  $h_1(t), h_2(t)$
belong to the same irreducible component
of $E_t$. We additionally assume that
$$
\Hom_{E_t}\bigl(\kV_{t}(h_2(t)), \kW_{t}(h_1(t))\bigr) = 0 =
\Ext^1_{E_t}\bigl(\kV_{t}(h_2(t)), \kW_{t}(h_1(t))\bigr).
$$
\end{itemize}

\noindent
The main result of this section is the following theorem.

\begin{theorem}\label{T:changedebase}
There exists an isomorphism of vector bundles on $S$
$$
\tilde{r}^{\kV, \kW}_{h_1,h_2}(\omega) = \tilde{r}^{\kV, \kW}_{h_1,h_2}: \quad
h_1^* {\mathcal Hom}_E(\kV, \kW) \lar h_2^* {\mathcal Hom}_E(\kV, \kW)
$$
such that for any base-change diagram
$$
\xymatrix
{ E' \ar[r]^f \ar[d]_{p'} & E  \ar[d]^p\\
  S' \ar[r]^g & S
}
$$
the following diagram
is commutative:
\[
\xymatrix
{
g^* h_1^* {\mathcal Hom}_E(\kV, \kW) \ar[d]_{\cong}
\ar[rrr]^-{g^*\bigl(\tilde{r}^{\kV, \kW}_{h_1, h_2}(\omega)\bigr)} & & &
g^* h_2^* {\mathcal Hom}_E(\kV, \kW) \ar[d]^{\cong} \\
h'^*_1 {\mathcal Hom}_{E'}(f^*\kV, f^*\kW)
\ar[rrr]^-{\tilde{r}^{f^*\kV, f^*\kW}_{h'_1, h'_2}(\omega')} & & &
h'^*_2 {\mathcal Hom}_{E'}(f^*\kV, f^*\kW)
}
\]
where $h'_1, h'_2: S' \lar E'$ are sections of $p'$ obtained as
pull-backs of $h_1$ and $h_2$.
The vertical arrows are canonical isomorphisms and the section $\omega' \in
H^0(\omega_{E'/S'})$ is defined via the commutative diagram
$$
\xymatrix
{f^* \kO_E \ar[r]^\cong  \ar[d]_{f^*(\omega)} & \kO_{E'} \ar[d]^{\omega'} \\
 f^*\omega_{E/S}  \ar[r]^\cong & \omega_{E'/S'}.
}
$$
Moreover, for any $s \in S$ the morphism $\tilde{r}^{\kV_s, \kW_s}_{h_1(s), h_2(s)}$
coincides with the morphism describing triple Massey products constructed in
Section \ref{S:PolConstr}.
\end{theorem}

\begin{proof} The construction of the morphism $\tilde{r}^{\kV, \kW}_{h_1, h_2}$
is the following. Let $D_i = h_i(S)$ and $D_{i}'=h_{i}'(S)$,
then the  exact sequence
\begin{equation}\label{E:residrel}
0 \lar \omega_{E/S} \lar \omega_{E/S}(D_1)
\xrightarrow{\underline{\res}_{D_1}}   \kO_{D_1} \lar 0
\end{equation}
induces an exact sequence
$
0 \lar \kW \otimes \omega_{E/S} \lar \kW \otimes \omega_{E/S}(D_1)
\lar  \kW \otimes \kO_{D_1} \lar 0.
$
Since ${\mathcal Ext}^1_{E}(\kV, \kW) = 0$ and
$\omega_{E/S} \cong \kO_E$, applying the functor ${\mathcal Hom}_E(\kV,\,-\,)$
we obtain the exact sequence
\begin{equation}\label{E:seqVW}
0 \to {\mathcal Hom}_{E}(\kV, \kW)
\to {\mathcal Hom}_{E}(\kV, \kW \otimes \omega_{E/S}(D_1))
\to {\mathcal Hom}_{E}(\kV, \kW \otimes \kO_{D_1})
\to 0.
\end{equation}

\medskip

\begin{lemma}\label{L:imagedirecte}
In the notations of the theorem we have:
$$
p_*\bigl({\mathcal Hom}_{E}(\kV, \kW)\bigr) \cong
R^1p_{*}\bigl({\mathcal Hom}_{E}(\kV, \kW)\bigr) = 0.
$$
\end{lemma}

\noindent
\emph{Proof of the lemma}. It suffices to show that
$\mathbb{R}p_{*}\bigl({\mathcal Hom}_{E}(\kV, \kW)\bigr) = 0$
viewed as an object
of the derived category of coherent sheaves $\Dbcoh(S)$.
Note that
a complex $\kF \in \Dbcoh(S)$ is zero if and only if
for all points
$t \in S$ we have: $\kF \stackrel{\mathbb{L}}\otimes \CC_t \cong  0$.
Since the morphism $p$ is flat, from  a  base-change isomorphism
it follows that
$$
\mathbb{R}p_{*}\bigl({\mathcal Hom}_{E}(\kV, \kW)\bigr)
\stackrel{\mathbb{L}}\otimes \CC_t \cong
\RHom_{E_t}(\kV_{t}, \kW_{t}) \cong  0,
$$
where the last equality follows from the assumption
$\Ext^i_{E_t}(\kV_{_t}, \kW_{_t}) = 0$ for all $i \in \mathbb{Z}$ and
$t\in S$.
\qed

\medskip
Hence, applying the left-exact functor $p_*$ to the exact Sequence
(\ref{E:seqVW}) we obtain an isomorphism
$
p_*{\mathcal Hom}_{E}\bigl(\kV, \kW \otimes \omega_{E/S}(D_1)\bigr)
\stackrel{\cong}\lar
p_*{\mathcal Hom}_{E}\bigl(\kV, \kW \otimes \kO_{D_1}\bigr).
$
Combining it with the  canonical isomorphisms
$$
{\mathcal Hom}_{E}\bigl(\kV, \kW \otimes \kO_{D_1}\bigr) \stackrel{\cong}\lar
{h_1}_*{\mathcal Hom}_{S}\bigl(h_1^*\kV, h_1^*\kW\bigr) \stackrel{\cong}\lar
{h_1}_* h_1^*{\mathcal Hom}_{E}\bigl(\kV, \kW\bigr)
$$
we obtain an isomorphism
$$
\underline{\res}^{\kV, \kW}_{h_1}: \;
p_*{\mathcal Hom}_{E}\bigl(\kV, \kW \otimes \omega_{E/S}(D_1)\bigr)
\stackrel{\cong}\lar
h_1^*{\mathcal Hom}_{E}\bigl(\kV, \kW\bigr).
$$
Moreover, the   choice of a global section
$\kO_E \stackrel{\omega}\lar \omega_{E/S}$ induces
an isomorphism
$$
\underline{\res}^{\kV, \kW}_{h_1}(\omega): \;
p_*{\mathcal Hom}_{E}\bigl(\kV, \kW(D_1)\bigr) \stackrel{\cong}\lar
h_1^*{\mathcal Hom}_{E}\bigl(\kV, \kW\bigr),
$$
which we shall frequently denote by $\underline{\res}^{\kV, \kW}_{h_1}$.

\begin{remark}\label{rem:sameres}
  If $S$ is a point, $E$ a curve and $D_{1}=\{x\}$, this isomorphism coincides
  with the map $\res_{x}^{\kV,\kW}(\omega)$ from Definition
  \ref{def:res-absolut}, if we identify
  $\Hom_{E}(\kV\otimes\CC_{x}, \kW\otimes\CC_{x})$ with the vector space
  ${\mathcal Hom}_{E}(\kV,\kW)|_{x}$.
  This follows from Remarks \ref{rem:thirdres}, \ref{rem:compareres} and
  \ref{rem:comparecl} by comparing the two constructions.
\end{remark}

The construction of another isomorphism
$$
\underline{\ev}^{\kV, \kW(D_{1})}_{h_2}: \;
p_*{\mathcal Hom}_{E}\bigl(\kV, \kW(D_1)\bigr) \stackrel{\cong}\lar
h_2^*{\mathcal Hom}_{E}\bigl(\kV, \kW\bigr)
$$
is similar. We start with the exact sequence
\begin{equation}\label{E:evseq}
0 \lar \kO_E(D_1 - D_2) \lar \kO_E(D_1) \lar \kO_{E}(D_1) \otimes \kO_{D_2}
\lar 0.
\end{equation}
For any Weil divisor $D \subset \breve{E}$ we view the line bundle
$\kO_E(D)$ as a subsheaf of the sheaf of  meromorphic functions $\kM_E$.
There exists   a canonical exact sequence
$$
0 \lar \kO_E \lar \kO_E(D_1)  \xrightarrow{\underline{\ev}_{D_1}}
\kO_{D_1}(D_1) \lar 0
$$
inducing an isomorphism $\kO_{D_2} \lar \kO_E(D_1)\otimes \kO_{D_2}$.
Tensoring the exact sequence (\ref{E:evseq}) with the vector bundle
$\kW$ and applying ${\mathcal Hom}_E(\kV,\,-\,)$ we obtain an
exact sequence
$$
0 \to
{\mathcal Hom}_E\bigl(\kV, \kW(D_1 - D_2)\bigr) \to
{\mathcal Hom}_E\bigl(\kV, \kW(D_1)\bigr) \to
{\mathcal Hom}_E(\kV, \kW \otimes \kO_{D_2}) \to 0.
$$
By the same argument as in Lemma \ref{L:imagedirecte} one can  show that
$$
p_*{\mathcal Hom}_E\bigl(\kV, \kW(D_1 - D_2)\bigr)
\cong  R^1p_*{\mathcal Hom}_E\bigl(\kV, \kW(D_1 - D_2)\bigr) = 0
$$
which implies that we obtain an  isomorphism of vector bundles on $S$
$$
\underline{\ev}^{\kV, \kW(D_{1})}_{h_2}: \;
p_*{\mathcal Hom}_E\bigl(\kV, \kW(D_1)\bigr)
\stackrel{\cong}\lar h_2^* {\mathcal Hom}_E(\kV, \kW).
$$

\begin{remark}\label{rem:sameev}
  If $S$ is a point, $E$ a curve, $D_{1}=\{x\}$ and $D_{2}=\{y\}$, this
  isomorphism coincides with the map $\ev_{y}^{\kV,\kW(x)}$ from Definition
  \ref{def:ev-absolut}. This follows easily by comparing the two
  constructions, using the canonical identification of
  $\Hom_{E}(\kV\otimes\CC_{x}, \kW\otimes\CC_{x})$ with
  ${\mathcal Hom}_{E}(\kV,\kW)|_{x}$.
\end{remark}

The isomorphism of vector bundles
$$\tilde{r}^{\kV, \kW}_{h_1, h_2}: \;
h_1^* {\mathcal Hom}_E(\kV, \kW) \lar h_2^* {\mathcal Hom}_E(\kV, \kW)
$$
is defined by the commutative diagram of vector bundles on $S$
$$
\xymatrix@C-C@R-3mm
{ & p_* {\mathcal Hom}_E\bigl(\kV, \kW(D_1)\bigr)
\ar[ldd]_{\underline{\res}^{\kV, \kW}_{h_1}(\omega)}
\ar[rdd]^{\underline{\ev}^{\kV, \kW(D_{1})}_{h_2}} & \\
 & & \\
h_1^*{\mathcal Hom}_E(\kV, \kW)
\ar[rr]^{\tilde{r}^{\kV, \kW}_{h_1, h_2}(\omega)} & &
h_2^*{\mathcal Hom}_E(\kV, \kW).
}
$$

\begin{remark}\label{rem:ok-without-vanishing}
  If we drop the assumption that
  $$
  \Hom_{E_t}\bigl(\kV_{t}(h_2(t)), \kW_{t}(h_1(t))\bigr) = 0 =
  \Ext^1_{E_t}\bigl(\kV_{t}(h_2(t)), \kW_{t}(h_1(t))\bigr)
  $$
  for all $t\in S$, we still get a morphism
  $\underline{\ev}^{\kV, \kW(D_{1})}_{h_2}$ but it may no longer be an
  isomorphism. This is the only change that occurs to the
  construction. Therefore, in this situation we still obtain a morphism of
  vector bundles
  $$\tilde{r}^{\kV, \kW}_{h_1, h_2} = \tilde{r}^{\kV, \kW}_{h_1, h_2}(\omega): \;
  h_1^* {\mathcal Hom}_E(\kV, \kW) \lar h_2^* {\mathcal Hom}_E(\kV, \kW).
  $$
\end{remark}

\begin{remark}\label{rem:dep-on-omega}
  Because $p:E\rightarrow S$ is proper, two nowhere vanishing global sections
  $\omega, \omega'\in H^{0}(\omega_{E|S})$ differ by a factor
  $\varphi=p^{\ast}(\psi)$ only, where $\psi\in H^{0}(\kO_{S}^{\ast})$. If
  $\omega'=\varphi\omega$, we obtain
  $\underline{\res}^{\kV, \kW}_{h_1}(\omega') =
  \varphi\cdot\underline{\res}^{\kV, \kW}_{h_1}(\omega)$ and
  $\tilde{r}^{\kV, \kW}_{h_1, h_2}(\omega)=
  \varphi\cdot\tilde{r}^{\kV, \kW}_{h_1, h_2}(\omega')$.
  In particular, if $S$ is a point, $\varphi$ is a constant factor.
\end{remark}

\medskip
Now let us prove the compatibility of $\tilde{r}^{\kV, \kW}_{h_1, h_2}$ with
respect to base-change.  We start with the commutative diagram of coherent
sheaves on $E'$
\begin{equation}
  \label{eq:twosqu}
  \xymatrix@C+C
  {
    f^{\ast}\kO_{E} \otimes f^{\ast}\kO_{E}(D_{1})
    \ar[r]^-{f^{\ast}(\omega)\otimes\id} \ar[d]_{\cong}^{\can}&
    f^*\bigl(\omega_{E/S}(D_1)\bigr) \ar[r]^-{f^*(\underline{\res}_{D_1})}
    \ar[d]_{\cong} &  f^*(\kO_{D_1})
    \ar[d]^{\cong}\\
    \kO_{E'}\otimes\kO_{E'}(D_{1}') \ar[r]^-{\omega'\otimes\id}&
    \omega_{E'/S'}(D_{1}') \ar[r]^-{\underline{\res}_{D'_1}}  &   \kO_{D'_1}
  }
\end{equation}
the right part of which was obtained in  Proposition
\ref{P:Resbasechange}. Next apply the functor
$$
p'_*{\mathcal Hom}_{E'}(f^*\kV, f^*\kW \otimes \,-\,):
\Coh(E') \lar \Coh(S')
$$
to the right square, which yields the commutative diagram
$$
\xymatrix
{p'_* {\mathcal Hom}_{E'}\bigl(f^*\kV,
f^*\kW \otimes f^*(\omega_{E/S}(D_1))\bigr) \ar[r] \ar[d]_\cong &
p'_* {\mathcal Hom}_{E'}\bigl(f^*\kV,
f^*\kW \otimes f^*\kO_{D_1}\bigr) \ar[d]^\cong  \\
p'_* {\mathcal Hom}_{E'}\bigl(f^*\kV,
f^*\kW \otimes \omega_{E'/S'}(D'_1)\bigr) \ar[r]  &
p'_* {\mathcal Hom}_{E'}\bigl(f^*\kV,
f^*\kW \otimes \kO_{D'_1}\bigr)
}
$$
in $\Coh(S')$.
There is an isomorphism of functors
$$
f^*{\mathcal Hom}_{E}\bigl(\kV, \kW \otimes \,-\,\bigr)
\lar
{\mathcal Hom}_{E'}\bigl(f^*\kV, f^*\kW \otimes f^*(\,-\,)\bigr)
$$
between the categories of coherent sheaves $\Coh(E)$ and $\Coh(E')$.
Composing these functors with $p'_*$ and then applying them
to the residue map
$\omega_{E/S}(D_1) \lar \kO_{D_1}$  we obtain a commutative diagram
$$
\xymatrix
{p'_*  f^*{\mathcal Hom}_{E}\bigl(\kV,
\kW \otimes \omega_{E/S}(D_1)\bigr) \ar[r] \ar[d]_\cong  &
p'_* f^*{\mathcal Hom}_{E'}\bigl(\kV,
\kW \otimes \kO_{D_1}\bigr) \ar[d]^\cong  \\
p'_* {\mathcal Hom}_{E'}\bigl(f^*\kV,
f^*\kW \otimes f^*(\omega_{E/S}(D_1))\bigr) \ar[r]  &
p'_* {\mathcal Hom}_{E'}\bigl(f^*\kV,  f^*\kW \otimes f^*\kO_{D_1}\bigr).
}
$$
Finally, there is a natural transformation of functors
$g^*p_* \lar p'_* f^*$ (base-change), which
is an isomorphism of functors on the category of $S$--flat coherent sheaves
on  $E$. Since both  sheaves
${\mathcal Hom}_E\bigl(\kV, \kW \otimes \omega_{E/S}(D_1)\bigr)$
and
${\mathcal Hom}_E(\kV, \kW \otimes \kO_{D_1}) \cong
{h_1}_{*} {\mathcal Hom}_{S}(h_1^*\kV, h_1^*\kW)$ are flat over $S$,
we  obtain a commutative diagram
$$
\xymatrix
{
g^* p_*{\mathcal Hom}_E\bigl(\kV, \kW \otimes \omega_{E/S}(D_1)\bigr)
\ar[r] \ar[d]_\cong  &
g^* p_*{\mathcal Hom}_E\bigl(\kV, \kW \otimes \kO_{D_1}\bigr)
 \ar[d]^\cong  \\
p'_* {\mathcal Hom}_{E'}\bigl(f^*\kV,
f^*\kW \otimes \omega_{E'/S'}(D'_1)\bigr) \ar[r]  &
p'_* {\mathcal Hom}_{E'}\bigl(f^*\kV,
f^*\kW \otimes \kO_{D'_1}\bigr).
}
$$
Using similar arguments one can show that the following  diagram is
commutative:
$$
\xymatrix
{g^* p_*{\mathcal Hom}_E\bigl(\kV, \kW \otimes \kO_{D_1}\bigr)
 \ar[d]_\cong  \ar[r]^-\cong  & g^* h_1^*{\mathcal Hom}_E(\kV, \kW) \ar[d]^\cong \\
p'_* {\mathcal Hom}_{E'}(f^*\kV, f^*\kW \otimes \kO_{D'_1}) \ar[r]^-\cong &
 {h'_1}^* {\mathcal Hom}_{E'}(f^*\kV, f^*\kW),
}
$$
in which all arrows are canonical isomorphisms.
Composinig the two previous diagrams, we obtain the compatibility of
$\underline{\res}^{\kV, \kW}_{h_1}$ with base change, i.e.\/ the commutative diagram
$$
\xymatrix
{
g^* p_*{\mathcal Hom}_E\bigl(\kV, \kW \otimes \omega_{E/S}(D_1)\bigr)
\ar[rrr]^-{g^*\left(\underline{\res}^{\kV, \kW}_{h_1}\right)} \ar[d]_\cong & & &
g^* h_1^*{\mathcal Hom}_E\bigl(\kV, \kW\bigr)
 \ar[d]^\cong \\
p'_* {\mathcal Hom}_{E'}\bigl(f^*\kV,
f^*\kW \otimes \omega_{E'/S'}(D'_1)\bigr)
\ar[rrr]^-{\underline{\res}^{f^*\kV, f^*\kW}_{h'_1}}  & & &
{h'_1}^* {\mathcal Hom}_{E'}\bigl(f^*\kV,
f^*\kW\bigr)
}
$$
in which the vertical arrows are compositions of the natural isomorphisms
constructed  above.
If we follow the same steps with the left square in diagram \eqref{eq:twosqu},
we obtain the compatibility with base change for
$\underline{\res}^{\kV, \kW}_{h_1}(\omega)$.

In an analogue way, we can show that $\underline{\ev}^{\kV, \kW}_{h_2}$ is
compatible with base change. This proves the base-change property for
$\tilde{r}^{\kV, \kW}_{h_1,h_2}(\omega)$.

It remains to show that, in case $S$ is a single point, the relative
construction yields the same result as the construction in Section
\ref{S:PolConstr}. This follows from Theorem \ref{T:GeomMassey},
Remark \ref{rem:sameres} and Remark \ref{rem:sameev}.
This finishes the proof of the theorem.
\end{proof}

\medskip
\noindent
Let $r^{\kV, \kW}_{h_1, h_2} = r^{\kV, \kW}_{h_1, h_2}(\omega)$ denote the image
of $\tilde{r}^{\kV, \kW}_{h_1, h_2}$ under the canonical isomorphism
$$
\begin{CD}
  \Hom_S\bigl(h_1^*{\mathcal Hom}_E(\kV, \kW),
              h_2^*{\mathcal Hom}_E(\kV, \kW)\bigr)\\
@VV{\cong}V\\
\Gamma\bigl(S,  h_1^*{\mathcal Hom}_E(\kW, \kV) \otimes
 h_2^*{\mathcal Hom}_E(\kV, \kW)\bigr).
\end{CD}
$$

\medskip
\noindent
From Theorem \ref{T:changedebase} we immediately obtain the following
corollary.

\begin{corollary}\label{C:changedebase}
In the notation of Theorem \ref{T:changedebase}
let
$\eta_{\kV, \kW}: g^*\bigl(h_1^*{\mathcal Hom}_E(\kW, \kV)
\otimes h_2^*{\mathcal Hom}(\kV, \kW)\bigr)
\lar
{h'_1}^* {\mathcal Hom}_{E'}(f^*\kW, f^*\kV) \otimes
{h'_2}^* {\mathcal Hom}_{E'}(f^*\kV, f^*\kW)$
be the canonical isomorphism of bifunctors. Then we have:
$$
\eta_{\kV, \kW}\left(g^*(r^{\kV, \kW}_{h_1, h_2})\right) =
r^{f^*\kV, f^*\kW}_{h'_1, h'_2}.
$$
\end{corollary}

\noindent
The following properties of the morphism $\tilde{r}^{\kV, \kW}_{h_1, h_2}$
are straightforward consequences of the naturality of all the morphisms
involved in the construction.

\begin{proposition}\label{P:multlb}
In the situation of Theorem \ref{T:changedebase}, the isomorphism
$\tilde{r}^{\kV, \kW}_{h_1, h_2}$ is functorial with respect to isomorphisms
$f:\kV \lar \kV'$ and $g:\kW \lar \kW'$, this means that
\[
\xymatrix
{
h_1^*{\mathcal Hom}(\kV,\kW)
\ar[d]_{h_1^*\bigl(\conj(f, g)\bigr)}
\ar[rrr]^{\tilde{r}^{\kV,\kW}_{h_1, h_2}}
& & &
h_2^*{\mathcal Hom}(\kV,\kW)
\ar[d]^{h_2^*\bigl(\conj(f, g)\bigr)} \\
h_1^*{\mathcal Hom}(\kV',\kW')
\ar[rrr]^{\tilde{r}^{\kV',\kW'}_{h_1, h_2}}
& & &
h_2^*{\mathcal Hom}(\kV',\kW')
}
\]
is commutative.
Moreover, for  any line bundle $\kL$ on $E$ the following diagram is
commutative
$$
\xymatrix
{
h_1^*{\mathcal Hom}_E(\kV, \kW) \ar[rrr]^-{\tilde{r}^{\kV, \kW}_{h_1, h_2}}
\ar[d]_\cong
& & &   h_2^*{\mathcal Hom}_E(\kV, \kW) \ar[d]^\cong\\
h_1^*{\mathcal Hom}_E(\kV \otimes \kL, \kW \otimes \kL)
\ar[rrr]^-{\tilde{r}^{\kV \otimes \kL, \kW\otimes \kL}_{h_1, h_2}} & & &
h_2^*{\mathcal Hom}_E(\kV\otimes \kL, \kW\otimes \kL),
}
$$
where the vertical arrows are induced by the canonical isomorphism
$${\mathcal Hom}_E(\kV, \kW) \stackrel{\cong}\lar
{\mathcal Hom}_E(\kV \otimes \kL, \kW\otimes \kL).
$$
\end{proposition}

\section{Geometric associative $r$--matrix}\label{S:GeomR}

The main goal of this section is to define the so-called geometric
associative $r$-matrix attached to a genus one fibration.
Throughout this section, we work either in the category $\Ans$  of complex
analytic spaces or in the category of algebraic schemes over an algebraically
closed field $\kk$ of characteristic zero.
We start with the following geometric data.

\begin{itemize}
\item  Let $E \stackrel{p}\lar T$ be a flat \emph{projective}
  morphism of relative dimension one between reduced complex spaces and
  denote by $\breve{E}$ the smooth locus of $p$.
\item   We assume there exists a section   $\imath: T \lar \breve{E}$ of
  $p$. Let $\Sigma := \imath(T) \subset E$ denote the corresponding Cartier
  divisor.
\item  Moreover, we assume that for all points
  $t \in T$ the fibre $E_t$ is a \emph{reduced} and \emph{irreducible}
  projective curve of \emph{arithmetic genus one}.
\item  The fibration $E \stackrel{p}\lar T$ is embeddable  into a smooth
  fibration of projective surfaces over $T$ and $\omega_{E/T} \cong \kO_E$.
\item For our applications it is convenient to assume  the  fibration $E
  \stackrel{p}\lar T$ is the \emph{analytification} of an algebraic fibration.
\end{itemize}

\subsection{The construction}
For a pair of coprime integers  $(n, d) \in \mathbb{N} \times \mathbb{Z}$
we denote by $\underline{\Mf}_{E/T}^{(n, d)}: \Ans_T \lar \Sets$ the moduli functor
of relatively stable vector bundles of rank $n$ and degree $d$.
In particular, $\underline{\Mf}_{E/T}^{(1, d)}  = \underline{\Pic}^d_{E/T}$
are the relative Picard functors  and $\underline{\Mf}_{E/T}^{(1, 0)}  =
\underline{\Pic}^0_{E/T}$  is the relative Jacobian.
The assumption that $p:E\lar T$ has a section is only needed to ascertain that
these functors have fine moduli spaces.

\begin{theorem}\label{T:ellfibrrepr}
In the above notation we have:
\begin{itemize}
\item The relative Jacobian $\underline{\Pic}^0_{E/T}$ is representable by the
  fibration $\breve{E} \stackrel{p}\lar T$ and the universal line bundle
  $\kL = \kO_{\breve{E} \times_T E}(-\Delta) \otimes \pi_2^*\kO_E(\Sigma)$, where
  $\Delta \subset \breve{E} \times_T E$ denotes the diagonal and
  $\pi_2: \breve{E} \times_T E\lar E$ is the natural projection.
\item The functors $\underline{\Pic}^d_{E/T}$ for all $d\in\ZZ$ are
  isomorphic to each other and these isomorphisms are induced by tensoring
  with the line bundle $\kO_E(\Sigma)$.
\item The natural transformation of functors
  $\underline{\det}: \underline{\Mf}_{E/T}^{(n, d)} \lar
  \underline{\Pic}^d_{E/T}$ is an isomorphism. In particular, the moduli
  functor $\underline{\Mf}_{E/T}^{(n, d)}$ is representable for all pairs
  of coprime integers $(n,d) \in \mathbb{N} \times \mathbb{Z}$.
\end{itemize}
\end{theorem}

\begin{proof}
The first part of this theorem can be found in \cite{AK} the  second statement
is trivial.
The third part seems to be well-known, see  for example \cite{Oda} for the
case of an elliptic curve and  \cite{BK5} for the proof of a more general
statement and further details.
\end{proof}

\medskip
\noindent
From now on, by $M \stackrel{l}\lar T$ we denote a fibration
which, together with a universal family
$\kP = \kP(n,d) \in \VB(M \times_T E)$, represents the
functor $\underline{\Mf}_{E/T}^{(n, d)}$.
For a closed point $t \in T$ we denote by $\kP_t \in \VB(M_t \times E_t)$ the
restriction of $\kP$ to $M_t \times E_t$.

\noindent
The morphisms $p$ and $l$ induce a  morphism
$C:= M \times_T M \times_T \breve{E} \times_T \breve{E} \stackrel{g}\lar  T$,
from which we obtain a Cartesian diagram of complex spaces:
\[
\xymatrix
{
M \times_T M \times_T \breve{E} \times_T \breve{E} \times_T E \ar[d]_q
\ar[rr]^-{f} & & E \ar[d]^p \\
C \ar[rr]^{g} & & T.
}
\]
Observe  that $q: M \times_T M \times_T \breve{E} \times_T \breve{E} \times_T E
\lar C$ is again a genus one fibration satisfying all the conditions listed at
the beginning of this section.

\begin{definition}\label{D:univgeom}
The diagonal map $\Delta: \breve{E} \lar \breve{E} \times_T E$ induces two
\emph{canonical} sections
$$
h_1, h_2: \quad  M \times_T M \times_T \breve{E} \times_T \breve{E} \lar
M \times_T M \times_T \breve{E} \times_T \breve{E} \times_T E
$$
of the morphism $q$,  given by the rule
$h_i(v_1, v_2; y_1, y_2) = (v_1, v_2; y_1, y_2, y_i)$
for $i = 1,2$. Let  $D_i$ be  the reduced image of $h_i$.
Next, consider the two projection maps
$$
\pi_i: M \times_T M \times_T \breve{E} \times_T \breve{E} \times_T E \lar
     M \times_T E,
$$
given by  $\pi_i(v_1, v_2; y_1, y_2, y) = (v_i, y)$ for $i = 1,2$.
For any base point
$x  = (v_1, v_2; y_1, y_2) \in
M \times_T M \times_T \breve{E} \times_T \breve{E}$ with $t = g(x)$ we denote:
$$
\kP^{v_i} := \pi_i^*\kP|_{q^{-1}(x)} \cong  \kP_{t}|_{\{v_i\} \times E_t}\in \VB(E_t).
$$
\end{definition}

\noindent
Consider the following closed subsets of the basis $C$:
$$
\Delta_M =
\bigl\{
(v_1, v_2; y_1, y_2) \in C \, | \,  v_1 =v_2
\bigr\}
\quad \mbox{\textrm{and}} \quad
\Delta_E =
\bigl\{
(v_1, v_2; y_1, y_2) \in C \, | \,  y_1 =y_2
\bigr\}
$$
and their complement
$B =  C \setminus \bigl(\Delta_M \cup \Delta_E\bigr)$. Then we have
 the  induced  genus one fibration:
$$
\xymatrix
{
X \ar[d]_{q|_X}  \ar@{^{(}->}[rr] & &  M \times_T M \times_T \breve{E}
\times_T \breve{E} \times_T E
\ar[rr]^-f \ar[d]_q & & E \ar[d]^p  \\
B \ar@{^{(}->}[rr] & & C
\ar[rr]^g & & T.
}
$$
Note that the images of the sections $h_1, h_2: B \lar X$ are disjoint and for any point
$x = (v_1, v_2; y_1, y_2) \in B$ we have:
$
\pi_1^*\kP|_{q^{-1}(x)} =   \kP^{v_1} \not\cong \kP^{v_2} = \pi_2^*\kP|_{q^{-1}(x)}.
$
Occasionally we shall use the abbreviation
$\kV_i = \pi_i^*\kP|_X \in \VB(X)$ for $i =1,2$.

\begin{lemma}
The set of points $\bar\Delta  =
\left\{x \in B \,\,  | \, \,
\kV_1(D_2)|_{q^{-1}(x)} \cong \kV_2(D_1)|_{q^{-1}(x)}\right\}$ is a closed analytic
subset of  $B$.
 \end{lemma}

\begin{proof}  Since the morphism $q$ is projective, the sheaf
$q_* {\mathcal Hom}\bigl(\kV_1(D_2), \kV_2(D_1)\bigr)$  is coherent.
Moreover, if $\kV$ and $\kW$ are two stable vector  bundles on an irreducible
projective curve $E_t$  of arithmetic  genus one having the same rank and
degree, then
$\Hom_{E_t}(\kV, \kW) \ne 0$ if and only if $\kV \cong \kW$.
Since the sheaf
${\mathcal Hom}\bigl(\kV_1(D_2), \kV_2(D_1)\bigr)$ is locally free,
it is flat over $B$. Therefore, the base-change formula implies that for a
point $x = (v_1, v_2; y_1, y_2)  \in B$ with $t = g(x)$, after identifying
$q^{-1}(x)$ with $E_{t}$, we have:
\[
q_*{\mathcal Hom}\bigl(\kV_1(D_2), \kV_2(D_1)\bigr)\big|_x \cong
\Hom_{E_t}\bigl(\kV_1|_{E_t}(y_2), \kV_2|_{E_t}(y_1)\bigr).
\]
Therefore, the set   $\bar\Delta$ coincides with the
 reduced support of $q_*{\mathcal Hom}\bigl(\kV_1(D_2), \kV_2(D_1)\bigr)$,
hence it is a closed analytic subset.
\end{proof}

\medskip

\begin{definition}
Let  $\omega \in H^0(\omega_{E/T})$ be a nowhere vanishing section
of the dualising sheaf $\omega_{E/T}$ and $f^*(\omega) \in H^0(\omega_{X/B})$
its pull-back to $X$. Theorem \ref{T:changedebase} provides us with a
canonical homomorphism of vector bundles on $B$ (see also Remark
\ref{rem:ok-without-vanishing}):
$$
\tilde{r} = \tilde{r}^{\kV_1, \kV_2}_{h_1, h_2}(\omega)  :=
\tilde{r}^{\kV_1, \kV_2}_{h_1, h_2}\bigl(f^*(\omega)\bigr): \quad
h_1^* {\mathcal Hom}_X(\kV_1, \kV_2) \lar h_2^*{\mathcal Hom}_X(\kV_1, \kV_2)
$$
and a canonical holomorphic section
$$
r = r^{\kV_1, \kV_2}_{h_1, h_2}(\omega)\in
H^0\bigl(B, h_1^* {\mathcal Hom}_X(\kV_2, \kV_1)
\otimes
h_2^*{\mathcal Hom}_X(\kV_1, \kV_2)\bigr).
$$
We call $\tilde{r}$ and $r$ the \emph{geometric associative $r$-matrix} of the
fibration $E\rightarrow T$.
\end{definition}

\noindent
Note that $r$ and $\tilde{r}$ depend on the pair of coprime integers $(n,d)$,
the fibration $E\rightarrow T$ and the section $\omega \in H^0(\omega_{E/T})$
only.
Two different choices of a universal bundle $\kP$ lead to
a canonical isomorphism between the corresponding section spaces
$H^{0}\bigl(B, h_1^* {\mathcal Hom}_X(\kV_2, \kV_1)\otimes
h_2^*{\mathcal Hom}_X(\kV_1, \kV_2)\bigr)$ under which the constructed
sections $r$ are identified. This isomorphism may involve an automorphism of
the moduli space $M$ and the tensor product with the pull back of a line
bundle on $M$ (see Prop.~\ref{P:multlb}).

To formulate the base-change property (Proposition \ref{P:chaba}) for the
geometric associative $r$-matrix $r$, let $\tilde{g}: T' \lar T$ be any
morphism between reduced analytic spaces and
\[
\xymatrix
{
E' \ar[r]^{\tilde f} \ar[d]_{p'} & E \ar[d]^p \\
T' \ar[r]^{\tilde g} & T
}
\]
the corresponding base-change diagram.
From the representability of the functors
$\underline{\Mf}_{E/T}^{(n, d)}$ and $\underline{\Mf}_{E'/T'}^{(n, d)}$ it
follows easily that $M':= M\times_{T} T'$ with
$\kP' := (u \times \tilde{f})^*\kP$ represents the functor
$\underline{\Mf}_{E'/T'}^{(n, d)}$. Here we denoted by $u:M\times_{T} T' \lar
M$ the first projection and $u \times \tilde{f} : M'\times_{T'} E' \lar
M\times_{T} E$ coincides with
$u\times \id_{E}: M' \times_{T} E \lar M\times_{T} E$ under the canonical
identification $M'\times_{T'} E' = M' \times_{T} E $.

By $X', B'$ we denote the spaces obtained from $E'\lar T'$ in the same way as
$X, B$ were obtained from $E\lar T$.
Using $\hat{g} = u \times u\times\tilde{f}\times\tilde{f}\times\tilde{f}:
X' \lar X$ and
$\tilde{g}_{B}=u \times u \times \tilde{f} \times \tilde{f} :B' \lar B$,
we obtain the Cartesian diagram
\[
\xymatrix
{
 X' \ar[r]^{\hat{g}} \ar[d]_{q'} & X \ar[d]^q \\
B' \ar[r]^{\tilde{g}_{B}} & B.
}
\]
Note that there exist canonical isomorphisms
$\phi_i: \hat{g}^*\kV_i = \hat{g}^* \pi_{i}^{\ast} \kP
\lar \kV'_i := {\pi'_i}^*\kP'$, using the
notation of Definition \ref{D:univgeom}.

\begin{proposition}\label{P:chaba}
Let $\omega \in H^0(\omega_{E/T})$ and
$\omega' = \tilde{f}^*(\omega) \in H^0(\omega_{E'/T'})$, then
the image of the section $r = r^{\kV_1, \kV_2}_{h_1, h_2}(\omega) \in
H^0\bigl(h_1^* {\mathcal Hom}_X(\kV_2, \kV_1) \otimes
h_2^* {\mathcal Hom}_X(\kV_1, \kV_2)\bigr)$
under the chain of canonical morphisms
$$
\xymatrix{
H^0\bigl(h_1^* {\mathcal Hom}_X(\kV_2, \kV_1) \otimes
h_2^* {\mathcal Hom}_X(\kV_1, \kV_2)\bigr)
\ar[d]^{H^{0}(\tilde{g}^{\ast}_{B})} \\
H^0\bigl({\tilde g}_{B}^*(h_1^* {\mathcal Hom}_X(\kV_2, \kV_1) \otimes
h_2^* {\mathcal Hom}_X(\kV_1, \kV_2))\bigr)
\ar[d]^{\eta_{\kV_{1},\kV_{2}}} \\
H^0\bigl({h'_1}^* {\mathcal Hom}_{X'}(\hat{g}^*\kV_2, \hat{g}^*\kV_1) \otimes
{h'_2}^* {\mathcal Hom}_{X'}(\hat{g}^*\kV_1, \hat{g}^*\kV_2)\bigr)
\ar[d]^{H^{0}\bigl({h'_{1}}^{\ast}\conj(\phi_{2},\phi_{1}) \otimes
{h'_{2}}^{\ast}\conj(\phi_{1},\phi_{2})\bigr)} \\
H^0\bigl({h'_1}^* {\mathcal Hom}_{X'}(\kV'_2, \kV'_1) \otimes
{h'_2}^* {\mathcal Hom}_{X'}(\kV'_1, \kV'_2)\bigr)
}
$$
is $r' = r^{\kV'_1, \kV'_2}_{h'_1, h'_2}(\omega')$, where
the first arrow  is induced by the functor $\tilde{g}_{B}^*$, the second by the
canonical isomorphisms of functors
$\tilde{g}_{B}^* h_i^* \cong {h'_i}^* \hat{g}^*$
and the third by the isomorphisms of vector bundles
$\phi_i: \hat{g}^*\kV_i \lar \kV_i'$, $i = 1,2$.
\end{proposition}

\noindent
\emph{Proof}. This proposition
is an immediate consequence of Corollary \ref{C:changedebase}.
\qed

\begin{corollary}\label{C:chaba}
  In the notations as above, let $x =(v_1, v_2; y_1, y_2) \in B$ and
  $t = g(x) \in T$.
  Let $\omega \in H^0(\omega_{E/T})$ be a nowhere vanishing section and
  $\omega_t$ be its restriction to $E_t$. Then the image of the section
  $r= r(\omega)$ under the chain of canonical morphisms
  $$
  \xymatrix{
    H^0\bigl(B, h_1^*{\mathcal Hom}_X(\pi_2^*\kP, \pi_1^*\kP) \otimes
    h_2^* {\mathcal Hom}_X(\pi_1^*\kP, \pi_2^*\kP)\bigr) \ar[d]^{\mathsf{can}} \\
    H^0\bigl(B,  h_1^*{\mathcal Hom}_X(\pi_2^*\kP, \pi_1^*\kP) \otimes
    h_2^* {\mathcal Hom}_X(\pi_1^*\kP, \pi_2^*\kP)
    \otimes \CC_x\bigr) \ar[d]^\cong  \\
    \Hom_{\CC}(\kP^{v_2}|_{y_1}, \kP^{v_1}|_{y_1}) \otimes
    \Hom_{\CC}(\kP^{v_1}|_{y_2}, \kP^{v_2}|_{y_2})
  }
  $$
  is the tensor $r^{\kP^{v_1}, \kP^{v_2}}_{y_1, y_2}(\omega_t)$ obtained by
  the construction in Section \ref{S:PolConstr} on the curve $E_{t}$.
  In particular, the section $r$ is non-degenerate on
  $B\setminus\bar{\Delta}$. Equivalently, the morphism of vector bundles
  $$
  \tilde{r}(\omega):
  h_1^* {\mathcal Hom}_X(\pi_1^*\kP, \pi_2^*\kP) \lar
  h_2^* {\mathcal Hom}_X(\pi_1^*\kP, \pi_2^*\kP)
  $$
  is an isomorphism over $B\setminus\bar{\Delta}$.
\end{corollary}

\medskip
\begin{remark}
Since we assume  the fibration $E \stackrel{p}\lar T$ is \emph{algebraic}, the
above construction yields a \emph{meromorphic} section $r(\omega)$ of the
vector bundle $$h_1^*{\mathcal Hom}(\pi_2^*\kP, \pi_1^*\kP) \otimes
h_2^* {\mathcal Hom}(\pi_1^*\kP, \pi_2^*\kP)$$ over $M \times_T M \times_T
 \breve{E} \times_T \breve{E}$, which is holomorphic on
$B =  M \times_T M \times_T \breve{E} \times_T \breve{E}
\setminus (\Delta_M \cup \Delta_E)$ and non-degenerate on
$B \setminus \bar\Delta$, see Remark \ref{rem:ok-without-vanishing}.
\end{remark}
\medskip
\noindent
Our next goal is to show that the constructed canonical
section $r= r(\omega)$ satisfies a version of the
\emph{associative Yang--Baxter equation}.
For this purpose  we need further notation.
Let $$p^{ij}_{kl}: M \times_T  M \times_T  M \times_T  \breve{E} \times_T
\breve{E} \times_T \breve{E}    \lar M \times_T  M \times_T  \breve{E} \times_T
\breve{E}$$ be the projection
$
p^{ij}_{kl}(v_1, v_2, v_3; y_1, y_2, y_3) = (v_i, v_j; y_k, y_l),
$
where $1 \le i \ne j \le 3$ and $1 \le k \ne l \le 3$.
We also denote by
$$\hat{\pi}_j: M \times_T  M \times_T  M \times_T  \breve{E} \times_T
\breve{E} \times_T \breve{E}  \times_T E \lar M \times_T E$$
the projection given by the formula
$\hat{\pi}_j(v_1, v_2, v_3;  y_1, y_2, y_3, y) = (v_j; y)$,
where $1\le j\le 3$. Similarly, we have three canonical sections
$$\hat{h}_i: M \times_T  M \times_T  M \times_T  \breve{E} \times_T
\breve{E} \times_T \breve{E}
\lar M \times_T  M \times_T  M \times_T  \breve{E} \times_T
\breve{E} \times_T \breve{E}  \times_T E$$
given by the formulae
$\hat{h}_i(v_1, v_2, v_3, y_1, y_2, y_3) =
(v_1, v_2, v_3, y_1, y_2, y_3, y_i)$,  $1 \le i \le 3$.

The set over which the Yang-Baxter relation will be defined is
\[
\hat{B}:=\bigcap_{i,j,k,l} \left(p_{kl}^{ij}\right)^{-1}(B) =
\left(p_{12}^{12}\right)^{-1}(B) \cap \left(p_{13}^{13}\right)^{-1}(B)
\cap \left(p_{23}^{23}\right)^{-1}(B).
\]
Let $r  \in H^0\bigl(B, h_1^*{\mathcal Hom}(\pi_2^*\kP, \pi_1^*\kP) \otimes
h_2^* {\mathcal Hom}(\pi_1^*\kP, \pi_2^*\kP)\bigr)$
be the canonical holomorphic section constructed above.
Observe that
\begin{multline*}
  (p^{ij}_{kl})^*\bigl(h_1^*{\mathcal Hom}(\pi_2^*\kP, \pi_1^*\kP) \otimes
  h_2^* {\mathcal Hom}(\pi_1^*\kP, \pi_2^*\kP)\bigr)
  \\
  \cong\hat{h}_k^*{\mathcal Hom}(\hat{\pi}_j^*\kP, \hat{\pi}_i^*\kP) \otimes
  \hat{h}_l^* {\mathcal Hom}(\hat{\pi}_i^*\kP, \hat{\pi}_j^*\kP).
\end{multline*}
Let
$
r^{ij}_{kl} = \left(p^{ij}_{kl}\right)^*r
$
be the pull-back to $\hat{B}\subset M \times_T M \times_T M \times_T \breve{E}
\times_T \breve{E} \times_T \breve{E}$. Then we have:
\begin{itemize}
\item $r^{32}_{12}$ is a meromorphic section of
  $\hat{h}_1^*{\mathcal Hom}(\hat{\pi}_2^*\kP, \hat{\pi}_3^*\kP) \otimes
  \hat{h}_2^* {\mathcal Hom}(\hat{\pi}_3^*\kP, \hat{\pi}_2^*\kP)$.
\item $r^{13}_{13}$ is a meromorphic section of
  $\hat{h}_1^*{\mathcal Hom}(\hat{\pi}_3^*\kP, \hat{\pi}_1^*\kP) \otimes
  \hat{h}_3^* {\mathcal Hom}(\hat{\pi}_1^*\kP, \hat{\pi}_3^*\kP)$.
\end{itemize}
Taking their composition, we obtain a meromorphic section
$(r^{13}_{13})^{13} (r^{32}_{12})^{12}$ of the holomorphic vector bundle
$\hat{h}_1^*{\mathcal Hom}(\hat{\pi}_2^*\kP, \hat{\pi}_1^*\kP) \otimes
\hat{h}_2^* {\mathcal Hom}(\hat{\pi}_3^*\kP, \hat{\pi}_2^*\kP)\otimes
\hat{h}_3^* {\mathcal Hom}(\hat{\pi}_1^*\kP, \hat{\pi}_3^*\kP).
$
In a similar way, two other meromorphic sections
$(r^{12}_{12})^{12} (r^{13}_{23})^{23}$
and $(r^{23}_{23})^{23} (r^{12}_{13})^{13}$ of this  vector bundle can be defined.
These sections are holomorphic over $\hat{B}$.

Let $x = (v_1, v_2, v_3; y_1,  y_2, y_3)$ be a point in $\hat{B}\subset  M
\times_T M \times_T M \times_T \breve{E}  \times_T  \breve{E} \times_T
\breve{E}$ lying over the point $t \in T$.
Because $\hat{h}_k^*{\mathcal Hom}(\hat{\pi}_i^*\kP, \hat{\pi}_j^*\kP)\big|_x$
is canonically isomorphic to
$\Hom_{\CC}(\kP^{v_i}|_{y_k}, \kP^{v_j}|_{y_k})$, we may consider
$r^{ij}_{kl}(x)$ as an element of the tensor product of vector spaces
$\Hom_{\CC}(\kP^{v_j}|_{y_k}, \kP^{v_i}|_{y_k}) \otimes
\Hom_{\CC}(\kP^{v_i}|_{y_l}, \kP^{v_j}|_{y_l})$ and we have a canonical
isomorphism of vector spaces
\begin{multline*}
  \bigl(
  \hat{h}_1^*{\mathcal Hom}(\hat{\pi}_2^*\kP, \hat{\pi}_1^*\kP) \otimes
  \hat{h}_2^* {\mathcal Hom}(\hat{\pi}_3^*\kP, \hat{\pi}_2^*\kP)\otimes
  \hat{h}_3^* {\mathcal Hom}(\hat{\pi}_1^*\kP, \hat{\pi}_3^*\kP)\bigr) \big|_x
  \cong
  \\
  \cong \Hom_{\CC}(\kP^{v_2}|_{y_1}, \kP^{v_1}|_{y_1}) \otimes
  \Hom_{\CC}(\kP^{v_3}|_{y_2}, \kP^{v_2}|_{y_2}) \otimes
  \Hom_{\CC}(\kP^{v_1}|_{y_3}, \kP^{v_3}|_{y_3}).
\end{multline*}

\begin{definition}\label{def:YBR}
  Assume $E \stackrel{p}\lar T$ is a genus one fibration which satisfies the
  conditions set out at the beginning of Section \ref{S:GeomR} and fix $(n,d)$
  and $\omega$ as before.
  We call
  \begin{equation}\label{E:AYBEweak}
    (r^{13}_{13})^{13} (r^{32}_{12})^{12} -
    (r^{12}_{12})^{12} (r^{13}_{23})^{23} +
    (r_{23}^{23})^{23} (r^{12}_{13})^{13} = 0
  \end{equation}
  the \emph{Yang-Baxter relation}. The left-hand side of this equation
  is a holomorphic section of the vector bundle
  $
  \hat{h}_1^*{\mathcal Hom}(\hat{\pi}_2^*\kP, \hat{\pi}_1^*\kP) \otimes
  \hat{h}_2^* {\mathcal Hom}(\hat{\pi}_3^*\kP, \hat{\pi}_2^*\kP)\otimes
  \hat{h}_3^* {\mathcal Hom}(\hat{\pi}_1^*\kP, \hat{\pi}_3^*\kP)
  $
  over $\hat{B}$.

  Let $\tau$ be the  involution  of
  $M \times_T M \times_T \breve{E} \times_T \breve{E}$ which is given by
  $\tau(v_1, v_2; y_1, y_2) = (v_2, v_1; y_2, y_1)$. We say that
  $r\in H^0\bigl(B, h_1^*{\mathcal Hom}(\pi_2^*\kP, \pi_1^*\kP) \otimes
  h_2^*{\mathcal Hom}(\pi_1^*\kP, \pi_2^*\kP)\bigr)$ is \emph{unitary}, if
  \begin{equation}\label{E:unitarysect}
    r^{\pi_1^*\kP, \pi_2^*\kP}_{h_1, h_2}(\omega) = -
    \tau^*\left(r^{\pi_2^*\kP, \pi_1^*\kP}_{h_2, h_1}(\omega)\right).
  \end{equation}
  This means that the section $\tau^{\ast}(r)$ is mapped to $-r$ under the
  composition of the canonical isomorphisms
  $$
  \tau^*\bigl(h_1^*{\mathcal Hom}(\pi_2^*\kP, \pi_1^*\kP) \otimes
  h_2^* {\mathcal Hom}(\pi_1^*\kP, \pi_2^*\kP)\bigr) \cong
  $$
  $$
  h_2^*{\mathcal Hom}(\pi_1^*\kP, \pi_2^*\kP) \otimes
  h_1^* {\mathcal Hom}(\pi_2^*\kP, \pi_1^*\kP) \cong
  h_1^*{\mathcal Hom}(\pi_2^*\kP, \pi_1^*\kP) \otimes
  h_2^* {\mathcal Hom}(\pi_1^*\kP, \pi_2^*\kP).
  $$
\end{definition}

\noindent
The purpose of the following lemma is to use the relations
\eqref{E:yb2} and \eqref{E:ybunit}, which were shown for tensors $r^{\kV_{1},
  \kV_{2}}_{y_1, y_2}$  on smooth curves in Section \ref{S:PolConstr}, in
order to prove that $r$ satisfies the Yang-Baxter relation
\eqref{E:AYBEweak} and is unitary \eqref{E:unitarysect} in the case of
elliptic fibrations with arbitrary fibers.

\begin{lemma}\label{L:YBE-local-global}
  Assume $E \stackrel{p}\lar T$ is a genus one fibration which satisfies the
  conditions set out at the beginning of Section \ref{S:GeomR} and fix $(n,d)$
  and $\omega$ as before. Denote $r$ by $r_{T}$ if constructed from the family
  $E\rightarrow T$. Let $U\subset T$ be a dense subset. Use the
  restriction of $\omega$ and the same pair of integers $(n,d)$ to
  construct the section $r_{U}$ and $r_{t}$ from the induced families
  $E|_{U}\rightarrow U$ and $E_{t}\rightarrow \{t\}$.
  Then the following conditions are equivalent.
  \begin{enumerate}
  \item The Yang-Baxter relation \eqref{E:AYBEweak} holds for $r_{U}$.
  \item The Yang-Baxter relation \eqref{E:AYBEweak} holds for $r_{T}$.
  \item The Yang-Baxter relation \eqref{E:AYBEweak} holds for $r_{t}$ for all
    $t\in T$.
  \item The Yang-Baxter relation \eqref{E:AYBEweak} holds for $r_{t}$ for all
    $t\in U$.
  \end{enumerate}
  If $E_{t}$ is smooth, the Yang-Baxter relation \eqref{E:AYBEweak} holds for
  $r_{t}$.

  \vspace*{1mm}
  \noindent
  Similar statements hold for unitarity \eqref{E:unitarysect}.
\end{lemma}

\begin{proof}
  Let $B_{U}:=B\times_{T}U, B_{t}:=B\times_{T}\{t\}$ and define $\hat{B}_{U},
  \hat{B}_{t}$ in a similar way. Note that $\hat{B}_{U}\subset \hat{B}$ is dense.
  If we apply Proposition \ref{P:chaba} to the base-change
  $U\subset T$, we obtain that $r|_{B_{U}}$ corresponds (under a certain
  canonical isomorphism) to the section $r_{U}$ obtained from the family
  $E_{U}\rightarrow U$. Similarly, $r|_{B_{t}}$ corresponds to the section
  $r_{t}$ obtained from the family $E_{t}\rightarrow \{t\}$.

  Let us denote the left hand side of the relation \eqref{E:AYBEweak} by
  $R_{T}$, if it is a relation for $r_{T}$. Because the projections
  $p_{kl}^{ij}$ are compatibel with restrictions to the subsets $B_{U}\subset
  B$ and $\hat{B}_{U}\subset \hat{B}$, we obtain from the above that
  $R_{T}|_{\hat{B}_{U}}$ corresponds to $R_{U}$ under a certain canonical
  isomorphism. Similarly, $R_{T}|_{\hat{B}_{t}}$ corresponds to $R_{t}$ for
  each $t\in T$. In particular, $R_{T}|_{\hat{B}_{U}}$ vanishes if and only if
  $R_{U}$ does so and the vanishing of $R_{T}|_{\hat{B}_{t}}$ is equivalent to
  the vanishing of $R_{t}$.

  As $T$ is reduced, $R_{T}=0$ is equivalent to $R_{T}(x)=0$ for all
  $x\in\hat{B}$. Similar statements hold for $R_{U}$ and $R_{t}$. Because the
  zero locus of $R_{T}$, which is a section of a coherent sheaf on $\hat{B}$,
  is a closed subset of $\hat{B}$, the equivalence of the statements (1)--(4)
  is now obvious.

  Corollary \ref{C:chaba} says that the restriction
  $r(\omega)|_{(v_{1},v_{2}; y_{1},y_{2})}$ corresponds to
  $$
  r^{\kP^{v_1}, \kP^{v_2}}_{y_1, y_2}(\omega_t) \in
  \Hom_{\CC}(\kP^{v_2}|_{y_1}, \kP^{v_1}|_{y_1}) \otimes
  \Hom_{\CC}(\kP^{v_1}|_{y_2}, \kP^{v_2}|_{y_2})
  $$
  under a certain canonical isomorphism.
  For each $x=(v_{1},v_{2},v_{3}; y_{1},y_{2},y_{3}) \in \hat{B}_{t}$
  this implies that $R_{t}(x)=0$ is equivalent to \eqref{E:yb2} which was
  shown to be true in Section \ref{S:PolConstr} for smooth curves.

  The proofs for unitarity are similar.
\end{proof}

\begin{theorem}\label{T:main}
  \begin{itemize}
  \item[(a)] For each Weierstra\ss{} curve $E$, smooth or singular, the
    section $r\in H^0\bigl(B, h_1^*{\mathcal Hom}(\pi_2^*\kP, \pi_1^*\kP) \otimes
    h_2^*{\mathcal Hom}(\pi_1^*\kP, \pi_2^*\kP)\bigr)$ from Definition
    \ref{D:univgeom} for the constant family $E\rightarrow\Spec(\CC)$
    satisfies the Yang--Baxter relation (\ref{E:AYBEweak}) and the unitarity
    condition (\ref{E:unitarysect}) for each choice of $(n,d), \kP$ and
    $\omega$.
  \item[(b)] Let $E \stackrel{p}\lar T$ be a genus one fibration satisfying
    the conditions from the beginning of this section.
    Let $\omega \in H^0(\omega_{E/T})$ be a no\-where vanishing differential
    form, $(n, d)$ a pair of coprime integers and  $\kP = \kP(n, d)$ a
    universal family. Then, the  universal section
    $$r  = r(\omega) \in H^0\bigl(B, h_1^*{\mathcal
      Hom}(\pi_2^*\kP, \pi_1^*\kP) \otimes
    h_2^*{\mathcal Hom}(\pi_1^*\kP, \pi_2^*\kP)\bigr)$$ satisfies the
    Yang--Baxter relation (\ref{E:AYBEweak}) and the unitarity condition
    (\ref{E:unitarysect}).
    Moreover,  $r$  depends holomorphically (and in particular, continuously)
    on the parameter $t \in T$.
  \end{itemize}
\end{theorem}

\begin{proof}
(a)
Let $E_{T} \subset  {\mathbb P}^2 \times \CC^2 \lar \CC^2 =:T$
be the elliptic  fibration given by the equation
$zy^2 = 4 x^3 -  g_2 xz^2 - g_3 z^3$ and let
$\Delta(g_2, g_3)  = g_2^3 - 27 g_3^2$ be the discriminant of this family.
This fibration has a section
$
(g_2, g_3) \mapsto \bigl((0:1:0), (g_2, g_3)\bigr)
$
and satisfies the condition $\omega_{E/T} \cong \kO_E$.
Let $\omega \in H^0(\omega_{E/T})$ be a no\-where vanishing differential form.

There exists $t\in T$, such that the given curve $E$ is isomorphic to the
fibre $E_{t}$. The chosen differential form on $E$ coincides with
$\omega_{t}$ up to a constant factor. The restriction $\kP|_{M_t \times E_t}$
of a universal family $\kP\in \VB(M \times_{T} E)$ is a universal family of
stable vector bundles of rank $n$ and degree $d$ on the curve $E_t$.

Using the open dense subset $U=T\setminus \Delta\subset T$ in
the equivalence of the statements (3) and (4) in Lemma
\ref{L:YBE-local-global} and the fact that \eqref{E:AYBEweak} and
\eqref{E:unitarysect} are satisfied for smooth fibres by Lemma
\ref{L:YBE-local-global}, we obtain the claim.

\noindent
(b) Because each fibre of $E\rightarrow T$ is isomorphic to a Weierstra\ss{}
curve, part (a) and Lemma \ref{L:YBE-local-global} (2) and (3) imply the claim.
\end{proof}

\subsection{Passing to Matrices}

Our next goal is to pass from the categorical version of the associative
Yang--Baxter equation (\ref{E:AYBEweak}) to the one which was studied in
Section \ref{S:YBE}.
Our construction is based on the choice of a trivialization of the universal
bundle $\kP$ and on the choice of local coordinates on $M$ and $E$. These two
choices can be made independently.

Let $o = (m_0, e_0) \in M \times_T \breve{E}\subset M \times_T E$ be an
arbitrary point, which lies over $t_0 \in T$.
Consider a small open neighbourhood $V \subset M\times_T \breve{E}$ of the
point $o$ such that there exists an isomorphism of vector bundles
$\xi: \kP|_V \stackrel{\cong}\lar  V \times \CC^n$.

Let $\varphi_{ij} := \pi_j \circ h_i:
M \times_T M \times_T \breve{E} \times_T \breve{E} \lar M \times_T E$,
$B_0 := \bigcap_{i,j =1}^2 \varphi^{-1}_{ij}(V)$ and $O = \kO_B(B_0\cap B)$ be
the ring of holomorphic functions on $B_0\cap B$.
The isomorphism  $\xi$ induces trivializations
$\varphi_{ij}^*(\xi): \varphi_{ij}^*\kP|_{B_0} \lar B_0 \times \CC^n$ from
which we obtain isomorphisms
$H^0\bigl(B_0\cap B, h_k^*{\mathcal Hom}(\pi_i^*\kP, \pi_j^*\kP)\bigr)
\stackrel{\cong}\lar \Mat_{n \times n}(O)$.
Under the induced isomorphism
$$
H^0\bigl(B_0\cap B, h_1^*{\mathcal Hom}(\pi_2^*\kP, \pi_1^*\kP) \otimes
h_2^* {\mathcal Hom}(\pi_1^*\kP, \pi_2^*\kP)\bigr)
\stackrel{\cong}\lar \Mat_{n \times n}(O) \otimes_O \Mat_{n\times n}(O),
$$
the section $r$  is mapped to a tensor
$$
r^\xi  = r^\xi(v_1, v_2; y_1, y_2) =
\sum_{\nu}  a_{\nu}(v_1, v_2; y_1, y_2) \otimes
b_{\nu}(v_1, v_2; y_1, y_2)
$$
in $\Mat_{n\times n}(O) \otimes_O \Mat_{n\times n}(O)$,
where
$a_\nu = a_{\nu}(v_1, v_2; y_1, y_2)$ and $b_\nu  = b_{\nu}(v_1, v_2; y_1, y_2)$
are holomorphic functions on $B_0\cap B$.

Because the fibration $p:E\rightarrow T$ is smooth at $e_{0}$ and so also
$M\rightarrow T$ at $m_{0}$, there exist an open neighbourhood $T_{0}$ of
$t_{0}\in T$ and an open disc $\DD\subset \CC$ such that there are
open neighbourhoods of $e_{0}\subset E$ and of $m_{0}\subset M$ which are
isomorphic to $T_{0}\times\DD$ and such that the following diagrams of complex
spaces are commutative:
$$
\xymatrix{
E    \ar[d]
&  T_0 \times \DD \ar@{_{(}->}[l] \ar[d]^{\mathsf{pr}}
&&
M   \ar[d]
&  T_0 \times \DD \ar@{_{(}->}[l] \ar[d]^{\mathsf{pr}}
\\
T
&  T_0 \ar@{_{(}->}[l]
&&
T
&  T_0. \ar@{_{(}->}[l]
}
$$
We assume that $V$ is isomorphic to the fibred product of open
neighbourhoods of the form $T_{0}\times\DD$, so that
$V\cong T_{0}\times\DD^{2}$.
With such $V$ we obtain an isomorphism $B_0 \cong T_{0}\times\DD^{4}$ and the
tensor $r^\xi  = r^\xi(v_1, v_2; y_1, y_2)$ will be written in such
coordinates as $r^\xi(t; v_1, v_2; y_1, y_2)$ with $t\in T_{0}$ and $v_1, v_2,
y_1, y_2 \in \DD$.
We also define
$\hat{B}_{0}:=\bigcap_{i,j,k,l} \left(p_{kl}^{ij}\right)^{-1}(B_{0})$
and obtain an isomorphism $\hat{B}_{0}\cong T_{0}\times\DD^{6}$.
In these coordinates, we have
$p_{kl}^{ij}(t; v_{1},v_{2},v_{3}; y_{1},y_{2},y_{3}) =
(t;v_{i},v_{j};y_{k},y_{l})$. This equation implies that the Yang-Baxter relation
\eqref{E:AYBEweak} and unitarity \eqref{E:unitarysect} translate into
\eqref{E:AYBEunitary} and \eqref{E:AYBE1} respectively. Therefore,
Theorem \ref{T:main} (b) implies the following corollary.

\begin{corollary}
Let $E \stackrel{p}\lar T$ be an elliptic fibration satisfying all the
conditions from the beginning of this section.
Let $\omega \in H^0(\omega_{E/T})$ be a nowhere vanishing section, $(n, d)$ be
coprime integers, $M = M_{E/T}^{(n,d)} \stackrel{l}\lar T$ be the moduli space
of relatively stable vector bundles of rank $n$ and degree $d$,
$\kP = \kP(n, d) \in \VB(M \times_T E)$ be a universal family.
Let $o = (m_0, e_0) \in M \times_T \breve{E}$ be an
arbitrary point lying over $t_0 \in T$ and choose coordinates around $e_{0}$
and $m_{0}$ as described above.
Finally, let $\xi: \kP|_{V} \lar \kO_V^n$ be a trivialization of the universal
family over a neighbourhood $V$ of the point $o$.

Then, the tensor $r^\xi(t; v_1, v_2; y_1, y_2)$ is unitary
(i.e.\ it fulfils \eqref{E:AYBEunitary})  and satisfies the associative
Yang-Baxter equation \eqref{E:AYBE1}.
\end{corollary}

Our next goal is to explain how the tensor $r^\xi$ depends on the choice of
the trivialization $\xi$. We do not need to choose coordinates here.
If $\kP|_V \stackrel{\zeta}\lar  V \times \CC^n$ is another trivialization
of $\kP$, we obtain a commutative diagram
$$
\xymatrix
{
 & \kP|_V \ar[ld]_{\xi} \ar[rd]^{\zeta} & \\
V \times \CC^n \ar[rr]^{\id \times \phi}& & V \times\CC^n,
}
$$
where $\phi = \phi(v, y): V \lar \GL_n(\CC)$ is a holomorphic function.

\begin{proposition}\label{P:gauge}
The solutions $r^\xi$ and $r^\zeta$ are gauge equivalent and such an
equivalence is given by the function $\phi$.
In other words, \emph{gauge transformations} of the solutions of the
associative Yang--Baxter equation, which are obtained from a geometric
associative $r$-matrix, correspond exactly to a
\emph{change of trivialization} of the universal family $\kP$.
\end{proposition}

\begin{proof}
With respect to the second trivialization $\zeta$, the section $r$ can be
written as a tensor
$$
r^\zeta  = r^\zeta(v_1, v_2; y_1, y_2) =
\sum_{\nu}  a'_{\nu}(v_1, v_2; y_1, y_2) \otimes
b'_{\nu}(v_1, v_2; y_1, y_2).
$$
The functions $a_\nu$ and $a'_\nu$ are related by the following commutative
diagram:
$$
\xymatrix@R-1mm@C-3mm
{
\CC^n \ar[rrr]^{a_\nu} \ar@/_30pt/[dd]_{\phi(v_2, y_1)} & & &  \CC^n
\ar@/^30pt/[dd]^{\phi(v_1, y_1)}\\
\kP^{v_2}|_{y_1} \ar[rrr] \ar[d]^{\zeta(v_2, y_1)}
\ar[u]_{\xi(v_2, y_1)} & & & \kP^{v_1}|_{y_1}
\ar[d]_{\zeta(v_1, y_1)}
\ar[u]^{\xi(v_1, y_1)}\\
\CC^n \ar[rrr]^{a'_\nu} & & &  \CC^n. \\
}
$$
Similarly for $b_\nu$ and $b'_\nu$ we have:
$$
\xymatrix@R-1mm@C-3mm
{
\CC^n \ar[rrr]^{b_\nu} \ar@/_30pt/[dd]_{\phi(v_1, y_2)} & & &  \CC^n
\ar@/^30pt/[dd]^{\phi(v_2, y_2)}\\
\kP^{v_1}|_{y_2} \ar[rrr] \ar[d]^{\zeta(v_1, y_2)}
\ar[u]_{\xi(v_1, y_2)} & & & \kP^{v_2}|_{y_2}
\ar[d]_{\zeta(v_2, y_2)}
\ar[u]^{\xi(v_2, y_2)}\\
\CC^n \ar[rrr]^{b'_\nu} & & &  \CC^n. \\
}
$$
These diagrams  imply   the following transformation rules:
$$
\begin{array}{lcl}
a'_\nu & = &  \phi(v_1, y_1)\; a_\nu\; \phi^{-1}(v_2, y_1), \\
b'_\nu  & =  & \phi(v_2, y_2)\; b_\nu\;  \phi^{-1}(v_1, y_2).
\end{array}
$$
This means that the choice of a different trivialization $\zeta$
of the universal bundle $\kP$ leads to the transformation rule
$$
r^\xi \mapsto
\bigl(\phi(v_1, y_1) \otimes \phi(v_2, y_2)\bigr)  \; r^\xi \;
\bigl(\phi(v_2, y_1)^{-1} \otimes \phi(v_1, y_2)^{-1}\bigr) = r^\zeta,
$$
which means that $r^{\xi}$ and $r^{\zeta}$ are gauge equivalent.
Conversely, if we start with $r^{\xi}$ and apply a gauge transformation $\phi$,
the same calculation shows that the result is $r^{\zeta}$ where $\zeta :=
(\id\times\phi)\circ\xi$.
\end{proof}

\medskip
\noindent
Summing up all results of this section, we get the following theorem,
which is one  of the main results of this  article.

\medskip
\begin{theorem}\label{C:MainCorol}
Let $E \stackrel{p}\lar T$ be a genus one fibration satisfying the conditions
from the beginning of this section.
$$
% duzze garna kartynka
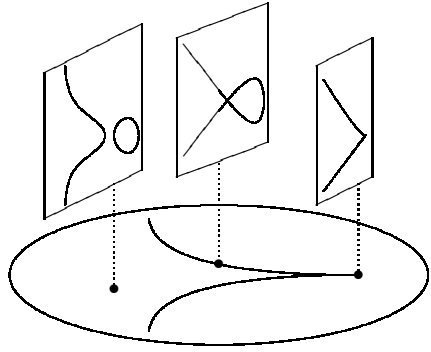
{\xy 0;/r0.22pc/:
\POS(40,0);
{(25,0)\ellipse(30,10){-}},%{\ellipse(,.75){}}
%
%elliptic curve
{(12,20)\ellipse(1.5,2.5)_,=:a(180){-}},
{(12,20)\ellipse(2,2.5)^,=:a(180){-}},
\POS(3,30);
@(,
\POS(3,25)@+,  \POS(8,22)@+, \POS(9,20)@+,
\POS(8,18)@+, \POS(3,15)@+, \POS(3,10)@+,
**\crvs{-}
,@i @);
{\ar@{-}(0,8);(0,29)},
{\ar@{-}(0,8);(14,15)},
{\ar@{-}(0,29);(14,36)},
{\ar@{-}(14,15);(14,36)},
%nodal curve
\POS(20,33);@-
{%(28,28)@+,
(20,30)@+,  (30,20)@+, (32,25)@+
,(30,30)@+, (20,20)@+, \POS(20,17)@-,
,**\qspline{}};  %@)
{\ar@{-}(19,14);(19,34)},
{\ar@{-}(19,14);(32,19)},
{\ar@{-}(19,34);(32,39)},
{\ar@{-}(32,19);(32,39)},
%cusp curve
\POS(40,28); \POS(46,20);
**\crv{(40,28);
(30,25); (42,22); (45,20) };
\POS(40,12); \POS(46,20);
**\crv{(40,18);
(45,5);   (42,8); (45,20); (46,20)};
%\POS(46.5,20)*{\bul}; \POS(48,20);
{\ar@{-}(39,10);(39,30)},
{\ar@{-}(39,10);(47,14)},
{\ar@{-}(39,30);(47,34)},
{\ar@{-}(47,14);(47,34)},
%%
%%cusp base
\POS(15,8); \POS(45,0)*{\bul} ;
**\crv{(16,0)};
\POS(15,-8); \POS(45,0);
**\crv{(16,0)};
\POS(15,8);
{\ar@{.}(45,13); (45,0) }
\POS(25,1.5)*{\bul};
{\ar@{.}(25,2); (25,16) }
\POS(10,-2)*{\bul};
{\ar@{.}(10,-2); (10,13) }
\endxy}
$$
Let $M = M_{E/T}^{(n,d)} \stackrel{l}\lar  T$  be the moduli space of relatively
stable vector bundles of rank $n$ and degree $d$
and $\kP = \kP(n, d) \in \VB(M \times_T E)$ be a
universal family on $M$.
We fix a nowhere vanishing differential form $\omega \in H^0(\omega_{E/T})$.
Let $e_{0}\in\breve{E}$ and $m_{0}\in M$ be arbitrary points and choose
coordinate neighbourhoods of the form $T_{0}\times\DD$ around these two
points.
Let $\xi$ be a trivialization of $\kP$ in the corresponding neighbourhood of
$o:=(m_{0}, e_{0})$.
Let $\hat{o} = (m_{0}, m_{0}, e_{0}, e_{0})
\in M \times_T M \times_T \breve{E} \times_T \breve{E}
\cong T_{0}\times \DD^{4}$.

\noindent
Then, we get the  germ of a meromorphic function
\[
r^\xi = \bigl(r^{(n,d)}_{E/T}(\omega)\bigr)^{\xi}:
\bigl(M\times_T  M\times_T\breve{E}\times_T\breve{E}, \hat{o}\bigr)
\lar \Mat_{n \times n}(\CC) \otimes \Mat_{n \times n}(\CC)
\]
which satisfies  the associative Yang-Baxter equation
\[
r^\xi(t; v_1, v_2; y_1, y_2)^{12} r^\xi(t; v_1, v_3; y_2, y_3)^{23} =
\]
\[
r^\xi(t; v_1, v_3; y_1, y_3)^{13} r^\xi(t; v_3, v_2; y_1, y_2)^{12} +
r^\xi(t; v_2, v_3; y_2, y_3)^{23} r^\xi(t; v_1, v_2; y_1, y_3)^{13}
\]
and its ``dual''
\clearpage

\[
r^\xi(t; v_2, v_3; y_1, y_2)^{23} r^\xi(t; v_1, v_3; y_1, y_2)^{12} =
\]
\[
r^\xi(t; v_1, v_2; y_1, y_2)^{12} r^\xi(t; v_2, v_3; y_1, y_3)^{13} +
r^\xi(t; v_1, v_3; y_1, y_3)^{13} r^\xi(t; v_2, v_1; y_2, y_3)^{23}.
\]
Moreover, it fulfills  the  unitarity condition
\[
r^\xi(t; v_1, v_2; y_1, y_2) = - \tau\bigl(r^\xi(t; v_2, v_1; y_2, y_1) \bigr),
\]
where $\tau(a \otimes b) = b \otimes a$.
The function $r^\xi(t; v_1, v_2; y_1, y_2)$ depends analytically on the
parameter $t \in T$ and its poles lie on the hypersurfaces $v_1 = v_2$ and
$y_1 = y_2$.

\noindent
Different choices of trivializations of the universal family $\kP$
lead to equivalent solutions: if $\zeta$ is another trivialization of $\kP$ and
$\phi = \zeta \, \circ \,  \xi^{-1}: (M \times_T E, o) \lar \GL_n(\CC)$ is the
corresponding holomorphic function, then we have:
$$
r^\zeta =
\bigl(\phi(t, v_1, y_1) \otimes \phi(t, v_2, y_2)\bigr)
 r^\xi\bigl(\phi(t, v_2, y_1)^{-1}
\otimes \phi(t, v_1, y_2)^{-1}\bigr).
$$
Finally, let $T' \stackrel{g}\lar  T$ be an arbitrary base change. Let $E' = E
\times_T T' \stackrel{p'}\lar T'$ be the induced genus one fibration. Then the
corresponding moduli space $M' = M_{E'/T'}^{(n,d)}$ of relatively stable
vector bundles of rank $n$ and degree $d$ is isomorphic to $M \times_T T'$.
Let $\omega' \in H^0(\omega_{E'/T'})$ be the pull-back of $\omega$
and $\xi'$ be the induced trivialization of the pull-back
$\kP' \in \VB(M' \times_{T'} E')$ of the universal family $\kP$.
Let $e_{0}'\in \breve{E}'$ and $m_{0}'\in M'$ be points which are mapped to
$e_{0}$ and $m_{0}$ respectively.
After choosing coordinate neighbourhoods $T_{0}\times \DD$ of $e_{0}\in
E$ and $m_{0}\in M$, we have induced neighbourhoods $T'_{0}\times \DD$
of $e_{0}'\in E'$ and $m_{0}'\in M'$ so that the morphism
between $E'$ and $E$ (and between $M'$ and $M$) is given by $g\times\id_{\DD}$.
Then we have
\[
r^\xi(g(t); v_1, v_2; y_1, y_2) = r^{\xi'}(t; v_1, v_2; y_1, y_2)
\]
for all $t\in T_{0}'$ and all $v_1, v_2, y_1, y_2\in\DD$.
In other words, the tensor $r^\xi(t; v_1, v_2; y_1, y_2)$ is compatible with
base change in the variable $t$.
\end{theorem}

\subsection{Comment on reducible curves}

The developed theory of the geometric $r$-matrices can be generalized
literally to the case of reduced but reducible curves of arithmetic genus one
with trivial dualizing sheaf.
In this subsection we discuss some necessary technical results which are not
yet available.

\medskip
\noindent
Throughout this  section, $E$ is  a reduced  projective curve with trivial
dualizing sheaf.

\begin{itemize}
\item If $E$ is smooth then it is isomorphic to an elliptic curve.
\item Assume  $E$ is singular  with singularities of
  embedding dimension \emph{equal to two}.  Then Kodaira's classification of
  degenerations of elliptic curves implies that $E$ is either a cycle of $m$
  projective lines (type I$_m$), a cuspidal cubic curve (type II), a tachnode
  curve (type III) or a configuration of three concurrent lines in a plane
  (type IV).
\item  However, the class of reduced genus one curve with trivial dualizing
  sheaf is larger. For example, a generic configuration of $m$ concurrent
  lines in $\PP^{m-1}$ for $m \ge 4$ is such a curve.
\end{itemize}
Let $\pi: \widetilde{E} \rightarrow E$ be the normalization of $E$ and
$\widetilde{E} = \widetilde{E}_1 \cup \dots \cup \widetilde{E}_m$ be the
decomposition into irreducible components. If $E$ was not smooth, then
$\widetilde{E}_i \cong \PP^{1}$ for all $i$.
For a vector bundle $\kF$ on $E$ we define its \emph{multi-degree} to be
\[
\underline{\deg}(\kF) = \bigl(d_1(\kF), \dots, d_m(\kF)\bigr) \in \ZZ^m
\]
where $d_i(\kF) = \deg\bigl(\pi^*(\kF)\big|_{\widetilde{E}_i}\bigr)$. The
following lemma can be shown in the same way as Lemma \ref{L:deg-of-vb}.

\begin{lemma}
  Let $\kF$ be a vector bundle on $E$. Then we have:
  \[
  \chi(\kF) = h^0(\kF) - h^{1}(\kF) = \deg(\kF) := d_1(\kF) + \dots + d_m(\kF).
  \]
\end{lemma}

\noindent
For  fixed $n \in \ZZ_+$ and $\mathbbm{d} \in \ZZ^m$  denote by
$\Spl^{(n, \mathbbm{d})}(E)$ the category of simple vector bundles of rank $n$
and multi-degree $\mathbbm{d}$. Then we have the following result.

\begin{theorem}\label{T:SimplBdl}
  Let $E$ be a reduced plane cubic curve.
  \begin{itemize}
  \item If  $\kF$ is  a simple vector bundle on $E$ then
    $\gcd\bigl(\rk(\kF), \chi(\kF)\bigr) = 1$;
  \item If $\mathbbm{d} = (d_1, d_2, \dots, d_m)$ satisfies
    $\gcd(n, d_1 + \dots + d_m) = 1$, then the functor
    \[\det: \Spl^{(n, \mathbbm{d})}(E) \lar \Pic^{\mathbbm{d}}(E)\]
    is an equivalence of categories.
    In particular, for any  pair $\kF', \kF'' \in \Spl^{(n, \mathbbm{d})}(E)$
    such that $\kF' \not\cong \kF''$ we have $\Hom_E(\kF', \kF'') = 0$.
  \end{itemize}
\end{theorem}

\noindent
The case of a smooth curve $E$ is due to Atiyah  \cite{Atiyah}. A proof for
the irreducible  singular Weierstra\ss{} cubic curves  can be found in
\cite{Burban1} (nodal case) and \cite{BodnarchukDrozd} (cuspidal case), see
also Section \ref{S:triples}.
The remaining cases, i.e.\ the Kodaira fibers of type I$_2$, I$_3$, III and
IV,  were considered in a recent paper of Bodnarchuk, Drozd and Greuel
\cite{BodDroGre}.

\begin{conjecture}\label{C:simplonCYcurve} Let $E$ be an arbitrary reduced
  projective curve with trivial dualizing sheaf and $m$ irreducible
  components. Then we have:
  \begin{enumerate}
  \item The description of simple vector bundles on $E$ given in Theorem
    \ref{T:SimplBdl} remains true:
    the rank and degree of a simple vector bundle are coprime;
    a simple vector bundle is determined by its rank and determinant;
    for given  $n$ and $\mathbbm{d} = (d_1, \dots, d_m)$  satisfying
    $\gcd(n, d_1 + \dots + d_m) =1$, the category $\Spl^{(n, \mathbbm{d})}(E)$
    is equivalent to $\Pic^{\mathbbm{d}}(E)$, in particular, it is non-empty.
  \item For any pair $(n, \mathbbm{d}) \in \ZZ_+ \times \ZZ^m$ as above,
    there exists an auto-equivalence
    $\FF \in \bigl\langle T_{\kO}, \,  \Pic(E), \, [1]\bigr\rangle$
    of the derived category $D^b\bigl(\Coh(E)\bigr)$ inducing an equivalence
    between $\Spl^{(n, \mathbbm{d})}(E)$ and the category of torsion
    sheaves of length one supported at the regular part of a single irreducible
    component of $E$.
  \item
    Consider the functor
    $\underline{\Mf}^{(n,\mathbbm{d})}_{E}: \Ans \lar  \Sets$ given by
    \[
    \underline{\Mf}^{(n,\mathbbm{d})}_{E}(T) =
    \left\{
      \kF \in \VB(E\times T)\left|
        \kF|_{E \times \{t\}} \in \Spl^{(n,\mathbbm{d})}(E)
        \quad \mbox{for all} \quad  t \in T
      \right.\right\}\Big/\sim
    \]
    where $\kF_1 \sim \kF_2$ if and only if there exists $\kL \in \Pic(T)$
    such that $\kF_1 \cong \kF_2 \otimes \pr_T^{*}(\kL)$.
    Then $\underline{\Mf}^{n, \mathbbm{d}}_E$ is isomorphic to
    $\Hom_{\Ans}(\;-\;,G)$, where $G = \CC^*$ for a cycle of projective lines
    and $G = \CC$ in the other cases.
  \end{enumerate}
\end{conjecture}

\noindent
A proof of the first part of this conjecture in the case of  cycles  of
projective lines was recently announced by Bodnarchuk and Drozd. Note that in
this case, a description of simple vector bundles in terms of \'etale
coverings is also known  \cite{OldSurvey}.

\medskip
\noindent
Let $E_{*}$ be an irreducible component of the curve $E$.
Provided a universal family $\kP(n, \mathbbm{d}) \in \VB(E \times G)$ of
simple vector bundles of rank $n$ and multi-degree $\mathbbm{d}$ does exist,
one can proceed in the same way as in the present section to end up with a
solution $r^{(n, \mathbbm{d}, *)}_E$ of the associative Yang--Baxter equation
with all the properties of Theorem \ref{C:MainCorol}.

\medskip
Conjecture \ref{C:simplonCYcurve} is closely related to the study of
moduli spaces of Simpson stable sheaves on degenerations of elliptic
curves, see recent papers  by Hern\'andez Ruip\'erez, L\'opez Mart\'in,
S\'anchez G\'omez and  Tejero Prieto \cite{Salamanca}, and Lowrey
\cite{Lowrey}.

\begin{conjecture}
  Let $E$ be as in Conjecture \ref{C:simplonCYcurve} and
  $(n, \mathbbm{d}) \in \ZZ_+ \times \ZZ^{m}$ be such that
  $\gcd(n, d_1 + \dots + d_m) = 1$.
  Then there exists a polarization $H$ of the curve $E$ and a Hilbert
  polynomial $p(t) \in \mathbb{Q}[t]$ such that all simple vector bundles of
  rank $n$ and multi-degree $\mathbbm{d}$ are Simpson-stable with Hilbert
  polynomial equal to $p$. Moreover, the moduli space $M_E(n, \mathbbm{d})$ is
  an open and dense subset of an irreducible component of the moduli space
  $M(p)$ of Simpson-stable sheaves with Hilbert polynomial $p$.
\end{conjecture}

Since the exact combinatorics of simple vector bundles on reduced
projective curves with trivial dualizing sheaf is still to be clarified and
their relationship to Simpson-stable sheaves is not completely clear, we
postpone a discussion of possible generalizations of the results of
Section \ref{S:GeomR} to a future publication.

\section{Action of the Jacobian and geometric associative $r$--matrices}
\label{S:GeomYB}

Let $E = V(zy^2 - 4 x^3 + g_2 xz^2 + g_3 z^3) \subseteq \mathbb{P}^2$ be a
Weierstra\ss{} cubic curve and denote by $\breve{E}$ the regular part of $E$.
Let $e\in\breve{E}$ be any point.
In this section we are dealing with a single curve, not with a family of
curves.

Fix a pair of coprime integers $(n, d) \in \mathbb{N} \times \mathbb{Z}$ and
let $(M,\kP)$ be a pair which represents the moduli functor
$\underline{\Mf}_E^{(n,d)}$.
In the previous section, we have shown how to construct a tensor
$r^\xi(v_1, v_2; y_1, y_2)$ satisfying the associative Yang-Baxter equation.
In order to do so, we had to choose a point $m\in M$, a trivialization
$\xi$ of $\kP$ over a neighbourhood of $(m,e)$ and coordinates around
$e\in E$ and around $m\in M$.

The main goal of this section is to show that it is possible to choose
coordinates on $M$ such that the associative $r$-matrix
$r^\xi$ is gauge equivalent to a solution $r^\zeta$ depending only on the
difference  $v = v_2 - v_1$ of the ``vector bundle'' spectral parameters.
More precisely, we are going to prove that there are coordinates on
$M$ and there exists a gauge transformation
$\phi: \bigl(\CC^2_{(v,y)}, 0\bigr) \to \GL_n(\CC)$  such that the function
\[
r^\zeta(v_1, v_2; y_1, y_2) =
\bigl(\phi(v_1, y_1) \otimes \phi(v_2, y_2)\bigr) r^\xi(v_1, v_2; y_1, y_2)
\bigl(\phi(v_2, y_1)^{-1} \otimes \phi(v_1, y_2)^{-1}\bigr)
\]
is invariant under transformations $(v_1, v_2) \mapsto (v_1 + v, v_2 + v)$.
In other words, there exists a trivialization $\zeta$ of the universal
family $\kP$ in a neighbourhood of the point $(m, e) \in M \times E$ such that
we have
\[
r^\zeta(v_1+v, v_2+v; y_1, y_2) = r^\zeta(v_1, v_2; y_1, y_2).
\]
For simplicity of notation, we shall  write  $r^\zeta(v_1, v_2; y_1, y_2) =
r^\zeta(v_1 - v_2; y_1, y_2) = r^\zeta(v; y_1, y_2)$, where $v = v_1 - v_2$.
As it was explained in Section \ref{S:YBE}, the tensor $r^\zeta(v; y_1, y_2)$
satisfies  the quantum Yang--Baxter equation and defines interesting first
order differential operators.

The key idea to find such a distinguished trivialization $\zeta$  is to study
the behaviour of the geometric $r$-matrix under the
\emph{action of the Jacobian} $J$ on $M$.
The coordinates on $M$ are obtained from an isomorphism $M\lar J$ and a
surjective homomorphism of groups $\CC\lar J$.
Using the isomorphism $M\lar J$ only, allows to give a coordinate free
description of the ``dependence on the difference of the $v_{i}$'', because
$J$ has the structure of a group. The coordinates are used because they link
the general discussion of this section with the explicit calculations of the
subsequent sections.

The functors $\underline{\Mf}^{(n,d)}_E$, $\underline{\Pic}_E^0$ and
$\underline{\Pic}_E^d$ are  representable. Let $(M,\kP)$, $(J,\kL)$ and
$(J^{d},\kL^{(d)})$ be spaces with universal families which represent these
functors, so that we have isomorphisms of functors as follows
\begin{gather*}
  \alpha: \underline{\Mf}^{(n,d)}_E   \lar \Mor(\, - \, ,M),\\
  \beta:   \underline{\Pic}_E^0  \lar \Mor(\, - \, ,J)
  \qquad \text{ and }\qquad
  \beta^{d}:   \underline{\Pic}_E^d  \lar \Mor(\, - \, ,J^{d}).
\end{gather*}
Recall that the product $F \times G: \Ans  \lar \Sets$ of
two functors $F,G: \Ans  \lar \Sets$ is defined by
$
(F \times G)(S) = F(S) \times G(S).
$
Because the tensor product of a stable vector bundle of rank $n$ and degree
$d$ with a line bundle of degree zero is again a stable vector bundle of rank
$n$ and degree $d$, the Jacobian acts on the moduli space. On the functorial
level, this action is described as a natural transformation of functors
$$
\underline{\tau}: \underline{\Pic}_E^{0} \times  \underline{\Mf}^{(n,d)}_E \lar
\underline{\Mf}^{(n,d)}_E,
$$
which is defined as $\underline{\tau}_S(\kN, \kF) = \kN \otimes \kF$ for any
complex space $S$, any line bundle $\kN \in \underline{\Pic}_E^{0}(S)$ and any
vector bundle $\kF \in \underline{\Mf}^{(n,d)}_E(S)$.
The natural transformation $\underline{\tau}$ induces a morphism  of complex
spaces $\tau: J \times M \lar M$ making the following  diagram commutative
\[
\xymatrix
{
\underline{\Pic}_E^0 \times \underline{\Mf}_E^{(r,d)}
\ar[rr]^(.4){\beta \times  \alpha}
\ar[d]_{\underline{\tau}}
& & \Mor(\, - \, ,J) \times \Mor(\,-\,,M)
\ar[rr] & & \Mor(\, - \, ,J \times M) \ar[d]^{\tau_{*}}\\
\underline{\Mf}_E^{(r,d)} \ar[rrrr]^{\alpha} & & & & \Mor(\,-\,,M).
}
\]
The unlabelled horizontal arrow is the isomorphism which results from the
universal property of the product of two complex spaces.
The morphism $J \times M \stackrel{\tau}\lar M$  corresponds to the natural
transformation $\tau_{\ast}$ using Yoneda's Lemma
\[
\Hom\bigl(\Mor(\, - \, ,J \times M),\,  \Mor(\,-\,,M)\bigr) \cong
\Mor(J \times M, M).
\]
The following lemma describes the equivalence class of vector bundles on
$J\times M$ which corresponds to the morphism $\tau\in\Mor(J\times M, M)$
under $\alpha_{J\times M}$.

\begin{lemma}\label{L:IsomLP}
Denote $\tilde\tau:=\tau\times\id_{E}:J \times M \times E \lar M \times E$
and let $p, q, t$ be the natural projections
\[
\xymatrix
{
& J \times M \times E \ar[ld]_{p} \ar[d]_{q} \ar[rd]^{t}& \\
M \times E & J \times E& J \times M
}
\]
Then $p^* \kP \otimes q^* \kL \sim \tilde\tau^*\kP$, i.e.\ there exists
a line bundle $\kN \in \Pic(J \times M)$ such that
\[
p^*\kP \otimes q^*\kL  \otimes t^*\kN \cong \tilde\tau^*\kP.
\]
\end{lemma}

\begin{proof} Note that we have the commutative diagram
\[
\xymatrix@C-=-5mm{
&  \Mor(J\times M, J) \times \Mor(J\times M, M) \ar[rd]^\can & \\
\underline{\Pic}_E^0(J\times M)  \times \underline{\Mf}_E^{(n,d)}(J\times M)
\ar[ru]^{\beta\times\alpha} \ar[d]_{\tau_{J\times M}}
&   & \Mor(J\times M, J \times M) \ar[d]^{\tau_{*}}\\
\underline{\Mf}_E^{(n,d)}(J\times M)
\ar[rr]^{\alpha_{J\times M}} & & \Mor(J\times M, M) \\
\underline{\Mf}_E^{(n,d)}(M) \ar[u]^{\tau^*}
\ar[rr]^{\alpha_{M}} & & \Mor(M, M). \ar[u]_{\tau^*}
}
\]
Since $\tau_*(\id_{J\times M}) = \tau^*(\id_{M}) = \tau$, we obtain:
$p^* \kP \otimes q^* \kL \sim (\tau\times \id_{E})^*\kP$.
\end{proof}

\noindent
Let $o$ denote the neutral element of $J$ and recall that we have chosen $e
\in E$ and $m \in M$.
Lemma \ref{L:IsomLP}  implies that there exist open neighbourhoods
$o \in J' \subseteq J$ and $m \in M' \subseteq M$ such that
$
q^*\kL \otimes p^*\kP|_{J' \times M' \times E} \cong
\tilde\tau^*\kP|_{J' \times M' \times E}.
$
The following proposition  is crucial.

\begin{proposition}\label{P:goodtriv}
Let $E$ be a Weierstra\ss{} cubic curve, $M = M_E^{(n,d)}$ be the moduli
space of stable vector bundles of rank $n$ and degree $d$ and
$J$ the Jacobian of $E$. Then there exist
open neighbourhoods $o \in J' \subseteq J$ and $m \in M' \subseteq M$,
trivializations $\xi: \kP|_{M' \times E'} \lar \kO_{M' \times E'}^n$,
$\eta: \kL|_{J' \times E'} \lar \kO_{J' \times E'}^n$
and an  isomorphism
\[
\varphi:  q^*\kL \otimes p^*\kP|_{J' \times M' \times E} \lar
\tilde\tau^*\kP|_{J' \times M' \times E},
\]
such that the following diagram in the category of vector bundles
$\VB(J' \times M' \times E')$ is commutative:
$$
\xymatrix
{
  \bigl(q^*\kL \otimes p^*\kP\bigr)|_{J' \times M' \times E'}
  \ar[rr]^-{\varphi|_{J' \times M' \times E'}} \ar[d] & &
  \tilde\tau^*\kP|_{J' \times M' \times E'} \ar[d]  \\
  q^*\bigl(\kL|_{J' \times E'}\bigr)
  \otimes p^*\bigl(\kP|_{M' \times E'}\bigr)
  \ar[d]_{q^*(\eta) \otimes p^*(\xi)}  & &
  \tilde\tau^*\bigl(\kP|_{M' \times E'}\bigr)
  \ar[dd]^{\tilde\tau^*(\xi)} \\
  \kO_{J' \times M' \times E'} \otimes  \kO_{J' \times M' \times E'}^n
  \ar[d]_{\mathsf{mult}} & & \\
  \kO_{J' \times M' \times E'}^n \ar[rr]^{I_n} & &
  \kO_{J' \times M' \times E'}^n,
}
$$
where $I_n$ is the identity morphism.
\end{proposition}

\begin{proof} This proposition follows from a case-by-case analysis made
below for each of the three types of Weierstra\ss{} cubic curves, see
Proposition \ref{P:goodtriell} and Theorem \ref{T:goodtrivsing}.
\end{proof}

\begin{remark}
Note that we require  the existence of a  \emph{global} isomorphism $\varphi$
on the whole space $J' \times M' \times E$, although the condition on $\varphi$
is \emph{local}. The necessity of these  assumptions is  explained in the
course of the proof of  Theorem \ref{T:diffinvbpar}, which is one of the main
results of this article.
\end{remark}

To describe the correct coordinates on the moduli space $M$, we use a
canonical isomorphism $M\lar J^{d}$ and a non-canonical one $J\lar J^{d}$,
which depends on the chosen point $e\in E$. To define the canonical
one, which is given by taking the determinant bundle, recall from
Theorem \ref{T:ellfibrrepr} that we have an isomorphism of functors
$\underline{\det}: \underline{\Mf}_E^{(r,d)} \lar \underline{\Pic}^d_E$.
Using $\alpha$ and $\beta^{d}$, it defines an isomorphism of complex spaces
$\det:M\lar J^{d}$, such that $\det^{\ast}\kL^{(d)} \sim \det(\kP)$.

\begin{lemma}\label{L:ActPic}
The following diagram is commutative
\[
\xymatrix
{
J \times M \ar[r]^-\tau \ar[d]_{\id_{J} \times \det}& M \ar[d]^\det  \\
J \times J^d \ar[r]^-{\sigma} & J^d,
}
\]
where $\sigma: J \times J^d \lar J^d$ corresponds (with the aid of $\beta$ and
$\beta^{d}$) to the natural transformation of functors
$\underline{\sigma}:  \underline{\Pic}_E^0 \times \underline{\Pic}_E^d \lar
\underline{\Pic}_E^d$
which sends $(\kL', \kL'')$ to ${\kL'}^{\otimes n} \otimes \kL''$.
\end{lemma}

\begin{proof}
This result follows from the isomorphism $\det(\kL \otimes \kF) \cong
\kL^{\otimes n} \otimes \det(\kF)$, where $\kF$ is a vector bundle of rank $n$ and
$\kL$ a line bundle.
\end{proof}

The second isomorphism, $t^{e}:J\lar J^{d}$, is defined by taking the tensor
product with $\kO_{E}(de)$. More precisely, on functorial level, it is given
by
\[\underline{t}^{e}_{S}(\kM) =
\kM\otimes\kO_{S\times E}\bigl(d(S\times\{e\})\bigr).\]
This gives us a commutative diagram
\begin{equation}\label{diag:jacob}
\xymatrix
{
J \times J \ar[r]^-{\sigma'} \ar[d]_{\id_{J} \times t^{e}}& J \ar[d]^{t^{e}}  \\
J \times J^{d} \ar[r]^-{\sigma} & J^d,
}
\end{equation}
in which $\sigma'$ is defined on functorial level by the same formula as
$\sigma$ in Lemma \ref{L:ActPic}.

Finally, recall that the Jacobian  $J = \Pic^0(E)$ has the following
description:
\[
J \cong
\begin{cases}
  \CC/\Gamma        & \text{if } \; $E$ \; \text{ is elliptic,}\\
  \CC^*\cong\CC/\ZZ & \text{if } \; $E$ \; \text{ is nodal,}\\
  \CC               & \text{if } \; $E$ \; \text{ is cuspidal.}
\end{cases}
\]
In particular, in all three cases we have a surjective homomorphism of
Lie groups $\pi_{J}:\CC\lar J$.
We combine $\pi_{J}$ with the two isomorphisms $\det:M\lar J^{d}$ and
$t^{e}:J\lar J^{d}$ to obtain a local isomorphism $\pi_{M}:\CC\lar M$, which
gives us local coordinates on $M$. This local isomorphism sits in the
commutative diagram
\[
\xymatrix
{
\CC \times \CC \ar[rrr]^-{\sigma}\ar[d]_{\pi_{J}\times\pi_{M}}&&&\CC\ar[d]^{\pi_{M}}&\\
J \times M  \ar[rrr]^-\tau & & & M,&\text{where $\sigma(a,b) = na + b$.}
}
\]
It is obtained by combining diagram \eqref{diag:jacob} with the one in Lemma
\ref{L:ActPic} and the local isomorphism $\pi_{J}:\CC\lar J$, which is a
homomorphism of groups.

\begin{theorem}\label{T:diffinvbpar}
Let $E$ be a Weierstra\ss{} cubic curve, $\breve{E}$ its smooth part,
$(n,d) \in \mathbb{N}\times {\mathbb Z}$ a pair of coprime integers,
$M = M_E^{(n, d)}$ denote the moduli space of stable vector bundles of rank
$n$ and degree $d$ and  $J = \Pic^0(E)$ be the Jacobian of $E$.
Let $\kP|_{M' \times E'} \stackrel{\xi}\lar \kO_{M' \times E'}^n$ be a
trivialization satisfying the conditions of Proposition \ref{P:goodtriv}.
Then there exist coordinates on $M'$ and on $E'$, such that the corresponding
associative $r$-matrix $r^\xi(v_1, v_2; y_1, y_2)$ satisfies
\[
r^\xi(v_1 + v, v_2+v; y_1, y_2) = r^\xi(v_1, v_2; y_1, y_2).
\]
\end{theorem}

\noindent
\begin{proof} Introduce  the following notation.
As in the previous section, for $i = 1,2$ let
\[
\pi_i: M \times M \times \breve{E} \times \breve{E} \times E
\lar M \times E
\]
be the canonical projections $\pi_i(v_1, v_2; y_1, y_2; y) = (v_i, y)$
and let
\[
h_i:  M \times M \times \breve{E} \times \breve{E} \lar
M \times M \times \breve{E} \times \breve{E} \times E
\]
be the canonical sections, given by
$h_i(v_1, v_2; y_1, y_2) = (v_1, v_2; y_1, y_2; y_i)$.
Let
\[
\bar\pi: J \times M \times M \times \breve{E} \times \breve{E} \lar
M \times M \times \breve{E} \times \breve{E}
\]
be the canonical projection and
\[
\bar\tau: J \times M \times M \times \breve{E} \times \breve{E} \lar
M \times M \times \breve{E} \times \breve{E}
\]
be the ``diagonal'' action of the Jacobian
$J$, which is given by $\bar\tau(g; v_1, v_2; y_1, y_2) :=
\bigl(\tau(g, v_1), \tau(g, v_2); y_1, y_2\bigr)$.
Let  $\hat\tau = \bar\tau \times \id_{E}$ and
$\hat\pi = \bar\pi \times \id_{E}$ be the morphisms
$$
\hat\tau, \hat\pi: \,
J \times M \times M \times \breve{E} \times \breve{E} \times E\lar
M \times M \times \breve{E} \times \breve{E} \times E,
$$
$
\hat{q}: J \times M \times M \times \breve{E} \times \breve{E} \times E
\lar J \times E
$
the canonical projection. Finally, for $i = 1,2$, define the projections
\[
\hat\pi_i:=\pi_{i}\circ\hat\pi:
J \times M \times M \times \breve{E} \times \breve{E} \times E
\lar M \times E
\]
and the sections
\[
\hat{h}_i:=\id_{J}\times h_{i}:
J \times M \times M \times \breve{E} \times \breve{E}
\lar
J \times M \times M \times \breve{E} \times \breve{E} \times E.
\]
Using Theorem \ref{T:changedebase} and Proposition \ref{P:multlb},
we obtain the following commutative diagram of vector bundles on the complex
manifold $U:=J' \times M' \times M' \times \breve{E} \times \breve{E}$
\[
\xymatrix
{
\bar\pi^*h_1^* {\mathcal Hom}(\pi_1^*\kP, \pi_2^*\kP) \ar[d]_{\can}
\ar[rrr]^{\bar\pi^*\bigl(\tilde{r}^{\pi_1^*\kP, \pi_2^*\kP}_{h_1, h_2}\bigr)} & & & \bar\pi^*h_2^* {\mathcal Hom}(\pi_1^*\kP, \pi_2^*\kP) \ar[d]^{\can} \\
\hat{h}_1^*{\mathcal Hom}(\hat\pi_1^*\kP, \hat\pi_2^*\kP) \ar[d]_{\can}
\ar[rrr]^{\tilde{r}^{\hat\pi_1^*\kP, \hat\pi_2^*\kP}_{\hat{h}_1, \hat{h}_2}} & & &
\hat{h}_2^*{\mathcal Hom}(\hat\pi_1^*\kP, \hat\pi_2^*\kP) \ar[d]^{\can} \\
\hat{h}_1^*{\mathcal Hom}(\hat\pi_1^*\kP \otimes \hat{q}^*\kL, \hat\pi_2^*\kP \otimes \hat{q}^*\kL)
\ar[d]_{\hat{h}_1^*\bigl(\conj(\hat\varphi_1, \hat\varphi_2)\bigr)}
\ar[rrr]^{\tilde{r}^{\hat\pi_1^*\kP \otimes \hat{q}^*\kL, \hat\pi_2^*\kP \otimes \hat{q}^*\kL}_{\hat{h}_1, \hat{h}_2}}
& & &
\hat{h}_2^*{\mathcal Hom}(\hat\pi_1^*\kP \otimes \hat{q}^*\kL, \hat\pi_2^*\kP\otimes \hat{q}^*\kL) \ar[d]^{\hat{h}_2^*\bigl(\conj(\hat\varphi_1, \hat\varphi_2)\bigr)} \\
\hat{h}_1^*{\mathcal Hom}(\hat\tau^*\pi_1^*\kP, \hat\tau^*\pi_2^*\kP) \ar[d]_{\can}
\ar[rrr]^{\tilde{r}^{\hat\tau^*\pi_1^*\kP, \hat\tau^*\pi_2^*\kP}_{\hat{h}_1, \hat{h}_2}} & & &
\hat{h}_2^*{\mathcal Hom}(\hat\tau^*\pi_1^*\kP, \hat\tau^*\pi_2^*\kP) \ar[d]^{\can}\\
\bar\tau^*h_1^* {\mathcal Hom}(\pi_1^*\kP, \pi_2^*\kP)
\ar[rrr]^{\bar\tau^*\bigl(\tilde{r}^{\pi_1^*\kP, \pi_2^*\kP}_{h_1, h_2}\bigr)} & &  &
\bar\tau^*h_2^* {\mathcal Hom}(\pi_1^*\kP, \pi_2^*\kP).
}
\]
In this diagram, the isomorphisms  of vector bundles
$
\hat\varphi_i: \hat{q}^*\kL \otimes \hat\pi_i^*\kP \lar \hat\tau^* \pi_i^*\kP
$
are defined as follows. For $i = 1,2$ let
\[
p_i: J'\times M' \times M' \times \breve{E} \times \breve{E} \times E
\lar J' \times M' \times E
\]
be the natural projection
$p_i(g; v_1, v_2; y_1, y_2; y) = (g, v_i, y)$.
Let $\xi: \kP|_{M' \times E'} \rightarrow \kO_{M' \times E'}^n$
be a  trivialization and
$\varphi:q^*\kL \otimes p^*\kP|_{J' \times M' \times E}
\lar
\tilde\tau^*\kP_{J' \times M' \times E}$ an isomorphism, both
satisfying the conditions of Proposition \ref{P:goodtriv}.
Then we set $\hat\varphi_i$ to be the composition of morphisms of vector
bundles on $U\times E$
\[
\hat{q}^*\kL \otimes \hat\pi_i^*\kP \stackrel{\can}\lar
p_i^*(q^*\kL \otimes p^*\kP) \stackrel{p_i^*(\varphi)}\lar
p_i^*\tilde\tau^*\kP \stackrel{\can}\lar  \hat\tau^* \pi_i^*\kP.
\]
Note that the commutative square which involves the $\hat\varphi_{i}$ is only
available if $\hat\varphi_{i}$ is defined on $U\times E$ and not only on
$U\times E'$. The reason is that $\underline{\res}$ in the definition of
$\tilde{r}$ would not be an isomorphism on $E'$ (see the proof of Theorem
\ref{T:changedebase}).

\medskip
\noindent
Let $O$ denote the ring of holomorphic functions on
$U' := J' \times M' \times M' \times E' \times E' \subset U$.
With the aid of the trivialization $\xi$, for $i = 1,2$ we get isomorphisms
\begin{align*}
H^{0}\left(U', \bar\pi^*h_i^* {\mathcal Hom}(\pi_1^*\kP, \pi_2^*\kP)\right)
&\cong
\Mat_{n \times n}(O), \\
H^{0}\left(U', \bar\tau^*h_i^* {\mathcal Hom}(\pi_1^*\kP, \pi_2^*\kP)\right)
&\cong
\Mat_{n \times n}(O).
\end{align*}
Under these identifications, we can write the morphisms
$H^{0}\left(\bar\pi^*\bigl(\tilde{r}^{\pi_1^*\kP, \pi_2^*\kP}_{h_1, h_2}\bigr)\right)$ and
$H^{0}\left(\bar\tau^*\bigl(\tilde{r}^{\pi_1^*\kP, \pi_2^*\kP}_{h_1,
    h_2}\bigr)\right)$ as matrices
$\tilde{\mathbbm{r}}$ and $\tilde{\mathbbm{r}}'$, such that the large diagram
we set up earlier in this proof boils down to
\[
\xymatrix
{
\Mat_{n \times n}(O) \ar[rr]^{\tilde{\mathbbm{r}}} \ar[d]_{\mathsf{Id}} & &
\Mat_{n \times n}(O) \ar[d]^{\mathsf{Id}} & \\
\Mat_{n \times n}(O) \ar[rr]^{\tilde{\mathbbm{r}}'} & &  \Mat_{n \times n}(O),&
\text{ hence $\tilde{\mathbbm{r}} = \tilde{\mathbbm{r}}'$.}
}
\]
That the vertical arrows are identities is a consequence of Proposition
\ref{P:goodtriv}.

If we choose arbitrary coordinates on $J', M'$ and $E'$, we have
$\tilde{\mathbbm{r}}(g; v_1, v_2; y_1, y_2) = \tilde{r}^\xi(v_1, v_2; y_1, y_2)$
and $\tilde{\mathbbm{r}}'(g; v_1, v_2; y_1, y_2) =
\tilde{r}^\xi\bigl(\tau(g, v_1), \tau(g, v_2); y_1, y_2\bigr)$.
If we use the special coordinates on $J'$ and $ M'$ introduce above with the
aid of $\pi_{J}:\CC \lar J$, we obtain
$\tilde{r}^\xi\bigl(\tau(g, v_1), \tau(g, v_2); y_1, y_2\bigr) =
\tilde{r}^\xi(v_1 + ng, v_2+ ng; y_1, y_2)$. This implies
\[
\tilde{r}^\xi(v_1 + ng, v_2+ ng; y_1, y_2) =
\tilde{\mathbbm{r}}'(g; v_1, v_2; y_1, y_2) =
\tilde{\mathbbm{r}}(g; v_1, v_2; y_1, y_2) =
\tilde{r}^\xi(v_1, v_2; y_1, y_2),
\]
which gives the desired property of the tensor
$\tilde{r}^\xi(v_1, v_2; y_1, y_2)$.
\end{proof}

\begin{remark}
  Proposition \ref{P:goodtriv} and Theorem \ref{T:diffinvbpar} remain valid if
  the open neighbourhoods $J'$ and $M'$ are replaced by the maps
  $\pi_{J}:\CC\lar J$ and $\pi_{M}:\CC\lar M$. Similarly, by
  identifying $\breve{E}$ with $J^{1}$ and then proceeding as in the case of
  $M$, we may define a map $\pi_{E}:\CC\lar \breve{E}$,
  which can be used instead of $E'$. The advantage of this point of view is
  that $v_{i}, v, y_{i}$ can be arbitrary complex numbers in the statement of
  Theorem \ref{T:diffinvbpar}, whereas, if small neighbourhoods $J', M'$ and
  $E'$ are used, we have to make sure that $v_{i}$ and $v_{i}+v$ are in $M'$
  and $y_{i}\in E'$.
\end{remark}

\begin{remark}
Unfortunately, we have not found a ``conceptual way'' to prove Proposition
\ref{P:goodtriv}, without going to a case-by-case analysis. As a consequence,
we do not know whether this result generalizes to the relative case, when we
replace a Weierstra\ss{} curve $E$ by the Weierstra\ss{} fibration
$zy^2 = 4 x^3 -  g_2 xz^2 - g_3 z^3$.

\medskip\noindent
Motivated by the corresponding result for the classical $r$--matrices
\cite{BelavinDrinfeld2}, it is natural to conjecture, that the statement of
Theorem \ref{T:diffinvbpar}  holds for the other pair of
spectral variables, the ``skyscraper'' variables $(y_1, y_2)$.
Namely, there should exist coordinates on $E$ and a trivialization $\xi$ of
the universal family $\kP$  such that we have
\[
r^\xi(v_1+v, v_2+v; y_1+y, y_2+y) = r^\xi(v_1, v_2; y_1, y_2).
\]
\end{remark}

\begin{definition}
Let $r(v; y_1, y_2) \in \Mat_{n \times n}(\CC) \otimes \Mat_{n \times n}(\CC)$
be a non-degenerate unitary solution of the associative Yang--Baxter equation
such that there exists the limit  $\bar{r}(y_1, y_2)
=\lim\limits_{v \to e} \bigl(\pr \otimes \pr \bigr) r(v; y_1, y_2)$.
We say that $r$ is of \emph{elliptic type} if $\bar{r}$ is an elliptic
classical $r$-matrix, of \emph{trigonometric type} if $\bar{r}$ is
trigonometric and of \emph{rational type} if $\bar{r}$ is rational.
\end{definition}

\noindent
It was shown by Polishchuk \cite{Polishchuk1, Polishchuk2} that
in the case of elliptic curves one always gets an associative $r$-matrix
of elliptic type and in the case of Kodaira cycles a solution of trigonometric
type.

\begin{remark}
It is natural to conjecture that for any pair of coprime integers $(n,d)$
the geometric $r$-matrix corresponding to a cuspidal cubic curve always is of
rational type.
\end{remark}

The goal of the following three sections is to get an explicit form of the
geometric $r$-matrix attached to the Weierstra\ss{} fibration
$zy^2 = 4 x^3 -  g_2 xz^2 - g_3 z^3$ and the pair $(n, d) = (2, 1)$
at any given point $(g_2, g_3) \in \CC^2$.

\section{Elliptic solutions of the associative Yang--Baxter equation}\label{S:elliptrm}

In  this section we are going  to compute the solution of the associative
Yang--Baxter equation and the corresponding classical $r$--matrix, obtained
from the universal family of stable vector bundles of rank two and degree one
on a smooth elliptic curve.
In \cite[Section 2]{Polishchuk1}, Polishchuk computed the corresponding
triple Massey products using homological mirror symmetry and formulae
for higher products in the Fukaya category of an elliptic curve.

It is very instructive, however, to carry out a direct computation of the
geometric triple Massey products for an elliptic curve, independent of
homological mirror symmetry.  This approach  allows us to express the
resulting associative $r$--matrix in terms of Jacobi's theta-functions and the
corresponding classical $r$--matrix in terms of the elliptic functions
$\mathrm{cn}(z)$, $\mathrm{sn}(z)$ and $\mathrm{dn}(z)$.

In order to proceed with  the necessary calculations we  recall some standard
results about holomorphic vector bundles on one-dimensional complex tori, a
description of morphisms between them in terms of theta-functions etc.

\subsection{Vector bundles on a one-dimensional complex torus}

We start with  some classical  results about  vector bundles on
smooth elliptic curves.

\begin{theorem}[Atiyah,   Theorem 7 in \cite{Atiyah}]\label{T:Atiyah}
Let $E$ be a smooth elliptic curve over $\CC$ and $\kV$ a vector bundle on $E$.
\begin{itemize}
\item If $\End_E(\kV) = \CC$ then $\gcd\bigl(\rk(\kV), \deg(\kV)\bigr) = 1$, $\kV$
  is \emph{stable} and is determined by $\bigl(\rk(\kV), \deg(\kV),
  \det(\kV)\bigr) \in \mathbb{N} \times \mathbb{Z} \times E$, where we use an
  isomorphism  $\Pic^d(E) \cong E$.
\item If $\kV$ is indecomposable and $m = \gcd\bigl(\rk(V), \deg(V)\bigr)$
  then there exists a unique stable  vector bundle $\kV'$ such that $\kV =
  \kV' \otimes \kA_m$, where $\kA_m$ is the  indecomposable vector bundle of
  rank $m$ and degree $0$ recursively defined by non-split the exact sequences
  $$0 \lar \kA_{m-1} \lar \kA_m \lar \kO \lar 0, \qquad\qquad \kA_1 = \kO.$$
\end{itemize}
\end{theorem}

\noindent
In the complex-analytic case, one can give an explicit description
of the stable holomorphic vector bundles on a one-dimensional complex torus.

\begin{theorem}[Oda,   Theorem 1.2 in \cite{Oda}]\label{T:Oda}
Let $E$ be an elliptic curve
and  $\pi_{n}: E' \to E$ be an  \'etale covering of degree $n$.
\begin{itemize}
\item If $\kV$ is a stable  vector bundle on $E$ of rank $n$ and degree
  $d$,  then there exists a line bundle $\kL \in \Pic^d(E')$ such that
  $\kV \cong \pi_{n*}(\kL)$. Conversely, if $\gcd(n,d) =1$, then for any $\kL
  \in \Pic^d(E')$ the vector bundle $\kV \cong \pi_{n*}(\kL)$ is
  stable of rank $n$ and degree $d$.
\item If $\kL, \kN \in \Pic^d(E')$ satisfy $\pi_{n\ast}(\kL) \cong
  \pi_{n\ast}(\kN)$, then $(\kL \otimes \kN^{-1})^{\otimes n} \cong
  \kO_{E'}.$
\end{itemize}
\end{theorem}

\noindent
A very convenient way  to carry out calculations with vector bundles on
complex tori is provided  by the theory of automorphy factors, see
\cite{Morimoto} or \cite[Section I.2]{MumfordAbelian}.

\begin{definition}
Let $\tau \in \CC$ be a complex number such that $\mathrm{Im}(\tau) > 0$.
The category of \emph{automorphy factors} $\AF_\tau$ is defined as follows.
Its objects are pairs $(\Phi, n)$ where $n \ge 0$ is an integer and
$\Phi: \mathbb{C} \to \GL_n(\mathbb{C})$ is a holomorphic function
such  that for all $z \in \CC$ we have: $\Phi(z+1) = \Phi(z)$.
Given two automorphy factors $(\Phi, n)$ and $(\Psi, m)$, we define
$$
\Hom_{\AF_\tau}\bigl((\Phi,n), (\Psi,m)\bigr) =
\left\{A: \CC \to \Mat_{m \times n}(\CC)
\left|
\begin{array}{l}
A \mbox{\textrm{\quad is holomorphic}} \\
A(z+1) = A(z) \\
A(z+\tau)\Phi(z) = \Psi(z) A(z)
\end{array}
\right.
  \right\}
$$
and the composition of morphisms in $\AF_\tau$ is given by the multiplication of matrices.
In what follows, we shall frequently
denote the object $(\Phi, n)$ of $\AF_\tau$ by
$\Phi$. Note that one can define an interior tensor product
in the category $\AF_\tau$ induced  by the tensor product of matrices.
\end{definition}

Let
$\Lambda = \Lambda_\tau = \mathbb{Z} \oplus \mathbb{Z}\tau \subset \CC^2$ be the lattice
defined by $\tau$,   $E = E_\tau = \CC/\Lambda_\tau$ the corresponding complex torus
 and $\pi:  \CC \to
E$   the  universal covering of $E$.
For an object  $(\Phi, n)$ of the category  $\AF_\tau$ we define the sheaf
$\kE(\Phi)$ of $\kO_E$--modules and an embedding of sheaves
$\mm_\Phi: \kE(\Phi) \to \pi_* \kO_\CC^n$ by the following rule.

The open subsets $U\subset E$ for which all connected components of
$\pi^{-1}(U)$ map isomorphically to $U$, form a basis of the topology of $E$.
For such $U$, we let $U_{0}$ be a connected component of $\pi^{-1}(U)$ and
denote $U_\gamma = \gamma + U_{0}$ for all $\gamma\in\Lambda_\tau$.
Then
$\pi_*\kO_\CC^n\bigl(U\bigr) = \prod_{\gamma \in \Lambda} \kO^{n}_\CC(U_\gamma)$
and we define
$$
\kE(\Phi)\bigl(U\bigr) :=
\left\{
(F_\gamma(z))_{\gamma \in \Lambda} \in \pi_*(\kO_\CC^n)\bigl(U\bigr)
\left|
\begin{array}{l}
F_{\gamma + 1}(z +1) = F_\gamma(z) \\
F_{\gamma + \tau}(z + \tau)  = \Phi(z)F_\gamma(z)  \\
\end{array}\!\!\!
\right.
  \right\}.
$$
By $\mm_\Phi$ we denote the canonical embedding
$\kE(\Phi)\subset \pi_* \kO_\CC^n$.
The next theorem plays a key role in our computation of elliptic $r$-matrices.

\begin{theorem}\label{T:factaft} In the notations as above the following
  properties are true.
\begin{itemize}
\item Let $(\Phi, n)$ be an object of $\AF_\tau$ then the corresponding sheaf
$\kE(\Phi)$ is locally free of rank $n$.
\item The map $(\Phi, n) \mapsto \kE(\Phi)$ extends to a functor
 $\FF: \AF_\tau \lar \VB(E_\tau)$ which is an equivalence of categories
\item The functor $\FF$ commutes with tensor products:
$\kE(\Phi \otimes \Psi) \cong  \kE(\Phi) \otimes \kE(\Psi)$.
\item Let $\mathsf{For}: \AF_\tau \lar \VB(\CC)$ be the forgetful functor, i.e.
$\mathsf{For}(\Phi, n) = \kO_\CC^n$ and $\mathsf{For}(A) = A$ for any
object $(\Phi, n)$ and any  morphism
$A$.
 Then there is an isomorphism of functors
$
\gamma: \pi^* \circ \FF \lar \mathsf{For},
$
where for an automorphy factor $\Phi$ we set $\gamma_\Phi$
to be the composition
$
\pi^* \kE(\Phi) \xrightarrow{\pi^*(\mm_\Phi)}
\pi^*\pi_* \kO_\CC^n \stackrel{\can}\lar \kO_\CC^n.
$
\item
  The natural transformation $\gamma$ is compatible with tensor products,
  i.e.\ for any pair $(\Phi, n), (\Psi, m)$ of automorphy factors we have a
  commutative diagram:
$$
\xymatrix
{
\pi^*\bigl(\kE(\Phi \otimes \Psi)\bigr)
\ar[rr]^-{\gamma_{\Phi \otimes \Psi}} \ar[d]_\cong  & &
\kO_{\CC}^{mn} \\
\pi^*\kE(\Phi) \otimes \pi^*\kE(\Psi) \ar[rr]^-{\gamma_{\Phi} \otimes \gamma_\Psi} & &
\kO_{\CC}^m \otimes \kO_{\CC}^n \ar[u]_{\mathsf{mult}}
}
$$
\item If $\pi_n: E_{n\tau} \to E_\tau$ is the \'etale covering given by the
  inclusion of lattices $\Lambda_{n\tau} \subset \Lambda_{\tau}$, then
$\pi_n^*\bigl(\kE(\Phi)\bigr) \cong  \kE(\widetilde\Phi)$, where
$\widetilde\Phi(z) :=
\Phi(z + (n-1)\tau)\cdot\ldots\cdot \Phi(z+\tau) \Phi(z).$
\item The direct image
$\pi_{n*}\bigl(\kE(\Phi, m)\bigr) \cong  \kE(\widetilde\Phi, mn)$
of a vector bundle $\kE(\Phi, m)$ is given by the automorphy
factor\footnote{we thank Oleksandr Iena for helping us at this point}
$$
\widetilde\Phi =
\left(
\begin{array}{ccccc}
0 & I_m & 0 & \dots & 0 \\
0 & 0 & I_m &\dots  & 0 \\
\vdots & \vdots & \vdots & \ddots & \vdots \\
0 & 0 & 0 & \dots & I_m \\
\Phi & 0 & 0 & \dots & 0
\end{array}
\right).
$$
\end{itemize}
\end{theorem}

\noindent
In particular, if $\Phi(z)$ is an automorphy factor and
$A: \CC \to \GL_n(\CC)$ is a holomorphic function  such that $A(z+1) = A(z)$,
then $\Psi(z) = A(z+\tau)^{-1} \Phi(z) A(z)$ defines an isomorphic locally free sheaf
$\kE(\Phi) \cong \kE(\Psi)$.

\begin{remark}\label{R:bundles-via-autom-fact}
There is another way to describe the functor $\FF$.
Let $(\Phi, n)$ be an automorphy factor then
the corresponding holomorphic vector bundle $\kE(\Phi)$ can be  defined as the
quotient $\CC \times \CC^n/\!\sim$, where the equivalence relation is generated
by  $(z, v) \sim (z+1, v) \sim \bigl(z+\tau, \Phi(z)v\bigr)$. Using this
description,  we have the following commutative  diagram of complex manifolds
$$
\xymatrix
{
\CC \times \CC^n \ar[r] \ar[d]_{\mathrm{pr}_1} & \kE(\Phi) = \CC \times \CC^n/\sim
\ar[d]\\
\CC \ar[r]^\pi  & E.
}
$$
The natural transformation $\gamma$ can be constructed using
the fact that this diagram is Cartesian.
\end{remark}

\medskip
Our next goal is to give a description of those automorphy factors which
correspond to indecomposable vector bundles on $E$.
To do this, recall that the
 holomorphic morphism  $\CC\lar\CC^{\ast}$ sending  $z$ to $\exp(2\pi iz)$
identifies $\CC^{\ast}$ with $\CC/\mathbb{Z}$ and maps  $\tau$ to $q^{2}$,
where $q=\exp(\pi i\tau)$. Hence, it induces an isomorphism $E \cong
\CC^{\ast}/q^{2}$, where the quotient is formed modulo the multiplicative
subgroup generated by $q^{2}$.

Note that in
 the case of line bundles, automorphy factors are holomorphic functions
$\varphi:\CC\lar\CC^{\ast}$ which satisfy $\varphi(z+1)=\varphi(z)$.
In what follows  we shall use the notation
 $\kL(\varphi) := \kE(\varphi)$. Line bundles of degree zero can be given by
 constant automorphy factors, for example $\kL(1) = \kO_{E}$.
Because the function $a(z) = \exp(2\pi iz)$ satisfies $a(z+1)=a(z)$
and $a(z+\tau) = q^2 a(z)$ with $q= \exp(\pi i\tau)$ as above, the constants
$\varphi\in\CC^{\ast}$ and $q^2 \varphi$ define isomorphic line bundles on
$E$. In fact, the map
$
E := \CC^*/q^2 \lar  \Pic^0(E)
$
sending  $\varphi \in \CC^*$ to $\kL(\varphi) \in \Pic^0(E)$, is an
isomorphism.

\medskip
To describe line bundles of non-zero degree, we denote
$p_0 = \frac{\displaystyle 1 + \tau}{\displaystyle 2} \in E = \CC/\Lambda$.
The automorphy factor
$$\varphi_0(z) = \exp(-\pi i\tau - 2\pi iz)$$
satisfies
$\kL(\varphi_0) = \kO_{E}(p_0)$. To see this, recall that, by definition,
$$
H^0\bigl(\kL(\varphi_0)\bigr)  \cong
\Hom\bigl(\kL(1), \kL(\varphi_0)\bigr) =
\left\{f: \CC \to \CC
\left|
\begin{array}{l}
f \mbox{\textrm{\quad is holomorphic}} \\
f(z+1) = f(z) \\
f(z+\tau) = \varphi_0(z) f(z)
\end{array}
\right.
  \right\}
$$
and that this  vector space is \emph{one-dimensional} and
 is generated by the \emph{basic theta
function}
$$
\theta(z|\tau) = \theta_3(z|\tau) =
\sum\limits_{n \in \mathbb{Z}} \exp(\pi i n^2\tau + 2\pi inz),
$$
see for example \cite{MumfordTheta}. It is well-known that $\theta(z|\tau)$
vanishes at $p_0 = \frac{\displaystyle 1 + \tau}{\displaystyle 2}$. Moreover,
this is the only zero in the fundamental parallelogram of
$\Lambda_\tau$. Hence, $H^0\bigl(\kL(\varphi_0)\bigr) \cong  \CC$ and
by the Riemann-Roch theorem   the line bundle $\kL(\varphi_0)$
has degree one. Moreover, if $\underline{\theta}(z|\tau)$ it its non-zero global section then $\mathsf{div}\bigl(\underline{\theta}(z|\tau)\bigr) =
[p_0]$ and hence $\kL(\varphi_0) \cong \kO_E(p_0)$.
Because $\theta\left(\left. z + \dfrac{1+\tau}{2} - x\right|\tau\right)$
has its unique zero at $x\in E$, we obtain
\begin{equation}\label{E:canidentell}
\kO_{E}(x) \cong
\kL\left(t^{\ast}_{\frac{1+\tau}{2} - x}\varphi_0\right).
\end{equation}
where $t_x^*\varphi_0(z) := \varphi_0(z+x)$.
This gives a complete description of $\Pic^{1}(E)$.

Finally, any line bundle of degree $d$ can be written as
$\kO_{E}\bigl([(d-1)p_0] + [p_0 - x]\bigr)$
for some point $x \in E$. To complete the description of
$\Pic(E)$, it remains to observe that
$$
\kO_{E}\bigl([p_0 - x]+(d-1)[p_0]\bigr) \cong
\kL\bigl(t_x^*\varphi_0\cdot  \varphi_0^{d-1}\bigr).
$$

Our next aim is to find an explicit family of automorphy factors describing
the set of  stable vector bundles of rank $n$ and degree $d$ on $E$, where
$\gcd(n,d) =1$. Interpreting Oda's description from  Theorem \ref{T:Oda} in
terms of automorphy factors we immediately obtain  the  following proposition.

\begin{proposition}\label{P:nicefamily}
Let $(n,d) \in \mathbb{N} \times \mathbb{Z}$ be coprime. For $x \in \CC/\langle 1,  \tau\rangle$
let
$\widetilde\varphi_{n,d}(z, x) := \exp(-\pi ind\tau - 2\pi idz - 2 \pi ix)$.
Then  the family of automorphy factors
$$
\widetilde\Phi_{n,d}(z, x) =
\left(
\begin{array}{ccccc}
0 & 1 & 0 & \dots & 0 \\
0 & 0 & 1 &\dots  & 0 \\
\vdots & \vdots & \vdots & \ddots & \vdots \\
0 & 0 & 0 & \dots & 1 \\
\widetilde\varphi_{n,d} & 0 & 0 & \dots & 0
\end{array}
\right)
$$
describes the set  of stable vector bundles of rank $n$ and degree $d$ on $E$.
\end{proposition}

\noindent
However, these automorphy factors are  not compatible with the action
of the Jacobian $\Pic^0(E)$. In order to overcome this problem, denote
$q_{\frac{x}{n}} = \exp\left(-\frac{\displaystyle 2\pi ix}{\displaystyle n}\right)$
and let
$$
A =
\left(
\begin{array}{ccccc}
1 & 0 & 0 & \dots & 0 \\
0 & q_{\frac{x}{n}} & 0 & \dots & 0 \\
0 & 0 & q_{\frac{x}{n}}^2 & \dots & 0 \\
\vdots & \vdots & \vdots & \ddots & \vdots \\
0 & 0 & 0 & \dots & q_{\frac{x}{n}}^{n-1}
\end{array}
\right).
$$
Then $A^{-1} \widetilde\Phi_{n,d} A =:  \Phi_{n,d}$ is the following matrix-valued function:
\begin{equation}\label{E:univfact}
\Phi_{n,d}(z, x) =
\left(
\begin{array}{ccccc}
0 & q_{\frac{x}{n}} & 0 & \dots & 0 \\
0 & 0  & q_{\frac{x}{n}} & \dots & 0 \\
\vdots & \vdots & \vdots & \ddots & \vdots \\
0 & 0 & 0 & \dots & q_{\frac{x}{n}} \\
q_{\frac{x}{n}}\varphi_{n,d} & 0 & 0 & \dots & 0
\end{array}
\right) =
q_{\frac{x}{n}}
\left(
\begin{array}{ccccc}
0 & 1 & 0 & \dots & 0 \\
0 & 0  & 1  & \dots & 0 \\
\vdots & \vdots & \vdots & \ddots & \vdots \\
0 & 0 & 0 & \dots & 1 \\
\varphi_{n}^d & 0 & 0 & \dots & 0
\end{array}
\right),
\end{equation}
where $\varphi_{n}(z) = \exp(- \pi i n\tau - 2 \pi iz)$ and $\varphi_{n,d} = \varphi_n^d$.
Note that we have the equality
$$  \exp(-2\pi iy)  \cdot  \Phi_{n,d}(z, x) = \Phi_{n,d}(z, x + ny).
$$

\begin{lemma}\label{L:moduli-and-autom-fact}
In the notations as above, for  two points $x, x' \in \CC/\langle 1, n\tau\rangle$ we have:
$\kE\bigl(\widetilde\Phi_{n,d}(z, x)\bigr)
\cong \kE\bigl(\widetilde\Phi_{n,d}(z, x')\bigr)$ if and only
if $x - x' \in \Lambda_\tau$.
\end{lemma}

\begin{proof}
By Theorem \ref{T:Oda}, the vector  bundle $\kE^x:= \kE\bigl(\widetilde\Phi_{n,d}(z, x)\bigr)$
is stable of rank $n$ and degree $d$  for any point $x \in E_{n\tau}$. By
\cite[Theorem 7]{Atiyah}  the Jacobian $\Pic^0(E)$ acts transitively on the moduli space $M^{(n, d)}_E$. Moreover,
  $\kE^x \cong \kE^x \otimes \kL$ for a line
bundle $\kL \in \Pic^0(E)$ if and only if  $\kL^{\otimes n} \cong \kO$. For  $a \in \CC$,  the line bundle
  $\kL^a :=  \kL\bigl(\exp(2\pi i a)\bigr)$ fulfills the property  $\kL^{\otimes n} \cong \kO$ precisely
  if   $n a \in \Lambda_\tau$. Observe that  $\kL^a \otimes \kE^x \cong \kE^{x + na}$. In particular,
it shows that $\kE^{x} \cong \kE^{x + a + b\tau}$ for any $x \in E_{n \tau}$ and $a, b \in \ZZ$.
\end{proof}

Our next goal is to explain how the language of
automorphy factors can be used to describe a universal family of stable
vector bundles of rank $n$ and degree $d$ on $E$  as well as
to construct a trivialization of it.
In order to do this, we need the following generalization
of Theorem  \ref{T:factaft}.

As usual,
let $\tau \in \CC$ be such that $\mathrm{Im}(\tau) > 0$ and  $M$ be a complex manifold. Then we  define
a category $\AF_\tau(M)$, whose objects are pairs $(\Phi, n)$, where
$n \ge 1$ is an integer and $\Phi$ is a holomorphic function
$\Phi: \CC \times M \lar \GL_n(\CC)$ such that $\Phi(z+1, y) = \Phi(z, y)$
for all $(z, y) \in \CC \times M$. For a pair of automorphy factors
$(\Phi, n)$ and $(\Psi, m)$ we set
$$
\Hom_{\AF_\tau(M)}\bigl(\Phi, \Psi\bigr) =
\left\{\CC \times M \xrightarrow{A} \Mat_{m \times n}(\CC)
\left|\!
\begin{array}{l}
A \mbox{\textrm{\quad is holomorphic}} \\
A(z+1, y) = A(z, y) \\
A(z+\tau, y)\Phi(z, y) = \Psi(z, y) A(z, y)
\end{array}\!\!
\right.
  \right\}
$$
and the composition of morphisms in $\AF_\tau(M)$ is given by the
multiplication of matrices.
As before, we have a fully faithful functor
$$\mathbb{F}_M: \AF_\tau(M) \lar \VB(E \times M)$$ mapping
an automorphy factor $(\Phi, n)$ to the subsheaf $\kE(\Phi)$ of the sheaf
$\pi_{M*}\kO_{\CC \times M}^n$, where
$\pi_M = \pi \times \id: \CC \times M \lar E \times M$.
The sheaf $\kE(\Phi)$ is defined exactly as in the absolute case.
Moreover, $\mathbb{F}_M$ is dense (hence an equivalence of categories) if,
 for example, $M \cong \Delta_1 \times \dots \times
\Delta_m \subseteq \CC^m$, where each $\Delta_i \subseteq \CC, 1 \le i \le m$  is either an open disc or
$\CC$ itself.
  This functor maps the tensor product of automorphy factors into the tensor product of the corresponding vector bundles. Next,
there is an isomorphism of functors
$
\gamma: \pi_M^* \circ \FF \lar \mathsf{For},
$
where $\mathsf{For}: \AF_\tau(M) \lar \VB(\CC \times M)$ is the forgetful
functor.  For an automorphy factor $\Phi$ the morphism  $\gamma_\Phi$
is  the composition
$$
\pi_M^* \kE(\Phi) \xrightarrow{\pi_M^*(\mm_\Phi)}
\pi_M^*\pi_{M*} \kO_{\CC\times M}^n \stackrel{\can}\lar \kO_{\CC\times M}^n.
$$
Let $U$ be an open subset of $E$ such that there exists
a connected component $\widetilde U$ of $\pi^{-1}(U)$ which maps
isomorphically to $U$.
Hence, $\pi: \widetilde{U} \lar U$ is an isomorphism of
Riemann surfaces and the morphism $\gamma_\Phi$ induces a trivialization
of the vector bundle $\kE(\Phi)|_{U \times M}$.

It is important to note  that the  natural transformation $\gamma$ is
compatible with tensor products, i.e.\ for any pair $(\Phi, n)$ and $(\Psi,
m)$ of automorphy factors we have a commutative diagram:
$$
\xymatrix
{
\pi^*\bigl(\kE(\Phi \otimes \Psi)\bigr)
\ar[rr]^-{\gamma_{\Phi \otimes \Psi}} \ar[d]_\cong  & &
\kO_{\CC \times M}^{mn} \\
\pi^*\kE(\Phi) \otimes \pi^*\kE(\Psi) \ar[rr]^-{\gamma_{\Phi} \otimes \gamma_\Psi} & &
\kO_{\CC \times M}^m \otimes \kO_{\CC \times M}^n. \ar[u]_{\mathsf{mult}}
}
$$
Moreover, any holomorphic map $f: M \lar N$ of
open domains  induces a  functor
$f^*: \AF_\tau(N) \lar \AF_\tau(M)$ mapping an automorphy factor
$\Phi: \CC \times N \to \GL_n(\CC)$ to the automorphy factor
$f^*(\Phi) : \CC \times M \xrightarrow{\id \times f } \CC \times N
\stackrel{\Phi}\lar \GL_n(\CC)$. In these terms,  we have the following diagram of functors
$$
\xymatrix
{
\AF_\tau(N) \ar[rr]^{f^*} \ar[d]_{\FF_N} & & \AF_\tau(M)  \ar[d]^{\FF_M} \\
\VB(E \times N)  \ar[rr]^{(1 \times f)^*}  & & \VB(E \times M),
}
$$
where both compositions $(\id \times f)^* \circ \FF_N$ and $\FF_M  \circ f^*$ are canonically isomorphic.

\medskip
Our next goal is to give an explicit description of a universal family of  stable vector
bundles of rank $n$ and degree $d$ on the  complex torus $E = E_\tau$.
In what follows, $M = E_\tau$ stands for the moduli space of such bundles.
Consider a pair
of matrix-valued functions  $\Phi, \Psi: \CC \times \CC \rightarrow \GL_n(\CC)$ given by the formulae:
\begin{equation}\label{E:new-def-of-ell-family}
\Phi(z, x) =
\exp(-2 \pi i x)
\left(
\begin{array}{ccccc}
0 & 1 & 0 & \dots & 0 \\
0 & 0  & 1  & \dots & 0 \\
\vdots & \vdots & \vdots & \ddots & \vdots \\
0 & 0 & 0 & \dots & 1 \\
\varphi_{n}^d & 0 & 0 & \dots & 0
\end{array}
\right)
\quad
\Psi(z, x) = \exp(-2 \pi i z) I_n,
\end{equation}
where $\varphi_n(z)$ is the same as in (\ref{E:univfact}).
As in Remark \ref{R:bundles-via-autom-fact}, we define the vector bundle
$\kE(\Phi, \Psi) \in \VB(E \times M)$
via the following commutative diagram of complex manifolds:
$$
\xymatrix
{
\CC \times \CC \times \CC^n \ar[rr] \ar[d]_{\mathrm{pr}_{1, 2}} &  & \kE(\Phi, \Psi) = \CC \times \CC \times \CC^n/\sim
\ar[d]\\
\CC \times \CC \ar[rr]^{\pi \times \pi}  & & E \times M,
}
$$
where the equivalence relation is given by the formulae:
\begin{align*}
  (z,x;v) &\sim (z+1,x;v) \sim (z,x+1;v), \\
  (z,x;v) &\sim \bigl(z+\tau,x;\Phi(z,x)v\bigr),\\
  (z,x;v) &\sim \bigl(z, x+\tau;\Psi(z,x)v\bigr).
\end{align*}
Note that the following equalities  are fulfilled:
$$
\Psi(z+\tau,x)\Phi(z,x)  = \Phi(z,x+\tau)\Psi(z,x)
$$
as well as
\begin{align*}
\Phi(z+1,x) &= \Phi(z,x)& \Phi(z, x+1) &= \Phi(z,x)\\
\Psi(z+1,x) &= \Psi(z,x)& \Psi(z, x+1) &= \Psi(z,x).
\end{align*}
Hence the equivalence relation $\sim$ is well-defined and
$\kE(\Phi, \Psi)$ is a holomorphic vector bundle on $E \times M$. Note that
for any point $x_0 \in M$ we have:
$$
\kE(\Phi, \Psi)\big|_{E \times \{x_0\}} \cong \kE\bigl(\Phi(\,-\,, x_0)\bigr).
$$
In particular, $\kE(\Phi, \Psi)$ is a family of stable vector bundles of rank $n$ and
degree $d$ on $E$ parameterized by the manifold $M$.
By Lemma \ref{L:moduli-and-autom-fact} we know that
\begin{equation}\label{E:remark-on-family}
\kE(\Phi, \Psi)\big|_{E \times \{x_0\}} \cong \kE(\Phi, \Psi)\big|_{E \times \{x_0'\}}
\, \, \Longleftrightarrow \, \,
n(x_0 - x_0') \in \Lambda_\tau.
\end{equation}

\begin{lemma}\label{L:some-family-of-stab-bundl}
Let $\mu_n: E_\tau \rightarrow E_\tau$ be an \'etale covering of degree $n^2$ given by the rule
$\mu_n(x) = n \cdot x$. Then there exists a universal family $\kP = \kP(n,d)$ of stable vector bundles
on $E \times M$ such that $\kE(\Phi, \Psi) \sim (1 \times \mu_n)^* \kP$.
\end{lemma}

\begin{proof} We know that $M = E_\tau$ is the moduli space of stable vector bundles of rank
$n$ and degree $d$.
Let $\kQ \in \VB(E \times M)$ be any universal family, then by the universal property there exists a unique morphism
$\nu: M \rightarrow M$ such that $\kE(\Phi, \Psi) \sim (1 \times \nu)^* \kQ$.
Since a morphism between two compact Riemann surfaces is either surjective or constant,
the morphism $\nu$ is surjective. From  the equality (\ref{E:remark-on-family}) we  obtain
that $\nu$ factorizes as
$$
 \xymatrix
 {
 M \ar[rr]^\nu \ar[rd]_{\mu_n} & & M  \\
 & M  \ar[ru]_{\hat\nu}&
 }
$$
and the induced map $\hat\nu$ is both injective and surjective, hence biholomorphic. Then the universal
family $\kP =  (1 \times \hat\nu)^* \kQ \in \VB(E \times M)$ is the one we are looking for.
\end{proof}

\begin{remark}\label{R:univ-fam-of-lin-bundl}
In a similar way, the  functions $\varphi, \psi: \CC \times \CC \rightarrow \CC^*$
given by $\varphi(z, x) = \exp(-2 \pi i x)$ and $\psi(z, x) = \exp(-2 \pi i z)$ define
a universal family $$\kL = \kL(\varphi, \psi) \in \Pic(E \times J),$$
 of degree zero line bundles on $E = E_\tau$ where $J = E_\tau$ is the Jacobian of $E$.
\end{remark}

\begin{lemma}
Let $\kP \in \VB(E \times M)$ and $\kL \in \Pic(E \times J)$ be as in
Lemma \ref{L:some-family-of-stab-bundl} and Remark \ref{R:univ-fam-of-lin-bundl}. Then
the following diagram is commutative:
\begin{equation}\label{E:imp-diag}
\xymatrix
{
\CC \times \CC \ar[rr]^{\sigma} \ar[d]_{\pi \times \bar\pi} & & \CC \ar[d]^{\bar\pi} \\
J \times M \ar[rr]^\tau & & M,
}
\end{equation}
where $\sigma(a, b) = a + b$, $\pi: \CC \rightarrow \CC/\langle 1, \tau\rangle$ is   the
canonical projection and $\bar\pi(z) = \pi(n\cdot z)$.
\end{lemma}

\begin{proof}
Recall that the morphism $\tau: J \times M \rightarrow M$
is uniquely determined by a choice of universal families
$\kP$ and $\kL$ by the following property: the isomorphism
$$
\kP\big|_{E \times \{\tau(a, b)\}} \cong \kL\big|_{E \times \{a\}} \otimes \kP|_{E \times \{b\}}
$$
holds for all points $(a, b) \in J \times E$. Consider the morphisms
$$
\tilde{p}: \quad E \times \CC \times \CC \xrightarrow{\pr_{1,3}} E \times \CC  \xrightarrow{1 \times \bar\pi}
E \times M,
$$
$$
\tilde{q}: \quad E \times \CC \times \CC \xrightarrow{\pr_{1,2}} E \times \CC  \xrightarrow{1 \times \pi}
E \times J
$$
and
$$
\tilde{\sigma}: \quad E \times \CC \times \CC \xrightarrow{1 \times \sigma} E \times \CC  \xrightarrow{1 \times \bar\pi}
E \times M.
$$
The commutativity of diagram (\ref{E:imp-diag}) is equivalent to
the fact that
$
\tilde{q}^* \kL \otimes \tilde{p}^*\kP  \sim \tilde{\sigma}^* \kP.
$
This is furthermore equivalent that for all $(a, b) \in \CC \times \CC$ we have:
$$
\tilde{q}^* \kL \otimes \tilde{p}^*\kP  \big|_{E \times \{a\} \times \{b\}}  \cong
\tilde{\sigma}^* \kP\big|_{E \times \{a\} \times \{b\}}.
$$
The last isomorphism can be rewritten as
$$
\kL\bigl(\varphi(a)\bigr) \otimes \kE\bigl(\Phi(\,-\,, b)\bigr)
\cong
\kE\bigl(\Phi(\,-\,, a+ b)\bigr).
 $$
Since for all $(a, b) \in \CC \times \CC$ we have
$\varphi(z, a) \cdot  \Phi(z, b) = \Phi(z, a+ b)$, the result follows.
\end{proof}

\medskip
\noindent
Consider the morphisms $\hat{p}, \, \hat{q}, \,  \hat{\tau}: E \times \CC \times \CC \rightarrow E \times \CC$,
where $\hat{p} = {\pr_{1,3}}$, $\hat{q} = \pr_{1,2}$ and $\hat{\sigma} = 1 \times \sigma$.
We define  the morphism $\hat{q}^* \kL(\varphi) \otimes \hat{p}^* \kE(\Phi) \stackrel{\phi}\lar
\hat{\sigma}^*\kE(\Phi)$  by the following commutative diagram:
$$
\xymatrix
{
\hat{q}^* \kL(\varphi) \otimes \hat{p}^* \kE(\Phi) \ar[rrr]^-{\phi}  \ar[d]_\cong& & &
\hat{\sigma}^*\kE(\Phi) \ar[d]^\cong \\
\kL(\hat{q}^*\varphi) \otimes  \kE(\hat{p}^*\Phi) \ar[rrr]  \ar[d]_\cong & & &
\kE(\hat{\sigma}^*\Phi) \ar[d]^=\\
\kE(\hat{q}^*\varphi \cdot \hat{p}^*\Phi) \ar[rrr]^-{\bar\phi} & & & \kE(\hat{\sigma}^*\Phi),
}
$$
where all vertical isomorphisms are canonical and $\bar\phi$ corresponds to the morphism in the
category $\AF_\tau(\CC \times \CC)$ given by  the identity matrix $I_n$. In particular,
the morphism $\phi$ is identity in the trivializations of $\hat{p}^* \kE(\Phi)$, $\hat{\sigma}^*\kE(\Phi)$
 and $\hat{q}^*\kL(\varphi)$ induced by
$\gamma_\Phi$ and $\gamma_\varphi$.
Summing up, we get the following proposition

\begin{proposition}\label{P:goodtriell}
  Let $\Phi:  \CC \times \CC  \lar \GL_n(\CC)$ be
  as in  (\ref{E:new-def-of-ell-family}),  $\varphi(y) = \exp(-2\pi y )$
  and $
  \hat{q}^* \kL(\varphi) \otimes \hat{p}^* \kE(\Phi) \stackrel{\phi}\lar
\hat{\sigma}^*\kE(\Phi)$
   be the morphism of vector bundles on $E \times \CC \times \CC$
  constructed above. Take small  local neighbourhoods in the moduli space $M$ and Jacobian $J$
  corresponding to small neighbourhoods of $0 \in \CC$  with respect to the
  diagram (\ref{E:imp-diag}). Then the induced trivializations
  $\xi^\kP$ and $\xi^\kL$  of the universal families $\kP = \kP(n,d )$ and $\kL = \kP(1,0)$
   and the morphism $\phi$ are the ones   we are looking for in Proposition
 \ref{P:goodtriv}. In particular,  Theorem \ref{T:diffinvbpar} is true in the elliptic case.
\end{proposition}

\begin{corollary}
Let $\tilde{r}_{\mathsf{ell}}^{(n,d)} = \tilde{r}^\xi_\tau(v_1, v_2; y_1, y_2)$ be the
 associative $r$--matrix obtained  from the universal family $\kP(n,d)$ of
stable vector bundles of rank $n$ and degree $d$ on the elliptic curve
$E_\tau$  using the trivialization $\xi = \gamma_\Phi$ described above.
Then we have:
$$\tilde{r}_{\mathsf{ell}}^{(n,d)}(v_1, v_2; y_1, y_2) =
\tilde{r}_{\mathsf{ell}}^{(n,d)}(v_1 + v, v_2 +v;
y_1, y_2)
$$
for all values $v_1, v_2, v$ and $y_1, y_2$ from a small neighbourhood of $0$.
\end{corollary}

\subsection{Rules to calculate the evaluation and the residue maps}
In this subsection, we consistently denote the vector space of complex linear
maps between two complex vector spaces $V,W$ by $\Lin(V,W)$. We reserve here
the notation $\Hom_{\CC}(\;\;,\;\;)$ for the vector space of morphisms of
sheaves on the complex manifold $\CC$.

Let $E = E_\tau$ be a complex torus, $\Omega_E$ the sheaf of
holomorphic differential one-forms and $\omega  = dz \in H^0(\Omega_E)$ a global
section. We fix a pair of coprime integers $(n, d) \in \mathbb{Z}^+ \times
\mathbb{Z}$ and let $M = M_E^{(n, d)}$ be the moduli space
of stable holomorphic vector bundles of rank $n$ and degree $d$ on $E$. By
$\kP = \kP(n, d) \in \VB(E \times M)$ we denote the universal family and by
$\xi^\kP: \kP|_{U \times M'} \lar \kO|_{U \times M'}^n$ a trivialization, as
constructed in the previous subsection.
Recall that these data define the germ of  a meromorphic function
$$
\tilde{r} = \tilde{r}^\xi : \bigl(M \times M \times E \times E, o\bigr) \lar
\Mat_{n \times n}(\CC) \otimes \Mat_{n \times n}(\CC),
$$ whose value at the  point $(v_1, v_2; y_1, y_2)$,
where $v_1 \ne v_1$ and $y_1 \ne y_2$,  is defined via the commutative diagram
\begin{equation}\label{E:maindiag}
\xymatrix
{ & \Hom_E\bigl(\kP^{v_1}, \kP^{v_2}(y_1)\bigr)
\ar[ldd]_{\res^{\kP^{v_1}, \kP^{v_2}}_{y_1}(\omega)}
\ar[rdd]^{\ev^{\kP^{v_1},\kP^{v_2}(y_1) }_{y_2}} &
\\
 & &
\\
\Lin(\kP^{v_1}|_{y_1}, \kP^{v_2}|_{y_1})
\ar[d]_{\conj(\bar\xi^{v_1}, \, \bar\xi^{v_2})}
\ar[rr]^{\tilde{r}^{\kP^{v_1}, \kP^{v_1}}_{y_1, y_2}(\omega)} & &
\Lin(\kP^{v_1}|_{y_2},\kP^{v_1}|_{y_2}) \ar[d]^{\conj(\bar\xi^{v_1}, \, \bar\xi^{v_2})}
\\
\Mat_{n \times n}(\CC) \ar[rr]^{\tilde{r}^\xi(v_1, v_2; y_1, y_2)} & &
\Mat_{n \times n}(\CC),
}
\end{equation}
where $\res^{\kP^{v_1},\kP^{v_2} }_{y_1}(\omega)$ and
$\ev^{\kP^{v_1},\kP^{v_2}(y_1)}_{y_2}$ are the maps defined in Section \ref{S:Massey} and
the vertical isomorphisms are induced by the trivialization $\xi$ of
the universal bundle.

Let $\pi: \CC \lar  \CC/\Lambda_\tau = E$
be the universal covering, $o = \pi(0) \in E$.
Take  an open neighbourhood $U$ of the point $o$ in $E$ such that there exists
a connected component $\widetilde U$ of $\pi^{-1}(U)$ which maps
isomorphically to $U$.
For the sake of convenience, we denote the preimage $\pi^{-1}(y) \in
\widetilde U$ of a point $y \in U$ by the same letter $y$.
By Theorem  \ref{C:Masseyprod} we have a commutative diagram
$$
\xymatrix
{
\Lin\bigl(\kP^{v_1}|_{y_1}, \kP^{v_2}|_{y_1}\bigr)
\ar[rr]^-{\pi^*} & &
\Lin\bigl(
\pi^*\kP^{v_1}|_{y_1}, \pi^*\kP^{v_2}|_{y_1}\bigr)
\\
 \Hom_E\bigl(\kP^{v_1}, \kP^{v_2}(y_1)\bigr)
\ar[u]^{\res_{y_1}^{\kP^{v_1}, \kP^{v_2}}(\omega)}
\ar[rr]^-{\pi^*}
\ar[d]_{\ev_{y_2}^{\kP^{v_1}, \kP^{v_2}(y_1)}} &  &
\Hom_\CC\bigl(\pi^*\kP^{v_1}, \pi^*\kP^{v_2}(D_1)\bigr)
\ar[u]_{\res_{y_1}^{\pi^*\kP^{v_1}, \pi^*\kP^{v_2}}(\tilde\omega)}
\ar[d]^{\ev_{y_2}^{\pi^*\kP^{v_1}, \pi^*\kP^{v_2}(D_1)}}
\\
\Lin\bigl(\kP^{v_1}|_{y_2}, \kP^{v_2}|_{y_2}\bigr)
\ar[rr]^-{\pi^*} & &
\Lin\bigl(\pi^*\kP^{v_1}|_{y_2}, \pi^*\kP^{v_2}|_{y_2}\bigr)
}
$$
where
$\kO_\CC(D_1) = \pi^* \kO_E(y_1)$ is the subsheaf of the sheaf $\kM_\CC$ of
meromorphic functions on $\CC$, whose local sections are meromorphic functions
having  at most simple poles along the infinite, but locally finite set
$D_1=\bigl\{y_1+\gamma\,\big|\,\gamma\in\Lambda_{\tau}\bigr\}=\pi^{-1}(y_{1})$.

Recall that the description of a universal family $\kP$ in terms
of automorphy factors yields an isomorphism of vector bundles
$
\gamma: \pi_M^* \kP \lar \kO_{\CC \times M}^n.
$
For any point $v \in M$ it induces an isomorphism
$\gamma^v: \pi^*\kP^v \lar \kO_{\CC}^n$.
If we apply Lemma \ref{L:conjofres} and Proposition \ref{L:functofev} to these
isomorphisms and use
the morphisms $\res_{y_1}(\tilde\omega)$ and $\ev_{y_2}$ from Lemma
\ref{L:resexplit} and Lemma \ref{L:evexplit}, which are given by
$\res_{y_1}(\tilde\omega)\bigl(F\bigr)=\res_{y_1}(F dz)$ and
$\ev_{y_2}(F)=F(y_2)$, we obtain the following commutative diagram
$$
\xymatrix
{
\Lin\bigl(\pi^*\kP^{v_1}|_{y_1}, \pi^*\kP^{v_2}|_{y_1}\bigr)
\ar[rrrr]^-{\conj(\bar\gamma^{v_1}, \, \bar\gamma^{v_2})}
& & & & \Mat_{n \times n}(\CC)
\\
\Hom_\CC\bigl(\pi^*\kP^{v_1}, \pi^*\kP^{v_2}(D_1)\bigr)
\ar[rrrr]^-{\conj\bigl(\gamma^{v_1}, \, \gamma^{v_2}(D_1)\bigr)}
\ar[u]^{\res_{y_1}^{\pi^*\kP^{v_1}, \pi^*\kP^{v_2}}(\tilde\omega)}
\ar[d]_{\ev_{y_2}^{\pi^*\kP^{v_1}, \pi^*\kP^{v_2}(D_1)}}
& & &   &  \Mat_{n \times n}\bigl(\kO_\CC(D_1)\bigr)
\ar[u]_{\res_{y_1}(\tilde\omega)}
\ar[d]^{\ev_{y_2}}
\\
\Lin\bigl(\pi^*\kP^{v_1}|_{y_2}, \pi^*\kP^{v_2}|_{y_2}\bigr)
\ar[rrrr]^-{\conj(\bar\gamma^{v_1}, \, \bar\gamma^{v_2})}
& &  & & \Mat_{n \times n}(\CC).
}
$$
The three previous diagrams in this subsection give us another one:
$$
\xymatrix
{
\Mat_{n \times n}(\CC)
\ar[rr]^{\tilde{r}^\xi(v_1, v_2; y_1, y_2)} & &
\Mat_{n \times n}(\CC)
\\
&
\Hom_{\CC}\bigl(\kO_\CC^n, \kO_{\CC}^n(D_1)\bigr)
\ar[lu]^{\res_{y_1}(\tilde\omega)}
\ar[ru]_{\ev_{y_2}}
&
\\
&
\Hom_E\bigl(\kP^{v_1}, \kP^{v_2}(y_1)\bigr)
\ar[ld]_{\res^{\kP^{v_1}, \kP^{v_2}}_{y_1}(\omega)}
\ar[rd]^{\ev^{\kP^{v_1},\kP^{v_2}(y_1) }_{y_2}}
\ar[u]^{\conj\bigl(\gamma^{v_1}\,,\gamma^{v_2}(D_1)\bigr)\,\circ\,\pi^*}
&
\\
\Lin(\kP^{v_1}|_{y_1}, \kP^{v_2}|_{y_1})
\ar[uuu]^{\conj(\bar\xi^{v_1}, \, \bar\xi^{v_2})}
\ar[rr]^{\tilde{r}^{\kP^{v_1}, \kP^{v_1}}_{y_1, y_2}(\omega)} & &
\ar[uuu]_{\conj(\bar\xi^{v_1}, \, \bar\xi^{v_2})}
\Lin(\kP^{v_1}|_{y_2},\kP^{v_1}|_{y_2}).
}
$$
Our next goal is to describe the image
of the morphism $\conj\bigl(\gamma^{v_1} \,,\gamma^{v_2}(D_1)\bigr) \circ \pi^*$.
Let $ \psi_y(z) = - \exp(-2\pi iz + 2\pi iy - 2\pi i\tau).$
It is easy to see that $\Hom_E\bigl(\kO, \kL(\psi_y)\bigr)$ is one-dimensional
and  generated by the section $\underline{\theta}_{y}$ corresponding
to the theta-function
$\theta_y(z)=\theta\left(\left.z+\frac{1+\tau}{2}-y\right|\tau\right)$.
Note that $\theta_y$ is a holomorphic function on $\CC$ having only one simple
zero at $y$ inside a fundamental parallelogram of $\Lambda_{\tau}$. Hence, we
have an isomorphism $\alpha: \kO_E(y) \lar \kL(\psi_y)$.

In order to be  more precise, recall
that $\kL(\psi_y)$ is a subsheaf of the sheaf $\pi_* \kO_\CC$.
In particular, we have:
$H^0\bigl(\kL(\psi_y)\bigr) = \CC \cdot \theta_y \subset H^0(\kO_\CC)$.
On the other hand, the sheaf $\kO_E(y)$ is a subsheaf  of the sheaf of
meromorphic functions $\kM_E$ and  $H^0\bigl(\kO_E(y)\bigr)  = \CC \cdot 1$.
Without loss of generality we may assume that $H^0(\alpha)(1) = \theta_y$.
This choice fixes the isomorphism $\alpha$.

Recall that for a point $x \in M$  we have:
$\kP^{v} = \kE(\Phi_v)$, where $\Phi_v = \Phi(-, v)$ is
the function defined by the equality  (\ref{E:univfact}).
Then we have an isomorphism
$$
\kP^{v_2}(y_1) \xrightarrow{\mathsf{id} \otimes \alpha}
\kP^{v_2} \otimes \kL(\psi_{y_1}) = \kE\bigl(\psi_{y_1} \cdot
\Phi_{v_2}\bigr)
$$
and the following diagram is commutative:
$$
\xymatrix@C-=3mm
{
\pi^* \kE(\psi_{y_1} \Phi_{v_2})
\ar[rr]^-{\gamma^{y_1, v_2}} & & \kO_\CC^n \\
\pi^* \kE(\psi_{y_1}) \otimes  \pi^*\kE(\Phi_{v_2})
\ar[rr]^-{\gamma^{y_1} \otimes \gamma^{v_2}} \ar[u]^{\can}
& & \kO_\CC \otimes \kO_\CC^n \ar[u]_{\mathsf{mult}}  \\
& \pi^* \kO_E(y_1) \otimes \pi^*\kE(\Phi_{v_2})
\ar[lu]^{\pi^*(\alpha) \otimes \mathsf{id}}
\ar[ru]_{\tilde\alpha \otimes \gamma^{v_2}}  &
}
$$
where
$\tilde\alpha: \pi^*\bigl(\kO_E(y_1)\bigr) = \kO_\CC(D_{1}) \lar \kO_\CC$
is defined to be $\tilde\alpha(f) = f \theta_{y_1}$,
$\gamma^{y_1, v_2}$ corresponds to the automorphy factor $\psi_{y_1}\cdot
\Phi_{v_2}$, $\gamma^{y_1}$ to $\psi_{y_1}$ and $\gamma^{v_2}$ to $\Phi_{v_2}$.
As a result, we get the following commutative diagram:
$$
\xymatrix@C-=-5mm
{
\Hom_\CC\Bigl(\pi^*\kP^{v_1},\pi^*\bigl(\kP^{v_2}\otimes\kO_E(y_1)\bigr)\Bigr)
\ar[rr]^-{\conj\bigl(\gamma^{v_1}, \, \gamma^{v_2}(D_1)\bigr)} & &
\Hom_\CC\bigl(\kO_\CC^n, \kO_\CC^n(D_{1})\bigr)
\\
\Hom_\CC\Bigl(\pi^*\kP^{v_1},
\pi^*\bigl(\kP^{v_2}\otimes\kL(\psi_{y_1})\bigr)\Bigr)
\ar[rr]^-{\conj\bigl(\gamma^{v_1}, \, \gamma^{v_2}(D_1)\bigr)}
\ar[u]^{(\mathsf{id} \otimes \alpha)_*} & &
\ar[u]_{\tilde\alpha_*}
\Hom_\CC\bigl(\kO_\CC^n, \kO_\CC^n\bigr)
\\
& \Hom_{\AF_\tau}(\Phi_{v_1}, \psi_{y_1} \Phi_{v_2}) \ar[lu]^{\pi^*}
\ar[ru]_{\mathsf{For}} &
}
$$
and $\tilde\alpha_*$ is given by the formula
$\tilde\alpha_*(F) = \frac{\displaystyle F}{\displaystyle \theta_{y_1}}$.

\begin{corollary}\label{C:ellcomp}
  Let $O = \Gamma(\CC, \kO_\CC)$,
  $O(D_{1}) = \Gamma\bigl(\CC, \kO_\CC(D_{1})\bigr)$ and
  $$
  \Pi^{v_1, v_2}_{y_1} =
  \mathsf{Im}\bigl(\Hom_{\AF_\tau}(\Phi_{v_1}, \psi_{y_1} \Phi_{v_2})
  \lar \Mat_{n \times n}(O)\bigr).
  $$
  Then the following diagram is commutative:
  $$
  \xymatrix
  {
    & \Pi^{v_1, v_2}_{y_1}
    \ar@/_13pt/[ddl]_{\overline{\res}_{y_1}}
    \ar@/^13pt/[ddr]^{\overline{\ev}_{y_2}}
    \ar[d]_{\tilde\alpha_*}  &
    \\
    & \Mat_{n \times n}\bigl(O(D_{1})\bigr)
    \ar[ld]_{\res_{y_1}(\tilde\omega)} \ar[rd]^{\ev_{y_2}} &
    \\
    \Mat_{n \times n}(\CC)
    \ar[rr]^{\tilde{r}^\xi(v_1, v_2; y_1, y_2)} & &
    \Mat_{n \times n}(\CC).
  }
  $$
  In particular, this gives the following algorithm  to compute the
  value of an associative geometric $r$--matrix of elliptic type at a point
  $(v_1, v_2; y_1, y_2) \in M \times M \times E \times E$ with respect
  to the trivialization $\xi$:
  \begin{enumerate}
  \item First describe the vector space
    $$
    \Pi^{v_1, v_2}_{y_1} = \mathsf{Im}\Bigl(
    \Hom_{\AF_\tau}\bigl(\Phi_{v_1}, \psi_{y_1}\cdot \Phi_{v_2}\bigr)  \lar
    \Mat_{n \times n}(O)
    \Bigr)
    $$
  \item The morphism
    $\overline{\res}_{y_1}:\Mat_{n\times n}(O)\rightarrow\Mat_{n\times n}(\CC)$
    is given by the formula
    $$
    F(z) \mapsto \overline{\res}_{y_1}
    \left(\frac{\displaystyle F(z)}{\theta_{y_1}(z)}dz\right) =
    \frac{\displaystyle F(y_1)}{\displaystyle \theta'_{y_1}(y_1)} =
    \frac{F(y_1)}{\theta'\left(\left.\frac{1 + \tau}{2}\right|\tau\right)}
    $$
    and the morphism
    $\overline{\ev}_{y_2}:\Mat_{n\times n}(O)\rightarrow\Mat_{n\times n}(\CC)$
    is given by the formula
    $$
    F(z) \mapsto \frac{\displaystyle F(y_2)}{\displaystyle \theta_{y_1}(y_2)} =
    \dfrac{F(y_2)}{\theta\left(\left.y + \frac{1+\tau}{2}\right|\tau\right)}.
    $$
  \item Compute $\tilde{r}^\xi(v_1, v_2; y_1, y_2)$ as the composition
    $$
    \Mat_{n \times n}(\CC)
    \xrightarrow{\bigl(\overline{\res}_{y_1}(\tilde\omega)\bigr)^{-1}}
    \Pi^{v_1, v_2}_{y_1}
    \xrightarrow{\overline{\ev}_{y_2}}  \Mat_{n \times n}(\CC) .
    $$
  \end{enumerate}
  Note that there is an ambiguity in choosing the morphism $\alpha$. Another
  choice of $\alpha$ corresponds to rescaling the section
  $\underline{\theta}_{y_1}$ by $\lambda \in \CC^*$ to
  $\lambda \underline{\theta}_{y_1}$. However, it is easy to see from the
  algorithm above, that the resulting linear map
  $\tilde{r}^\xi(v_1, v_2; y_1, y_2)$ does not depend on this choice.
\end{corollary}

\medskip

\subsection{Calculation of the elliptic $r$--matrix corresponding to $M_E^{(2,1)}$.}
Let $\kP^{x_1}$ and $\kP^{x_i}$ be a pair of  non-isomorphic
simple vector bundles of rank two and degree
one on the elliptic curve $E = E_\tau$, $y_1$ and $y_2$ two distinct points.
In what follows we denote $q = \exp(\pi i\tau)$, $q_x = \exp(-\pi ix)$,
$e(z) = \exp(-2\pi i z)$, $\varphi(z) = e(z+\tau)$,
$x = x_2 - x_1$ and $y = y_2 - y_1$.

\noindent
As we have seen in the previous subsection, one can write
$\kP^{x_i} =  \kE\bigl(\CC^2, \Phi_{2,1,x_i}(z)\bigr)$, where
$$
\Phi_{2,1,x_i}(z) =
q_{x_i}\left(
\begin{array}{cc}
0 & 1 \\
\varphi(z) & 0
\end{array}
\right) =:
q_{x_i} \Phi(z),
$$
and the  line bundle $\kO_E(y_1)$ corresponds to the automorphy factor
$$\psi_{y_1}(z) = - e(z+ \tau - y_1).$$
Recall that  $\Pi^{x_1, x_2}_{y_2}  = $
$$
\left\{\left.
A(z) =
\left(
\begin{array}{cc}
u(z) & v(z) \\
w(z) & t(z)
\end{array}
\right)
\right|
\begin{array}{c}
A(z+1) = A(z), \;
A(z+\tau) \Phi(z)  =
q_x \psi_{y_1}(z)
\Phi(z) A(z)
\end{array}
\right\}
$$
This leads to two systems of functional equations
$$
\left\{
\begin{array}{ccc}
u(z + \tau) & = &  q_x \psi_{y_1}(z) t(z) \\[2mm]
t(z+\tau)  & = & q_x  \psi_{y_1}(z) u(z)
\end{array}
\right.
\quad \mbox{\textrm{and}} \quad
\left\{
\begin{array}{ccc}
\varphi(z) v(z + \tau) & = &  q_x \psi_{y_1}(z) w(z) \\[2mm]
w(z+\tau)  & = & q_x  \varphi(z) \psi_{y_1}(z) v(z)
\end{array}
\right.
$$
which are equivalent to
$$
\left\{
\begin{array}{ccc}
u(z + 2 \tau) & = & a(z) u(z) \\
u(z+1)        & =  & u(z) \\
t(z)          & = &
\frac{\displaystyle u(z+\tau)}{\displaystyle q_x \psi_{y_1}(z)}
\end{array}
\right.
\quad \mbox{\textrm{and}} \quad
\left\{
\begin{array}{ccl}
v(z + 2 \tau) & = & b(z) v(z) \\
v(z+1)        & =  & v(z) \\
w(z)          & = &
\frac{\displaystyle \varphi(z)}{\displaystyle q_x \psi_{y_1}(z)}
v(z +\tau)
\end{array}
\right.
$$
where
\begin{align*}
  a(z) &=
\exp\left(- 2\pi i \tau - 2\pi i\left(z+\frac{x+\tau}{2}-y_1\right)\right)^2
\quad\text{ and }\\
b(z)&=\exp\left(-2\pi i \tau -2\pi i\left(z+\frac{x}{2} - y_1\right)\right)^2.
\end{align*}

\begin{lemma}\label{L:basis}
Let $E = E_\tau$ be an elliptic curve, $\varphi_0(z) = \exp(-\pi i \tau - 2\pi
iz)$, $l \in \mathbb{N}$. Then
$
H^0\bigl(\kL(\varphi_0^l)\bigr)$ has a basis
$\left.\left\{\theta\left[\tfrac{a}{l}, 0\right](lz|l\tau) \right|0 \le a < l,
 a\in\mathbb{Z}\right\},
$
where we use Mumford's notation
$$
\theta[a,b](z|\tau)  =
\sum\limits_{n \in \mathbb{Z}} \exp\bigl(\pi i(n+a)^2\tau + 2\pi i(n+a)(z+b)\bigr).
$$
\end{lemma}
In the particular case of bundles of rank two and degree one it is convenient
to use instead the four classical theta-functions of Jacobi:
$$
\begin{array}{l}
\theta_1(z|\tau) = 2 q^{\frac{1}{4}} \sum\limits_{n=0}^\infty (-1)^n q^{n(n+1)}
\sin\bigl((2n+1)\pi z\bigr), \\
\theta_2(z|\tau) = 2 q^{\frac{1}{4}} \sum\limits_{n=0}^\infty  q^{n(n+1)}
\cos\bigl((2n+1)\pi z\bigr), \\
\theta_3(z|\tau) = 1+  2 \sum\limits_{n=1}^\infty  q^{n^2}
\cos(2\pi nz), \\
\theta_4(z|\tau) = 1+  2 \sum\limits_{n=1}^\infty  (-1)^n q^{n^2}
\cos(2\pi nz).
\end{array}
$$
\begin{remark}\label{R:MumfordJacobi}
In Mumford's notation it holds:
\begin{align*}
  \theta_1(z|\tau) &= -\theta\left[\tfrac{1}{2}, \tfrac{1}{2}\right](z|\tau)&
  \theta_2(z|\tau) &= \theta\left[\tfrac{1}{2}, 0\right](z|\tau)\\
  \theta_3(z|\tau) &= \phantom{-}\theta\left[0,0\right](z|\tau)&
  \theta_4(z|\tau) &= \theta\left[0,\tfrac{1}{2}\right](z|\tau).
\end{align*}
\end{remark}

\medskip
\noindent
In what follows we shall express all our computations in terms of Jacobi's
theta-functions. From Lemma \ref{L:basis} and Remark \ref{R:MumfordJacobi}
we immediately obtain:
\begin{corollary}\label{C:ellbasis}
If we let
\begin{align*}
  u_1(z) &= \theta_3\left(2\left.\left(z-y_1+\dfrac{x +\tau}{2}\right)\right|4\tau\right)&
  v_1(z) &= \theta_3\left(2\left.\left(z-y_1+\dfrac{x}{2}      \right)\right|4\tau\right) \\
  u_2(z) &= \theta_2\left(2\left.\left(z-y_1+\dfrac{x +\tau}{2}\right)\right|4\tau\right)&
  v_2(z) &= \theta_2\left(2\left.\left(z-y_1+\dfrac{x}{2}      \right)\right|4\tau\right)
\end{align*}
and
$$
F_k(z) =
\begin{pmatrix}
u_k(z) & 0 \\
0 & \dfrac{u_k(z+\tau)}{q_x \psi_{y_1}(z)}
\end{pmatrix},
G_k(z) =
\begin{pmatrix}
0  & v_k(z) \\
\dfrac{\varphi(z)}{q_x  \psi_{y_1}(z)} v_k(z) & 0
\end{pmatrix}, \quad k =1,2,
$$
then $F_1(z), F_2(z), G_1(z), G_2(z)$ is a basis of
$\Pi^{x_1, x_2}_{y_1}$.
\end{corollary}

The following proposition summarises the main properties of Jacobi's
theta-functions which we need in our calculation of the associative
$r$--matrix corresponding to the universal family of stable vector bundles of
rank two and degree one.

\begin{proposition}[see \cite{BatErd} and  Section I.4 in \cite{Lawden}]
\label{P:theta}
The transformation rules for shifts of theta-functions are given by the table
$$
\begin{array}{|c||c|c|c|c|c|c|}
\hline
\theta(z) & \theta(-z) & \theta(z+1) &
\theta(z+\tau) & \theta(z+1+\tau) & \theta(z + \frac{1}{2})
& \theta(z + \frac{\tau}{2}) \\
\hline
\hline
\theta_1(z) & -\theta_1(z) & -\theta_1(z) & -p(z) \theta_1(z) & p(z)\theta_1(z)&
\theta_2(z) & i q(z) \theta_4(z) \\
\hline
\theta_2(z) & \theta_2(z) & -\theta_2(z) & p(z) \theta_2(z) & -p(z) \theta_2(z)&
- \theta_1(z) & q(z) \theta_3(z) \\
\hline
\theta_3(z) & \theta_3(z) & \theta_3(z) & p(z) \theta_3(z) & p(z) \theta_3(z)&
\theta_4(z) & q(z) \theta_2(z) \\
\hline
\theta_4(z) & \theta_4(z) & \theta_4(z) & -p(z) \theta_4(z) & -p(z) \theta_4(z)&
\theta_3(z) & i q(z) \theta_1(z) \\
\hline
\end{array}
$$
where $p(z) = \exp\bigl(-\pi i (2z + \tau)\bigr)$ and
$q(z) =
\exp\left(-\pi i\left(z + \dfrac{\tau}{4}\right)\right)$.
Moreover, Jacobi's theta-functions satisfy  the so-called Watson's
determinantal identities:
$$
\begin{array}{l}
\theta_3(2x|2\tau) \theta_2(2y|2\tau) -
\theta_3(2y|2\tau) \theta_2(2x|2\tau) = \theta_1(x+y|\tau)\theta_1(x-y|\tau), \\
\theta_1(2x|2\tau) \theta_4(2y|2\tau) -
\theta_1(2y|2\tau) \theta_4(2x|2\tau) = \theta_2(x+y|\tau)\theta_1(x-y|\tau), \\
\theta_1(2x|2\tau) \theta_4(2y|2\tau) +
\theta_1(2y|2\tau) \theta_4(2x|2\tau) = \theta_1(x+y|\tau)\theta_2(x-y|\tau), \\
\theta_4(2x|2\tau) \theta_4(2y|2\tau) -
\theta_1(2y|2\tau) \theta_1(2x|2\tau) = \theta_3(x+y|\tau)\theta_4(x-y|\tau), \\
\theta_4(2x|2\tau) \theta_4(2y|2\tau) +
\theta_1(2y|2\tau) \theta_1(2x|2\tau) = \theta_4(x+y|\tau)\theta_3(x-y|\tau). \\
\end{array}
$$
\end{proposition}

\medskip
\noindent
By Corollary \ref{C:ellbasis}, any element of
$\Pi^{x_1, x_2}_{y_1}$ can be written as a sum
$$
A(z) = \alpha F_1(z) + \beta F_2(z) + \gamma G_1(z) + \delta G_2(z)
$$
for some $\alpha, \beta, \gamma, \delta \in \CC$. In order to calculate the
geometric associative $r$--matrix $r^\xi(x_1, x_2; y_1, y_2)$
we have to solve the system of linear equations
$$
\overline{\res}_{y_1}\bigl(A(z)\bigr)  =
\left(
\begin{array}{cc}
a & b \\
c & d
\end{array}
\right).
$$
Then the  linear map $\tilde{r}^\xi(x_1, x_2; y_1, y_2): \Mat_{2 \times 2}(\CC) \lar
\Mat_{2 \times 2}(\CC)$ is given by the rule
$$
\left(
\begin{array}{cc}
a & b \\
c & d
\end{array}
\right)
\stackrel{\res_{y_1}^{-1}}\longmapsto
A(z)
\stackrel{\ev_{y_2}}\longmapsto
\frac{1}{\theta_3(y+ \frac{1+\tau}{2}|\tau)} A(y_2).
$$
It is easy to see that the system of linear  equations
$$
\overline{\res}_{y_1}\bigl(\alpha F_1(z) + \beta F_2(z) + \gamma G_1(z) + \delta G_2(z)\bigr)
=
\left(
\begin{array}{cc}
a & b \\
c & d
\end{array}
\right)
$$
splits into two independent systems
$$
\overline{\res}_{y_1}\bigl(F(z)\bigr) =
\left(
\begin{array}{cc}
a & 0 \\
0 & d
\end{array}
\right) \quad \mbox{\textrm{and}} \quad
\overline{\res}_{y_1}\bigl(G(z)\bigr)
=
\left(
\begin{array}{cc}
0 & b \\
c & 0
\end{array}
\right),
$$
where
$
F(z) = \alpha F_1(z) + \beta F_2(z)  \quad \mbox{\textrm{and}} \quad
G(z) = \gamma  G_1(z) + \delta G_2(z).
$

\medskip
\noindent
\textbf{Computation of the ``diagonal terms''}.
The system of linear equations
$$\overline{\res}_{y_1}\bigl(F(z)\bigr) :=
\frac{1}{\theta'_3(\frac{1+\tau}{2}|\tau)} F(y_1) = \left(
\begin{array}{cc}
a & 0 \\
0  & d
\end{array}
\right)
$$ reads as
$$
\left\{
\begin{array}{ccl}
\theta_3(x+\tau|4\tau)\alpha + \theta_2(x+\tau|4\tau) \beta & = &
\theta'_3(\frac{1+\tau}{2}|\tau) a \\
\theta_3(x+ 3\tau|4\tau)\alpha + \theta_2(x+ 3\tau|4\tau) \beta & = &
- e(\tau + \frac{x}{2}) \theta'_3(\frac{1+\tau}{2}|\tau) d.
\end{array}
\right.
$$
By Watson's identity, the determinant of this system is
$$
\Delta_1 = \left|
\begin{array}{cc}
\theta_3(x+\tau|4\tau) &  \theta_2(x+\tau|4\tau) \\
\theta_3(x+ 3\tau|4\tau) &  \theta_2(x+ 3\tau|4\tau)
\end{array}
\right| =
\theta_1(x+ 2\tau|2\tau) \theta_1(-\tau|2\tau) =
$$
$$
= e(x+\tau)
\theta_1(x|2\tau) \theta_1(\tau|2\tau)
$$
and we obtain:
$$
\left\{
\begin{array}{ccl}
\alpha & = & \frac{\theta'_3(\frac{1+\tau}{2}|\tau)}{
\Delta_1}
\bigl(\theta_2(x+3 \tau|4\tau)a + e(\tau + \frac{x}{2})
\theta_2(x+\tau|4\tau)d\bigr) \\
\beta  & = & -\frac{\theta'_3(\frac{1+\tau}{2}|\tau)}{
\Delta_1}
\bigl(\theta_3(x+3 \tau|4\tau)a + e(\tau + \frac{x}{2})
\theta_3(x+\tau|4\tau)d\bigr).
\end{array}
\right.
$$
This implies:
$$
\tilde{r}^\xi(x_1, x_2;y_1, y_2)\left[\left(\begin{array}{cc} a & 0 \\
0 & d\end{array}
\right)\right]  =
$$
$$
\frac{\theta'_3(\frac{1+\tau}{2}|\tau)}{\theta_3(y+ \frac{1+\tau}{2}|\tau)
\Delta_1}
\times
\left[ p_1(z)
\left(
\begin{array}{cc}
\theta_3(2y + x +\tau|4\tau) & 0 \\
0 & -
\frac{\displaystyle \theta_3(2y +
x + 3\tau|4\tau)}{\displaystyle e(x/2 + y + \tau)}
\end{array}
\right)
\right.
$$
$$
\left.
- p_2(z)
\left(
\begin{array}{cc}
\theta_2(2y + x +\tau|4\tau) & 0 \\
0 & -
\frac{\displaystyle \theta_2(2y +
x + 3\tau|4\tau)}{\displaystyle e(x/2 + y + \tau)}
\end{array}
\right)
\right].
$$
where
$$
\begin{array}{l}
p_1(z) =
\theta_2(x+3\tau|4\tau)a + e(x/2 +\tau) \theta_2(x+\tau|4\tau)d, \\
  p_2(z) =
\theta_3(x+3\tau|4\tau)a + e(x/2 +\tau) \theta_3(x+\tau|4\tau)d.
\end{array}
$$
In order to calculate  the ``diagonal part'' of the corresponding tensor
$r(x_1, x_2; y_1, y_2)$
 we use the inverse of the canonical isomorphism
$$
\Mat_{2 \times 2}(\CC) \otimes \Mat_{2 \times 2}(\CC) \lar
\Lin\bigl(\Mat_{2 \times 2}(\CC), \Mat_{2 \times 2}(\CC)\bigr)
$$
given by the formula $X \otimes Y \mapsto \tr(X \circ \,-) Y$.
It is easy to see that under the map
$$
\Lin\bigl(\Mat_{2 \times 2}(\CC), \Mat_{2 \times 2}(\CC)\bigr)
\lar \Mat_{2 \times 2}(\CC) \otimes \Mat_{2 \times 2}(\CC)
$$
a  linear function $e_{ij} \mapsto \alpha_{ij}^{kl} e_{kl}, \alpha_{ij}^{kl}
\in \CC^*$ corresponds to the tensor
$\alpha_{ij}^{kl}e_{ji} \otimes e_{kl}$.

\medskip
\noindent
Again,  Watson's identities imply:

\medskip
\noindent
$\bullet$  The coefficient at $e_{11} \otimes e_{11}$ is
$$
\theta_3(2y +x + \tau|4\tau) \theta_2(x+3\tau|4\tau) -
\theta_2(2y +x + \tau|4\tau) \theta_3(x+3\tau|4\tau) =
$$
$$
\theta_1(x +y + 2\tau|2\tau) \theta_1(y-\tau|2\tau).
$$

\medskip
\noindent
$\bullet$ The coefficient at  $e_{22} \otimes e_{22}$ is
$$
e(-y) \bigl(\theta_3(x + \tau|4\tau) \theta_2(2y + x+3\tau|4\tau) -
\theta_2(x + \tau|4\tau) \theta_3(2y + x+3\tau|4\tau)\bigr) =
$$
$$
= \theta_1(x +y + 2\tau|2\tau) \theta_1(y-\tau|2\tau).
$$

\medskip
\noindent
$\bullet$  The coefficient at  $e_{22} \otimes e_{11}$ is
$$
e(x/2 + \tau)\bigl(
\theta_2(x+\tau|4\tau)\theta_3(2y+x+\tau) - \theta_3(x+\tau|4\tau)\theta_2(2y+x+\tau)
\bigr)
=
$$
$$
= e(x/2 + \tau) \theta_1(y+x+\tau|2\tau) \theta_1(y|2\tau).
$$

\medskip
\noindent
$\bullet$ The coefficient at  $e_{11} \otimes e_{22}$ is
$$
e(- y - x/2 - \tau)\bigl(
\theta_2(2y + x + 3\tau|4\tau) \theta_3(x + 3\tau|4\tau) -
\theta_3(2y + x + 3\tau|4\tau) \theta_3(x + 3\tau|4\tau)
\bigr)
$$
$$
= e(x/2 + \tau) \theta_1(y+x+\tau|2\tau) \theta_1(y|2\tau).
$$
Now observe that
$$
\theta_1(x +y + 2\tau|2\tau) \theta_1(y-\tau|2\tau) =
i e(x + y/2 + 5\tau/4) \theta_1(x+y|2\tau) \theta_4(y|2\tau)
$$
and
$$
e(x/2 + \tau) \theta_1(y+x+\tau|2\tau) \theta_1(y|2\tau)
= i e(x + y/2 +  5\tau/4) \theta_4(x+y|2\tau) \theta_1(y|2\tau).
$$

\medskip
\noindent
Hence, the ``diagonal part''  of   $r(x_1, x_2; y_1, y_2)$ is
$$
C \bigl[\theta_1(x+ y|2\tau)\theta_4(y|2\tau)(e_{11}\otimes e_{11} +
e_{22}\otimes e_{22}) +
\theta_4(x+ y|2\tau)\theta_1(y|2\tau)(e_{11}\otimes e_{22} +
e_{22}\otimes e_{11})\bigr],
$$
where
$$
C = \frac{\theta_3'(\frac{1+\tau}{2}|\tau) e(x + \frac{y}{2} +
\frac{5\tau}{4})}{\theta_3(y+ \frac{1+\tau}{2}|\tau) \Delta_1}.
$$
From the identities  $\theta_3(y+ \frac{1+\tau}{2}|\tau) =
i \exp(-\pi i (y + \tau/4))\theta_1(y|\tau)$ and
 $\theta_1(0|\tau) = 0$ it follows:
$\theta_3'(\frac{1+\tau}{2}|\tau) = i e(\frac{\tau}{8}) \theta'(0|\tau)$.
Using the transformation rules from Proposition \ref{P:theta} we get:
$$
C = \frac{\displaystyle\theta_1'(0|\tau)}{\displaystyle
\theta_4(0|2\tau)\theta_1(x|2\tau) \theta_1(y|\tau)} =
\frac{\displaystyle\theta_1'(0|\tau)}{\displaystyle
\theta_1\left(\left.\frac{x}{2}\right|\tau\right)
\theta_2\left(\left.\frac{x}{2}\right|\tau\right)
\theta_1\left(y|\tau\right)},
$$
where we have used Landen's transform
$$
\theta_4(0|2\tau) \theta_1(2x|2\tau) = \theta_1(x|\tau) \theta_2(x|\tau).
$$

\medskip
\noindent
It remains to observe that
$
A(e_{11}\otimes e_{11} + e_{22}\otimes e_{22}) +
B(e_{11}\otimes e_{22} + e_{22}\otimes e_{11}) =
$
$$
\frac{1}{2}(A + B) (e_{11} + e_{22}) \otimes  (e_{11} + e_{22}) -
\frac{1}{2}(A - B) (e_{11} - e_{22}) \otimes  (e_{11} -  e_{22}),
$$
and that by Watson's identities we have
$$
\theta_1(x+y|2\tau) \theta_4(y|2\tau) + \theta_4(x+y|2\tau) \theta_1(y|2\tau) =
\theta_1\left(\left.y + \frac{x}{2}\right|\tau\right)
\theta_2\left(\left.\frac{x}{2}\right|\tau\right)
$$
and
$$
\theta_1(x+y|2\tau) \theta_4(y|2\tau) -  \theta_4(x+y|2\tau) \theta_1(y|2\tau) =
\theta_2\left(\left.y + \frac{x}{2}\right|\tau\right)
\theta_1\left(\left.\frac{x}{2}\right|\tau\right),
$$
so the contribution of the ``diagonal terms'' is
$$
\frac{1}{2} \frac{\theta'_1(0|\tau)}{\theta_1(y|\tau)}
\frac{\theta_1(y + \frac{x}{2}|\tau)}{\theta_1(\frac{x}{2}|\tau)}
\mathbbm{1} \otimes \mathbbm{1} +
\frac{1}{2} \frac{\theta'_1(0|\tau)}{\theta_1(y|\tau)}
\frac{\theta_2(y + \frac{x}{2}|\tau)}{\theta_2(\frac{x}{2}|\tau)}
h \otimes h.
$$

\medskip
\noindent
\textbf{Contribution of the ``skew terms''}.
We have to solve the system of linear equations
$$
\overline{\res}_{y_1}\bigl(\gamma G_1(z) + \delta G_2(z)\bigr) =
\left(
\begin{array}{cc}
0 & b \\
c & 0
\end{array}
\right).
$$
In explicit form this system reads as
$$
\left\{
\begin{array}{ccl}
\theta_3(x|4\tau)\gamma + \theta_2(x|4\tau)\delta & = &
\theta'_3(\frac{1+\tau}{2}|\tau)b \\
\theta_3(x + 2\tau|4\tau)\gamma + \theta_2(x+2\tau|4\tau)\delta & = &
- e(x/2 - y_1) \theta'_3(\frac{1+\tau}{2}|\tau)c.
\end{array}
\right.
$$
By Watson's formulas the determinant of this system is
$$
\Delta_2 =
\left|
\begin{array}{cc}
\theta_3(x|4\tau) & \theta_2(x|4\tau) \\
\theta_3(x + 2\tau|4\tau) & \theta_2(x+2\tau|4\tau)
\end{array}
\right| =
-\theta_1(x+\tau|2\tau) \theta_1(\tau|2\tau).
$$

\noindent
Hence, the  solution of this system of equations  is
$$
\left\{
\begin{array}{ccl}
\gamma & = & \theta'_3(\frac{1+\tau}{2}|\tau)\bigl(\theta_2(x+2\tau|4\tau)b +
e(x/2 - y_1) \theta_2(x|4\tau)c\bigr) \\
\delta  & = & \theta'_3(\frac{1+\tau}{2}|\tau)\bigl(\theta_3(x+2\tau|4\tau)b +
e(x/2 - y_1) \theta_3(x|4\tau)c.
\bigr)
\end{array}
\right.
$$
As a result,  we obtain:
$$
\tilde{r}^\xi(x_1, x_2;y_1, y_2)\left[\left(\begin{array}{cc} 0 & b \\
c & 0 \end{array}
\right)\right] =
$$
$$
\frac{\theta'_3(\frac{1+\tau}{2}|\tau)}{\theta_3(y+ \frac{1+\tau}{2}|\tau)
\Delta_2}
\times
\left[
q_1(z)
\left(
\begin{array}{cc}
0 &  \theta_3(2y + x|4\tau) \\
- e(y_1 - x/2) \theta_3(2y + x + 2\tau|4\tau) & 0
\end{array}
\right)
\right.
$$
$$
\left.
-  q_2(z)
\left(
\begin{array}{cc}
0 &  \theta_2(2y + x|4\tau) \\
- e(y_1 - x/2) \theta_2(2y + x + 2\tau|4\tau) & 0
\end{array}
\right)
\right],
$$
where
$$
\begin{array}{l}
q_1(z) = \theta_2(x+2\tau|4\tau)b + e(x/2 - y_1) \theta_2(x|4\tau)c, \\
 q_2(z) = \theta_3(x+2\tau|4\tau)b + e(x/2 - y_1) \theta_3(x|4\tau)c.
\end{array}
$$

\medskip
\noindent
Again,  Watson's identities imply:

\medskip
\noindent
$\bullet$  The coefficient at $e_{21} \otimes e_{12}$ is
$$
\theta_2(x+2\tau|4\tau) \theta_3(2y +x|4\tau) -
\theta_3(x+2\tau|4\tau) \theta_2(2y +x|4\tau) =
$$
$$ =
e\left(\frac{1}{2}(x+\tau)\right) \theta_4(x+y|2\tau) \theta_4(y|2\tau).
$$

\medskip
\noindent
$\bullet$ The coefficient at $e_{12} \otimes e_{21}$ is
$$
\theta_3(x|4\tau) \theta_2(2y+x+2\tau|4\tau) - \theta_2(x|4\tau)
\theta_3(2y+x+2\tau|4\tau) = $$
$$ =
e\left(\frac{1}{2}(x+\tau) +y\right) \theta_4(x+y|2\tau) \theta_4(y|2\tau).
$$

\medskip
\noindent
$\bullet$
 The coefficient at $e_{12} \otimes e_{12}$ is
$$
e(x/2 - y_1)\bigl(
\theta_2(x|4\tau) \theta_3(2y+x+2\tau|4\tau) - \theta_3(x|4\tau)
\theta_2(2y+x+2\tau|4\tau)\bigr) =
$$
$$
=
e(x/2 - y_1) \theta_1(x+y|2\tau) \theta_1(y|2\tau).
$$

\medskip
\noindent
$\bullet$ The coefficient at $e_{21} \otimes e_{21}$ is
$$
e(y_1 - x/2)\bigl(\theta_3(x+2\tau|4\tau) \theta_2(2y+x+2\tau|4\tau) -
\theta_2(x+2\tau|4\tau) \theta_3(2y+x+2\tau|4\tau)\bigr) =
$$
$$
=
e(y_2 + x/2 +\tau) \theta_1(y+x|2\tau) \theta_1(y|2\tau).
$$

\medskip
\noindent
Note that the  coefficients of the tensors $e_{12} \otimes e_{12}$ and
$e_{21} \otimes e_{21}$ are not functions of $y = y_2 - y_1$.
In order to overcome this problem we take
$\phi(y) = \left(\begin{array}{cc} e(y/2) & 0 \\
0 & e(-\tau/4)\end{array}\right)$ and consider the gauge transformation
$$
r(x;y_1, y_2) \mapsto \bigl(\phi(y_1) \otimes \phi(y_2)\bigr)
r(x;y_1, y_2) \bigl(\phi^{-1}(y_1) \otimes \phi^{-1}(y_2)\bigr).
$$
It is easy to see  that the ``diagonal tensors''
$e_{kk}\otimes e_{ll} (k,l = 1,2)$ remain unchanged (and, in particular,  this
gauge transformation does not influence the final answer for  the ``diagonal
terms'' obtained  before) and the transformation rule for the ``skew tensors''
is the following:
$$
\begin{array}{ccl}
e_{12} \otimes e_{12}
& \mapsto & e(\frac{\tau}{2}) e(\frac{y_1 + y_2}{2}) e_{12} \otimes e_{12}, \\
e_{21} \otimes e_{21}
& \mapsto & e(-\frac{\tau}{2}) e(-\frac{y_1 + y_2}{2}) e_{21} \otimes e_{21}, \\
e_{12} \otimes e_{21}
& \mapsto & e(-\frac{y}{2}) e_{12} \otimes e_{21}, \\
e_{21} \otimes e_{12}
& \mapsto & e(\frac{y}{2}) e_{21} \otimes e_{12}.
\end{array}
$$

\medskip
\noindent
Hence, the new tensor of ``skew terms'' is
$$
C\bigl[\theta_4(x+y|2\tau) \theta_4(y|2\tau)(e_{21} \otimes e_{12} +
e_{12} \otimes e_{21}) +
\theta_1(x+y|2\tau) \theta_1(y|2\tau)(e_{12} \otimes e_{12} +
e_{21} \otimes e_{21})
\bigr],
$$
where

\begin{align*}
  C = \frac{\theta'_3(\frac{1+\tau}{2}|\tau) e(\frac{1}{2}(x+y+\tau))}{
\theta_3(y+ \frac{1+\tau}{2}|\tau) \Delta_2} &=
\frac{\theta'_1(0|\tau)}{\theta_4(0|2\tau)
  \theta_4(x|2\tau)\theta_1(y|\tau)}\\
  &=
\frac{\theta'_1(0|\tau)}{\theta_3(x/2|\tau)\theta_4(x/2|\tau) \theta_1(y|\tau)}.
\end{align*}

\medskip
\noindent
Using the equality
$$
A \left(
\begin{array}{cc}
0 & 1 \\
1 & 0
\end{array}
\right)
\otimes
 \left(
\begin{array}{cc}
0 & 1 \\
1 & 0
\end{array}
\right) +
B \left(
\begin{array}{cc}
0 & -i \\
i & 0
\end{array}
\right)
\otimes
 \left(
\begin{array}{cc}
0 & -i \\
i & 0
\end{array}
\right) =
$$
$$
= (A+ B)(e_{12} \otimes e_{21} + e_{21} \otimes e_{12}) + (A-B)
(e_{21} \otimes e_{21} + e_{12} \otimes e_{12})
$$
and Watson's identities
$$
\theta_4(y+x|2\tau) \theta_4(y|2\tau) + \theta_1(y+x|2\tau) \theta_1(y|2\tau) =
\theta_4\left(\left.y + \frac{x}{2}\right|\tau\right)
\theta_3\left(\left.\frac{x}{2}\right|\tau\right)
$$
$$
\theta_4(y+x|2\tau) \theta_4(y|2\tau) - \theta_1(y+x|2\tau) \theta_1(y|2\tau) =
\theta_3\left(\left.y + \frac{x}{2}\right|\tau\right)
\theta_4\left(\left.\frac{x}{2}\right|\tau\right)
$$
it follows  that the contribution of the ``skew terms'' is
$$
\frac{1}{2} \frac{\theta'_1(0|\tau)}{\theta_1(y|\tau)}
\left(
\frac{\theta_3(y + \frac{x}{2}|\tau)}{\theta_3(\frac{x}{2}|\tau)}
\sigma \otimes \sigma +
\frac{\theta_4(y + \frac{x}{2}|\tau)}{\theta_4(\frac{x}{2}|\tau)}
\gamma \otimes \gamma
\right),
$$
where
$$
\sigma =
\left(
\begin{array}{cc}
0 & - i \\
i & 0
\end{array}
\right),
\qquad
\gamma  =
\left(
\begin{array}{cc}
0 & 1 \\
1 & 0
\end{array}
\right).
$$

\noindent
In summary, we obtain the following theorem.
\begin{theorem}
The universal family of stable vector bundles of rank two and degree one on
an elliptic curve $E_\tau$ gives the following solution of the associative
Yang--Baxter equation:
\begin{align*}
r^{(2,1)}_{\mathsf{ell}}(x;y) =
\frac{1}{2} \frac{\theta'_1(0|\tau)}{\theta_1(y|\tau)}
\left(
\frac{\theta_1(y + \frac{x}{2}|\tau)}{\theta_1(\frac{x}{2}|\tau)}
\mathbbm{1} \otimes \mathbbm{1} \right.&+
\frac{\theta_2(y + \frac{x}{2}|\tau)}{\theta_2(\frac{x}{2}|\tau)}
h  \otimes h + \\
+
\frac{\theta_3(y + \frac{x}{2}|\tau)}{\theta_3(\frac{x}{2}|\tau)}
\sigma  \otimes \sigma &\left.+
\frac{\theta_4(y + \frac{x}{2}|\tau)}{\theta_4(\frac{x}{2}|\tau)}
\gamma \otimes \gamma
\right).
\end{align*}
\end{theorem}

\noindent
Recall that
$$
\cn(z) =
\frac{\theta_4(0|\tau) \theta_2(z|\tau)}{\theta_2(0|\tau) \theta_4(z|\tau)},
\quad
\sn(z) =
\frac{\theta_3(0|\tau) \theta_1(z|\tau)}{\theta_2(0|\tau) \theta_4(z|\tau)},
\quad
\dn(z) =
\frac{\theta_4(0|\tau) \theta_3(z|\tau)}{\theta_3(0|\tau) \theta_4(z|\tau)}
$$
and
$$
\theta'_1(0|\tau) = \theta_2(0|\tau) \theta_3(0|\tau) \theta_4(0|\tau),
$$
see \cite[Sections I.5 and Section II.1]{Lawden}. Let
$\bar{r}(y) = \lim\limits_{x\to 0} (\pr \otimes \pr)r(x;y)$ then we have:

\begin{theorem}
The solution of the classical Yang--Baxter equation obtained from
the universal family of stable vector bundles of rank two and degree one on
a complex torus $E_\tau$ is
$$
\bar{r}(y) =
\frac{1}{2}\left(\frac{\cn(y)}{\sn(y)} h\otimes h +
\frac{1}{\sn(y)} \gamma \otimes \gamma +
 \frac{\dn(y)}{\sn(y)} \sigma \otimes \sigma\right).
$$
\end{theorem}

\begin{remark}
Note that
$\res_x\bigl(r(x;y)\bigr)  = \frac{1}{4} \mathbbm{1} \otimes \mathbbm{1}$,
hence the tensor $r_x(y) := r(x;y)$ also satisfies the quantum Yang--Baxter
equation for $x \ne 0$. In fact, it is the well-known  solution of the QYBE
which was found and  studied by Baxter.
\end{remark}

\begin{remark}[see for example Section VII.3 in \cite{Chandrasekharan}]
Let
$$
\wp(z) =
\frac{1}{z^2} + \sum\limits_{(n,m) \in \mathbb{Z}^2\setminus  \{0, 0\}}
\left(
\frac{1}{(z - n\tau - m)^2} - \frac{1}{(n\tau + m)^2}
\right)
$$
be the Weierstra\ss{} $\wp$--function.
Then $\wp'(\frac{1}{2}) = \wp'(\frac{\tau}{2}) = \wp'(\frac{1+ \tau}{2}) = 0$
and $\frac{1}{2}, \frac{\tau}{2}$ and $\frac{1+\tau}{2}$ are the only branch
points of $\wp(z)$ in the fundamental parallelogram of $\Lambda_\tau$.
Denote $e_1 = \wp(\frac{1}{2})$, $e_2 =  \wp(\frac{\tau}{2})$ and
$e_3 =  \wp(\frac{1+ \tau}{2})$. Then we have:
$$
\wp(z) - e_1 = \left(\frac{\cn(z)}{\sn(z)}\right)^2, \quad
\wp(z) - e_2 = \left(\frac{1}{\sn(z)}\right)^2, \quad
\wp(z) - e_3 = \left(\frac{\dn(z)}{\sn(z)}\right)^2.
$$
\end{remark}

\medskip

\section{Vector bundles on singular cubic curves}\label{S:triples}

To compute the associative $r$-matrices coming from the nodal and cuspidal
Weierstra\ss{} cubic curves, we use a description of vector bundles on
singular projective curves via the formalism of matrix problems
\cite{DrozdGreuel}, see also \cite{Thesis, Survey}. The purpose of this
section is to set up a clear language and provide the core technical tools
necessary for our applications.

\subsection{Description of vector bundles on singular curves}
\label{SS:catoftriples}

We start with recalling the general approach of Drozd and Greuel to study
torsion free sheaves on singular projective curves \cite{DrozdGreuel}. In our
applications, the normalization of the curve will always be rational.

\medskip
Let $X$ be a reduced singular (projective)  curve, $\pi:  \widetilde{X} \to  X$
its normalization, $\kI :=
{\mathcal Hom}_\kO\bigl(\pi_*(\kO_{\widetilde{X}}), \kO\bigr) =
{\mathcal A}nn_\kO\bigl(\pi_*(\kO_{\widetilde{X}})/\kO\bigr)$
the conductor ideal sheaf.
Denote  by $\eta: Z = V(\kI) \lar X$ the
closed artinian subspace  defined by $\kI$
(its topological support is precisely the singular locus of $X$)  and by
$\tilde\eta: \widetilde{Z} \lar \widetilde{X}$ its preimage in
$\widetilde{X}$, defined by the Cartesian  diagram
\begin{equation}\label{E:keydiag}
\xymatrix
{\widetilde{Z} \ar[r]^{\tilde{\eta}} \ar[d]_{\tilde{\pi}}
& \widetilde{X} \ar[d]^\pi \\
Z \ar[r]^\eta & X.
}
\end{equation}
The proof of the following lemma is straightforward.

\begin{lemma}
The  diagram (\ref{E:keydiag}) is also  a \emph{push-down} diagram.
Moreover,  denote $\nu = \eta \tilde\pi = \pi \tilde\eta$ and consider the
following natural transformations of functors:
$$
\left\{
\begin{array}{cccl}
\mathsf{j}: & \mathbbm{1}_X &  \lar &  \pi_* \pi^*, \\
\mathsf{q}: &  \mathbbm{1}_X &  \lar &  \eta_* \eta^*, \\
\mathsf{c}: &  \pi_* \pi^*   &  \lar &  \pi_* \tilde\eta_* \tilde\eta^*\pi^* \stackrel{\sim}\lar    \nu_* \nu^*, \\
\mm: &   \eta_* \eta^* &  \lar  &  \eta_* \tilde\pi_* \tilde\pi^*\eta^* \stackrel{\sim}\lar  \nu_* \nu^*.
\end{array}
\right.
$$
Then for any vector bundle $\kV$ on $X$ we have a short  exact sequence
$$
0 \lar \kV
\xrightarrow{
\left(
\begin{array}{r}
-\jj_\kV \\
\qq_\kV
\end{array}
\right)
}
\pi_* \pi^*\kV \oplus \eta_* \eta^*\kV
\xrightarrow{
\left(
\begin{array}{cc}
\cc_\kV  &
\mm_\kV
\end{array}
\right)
}
\nu_* \nu^*\kV \lar 0.
$$
\end{lemma}

\noindent
In  order to relate vector bundles on $X$ and $\widetilde{X}$ we use the
following definition.

\begin{definition}
  The category $\Tri(X)$ is defined as follows.
  \begin{itemize}
  \item Its objects are triples
    $\bigl(\widetilde\kV, \kN, \widetilde \mm\bigr)$, where
    $\widetilde\kV \in \VB(\widetilde{X})$, $\kN \in \VB(Z)$ and
    $$\widetilde \mm: \tilde{\pi}^*\kN \lar \tilde{\eta}^*\widetilde\kV$$ is an
    isomorphism of $\kO_{\widetilde{Z}}$-modules, called the \emph{gluing map}.
  \item The set of morphisms
    $\Hom_{\Tri(X)}\bigl((\widetilde\kV_1, \kN_1, \widetilde{\mm}_1),
    (\widetilde\kV_2, \kN_2, \widetilde{\mm}_2)\bigr)$ consists of all pairs
    $(F,f)$, where $F: \widetilde\kV_1 \to  \widetilde\kV_2$ and
    $f: \kN_1 \to \kN_2$
    are morphisms of vector bundles such that the following diagram is
    commutative
    $$
    \xymatrix
    {\tilde{\pi}^*\kN_1 \ar[r]^{\widetilde{\mm}_1}\ar[d]_{\tilde{\pi}^*(f)} &
      \tilde{\eta}^*\widetilde\kV_1
      \ar[d]^{\tilde{\eta}^*(F)} \\
      \tilde{\pi}^*\kN_2 \ar[r]^{\widetilde{\mm}_2} & \tilde{\eta}^*\widetilde\kV_2.
    }
    $$
  \end{itemize}
\end{definition}

\begin{remark}\label{R:tenstripel}
  The category $\Tri(X)$ is endowed with an interior tensor product:
  $$(\widetilde\kV_1, \kN_1, \widetilde{\mm}_1) \otimes (\widetilde\kV_2,
  \kN_2, \widetilde{\mm}_2) =
  (\widetilde\kV_1 \otimes \widetilde\kV_2, \kN_1 \otimes \kN_2,
  \widetilde{\mm}),$$
  where $\widetilde \mm$ is defined to be  the composition
  $$
  \tilde\pi^*(\kN_1 \otimes \kN_2) \stackrel{\cong}\lar
  \tilde\pi^*\kN_1 \otimes \tilde\pi^*\kN_2
  \xrightarrow{\widetilde{\mm}_1 \otimes \widetilde{\mm}_2}
  \tilde\eta^*\widetilde\kV_1 \otimes \tilde\eta^*\widetilde\kV_2
  \stackrel{\cong}\lar
  \tilde\eta^*(\widetilde\kV_1 \otimes \widetilde\kV_2).
  $$
  Similarly,  we define the functor
  $\det: \Tri_{n}(X) \lar \Tri_{1}(X)$, where $\Tri_{n}(X)$ denotes the full
  subcategory of $\Tri(X)$ whose objects
  $(\widetilde\kV, \kN, \widetilde{\mm})$ satisfy
  $\rk(\widetilde\kV)=\rk(\kN)=n$.
\end{remark}

\noindent
The following theorem summarizes main results about the category $\Tri(X)$ and
its relations with the category of vector bundles $\VB(X)$.

\begin{theorem}[Lemma 2.4 in \cite{DrozdGreuel} and
  Theorem 1.3 in \cite{Thesis}]\label{thm:Drozd-Greuel}
  Let $X$ be a reduced curve.

\vspace{1mm}

\noindent $\bullet$
Let  $\mathbb{F}: \VB(X) \lar \Tri(X)$ be the functor assigning  to a  vector
bundle $\kV$ the triple  $(\pi^*\kV, \eta^*\kV, \widetilde{\mm}_{\kV})$, where
$\widetilde{\mm}_{\kV}: \tilde{\pi}^*(\eta^*\kV) \lar \tilde\eta^*(\pi^*\kV)$ is
the canonical isomorphism. Then $\mathbb{F}$ is  an equivalence of categories.

\vspace{1mm}

\noindent $\bullet$ The functor  $\mathbb{F}$ commutes with tensor products:
we have a bifunctorial isomorphism
$$
\mathbb{F}(\kV_1 \otimes \kV_2) \stackrel{\cong}\lar \mathbb{F}(\kV_1) \otimes
\mathbb{F}(\kV_2).
$$
Moreover, we have an isomorphism  $\mathbb{F} \circ \det \stackrel{\cong}\lar
\det \circ \mathbb{F}$ of functors $\VB_{n}(E)\rightarrow\Tri_{1}(E)$, where
$\VB_{n}(E)$ denotes the category of vector bundles of fixed rank $n$.

\vspace{1mm}

\noindent $\bullet$ Let $\mathbb{G}: \Tri(X) \lar \Coh(X)$ be the functor
assigning to  a triple
$(\widetilde\kV, \kN, \widetilde{\mm})$ the coherent sheaf
$$
\kV := \ker\bigl(\pi_*\widetilde\kV \oplus  \eta_*\kN \xrightarrow{
(\cc \, \, \mm)
}
\nu_* \tilde\eta^*\widetilde\kV\bigr),
$$
where $\cc = \cc^{\widetilde\kV}$ is the canonical morphism
$
\pi_*\widetilde\kV \lar \pi_* \tilde\eta_* \tilde\eta^*\tilde\kV
= \nu_* \tilde\eta^*\widetilde\kV
$
and $\mm$ is the composition
$
\eta_*\kN \stackrel{\can}\lar \eta_* \tilde\pi_* \tilde\pi^*\kN \stackrel{=}\lar
\nu_* \tilde\pi^*\kN
\xrightarrow{\nu_*(\widetilde{\mm})} \nu_* \tilde\eta^*\widetilde\kV.
$
Then the coherent sheaf $\kV$ is locally free. Moreover, the functor
$\mathbb{G}$ is quasi-inverse to $\mathbb{F}$.
\item Being more precise, let $\kV
\xrightarrow{
  \left(\begin{smallmatrix}
    -\pp \\ \qq
  \end{smallmatrix}\right)
}
\pi_* \widetilde\kV \oplus \eta_* \kN$ be the canonical inclusion.
Then the morphisms $\widetilde{\pp}:
\pi^*\kV \stackrel{\pi^*(\pp)}\lar \pi^* \pi_* \widetilde\kV
\stackrel{\can}\lar \widetilde\kV$ and
$\widetilde{\qq}: \eta^*\kV \stackrel{\eta^*(\qq)}\lar
\eta^* \eta_* \kN \stackrel{\can}\lar \kN$ are isomorphisms and
$(\pi^*\kV, \eta^*\kV, \widetilde{\mm}_\kV)
\xrightarrow{(\widetilde{\pp},  \, \widetilde{\qq})} (\widetilde\kV, \kN,
\widetilde{\mm})$ is  an isomorphism in the category $\Tri(X)$.

\vspace{1mm}

\noindent $\bullet$
Let $\kT_i = (\widetilde\kV_i, \kN_i, \widetilde{\mm}_i)$, $i = 1,2$ be
objects of $\Tri(X)$ and $\kV_i = \mathbb{G}(\kT_i)$. Consider the short exact
sequences defining $\kV_i = \mathbb{G}(\kT_i)$:
$$
0 \lar \kV_i
\xrightarrow{
\left(
  \begin{smallmatrix}
    -\pp_i \\
    \qq_i
  \end{smallmatrix}
\right)
}
\pi_* \widetilde\kV_i \oplus \eta_* \kN_i
\xrightarrow{
\left(
  \begin{smallmatrix}
    \cc_i  &
    \mm_i
  \end{smallmatrix}
\right)
}
\nu_* \tilde\eta^*\widetilde\kV_i \lar 0.
$$
Then the sequence
$$
0 \lar \kV_1 \otimes \kV_2
\xrightarrow{
\left(
  \begin{smallmatrix}
    -\pp \\
    \qq
  \end{smallmatrix}
\right)
}
\pi_* (\widetilde\kV_1 \otimes \widetilde\kV_2)
\oplus \eta_*(\kN_1 \otimes \kN_2)
\xrightarrow{
\left(
  \begin{smallmatrix}
    \cc  &
    \mm
  \end{smallmatrix}
\right)
}
\nu_* \tilde\eta^*(\widetilde\kV_1 \otimes  \widetilde\kV_2)  \lar 0
$$
is exact, where
$$
\left\{
\begin{array}{cl}
  \pp: & \kV_1 \otimes \kV_2
  \xrightarrow{\pp_1 \otimes \pp_2}
  \pi_* \widetilde\kV_1 \otimes \pi_* \widetilde\kV_2
  \stackrel{\can}\lar  \pi_*(\widetilde\kV_1 \otimes \widetilde\kV_2) \\
  \qq: & \kV_1 \otimes \kV_2
  \xrightarrow{\qq_1 \otimes \qq_2}
  \eta_*\kN_1 \otimes \eta_* \kN_2
  \stackrel{\can}\lar \eta_*(\kN_1 \otimes \kN_2) \\
  \cc: & \pi_*(\widetilde\kV_1 \otimes \widetilde\kV_2)
  \xrightarrow{\pi_*(\can)}
  \pi_* \tilde\eta_* \tilde\eta^*(\widetilde\kV_1 \otimes \widetilde\kV_2)
  = \nu_* \tilde\eta^*(\widetilde\kV_1 \otimes \widetilde\kV_2) \\
  \mm: & \eta_*(\kN_1 \otimes \kN_2)
  \xrightarrow{\eta_*(\can)}
  \eta_* \tilde\pi_* \tilde\pi^*(\kN_1 \otimes \kN_2)
  \xrightarrow{\nu_*(\widetilde{\mm})}
  \nu_* \eta^*(\kV_1 \otimes \kV_2)\;,
\end{array}
\right.
$$
using $\widetilde{\mm}$ from Remark \ref{R:tenstripel}.
This means that
$\left(\begin{smallmatrix}-\pp\\ \qq\end{smallmatrix}\right)$ gives us a
bifunctorial isomorphism
\[
\alpha_{\kT_{1},\kT_{2}}:\mathbb{G}(\kT_{1})\otimes\mathbb{G}(\kT_{2})
\stackrel{\sim}{\lar} \mathbb{G}(\kT_{1}\otimes\kT_{2}).
\]

\vspace{1mm}

\noindent $\bullet$ Let $\For: \Tri(X) \lar \VB(\widetilde{X})$ be the
forgetful functor mapping a triple
$\kT = (\widetilde\kV, \kN, \widetilde{\mm})$ to $\widetilde\kV$.
Let $\kV = \mathbb{G}(\kT)$ and
$\gamma_\kT = \widetilde \pp: \pi^*\kV \lar \widetilde\kV$ be the isomorphism
introduced above.
Then we obtain   an isomorphism of functors
$\gamma: \pi^* \circ \mathbb{G} \lar \For$,
In particular, we have a commutative diagram:
$$
\xymatrix
{
  \Hom_{\Tri(X)}(\kT_1, \kT_2) \ar[rr]^{\mathbb{G}} \ar[d]_{\For} & &
  \Hom_X(\kV_1, \kV_2) \ar[d]^{\pi^*}\\
  \Hom_{\widetilde{X}}(\widetilde\kV_1, \widetilde\kV_2)  & &
  \Hom_{\widetilde{X}}(\pi^*\kV_1, \pi^*\kV_2). \ar[ll]_{\conj(\gamma_{\kT_1}, \gamma_{\kT_2})}
}
$$
Moreover, $\gamma$ is compatible with tensor products:
for $\kT_i = (\tilde\kV_i, \kN_i, \widetilde{\mm}_i)$ and
$\kV_i = \mathbb{G}(\kT_i)$ \, ($i = 1,2$) the diagram
$$
\xymatrix
{
  \pi^* \mathbb{G}(\kT_1) \otimes \pi^* \mathbb{G}(\kT_2)
  \ar[rr]^{\pi^{\ast}\left(\alpha_{\kT_{1},\kT_{2}}\right)}
  \ar[rd]_{\gamma_{\kT_1} \otimes \gamma_{\kT_2}} & &
  \pi^* \mathbb{G}(\kT_1 \otimes \kT_2)
  \ar[ld]^{\gamma_{\kT_1 \otimes \kT_2}}\\
  & \widetilde\kV_1 \otimes \widetilde\kV_2 &
}
$$
is commutative.

\end{theorem}

Our next goal is to obtain an explicit  description of stable vector bundles
on a singular Weierstra\ss{} curve $E$. In this context, we replace $X$ by $E$
and note that the normalisation is $\widetilde{E}\cong  \PP^1$.
In order to obtain a clearer description of objects of $\Tri(E)$, recall the
following well-known theorem.

\begin{theorem}[Birkhoff-Grothendieck]
On the projective line $\PP^1$, taking the degree gives an isomorphism
$\Pic(\PP^1) \cong \mathbb{Z}$. Any vector bundle
$\kE$ on $\PP^1$   splits into a direct sum of line  bundles:
$
\kE \cong \oplus_{n \in \mathbb{Z}} \kO_{\PP^1}(n)^{m_n}.
$
\end{theorem}

\noindent
This implies that if $(\widetilde{\kV}, \kN, \widetilde{\mm})$ is an
object of $\Tri(E)$ with $\rk(\widetilde{\kV}) = n$, we have
\[\widetilde{\kV} = \bigoplus\limits_{l \in \mathbb{Z}} \kO_{\PP^1}(l)^{k_l}
\quad \text{ and } \quad \quad \kN \cong \kO_{Z}^n, \quad
\text { where } \sum\limits_{l \in \mathbb{Z}} k_l = n.\]
Note that $\kN$ is in fact free, because $Z$ is artinian. From now on we shall
always fix a decomposition of $\widetilde{\kV}$ as above.

An explicit description of morphisms between objects in $\Tri(E)$ requires to
choose coordinates on $\PP^{1}$.
Let $(z_0, z_1)$ be coordinates on $V = \CC^2$. They induce homogeneous
coordinates $(z_0:z_1)$ on the projective line
$\PP^1(V) = (V \setminus\{0\})/\sim$, where $v \sim \lambda v$ for all
$\lambda \in \CC^*$.

We set
$U_0 =\{(z_0: z_1)| z_0 \ne 0\}$ and $U_\infty  =\{(z_0: z_1)| z_1 \ne 0\}$
and put $0 := (1: 0)$, $\infty := (0: 1)$, $z = z_1/z_0$ and $w = z_0/z_1$.
So, $z$ is a   coordinate in a  neighbourhood of $0$.
If $U = U_0 \cap U_\infty$ and $w=1/z$ is used as a coordinate on
$U_{\infty}$, then the transition function of the line bundle $\kO_{\PP^1}(n)$ is
$$
 U_0 \times \mathbb{C} \supset  U \times \mathbb{C}
\xrightarrow{(z,v) \mapsto \left(\frac{1}{z}, \frac{v}{z^n}\right)}
U \times \mathbb{C} \subset  U_\infty \times \mathbb{C}.
$$
The vector bundle $\kO_{\PP^1}(-1)$  is isomorphic to the sheaf of sections
of  the so-called tautological line bundle
$$\bigl\{(l, v) | \, v \in l \bigr\} \subset \PP^1(V) \times V = \kO_{\PP^1}^2.$$
The choice of coordinates on $\PP^1$ fixes two distinguished elements, $z_0$
and $z_1$, in the space $\Hom_{\PP^1}\bigl(\kO_{\PP^1}(-1), \kO_{\PP^1}\bigr)$:
$$
\xymatrix
{ \PP^1 \times \mathbb{C}^2 \ar[rd] \ar@{<-^{)}}[r] &
                   \kO_{\PP^1}(-1) \ar[d] \ar[r]^{z_i} &
                   \PP^1 \times \mathbb{C} \ar[ld] \\
& \PP^1 &
}
$$
where $z_i$ maps $\bigl(l, (v_0, v_1)\bigr)$ to $(l, v_i)$ for  $i = 0,1$.
It is clear that the section $z_0$ vanishes at $\infty$ and $z_1$
vanishes at $0$.  After having made this choice,  we may write
$$
\Hom_{\PP^1}\bigl(\kO_{\PP^1}(n), \kO_{\PP^1}(m)\bigr) =
 \CC[z_0, z_1]_{m-n} :=
\bigl\langle z_0^{m-n}, z_0^{m-n -1}z_1,\dots,z_1^{m-n}\bigr\rangle_\CC.
$$

\hspace{2mm}

\begin{lemma}\label{L:deg-of-vb}
Let $E$ be a singular Weierstra\ss{} curve,
$\pi: \PP^1 \to E$ its normalization and  $\kV$ a vector bundle on $E$. Then
$\deg_E(\kV) = \deg_{\PP^1}(\pi^*\kV)$.
\end{lemma}

\noindent
\begin{proof} If $n = \rk(\kV)$ then $\pi^*\kV$ is a vector bundle of rank $n$
on $\PP^1$. The canonical morphism  $g: \kV \rightarrow   \pi_{*} \pi^* \kV$
is generically injective and $\kV$ is torsion free, hence $\ker(g) = 0$ and we
have an exact sequence
$$
0 \lar \kV \stackrel{g}\lar \pi_{*} \pi^* \kV \lar \kS \lar 0,
$$
where $\kS$ is a torsion sheaf supported at the singular point $s$ of the
curve $E$.
Since $g$ commutes with restrictions to an open set, we have
$$\kS \cong \bigl(\coker(\kO_E \to \pi_*(\kO_{\PP^1}))\bigr)^n.$$
Because $s$ is either a node or a cusp, we obtain $h^0(\kS) = n$.
Using the Riemann-Roch formula, this implies
$$
\deg_E(\kV) = \chi(\kV) = \chi(\pi_{*} \pi^* \kV) - \chi(\kS) =
\chi(\pi^*\kV) - n = \deg_{\PP^1}(\pi^*\kV).
$$
\end{proof}

\begin{lemma}\label{L:normofsimple}
Let $E$ be a singular Weierstra\ss{} curve,
$\pi: \PP^1 \to E$ its normalization and  $\kV$ a simple vector bundle of rank
$n$ on $E$. Then
\begin{itemize}
\item $\kV$ is stable.
\item $\pi^*\kV \cong  \kO_{\PP^1}(c)^{n_1} \oplus \kO_{\PP^1}(c+1)^{n_2}$ for some
integer $c \in \mathbb{Z}$ and some non-negative integers $n_1, n_2$ which
satisfy $n = n_1 + n_2$.
\end{itemize}
\end{lemma}

\begin{proof} For the first statement see for example
\cite[Corollary 4.5]{BK3}.
To prove the second part, let $\mathbb{F}(\kV) = (\widetilde\kV, \kO_{Z}^n,
\widetilde{\mm})$ and assume
$$
\pi^* \kV \cong \widetilde\kV = \kO_{\PP^1}(c) \oplus \kO_{\PP^1}(d) \oplus
\widetilde\kV'',
$$
where $ d - c \ge 2$.
Because the length of $\widetilde{Z}$ is two, we can find a non-zero
homogeneous form
$p = p(z_0, z_1) \in \Hom_{\PP^1}\bigl(\kO_{\PP^1}(c), \kO_{\PP^1}(d)\bigr)$
such that $\tilde\eta^*(p) = 0$.
This gives us  a non-scalar endomorphism of $\kV$  corresponding to the
endomorphism  $(F, f)$ of the triple $\mathbb{F}(\kV)$ given by $f = \id$ and
$$
F =
\begin{pmatrix}
1 & 0 & 0 \\
p & 1 & 0 \\
0 & 0 & 1
\end{pmatrix}.
$$
This contradicts our assumption that $\kV$ was simple.
\end{proof}

\medskip

An explicit description of the morphism $\widetilde{\mm}$ by a matrix requires
to fix isomorphisms $\zeta_l: \tilde\eta^{*}\kO_{\PP^1}(l)  \to  \kO_{\widetilde{Z}}$.
In order to keep compatibility with tensor products in our description of
vector bundles, we have to ensure that for all $k,l \in \mathbb{Z}$
the following diagram is commutative:
\begin{equation}
  \label{eq:compatiblezeta}
  \xymatrix
  {
    \tilde\eta^{*}\kO_{\PP^1}(k) \otimes \tilde\eta^{*}\kO_{\PP^1}(l) \ar[rr]^\can
    \ar[d]_{\zeta_k \otimes \zeta_l} & & \tilde\eta^{*}\kO_{\PP^1}(k+l)
    \ar[d]^{\zeta_{k+l}}\\
    \kO_{\widetilde{Z}} \otimes \kO_{\widetilde{Z}} \ar[rr]^{\mathrm{mult}} & &
    \kO_{\widetilde{Z}}.
  }
\end{equation}
In the case of a nodal or cuspidal Weierstra\ss{}
cubic curve  we shall explicitly give our choice of these
isomorphisms.

\begin{remark}\label{R:singvssmooth}
It is natural to assume $\zeta_0 = \mathsf{id}$. Then such a family of
isomorphisms $\{\zeta_l\}_{l \in \mathbb{Z}}$ is uniquely determined by
$\zeta = \zeta_1: \tilde\eta^*\bigl(\kO_{\mathbb{P}^1}(1)\bigr)
\lar \kO_{\widetilde{Z}}$.
Moreover, the choice of a global section
$p = p_\zeta = a z_0 + b z_1 \in H^0\bigl(\kO_{\mathbb{P}^1}(1)\bigr)$, which
does not vanish on $\widetilde{Z}$, determines $\zeta$ as follows:
$\zeta(s)=\left.\frac{s}{p}\right|_{\widetilde{Z}}$. Modulo automorphisms of
$\mathbb{P}^1$ such a section $p$ is determined by its unique zero, which
should belong to $\mathbb{P}^1 \setminus \widetilde{Z} \cong \breve{E}$. In
other words, our choice of a set of trivializations
$\{\zeta_l\}_{l \in \mathbb{Z}}$  corresponds to the choice of a smooth point in
Atiyah's classification of vector bundles on an elliptic curve \cite{Atiyah}.
\end{remark}

Note that, because we have fixed a decomposition
$\widetilde{\kV}=\bigoplus_{l \in \mathbb{Z}} \kO_{\PP^1}(l)^{k_l}$, a
family $\{\zeta_l\}_{l \in \mathbb{Z}}$ induces  an isomorphism
$\zeta^{\widetilde\kV}: \tilde\eta^{*}\widetilde\kV \to \kO_{\widetilde{Z}}^n$.
Because $\kN \cong \kO_Z^n$, we also get an isomorphism
$\tilde{\pi}^* \kN \cong \kO_{\widetilde{Z}}^n$.
This allows us to describe
the map $\widetilde{\mm}: \tilde{\pi}^*\kN \lar  \tilde\eta^*\widetilde\kV$ as
a matrix in $\GL_n(\kO_{\widetilde{Z}})$.

\begin{corollary}\label{C:choiceoftriv}
Let $\Mat_{\widetilde{Z}}$ be the category of square matrices over the ring
$\kO_{\widetilde{Z}}$.
The choice of isomorphisms $\{\zeta_l\}_{l \in \mathbb{Z}}$ yields a functor
$
\mathbb{P}^\zeta: \Tri(E) \lar \Mat_{\widetilde{Z}},
$
assigning to a triple $(\widetilde\kV, \kO_{Z}^n, \widetilde{\mm})$ the matrix
of the $\kO_{\widetilde{Z}}$-linear map
$$
\kO_{\widetilde{Z}}^n \stackrel{\widetilde{\mm}}\lar \tilde\eta^*\widetilde\kV
\stackrel{\zeta^{\widetilde\kV}}\lar \kO_{\widetilde{Z}}^n.
$$
Moreover, let $\mathbb{H}^\zeta
= \mathbb{P}^\zeta \circ \mathbb{F}: \VB(E) \lar \Mat_{\widetilde{Z}}$.
Using \eqref{eq:compatiblezeta}, for any
$\kL \in \Pic(E)$ and $\kV \in \VB(E)$ we obtain:
$$
\mathbb{H}^\zeta(\kL \otimes \kV)
= \mathbb{H}^\zeta(\kL) \cdot \mathbb{H}^\zeta(\kV)
\quad \mbox{\textrm{and}} \quad \mathbb{H}^\zeta\bigl(\det(\kV)\bigr) =
\det\bigl(\mathbb{H}^\zeta(\kV)\bigr).
$$
\end{corollary}

\noindent
Let $(\widetilde\kV, \kN, \widetilde\mm)$ be an object of $\Tri(E)$.
We have a natural action of the group
$\Aut_{\PP^1}(\widetilde\kV) \times \Aut_{Z}(\kN)$ on the vector space
$\Hom_{\tilde{Z}}(\tilde{\pi}^*\kN, \tilde\eta^*\widetilde\kV)$. The orbits of
this action correspond precisely to the points in the fibre of the functor
$\pi^*: \VB(E) \to \VB(\PP^{1})$ over $\widetilde\kV$.
In what follows, we shall use this action to find a normal form for
$\widetilde{\mm}$. A description of the matrix problem describing \emph{all}
vector bundles on an irreducible Weierstra\ss{} cubic curve, can be found in
\cite{DrozdGreuel} and \cite{Survey}.  In this article, we are mainly
interested in a description of \emph{simple} vector bundles. Having in mind
Lemma \ref{L:normofsimple}, we introduce the following notation.

In order to recover a vector bundle $\kV$ from the matrix
$\mathbb{H}^{\zeta}(\kV)$, we need to specify $\widetilde{\kV}$.
For a singular Weierstra\ss{} cubic curve $E$,  let $\VB^{(0,1)}(E)$ be the
full subcategory of $\VB(E)$ consisting of vector bundles $\kV$ such that
$\pi^*\kV \cong  \kO_{\PP^1}^{n_1} \oplus \kO_{\PP^1}(1)^{n_2}$ for some
non-negative  integers  $n_1, n_2  \in \mathbb{Z}$. In a similar way, let
$\Tri^{(0,1)}(E)$ be the corresponding subcategory of the category $\Tri(E)$.

\begin{definition}\label{D:categoryBM}
  Let $E$ be a Weierstra\ss{}  cubic curve $E$. Consider the following
  category $\BM(E)$ of ``block matrices'':
  \begin{itemize}
  \item Its objects are invertible  matrices over the ring
    $\kO_{\widetilde{Z}}$ with a block structure:
    $$
    M =
    \left(
      \begin{array}{c|c}
        M_{00} & M_{01} \\
        \hline
        M_{10} & M_{11}
      \end{array}
    \right),
    $$
    where $M_{00}$ and $M_{11}$ are square matrices, possibly of size zero.

  \item Let $M$ and $N$ be two objects of $\BM(E)$ of sizes $m = m_0 + m_1$ and
    $n = n_0 + n_1$ respectively, where the block $M_{ij}$ has size
    $m_i \times m_j$ etc.
    Then a morphism from $M$ to $N$ in the category $\BM(E)$
    is given by a pair of matrices $(F, f)$, where
    $f \in \Mat_{n \times m}(\kO_Z)$ and
    $$
    F =
    \left(
      \begin{array}{c|c}
        F_{00} & 0 \\
        \hline
        F_{10} & F_{11}
      \end{array}
    \right)
    $$
    has blocks
    $F_{00} \in \Mat_{n_0 \times m_0}(\CC)$, $F_{11} \in \Mat_{n_1 \times m_1}(\CC)$
    and $F_{10} \in \Mat_{n_1 \times m_0}(\kO_{\widetilde Z})$, such that
    $F M = N \tilde{f}.$ Here $\tilde{f}$ is the image of the matrix $f$ under
    the morphism
    $\Mat_{n \times m}(\kO_Z) \lar \Mat_{n \times m}(\kO_{\widetilde{Z}})$
    induced by the ring homomorphism  $\kO_Z \lar \kO_{\widetilde{Z}}$.
  \item The composition of morphisms in $\BM(E)$ is given  by the matrix
    product.
  \end{itemize}
\end{definition}

\begin{proposition}\label{P:firstred}
  Take some isomorphism   $\zeta: \kO_{\mathbb{P}^1}(1)|_{\widetilde{Z}} \lar
  \kO_{\widetilde{Z}}$. Then in the notation of Remark
  \ref{R:singvssmooth} and Corollary \ref{C:choiceoftriv},  we have
  equivalences of categories:
  $$
  \VB^{(0,1)}(E) \stackrel{\mathbb{F}}\lar \Tri^{(0,1)}(E)
  \stackrel{\mathbb{P}^\zeta}\lar \BM(E)\;,
  $$
  with block structure on $\mathbb{P}^\zeta\bigl(\widetilde\kV,
    \kO_{Z}^{n}, \widetilde{\mm}\bigr)$ coming from the decomposition
  $\widetilde\kV = \kO_{\PP^{1}}^{n_{1}} \oplus \kO_{\PP^{1}}(1)^{n_{2}}$.
  Moreover, the functor $\mathbb{P}^\zeta \circ \mathbb{F}$ sends
  $\det(\kV) \in \VB(E)$ to the determinant of the
  corresponding matrix $\mathbb{P}^\zeta \bigl(\mathbb{F}(\kV)\bigr)$.
\end{proposition}

\begin{proof}
This result follows from Theorem \ref{thm:Drozd-Greuel} and the observation
that the map
$$
\tilde\eta^*:
 \Hom_{\PP^1}\bigl(\kO_{\PP^1}, \kO_{\PP^1}(1)\bigr) \lar
 \Hom_{\widetilde Z}(\kO_{\widetilde Z}, \kO_{\widetilde Z})
$$
is an isomorphism both for a nodal and a cuspidal cubic curve.
\end{proof}

\medskip

\subsection{Simple vector bundles on a nodal Weierstra\ss{} curve}\label{SS:node}
The main aim of this subsection is an explicit description of those objects in
$\Tri(E)$ which correspond to simple vector bundles on a nodal Weierstra\ss{}
curve $E$.
We give an algorithm which produces some kind of normal form of such triples
for each given rank and degree. Crucial for our application to the Yang--Baxter
equation is a description of the family of all simple vector bundles with
fixed rank and degree in a way which is compatible with the action of the
Jacobian.
We shall also see that rank and degree of a simple vector bundle on a nodal
Weierstra{\ss} curve are always coprime.

Let $E$ be a nodal Weierstra\ss{} curve, e.g.\ given by $zy^2 = x^3 + x^2 z$,
$s = (0:0:1)$ the singular point and $\pi: \PP^1 \lar E$ its
normalization.
Choose homogeneous coordinates $(z_0: z_1)$ on $\PP^1$ in such a way that
$\pi^{-1}(s) = \{0,\infty \}$.
Then, in notations  of the previous subsection, $Z$ and $\widetilde{Z}$ are
reduced complex spaces as follows
$$Z =\{s\} \quad \mbox{\textrm{and}} \quad
\widetilde{Z} = \{0\} \cup \{\infty \}.$$
Hence, for $\bigl(\tilde\kV, \kN, \widetilde{\mm}\bigr)\in\Tri(E)$ the map
$\widetilde{\mm}$ is just an isomorphism of $\CC\times \CC$-modules, i.e. it
is given by a pair of invertible matrices $\mm(0)$ and $\mm(\infty)$.

The linear form
$p = p_\zeta(z_0, z_1) = z_1 - z_0 \in H^0\bigl(\kO_{\PP^1}(1)\bigr)$ does not
vanish on $\widetilde{Z}$. Following the recipe from Remark
\ref{R:singvssmooth}, we  consider  the collection  of isomorphisms
$
\zeta_l:  \tilde\eta^* \kO_{\PP^1}(l) \lar \kO_{\widetilde Z}
$
given by the formula
$\left.\zeta_l(s) = \dfrac{s}{p^l}\right|_{\widetilde{Z}}$
for each open subset $V \subset \PP^1$ not containing $(1:1)$,  $l \in
\mathbb{Z}$ and any $s \in \Gamma\bigl(V, \kO_{\PP^1}(l)\bigr)$.

This implies the following evaluation rule for morphisms of vector bundles on
$\PP^1$:
if $q = a_0 z_0^{m-n} + a_1 z_0^{m-n-1} z_1 + \dots +  a_{m-n} z_1^{m-n} \in
\Hom_{\PP^1}\bigl(\kO_{\PP^1}(n), \kO_{\PP^1}(m)\bigr)$ then  we have a
commutative diagram
$$
\xymatrix
{
\tilde\eta^*\kO_{\PP^1}(n)
\ar[rrrr]^{\tilde\eta^*(q)} \ar[d]_{\zeta_n} & & & &
\tilde\eta^*\kO_{\PP^1}(m) \ar[d]^{\zeta_m}\\
\CC_0 \oplus \CC_\infty
\ar[rrrr]^{
\left(
\begin{smallmatrix}(-1)^{m-n} a_{0}&0\\0&a_{m-n}\end{smallmatrix}
\right)
}
& & & & \CC_0 \oplus \CC_\infty.
}
$$
If the family $\{\zeta_{l}\}$ is understood, we shall often write
$\widetilde\eta^{\ast}(q) =
\left(\bigl((-1)^{m-n}a_{0}\bigr), \bigl(a_{m-n}\bigr)\right)$.
\medskip

Our next goal is to describe the category $\BM(E)$ from Definition
\ref{D:categoryBM}. An object of $\BM(E)$ is a pair of matrices  $\mm(0)$ and
$\mm(\infty)$ simultaneously divided into blocks
$$
\mm(0) =
\left(
\begin{array}{c|c}
M_{00}(0) & M_{01}(0) \\
\hline
M_{10}(0) & M_{11}(0)
\end{array}
\right), \quad
\mm(\infty) =
\left(
\begin{array}{c|c}
M_{00}(\infty) & M_{01}(\infty) \\
\hline
M_{10}(\infty) & M_{11}(\infty)
\end{array}
\right).
$$
Two objects $\bigl(\mm(0), \mm(\infty)\bigr)$ and
$\bigl(\mm'(0), \mm'(\infty)\bigr)$ of  $\BM(E)$ are isomorphic if and only if
the corresponding blocks have the same sizes and there exist matrices
$$
F(0) =     \left(
  \begin{array}{c|c}
    F_{11}       & 0 \\
    \hline
    F_{21}(0) & F_{22}
  \end{array}
\right), \quad
 F(\infty) =   \left(
   \begin{array}{c|c}
     F_{11} & 0 \\
     \hline
     F_{21}(\infty) & F_{22}
   \end{array}
 \right)
$$
and $f$ such that
$$
F(0) \mm(0) = \mm'(0) f, \quad F(\infty) \mm(\infty) = \mm'(\infty) f.
$$
In particular, we have the following  isomorphism in the category $\BM(E)$:
$$
\bigl(\mm(0), \mm(\infty)\bigr) \cong
\bigl(\mm(0) \mm(\infty)^{-1}, \mathsf{id}\bigr)
$$
i.e.\ without loss of generality we may assume that the second matrix
$\mm(\infty)$ is the identity matrix.

To illustrate how the explicit identification of vector bundles on $E$ and
objects of $\BM(E)$ works in practice, we
shall now consider the simplest interesting case: the description of
$\Pic^{1}(E)$. We explicitly determine for each $y \in \breve{E}$ the
object in $\Tri(E)$ which corresponds to the line bundle $\kO_E(y)$.

The chosen coordinates provide us with an isomorphism
$\CC^{\ast} \cong U :=\PP^{1}\setminus\{0,\infty\}$ mapping $y \in \CC^{\ast}$ to
$(1: y) \in U$.
As $E$ is nodal, the normalization restricts to an isomorphism
$\pi:\PP^{1}\setminus\{0,\infty\} \rightarrow \breve{E}$.
Together, this gives us an identification $\breve{E} \cong \CC^*$,
under which $y\in\CC^{\ast}$ corresponds to $\tilde{y} := \pi^{-1}(y) =
(1:y)\in\PP^{1}$.
Obviously, $\pi^*\bigl(\kO_E(y)\bigr) = \kO_{\PP^1}(\tilde{y}) \cong
\kO_{\PP^1}(1)$ and the following lemma is true.

\begin{lemma}\label{L:idenonnode}
For the given choice of homogeneous coordinates on $\PP^1$ and
the set of trivializations  $\{\zeta_l\}_{l \in \mathbb{Z}}$ described above, we
obtain for all $y\in \breve{E} \cong \CC^*$
$$
\mathbb{F}\bigl(\kO_E(y)\bigr) =
\bigl(\kO_{\PP^1}(1), \CC_s, \bigl((y), (1)\bigr)\bigr).
$$
\end{lemma}

\noindent
\begin{proof}
Assume
$\kT_{\kO_E(y)} := \mathbb{F}\bigl(\kO_E(y)\bigr)
=  \bigl(\kO_{\PP^1}(1), \CC_s, \bigl((\lambda), (1)\bigr)\bigr)$.
It is clear that
$\kT_{\kO_E} := \mathbb{F}(\kO_E)
=  \bigl(\kO_{\PP^1}, \CC_s, \bigl((1), (1)\bigr)\bigr)$. Moreover, by
Theorem \ref{thm:Drozd-Greuel} we have a commutative diagram
$$
 \xymatrix
 {
 \Hom_{\Tri(E)}(\kT_{\kO_{E}}, \kT_{\kO_E(y)}) \ar[rrr]^{\mathbb{G}}
 \ar[d]_{\For} & & &
 \Hom_{E}\bigl(\kO_E,  \kO_E(y)\bigr) \ar[d]^{\pi^*}\\
 \Hom_{\PP^1}\bigl(\kO_{\PP^1},  \kO_{\PP^1}(1)\bigr)  & & &
 \Hom_{\PP^1}\bigl(\kO_{\PP^1},  \kO_{\PP^1}(\tilde y)\bigr).
 \ar[lll]_{\conj\left(\gamma_{\kT_\kO}, \gamma_{\kT_{\kO(y)}}\right)}
 }
$$
The section
$z_1 - yz_0 \in \Hom_{\PP^1}\bigl(\kO_{\PP^1}, \kO_{\PP^1}(1)\bigr)$
generates the image of $\pi^{\ast}$, hence belongs to the image of
$\For$. Using the description of morphisms in the category $\Tri(E)$  and the
evaluation rule  $\tilde\eta^*(z_1 - y z_0) = \bigl((y), (1)\bigr)$,  this is
equivalent to  the existence of  a constant $c \in \CC^*$ making  the
following diagram commutative:
$$
\xymatrix
{\CC_0 \oplus \CC_\infty
\ar[rr]^{\begin{pmatrix} c & 0 \\ 0 & c \end{pmatrix}}
\ar[dd]_{\begin{pmatrix} 1 & 0 \\ 0 & 1 \end{pmatrix}} & &
\CC_0 \oplus \CC_\infty
\ar[dd]^{\begin{pmatrix} \lambda & 0  \\ 0  & 1 \end{pmatrix}} \\
 & & \\
\CC_0 \oplus \CC_\infty
\ar[rr]^{\begin{pmatrix} y & 0  \\ 0   & 1 \end{pmatrix}} & &
\CC_0 \oplus \CC_\infty.
}
$$
This implies  that $\lambda = y$ and
$\mathbb{F}\bigl(\kO_E(y)\bigr)
\cong \bigl(\kO_{\PP^1}(1),\CC_s,\bigl((y),(1)\bigr)\bigr)$.
\end{proof}

\noindent
Our next goal is to describe the so-called \emph{Atiyah bundles}.

\begin{lemma}\label{lem:semi-stable-on-nodal}
Let $E$ be a nodal Weierstra\ss{} curve. Then there exists a unique
indecomposable semi-stable vector bundle $\kA_n$ of rank $n$ and degree $0$
such that all its  Jordan-H\"older factors are isomorphic to $\kO_E$. This
vector bundle is called the Atiyah bundle of rank $n$  and is given by the
triple $\bigl(\kO_{\PP^1}^n, \CC^n_s, \widetilde{\mm}\bigr),$ where
$$
\mm(0) = J_m(1) =
\left(
\begin{array}{ccccc}
1 & 1 & 0 & \dots & 0 \\
0 & 1 & 1 & \dots & 0 \\
\vdots & \vdots & \ddots & \ddots & \vdots  \\
0 & 0 & \dots & 1 & 1 \\
0 & 0 & \dots & 0 & 1 \\
\end{array}
\right)
\quad
\mm(\infty) = I_m =
\left(
\begin{array}{ccccc}
1 & 0 & 0 & \dots & 0 \\
0 & 1 & 0 & \dots & 0 \\
\vdots & \vdots & \ddots & \ddots & \vdots  \\
0 & 0 & \dots & 1 & 0 \\
0 & 0 & \dots & 0 & 1 \\
\end{array}
\right).
$$
\end{lemma}

\begin{proof} The category of semi-stable vector bundles with the
Jordan-H\"older factor
$\kO_E$ is equivalent to the category of finite-dimensional modules
over $\CC[[t]]$, see for example
\cite[Theorem 1.1 and Lemma 1.7]{FMW}.
Therefore,  there exists a unique indecomposable vector bundle
$\kA_n$ of rank $n$ recursively defined by the non-split exact sequences
$$
0 \lar \kA_n \lar \kA_{n+1} \lar \kO_E \lar 0 \quad\text{ and }\quad \kA_1 = \kO_E.
$$
In order to get a description of $\kA_n$ in terms of triples, first observe
that $\pi^*\kA_n \cong \kO_{\PP^1}^n$, hence  $\mathbb{F}(\kA_n) =
\bigl(\kO_{\PP^1}^n, \CC_s^n, \widetilde{\mm}\bigr)$.
The morphism $\widetilde{\mm}$ is given by two invertible
matrices $\mm(0), \mm(\infty) \in \GL_n(\CC)$. If $\widetilde{\mm}' =
\bigl(\mm'(0), \mm'(\infty)\bigr)$ is another pair such that
$$
\mm'(0) = S^{-1} \mm(0) T, \quad \mm'(\infty) = S^{-1} \mm(\infty) T
$$
with $S,T\in\GL_n(\CC)$,
then $\bigl(\kO_{\PP^1}^n, \CC_s^n, \widetilde{\mm}'\bigr)$ and $\bigl(\kO_{\PP^1}^n,
\CC_s^n, \widetilde{\mm}\bigr)$ define isomorphic  vector bundles on $E$.
We may, therefore, assume $\mm(\infty) = I_n$. Keeping $\mm(\infty) = I_n$ unchanged, the matrix
$\mm(0)$ can still be transformed to $S^{-1} \mm(0) S$. Hence,  $\mm(0)$ splits
into a direct sum of Jordan blocks. Since the vector bundle  $\kA_n$ is
indecomposable, $\mm(0) \sim J_n(\lambda)$ for some $\lambda \in \CC^*$. From
the condition $\Hom_E(\kA_n, \kO) = \CC$ one can  easily deduce $\lambda = 1$.
\end{proof}
\medskip

Now we start to focus on simple vector bundles on a nodal Weierstra\ss{} curve
$E$. We aim at giving a canonical form for those elements
$(\widetilde\kV, \kN, \widetilde{\mm})$ in $\Tri(E)$ which correspond to simple vector
bundles on $E$ under the functor $\mathbb{F}$ from Theorem
\ref{thm:Drozd-Greuel}.

\begin{definition}\label{D:MPnode}
Let $E$ be a nodal cubic curve and $n_{1}>0$, $n_{2}\ge0$ integers.
The category  $\MPnode(n_{1},n_{2})$  is defined as follows.
\begin{itemize}
\item Its objects are invertible matrices with blocks
  $M_{ij}\in\Mat_{n_{i}\times n_{j}}(\CC)$
  $$
  M =
  \left(
    \begin{array}{c|c}
      M_{11} & M_{12} \\
      \hline
      M_{21} & M_{22}
    \end{array}
  \right).
  $$
\item Morphisms are pairs of block matrices
  $$
  \Hom_{\MPnode(n_{1},n_{2})}(M, M') =
  \left\{ (S,T) \mid S M = M'  T \right\},
  $$
  with obvious composition and such that
  $$
  S =
  \left(
    \begin{array}{cc}
      A & 0 \\
      C' & B
    \end{array}
  \right) \quad\text{ and }\quad
  T =
  \left(
    \begin{array}{cc}
      A & 0 \\
      C'' & B
    \end{array}
  \right)
  $$
  have blocks of the same size as the blocks of $M$ and $M'$.
\end{itemize}
By $\MPnode^s(n_{1},n_{2})$ we denote the full subcategory of
\emph{simple}\footnote{Simple objects of $\MPnode(n_{1},n_{2})$
are by  definition the objects having  only scalar endomorphisms.
They are sometimes  called \emph{Schurian} objects or \emph{bricks}.}
objects of $\MPnode(n_{1},n_{2})$.
\end{definition}

\noindent
The proof of the following lemma is straightforward.
\begin{lemma}\label{L:VBequivMP}
Let $\VB^{(0,1)}_{n_1, n_2}(E)$ be the category of vector bundles $\kV \in \VB(E)$ such that
$\pi^*\kV \cong \kO_{\PP^1}^{n_1} \oplus \kO_{\PP^1}(1)^{n_2}$. Then $\VB^{(0,1)}_{n_1, n_2}(E)$
and $\MPnode(n_{1},n_{2})$ are equivalent.
\end{lemma}

\begin{proof} Let $\BM_{n_1, n_2}(E)$ be the full subcategory of $\BM(E)$
  consisting of matrices, whose diagonal blocks have sizes $n_1 \times n_1$
  and $n_2 \times n_2$. By Proposition \ref{P:firstred} the categories
  $\VB^{(0,1)}_{n_1, n_2}(E)$ and $\BM_{n_1, n_2}(E)$ are equivalent.
  But it is easy to see that sending and  $M\in\MPnode(n_{1},n_{2})$ to
  $(M,\id)\in\GL_{n_{1}+n_{2}}(\kO_{\widetilde{Z}})$ as an object in
  $\BM_{n_1, n_2}(E)$ with same block structure, is an equivalence.
\end{proof}

\begin{remark}\label{R:degree-zero}
  If $n_{2}=0$ the block structure becomes invisible and we end up in a
  situation of elementary linear algebra. Indeed, we have
  $\MPnode(n,0)=\GL_{n}(\CC)$ and
  $\Hom_{\MPnode(n,0)}(M, M') = \left\{ S \mid S M = M' S\right\}$.
  The indecomposable objects in this category are precisely
  those which are isomorphic to a Jordan block $J_{n}(\lambda)$,
  $\lambda\in\CC^{\ast}$. The endomorphism ring of $J_{n}(\lambda)$ is
  isomorphic to $\CC[t]/t^{n}$.
  Hence, $\MPnode^{s}(n,0)=\emptyset$ if $n>1$ and
  $\MPnode^{s}(1,0)=\GL_{1}(\CC)=\CC^{\ast}$.
\end{remark}

We aim now at finding a canonical form for objects
$M\in\MPnode^{s}(n_{1},n_{2})$. This means that we wish to find  in each
isomorphism class of $\MPnode^{s}(n_{1},n_{2})$ a unique object with a
particularly ``simple'' structure.
We shall often say that we can ``reduce'' a matrix $M$ to a matrix $N$ if
$M$ and $N$ are isomorphic in $\MPnode^{s}(n_{1},n_{2})$ and $N$ has a ``simpler''
form than $M$.
The reduction procedure described below is based on the following easy lemma.
\begin{lemma}\label{L:invblock}
  The block $M_{12}$ has  full rank, if $M \in \MPnode^s(n_1, n_2)$ is simple.
\end{lemma}

\begin{proof}
  If the matrix $M_{12}$ does not have full rank, $M$ can be reduced to the form
$$
M =
\left(
\begin{array}{cc|cc}
M_1 & M_2 & 0 & 0 \\
0   &  0  & I & 0 \\
\hline
M_3 & M_4 & 0 & M_5 \\
M_6 & M_7 & 0 & M_8
\end{array}
\right),
$$
where $M_{11}, M_{21}$ and $M_{22}$ are split into blocks such that $M_1$
and $M_8$ are square matrices. As an object of $\MPnode(n_{1},n_{2})$, such a
matrix has an endomorphism $(S,T)$ with
$$
S =
\left(
\begin{array}{cc|cc}
I  & 0  & 0 & 0 \\
0   &  I   & 0 & 0 \\
\hline
M_5 & 0  & I  & 0   \\
M_8 & 0 & 0 & I
\end{array}
\right) \quad
T =
\left(
\begin{array}{cc|cc}
I  & 0  & 0 & 0 \\
0   &  I   & 0 & 0 \\
\hline
0 & 0 & I & 0 \\
M_1 & M_2  & 0  & I   \\
\end{array}
\right).
$$
Since $M$ is invertible, at least one of the matrices $M_1$ and $M_2$ is not
the zero matrix, hence $(S,T)$ is not a scalar multiple of the identity. This
implies that $M$ was not simple.
\end{proof}

\medskip

\begin{example}\label{E:rank2}
  The triple $\bigl(\kO_{\PP^1}\oplus \kO_{\PP^1}(1), \CC^2_s, \widetilde{\mm}\bigr)$ with
  $$
  \widetilde{\mm}(0) =
  \left(
    \begin{array}{c|c}
      0 & 1 \\
      \hline
      \lambda & 0
    \end{array}
  \right) \quad \text{ and } \quad
  \widetilde{\mm}(\infty) =
  \left(
    \begin{array}{c|c}
      1      & 0 \\
      \hline
      0 & 1
    \end{array}
  \right),
  $$
  defines for any $\lambda \in \CC^*$ a simple  vector bundle of rank
  $2$ and degree $1$ on $E$. The corresponding matrix for this vector bundle
  is
  $M_{1,1}(\lambda) :=
  \left(\begin{smallmatrix}0&1\\
      \lambda&0\end{smallmatrix}\right)
  \in \MPnode^s(1,1)$.
  \qed
\end{example}

\begin{theorem}\label{T:nodecurve}
Let $E$ be a  nodal Weierstra\ss{} curve
and denote by $\Spl^{(n,d)}(E)$ the set of all isomorphism classes of simple
vector bundles of rank $n$ and degree $d$ on $E$.
If $\gcd(n,d)=1$ the map
$\det:\Spl^{(n,d)}(E) \rightarrow \Pic^{d}(E)\cong\CC^*$ is
bijective. If $\gcd(n,d)>1$ we have $\Spl^{(n,d)}(E) = \emptyset$.
\end{theorem}

\noindent
This  result can be proven by various methods, see for example
\cite[Theorem 3.6]{Burban1} for a description of simple vector bundles on $E$
in terms of \'etale coverings.
For the  reader's convenience we shall outline another proof\footnote{This
  proof is due to Lesya Bodnarchuk.}, which is parallel
to the case of a cuspidal cubic curve \cite{BodnarchukDrozd}.

\begin{proof}
First note that, without loss of generality, we may assume $0 \le d < n$.
If $\Spl^{(n,d)}(E) \ne \emptyset$ and $\kV$ is a non-zero element of
$\Spl^{(n,d)}(E)$ then, by Lemma \ref{L:normofsimple},
$$
\mathbb{F}(\kV) \cong
\Bigl(\kO_{\PP^1}^{n_1} \oplus \kO_{\PP^1}(1)^{n_2},
\CC^{n_1 + n_2}_s, \bigl(M, \mathsf{id}\bigr)\Bigr),
$$
where $n_2 = d$ and $n_1 = n-d$.
By Lemma \ref{L:VBequivMP} we have an equivalence
$\Spl^{(n,d)}(E) \cong \MPnode^{s}(n-d,d)$.

Assume first $n_{2}=0$ and $n_1 > 1$.
In this case, we have seen  in Remark \ref{R:degree-zero} that
$\MPnode^s(n,0)=\emptyset$. This implies $\Spl^{(n,0)}(E)=\emptyset$ for
$n > 1$. On the other hand,  $\Spl^{(1,0)}(E)=\Pic^{0}(E) \cong \breve{E}$.

For the rest of this proof we assume $n_{2}>0$.
By Lemma \ref{L:invblock}, the block $M_{12}$ of
$M\in\MPnode^{s}(n_{1},n_{2})$ has maximal rank.
If $n_1 = n_2$ this means that $M_{12}$ is invertible and $M$ can be reduced
to the form
$$
M =
\left(
\begin{array}{c|c}
0 & I \\
\hline
X & 0
\end{array}
\right)
$$
where $X$ splits into a direct sum of Jordan blocks with non-zero eigenvalues.
It is easy to see that $M$ is decomposable in $\MPnode(n_{1},n_{2})$ unless
$n_1 = n_2 = 1$. Hence, $\MPnode^{s}(m,m)=\emptyset$ if $m>1$.

On the other hand, if $n_{1}\ne n_{2}$, we can reduce $M$ to the form
(because both $n_{i}>0$)
$$
\left(
\begin{array}{c|cc}
0 & I & 0 \\
\hline
M'_{11} & 0 & M'_{12} \\
M'_{21} & 0 & M'_{22}
\end{array}
\right)
\text{ if }n_{2}> n_{1},\; \text{ or to }\quad
\left(
\begin{array}{cc|c}
M'_{11} & M'_{12} & 0 \\
0       & 0       & I \\
\hline
M'_{21} &  M'_{22}  & 0\\
\end{array}
\right)\text{ if } n_{1}> n_{2}.
$$
In both cases, the additional split of the blocks is made in such a way that
$M'_{11}$ and $M'_{22}$ are square matrices.
A straightforward calculation shows that
$$
M' =
\left(
\begin{array}{cc}
M'_{11} &  M'_{12} \\
M'_{21} &  M'_{22}
\end{array}
\right)
$$
is an object of $\MPnode^s(n_1, n_2 - n_1)$ or $\MPnode^s(n_1 - n_2, n_2)$
respectively.
This implies that, in case $n_2>n_1$, the fully faithful functor
$$\MPnode^s(n_1,n_2-n_1)\lar\MPnode^s(n_1,n_2)$$
which is defined on objects by sending
$$
N'=\left(
\begin{array}{cc}
N'_{11} &  N'_{12} \\
N'_{21} &  N'_{22}
\end{array}
\right)
\quad\text{ to }\quad
N=\left(
\begin{array}{c|cc}
0 & I & 0 \\
\hline
N'_{11} & 0 & N'_{12} \\
N'_{21} & 0 & N'_{22}
\end{array}
\right)\in\MPnode^s(n_1,n_2)
$$
and on morphisms by sending
$$
\left(
    \begin{pmatrix}
      A&0\\D'&E
    \end{pmatrix},
    \begin{pmatrix}
      A&0\\D''&E
    \end{pmatrix}
\right)
\in\Hom_{\MPnode^s(n_1,n_2-n_1)}(M',N')
$$
to
$$
\left(\left(
\begin{array}{c|cc}
A & 0 & 0 \\
\hline
N'_{12}D' & A  & 0 \\
N'_{22}D' & D' & E
\end{array}
\right),
\left(
\begin{array}{c|cc}
A & 0 & 0 \\
\hline
0   & A  & 0 \\
D'' & D' & E
\end{array}
\right)\right)
\in\Hom_{\MPnode^s(n_1,n_2)}(M,N)
$$
is in fact an equivalence of categories. Similarly, if $n_1>n_2$, we obtain an
equivalence
$$
\MPnode^s(n_1-n_2,n_2)\lar\MPnode^s(n_1,n_2).
$$
If we start with any pair of positive integers $n_{1}\ne n_{2}$ and continue
to reduce the size of the matrix in the way described above, we obtain an
equivalence of categories $\MPnode^{s}(m,m) \rightarrow \MPnode^{s}(n_{1},n_{2})$,
where $m=\gcd(n_{1},n_{2})$.
Our assumption $\Spl^{(n,d)}(E) \ne \emptyset$ implies now
$\gcd(n, d) = \gcd(n_1, n_2) = 1$ and our construction gives us an equivalence
$$\MPnode^{s}(1,1) \rightarrow \MPnode^{s}(n_{1},n_{2}).$$
Using Lemma \ref{L:invblock} we see that each object in $\MPnode^{s}(1,1)$ is
isomorphic to
$$
M_{1,1}(\lambda) =
\left(
\begin{array}{c|c}
0 & 1 \\
\hline
\lambda & 0
\end{array}
\right)
$$
for some $\lambda\in\CC^{\ast}$. We consider $M_{1,1}(\lambda)$ to be the
canonical form for objects in $\MPnode^{s}(1,1)$. Its image under the
equivalence $\MPnode^{s}(1,1) \rightarrow \MPnode^{s}(n_{1},n_{2})$
constructed above will be denoted by $M_{n_1, n_2}(\lambda)$. This is a
\emph{canonical form} for objects of the category $\MPnode^{s}(n_{1},n_{2})$.
An explicit description of $M_{n_1, n_2}(\lambda)$ is given in Algorithm
\ref{A:node} below.

Observe now that isomorphic objects of $\MPnode^{s}(n_{1},n_{2})$ have the same
determinant and that the functor
$\MPnode^{s}(1,1) \rightarrow \MPnode^{s}(n_{1},n_{2})$
respects the determinant up to a sign which depends on $(n_{1},n_{2})$ only.
As a consequence, we see that
$M_{n_1, n_2}(\lambda) \cong  M_{n_1, n_2}(\lambda')$ if and only if
$\det\bigl(M_{n_1, n_2}(\lambda)\bigr)=\det\bigl(M_{n_1,
  n_2}(\lambda')\bigr)$, which is equivalent to $\lambda = \lambda'$.
Because we have $$\mathbb{F}\bigl(\det(\kV)\bigr)
= \bigl(\kO_{\PP^1}(n_2), \CC_s, (\det(M),1)\bigr)$$
for any $\kV\in\Spl^{(n,d)}(E)$, we see now that
$\det:\Spl^{(n,d)}(E) \rightarrow \Pic^{d}(E)$ is bijective if $\gcd(n,d)=1$.
\end{proof}

\begin{algorithm}\label{A:node}
For any pair of positive coprime integers $(n_1, n_2)$, the simple
objects $M_{n_1, n_2}(\lambda) \in \MPnode^s(n_1, n_2)$ are described in the
following way.
\begin{enumerate}
\item First, we produce a sequence of pairs of coprime integers by replacing
  at each step a pair $(n_1, n_2)$ by $(n_1 - n_2, n_2)$ if $n_1 > n_2$ and
  by $(n_1, n_2 - n_1)$ if $n_2 > n_1$. We continue until we arrive at $(1, 1)$.
\item
  Starting with the matrix $M_{1,1}(\lambda) \in \MPnode^s$ from Example
  \ref{E:rank2} we recursively construct the matrix $M_{n_1, n_2}(\lambda)$ as
  follows. We follow the sequence constructed in part (1) in reverse order and
\begin{itemize}
\item if we go from $(m_1, m_2)$ to $(m_1 + m_2, m_2)$ we proceed as follows
$$
M_{m_1, m_2}(\lambda)=\left(
\begin{array}{c|c}
X & Y \\
\hline
Z & W
\end{array}
\right)
\Rightarrow
M_{m_1 + m_2, m_2}(\lambda)=\left(
\begin{array}{cc||c}
X & Y & 0 \\
\hline
0 & 0 & I_{m_2} \\
\hline
\hline
Z & W & 0
\end{array}
\right).
$$
\item and similarly, if we go from $(m_1, m_2)$ to $(m_1, m_1 +m_2)$ we set
$$
M_{m_1, m_2}(\lambda)=\left(
\begin{array}{c|c}
X & Y \\
\hline
Z & W
\end{array}
\right)
\Rightarrow
M_{m_1, m_1 +m_2}(\lambda)=
\left(
\begin{array}{c||c|c}
0 & I_{m_1} & 0 \\
\hline
\hline
X & 0 & Y \\
Z & 0 & W
\end{array}
\right).
$$
\end{itemize}
\end{enumerate}
\end{algorithm}

\begin{remark}\label{R:structofcanmatrix}
  From the construction it is clear that $M_{n_1, n_2}(\lambda)$ is an $n\times
  n$-matrix having exactly one non-zero  entry in each column and each row.
  Exactly one of these non-zero entries is equal to $\lambda\in\CC^{\ast}$ and
  this entry is found in the last row. All the other non-zero entries are equal
  to $1$.
\end{remark}

\begin{remark}\label{R:detvbnode}
  Because simple vector bundles on Weierstra{\ss} curves are stable
  (\cite[Cor.\ 4.5]{BK3}), we have $\Spl^{(n,d)}(E)=M_{E}^{(n,d)}$ and
  Theorem \ref{T:nodecurve} provides another proof of the part of
  Theorem \ref{T:ellfibrrepr} which says that two stable vector bundles
  $\kV_1$ and $\kV_2$ of the same rank on a nodal Weierstra\ss{} curve are
  isomorphic if and only if $\det(\kV_1) \cong \det(\kV_2)$.
\end{remark}

\begin{example}
Let us apply Algorithm \ref{A:node} to describe in terms of triples all simple
vector bundles on $E$ of rank $5$ and degree $12$. From earlier calculations
we see that the normalisation of such a bundle is
$\widetilde\kV = \kO_{\PP^1}(2)^3 \oplus \kO_{\PP^1}(3)^2$, in particular
$(n_1, n_2) = (3, 2)$.
The sequence of reductions for sizes of matrices from the category $\MPnode$
is:
$$
(3, 2) \lar (1, 2) \lar (1, 1).
$$
This induces a reverse sequence of functors
$$
\MPnode^s(1,1) \lar \MPnode^s(1,2) \lar \MPnode^s(3,2),
$$
giving the following sequence of canonical forms:
$$
\left(
\begin{array}{c|c}
0 & 1 \\
\hline
\lambda & 0
\end{array}
\right)
\mapsto
\left(
\begin{array}{c|cc}
0 & 1  & 0\\
\hline
0  & 0 & 1 \\
\lambda & 0 & 0
\end{array}
\right)
\mapsto
\left(
\begin{array}{ccc|cc}
0 & 1  & 0 & 0 & 0 \\
0 & 0  & 0 & 1 & 0 \\
0 & 0  & 0 & 0 & 1 \\
\hline
0  & 0 & 1  & 0 & 0 \\
\lambda & 0 & 0 & 0 & 0
\end{array}
\right).
$$
Therefore, the set of stable vector bundles of rank $5$ and degree $12$ is
described by the family of triples
$\bigl(\kO_{\PP^1}(2)^3 \oplus \kO_{\PP^1}(3)^2, \CC_s^5,
\widetilde{\mm}\bigr)$, where
$$
\mm(0) = \left(
\begin{array}{ccc|cc}
0 & 1  & 0 & 0 & 0 \\
0 & 0  & 0 & 1 & 0 \\
0 & 0  & 0 & 0 & 1 \\
\hline
0  & 0 & 1  & 0 & 0 \\
\lambda & 0 & 0 & 0 & 0
\end{array}
\right), \quad
\mm(\infty) = I_5.
$$
\end{example}

\begin{example}
The indecomposable semi-stable vector bundles $\kV$ of rank two and degree
zero, whose   Jordan-H\"older factors are locally free,  are  of the form
$\kL \otimes \kA_2$, where $\kL \in \Pic^0(E)$. Using Lemma
\ref{lem:semi-stable-on-nodal} and compatibility with tensor products, we see
that they are described by  the triples
$\bigl(\kO_{\PP^1}^2, \CC_s^2, \widetilde{\mm}\bigr)\in\Tri(E)$, where
$$
\mm(0) = \lambda J_{2}(1) =
\begin{pmatrix}
\lambda & \lambda \\
0       & \lambda
\end{pmatrix},  \,  \lambda \in \CC^* \quad \text{ and } \quad
\mm(\infty) = I_{2} =
\begin{pmatrix}
1 & 0 \\
0       & 1
\end{pmatrix}.
$$
\end{example}

\begin{lemma}\label{L:mult-diag-node}
Let $\lambda\in\CC^{\ast}$, $X = \mathrm{diag}(\alpha_1,
\alpha_2,\dots,\alpha_n)$ with $\alpha_{i}\in\CC^{\ast},\; i=1,\ldots,n$ and
$n=n_{1}+n_{2}$, then  $X M_{n_{1},n_{2}}(\lambda) \cong
M_{n_{1},n_{2}}(\lambda\cdot \alpha_1\cdot  \ldots \cdot\alpha_n)$
as objects of $\MPnode$.
\end{lemma}

\begin{proof}
  It is not hard to see that $X M_{n_{1},n_{2}}(\lambda)\in\MPnode$ is again
  simple, hence it is isomorphic to a canonical form
  $M_{n_{1},n_{2}}(\lambda')$. This means that there exist invertible matrices
  $S =
  \left(
    \begin{smallmatrix}A & 0 \\C' & B\end{smallmatrix}
  \right)$ and
  $T =
  \left(
    \begin{smallmatrix}A & 0 \\C'' & B\end{smallmatrix}
  \right)$
  such that $M_{n_{1},n_{2}}(\lambda') = S^{-1} X M_{n_{1},n_{2}}(\lambda)
  T$. Because
  $\det(S)=\det(A)\det(B)=\det(T)$, we obtain $\lambda'=\lambda\det(X)$.
\end{proof}

\begin{remark}\label{R:noderank2}
If $n_{1}, n_{2}$ are fixed and for each $\lambda\in\CC$ there is given an
object $M(\lambda)\in\MPnode^{s}(n_1, n_2)$, we obtain simple vector bundles
$\kV_{\lambda}$ of rank $n=n_{1}+n_{2}$ on $E$ such that
$\mathbb{F}(\kV_{\lambda}) \cong \bigl(\kO_{\PP^1}^{n_1} \oplus
\kO_{\PP^1}(1)^{n_2}, \CC^{n_{1}+n_{2}}_s,(M(\lambda),\id)\bigr)$. Similarly, for
each $\beta\in\CC^{\ast}$ there exists a unique line bundle
$\mathcal{L}_{\beta}$ on $E$ such that
$\mathbb{F}(\mathcal{L}_{\beta})
\cong \bigl(\kO_{\PP^1}, \CC,(\beta,1)\bigr)$.
We say that \textsl{the family $M(\lambda)$ is compatible with the action of
  the Jacobian}, if for all $\lambda,\beta\in\CC^{\ast}$ we have
$M(\beta^{n}\lambda) = \beta M(\lambda)$. This implies, but is stronger than
$\kV_{\beta^{n}\lambda} \cong \mathcal{L}_{\beta} \otimes \kV_{\lambda}$.

Lemma \ref{L:mult-diag-node} implies that the vector bundles $\kV_{\lambda}$
given by the family $M_{n_1, n_2}(\lambda)$ satisfy $\kV_{\beta^{n}\lambda}
\cong \mathcal{L}_{\beta} \otimes \kV_{\lambda}$, but $M_{n_1, n_2}(\lambda)$
is not compatible with the action of the Jacobian.
However, if we replace all non-zero entries of  $M_{n_1, n_2}(\lambda)$
by $\sqrt[n]{\lambda}$,  some fixed $n$-th root of $\lambda$, we obtain an
object $N_{n_1, n_2}(\lambda) \in \MPnode(n_1, n_2)$ which is isomorphic to
$M_{n_1, n_2}(\lambda)$ and which is compatible with the action of the Jacobian:
$N_{n_1, n_2}(\beta^{n}\lambda) = \beta N_{n_1, n_2}(\lambda)$.
Another choice of $\sqrt[n]{\lambda}$ gives an isomorphic vector bundle.

Because there is no global choice of an $n$-th root, we define
$\widetilde{N}_{n_1, n_2}(t) := N_{n_1, n_2}(t^{n})$ for all $t\in\CC^{\ast}$
in order to globalize this construction. Compatibility for this family has now
the form $\widetilde{N}_{n_1, n_2}(\beta t) = \beta \widetilde{N}_{n_1, n_2}(t)$.
As we shall see in Subsection \ref{SS:spectriv}, this will give us
a trivialization of a universal family of stable vector bundles compatible
with the action of Jacobian, necessary to construct a gauge transformation
of the geometric associative $r$-matrix depending on the difference of the
vector bundle spectral parameters.

If we apply this construction to the family of triples from Example
\ref{E:rank2}, which describe the simple vector bundles of
rank $2$ and degree $1$, we obtain a family of vector bundles
described by triples
$(\mathcal{O}_{\PP^{1}}\oplus \mathcal{O}_{\PP^{1}}(1),\CC^{2}_{s}, \widetilde{\mm})$
with
$$
\mm(0) = \widetilde{N}_{1,1}(\lambda) =
\left(
\begin{array}{c|c}
0 & \lambda \\
\hline
\lambda & 0 \\
\end{array}
\right),   \lambda \in \CC^* \quad \mbox{\textrm{and}} \quad
\mm(\infty) =
\left(
\begin{array}{c|c}
1     & 0 \\
\hline
0 & 1
\end{array}
\right).
$$
This family of matrices is compatible with the action of $\Pic^0(E)$. But for
$\lambda$ and $-\lambda$ we always obtain isomorphic vector bundles.
\end{remark}

\medskip
\medskip

\subsection{Vector bundles on a cuspidal Weierstra\ss{} curve}\label{SS:cusp}
We now recall  an explicit description of those objects in $\Tri(E)$ which
correspond to simple vector bundles on a cuspidal Weierstra\ss{} curve $E$,
following the approach of \cite{BodnarchukDrozd}.
The exposition is very similar to the nodal case and includes an algorithm
producing  a normal form for such triples for each given rank and
degree. Rank and degree of a simple vector bundle on a cuspidal
Weierstra{\ss} curve also turn out to be coprime.

Let $E$ be the cuspidal cubic curve, given by the equation
$zy^2 = x^3$.  Its normalisation   $\pi: \PP^1 \lar E$ is given by
$\pi(z_0:z_1) = (z_0^2 z_1: z_0^3: z_1^3)$.
With these  coordinates on $\PP^1$  the preimage of the singular point $s =
(0:0:1) \in E$ is $\pi^{-1}(s)  = (0:1) = \infty$. Then $Z$ is the reduced
point $s \in E$ with the structure sheaf  $\CC$. Moreover, $\widetilde{Z}$ is
non-reduced with support at $\infty = (0:1) \in \PP^1$ and structure sheaf
$\fR = \CC[\varepsilon]/\varepsilon^2$. The morphism  $\tilde\pi:
\widetilde{Z} \to Z$ corresponds to the canonical ring homomorphism  $\CC \to
\fR$.

\medskip
Recall that $w = z_0/z_1$ is a  coordinate in the neighbourhood $U_\infty$
of the point $(0:1)$. The morphism $\tilde\eta: \widetilde{Z} \to \PP^1$
corresponds to the map
$\ev_{U_\infty}: \kO_{\PP^1}(U_\infty) \to
\mathcal{O}_{\widetilde{Z}}(U_\infty) = \fR$, given
by $\ev_{U_\infty}(w)=\varepsilon$.
Next, following the recipe of Remark \ref{R:singvssmooth} we use the section
$p_\zeta(z_0, z_1) = z_{1}\in H^0\bigl(\kO_{\PP^1}(1)\bigr)$ to define the
collection  of isomorphisms
$
\zeta_l:  \tilde\eta^* \kO_{\PP^1}(l) \lar \kO_{\widetilde Z}.
$
They are given by  the formula
$\zeta_l(s) = \ev_{V}\left(\dfrac{s}{z_1^l}\right)$
for each open set $V \subset U_\infty$, all $l \in \mathbb{Z}$ and any
$s \in \Gamma\bigl(V, \kO_{\PP^1}(l)\bigr)$.
A morphism
$$q = q(z_0, z_1) = a_0 z_0^{m-n} + a_1 z_0^{m-n-1} z_1 + \dots  +
a_{m-n} z_1^{m-n} \in \Hom_{\PP^1}\bigl(\kO_{\PP^1}(n),\kO_{\PP^1}(m)\bigr)$$
is therefore evaluated according to the rule
$$
\xymatrix
{\tilde\eta^*\kO_{\PP^1}(n) \ar[d]_{\zeta_n} \ar[rr]^{\tilde\eta^*(q)} & &
\tilde\eta^*\kO_{\PP^1}(m) \ar[d]^{\zeta_m} \\
\fR \ar[rr]^{a_{m-n}+a_{m- n-1}\varepsilon} & &\fR.
}
$$
The following lemma shows how the explicit identification of $\widetilde{\mm}$
with a matrix is carried out for line bundles of degree one.

The chosen coordinates provide us with an isomorphism
$\CC \cong U_{\infty}=\PP^{1}\setminus\{\infty\}$ mapping $y \in \CC$ to
$(1:y) \in U_\infty$.
In the cuspidal case, the normalization restricts to an isomorphism
$\pi:\PP^{1}\setminus\{\infty\} \lar \breve{E}$.
Together, this gives us an identification $\breve{E} \cong \CC$,
under which $y\in\CC$ corresponds to $\tilde{y} := \pi^{-1}(y) =
(1:y)\in\PP^{1}$.

\begin{lemma}\label{L:idenoncusp}
With respect to the given choice of homogeneous coordinates on $\PP^1$ and the
set of trivializations $\{\zeta_l\}_{l \in \mathbb{Z}}$ described above, we have
for all $y \in \breve{E} \cong \CC$
$$
\mathbb{F}\bigl(\kO_E(y)\bigr)
\cong \bigl(\kO_{\PP^1}(1), \CC_s, 1 - y\varepsilon\bigr).
$$
\end{lemma}

\begin{proof}
As in the case of a nodal cubic curve,
because $\kO_E(y)$ is a line bundle of degree one, we know
$\kT_{\kO_E(y)} := \mathbb{F}\bigl(\kO_E(y)\bigr) = \bigl(\kO_{\PP^1}(1), \CC_s, (1 +
\lambda \varepsilon)\bigr)
$
for some $\lambda \in \CC$ and $\kT_{\kO_E} := \mathbb{F}(\kO_E)
=  \bigl(\kO_{\PP^1}, \CC_s, (1)\bigr)$.
By
Theorem \ref{thm:Drozd-Greuel} we have a commutative diagram
$$
 \xymatrix
 {
   \Hom_{\Tri(E)}\bigl(\kT_\kO, \kT_{\kO_E(y)}\bigr)
   \ar[rrr]^{\mathbb{G}} \ar[d]_{\For} & & &
   \Hom_{E}\bigl(\kO_E,  \kO_E(y)\bigr) \ar[d]^{\pi^*}\\
   \Hom_{\PP^1}\bigl(\kO_{\PP^1},  \kO_{\PP^1}(1)\bigr)  & & &
   \Hom_{\PP^1}\bigl(\kO_{\PP^1},  \kO_{\PP^1}(\tilde y)\bigr).
   \ar[lll]_{\conj\left(\gamma_{\kT_\kO}, \gamma_{\kT_{\kO(y)}}\right)}
 }
 $$
The section $z_1 - yz_0 \in \Hom_{\PP^1}\bigl(\kO_{\PP^1},  \kO_{\PP^1}(1)\bigr)$
generates the image of $\pi^{\ast}$, hence belongs to the image of $\For$ and
there exists  $c \in \CC$ such that the following diagram is commutative:
$$
\xymatrix{
\fR \ar[r]^c \ar[d]_1 & \fR \ar[d]^{1+\lambda \varepsilon} \\
\fR  \ar[r]^{1 - y\varepsilon} & \fR. \\
}
$$
This implies  $c=1$, $\lambda = -y$ and  $\mathbb{F}\bigl(\kO(y)\bigr) =
\bigl(\kO_{\PP^1}(1), \CC_s, (1- y  \varepsilon)\bigr)$.
\end{proof}

We aim now at giving a canonical form for those elements in $\Tri(E)$ which
correspond to simple vector bundles on $E$ having rank $n$ and degree $d$.
Just as in the nodal case, without loss of generality, we  may assume
$0 \le d < n$. Recall that $\VB^{(0,1)}_{n_{1},n_{2}}(E)$ is  the full
subcategory of $\VB^{(0,1)}(E)$,  whose objects are vector bundles $\kV$ with
fixed normalization
$
\pi^*\kV  \cong \widetilde\kV := \kO_{\PP^1}^{n_1} \oplus \kO_{\PP^1}(1)^{n_2}.
$
The category $\VB^{(0,1)}_{n_{1},n_{2}}(E)$ is equivalent to the full subcategory
of $\Tri(E)$ whose objects $(\widetilde{\kV},\mathcal{N},\widetilde{\mm})$
satisfy $\mathcal{N}\cong \CC_{s}^{n_{1}+n_{2}}$ and
$\widetilde\kV \cong \kO_{\PP^1}^{n_1} \oplus \kO_{\PP^1}(1)^{n_2}$.
Hence, these objects are described by an invertible
matrix $\widetilde{\mm}= \mm_{0}+\varepsilon \mm_{\varepsilon}\in\GL_{n_1 + n_2}(\fR)$.
Note the following easy lemma.

\begin{lemma}
Let $\widetilde{\mm} = \mm_0 + \varepsilon \mm_{\varepsilon}$ be an element
of  $\Mat_{n \times n}(\fR)$, where
$\mm_0, \mm_\varepsilon \in \Mat_{n \times n}(\CC)$. Then the matrix $\mm$ is
invertible if and only if $\mm_0$ is. Moreover, in the latter case we have:
$
\det(\widetilde{\mm})
= \det(\mm_0)\left(1 + \varepsilon \, \tr(\mm_0^{-1} \mm_\varepsilon)\right).
$
\end{lemma}

\noindent
Next, two such matrices
$\widetilde{\mm}, \widetilde{\mm}' \in \GL_{n_1 + n_2}(\fR)$ correspond to
isomorphic vector bundles if and only if there exist a matrix
$f\in\GL_{n_1+n_2}(\CC)$ and an automorphism $F$ of
$\kO_{\PP^1}^{n_1} \oplus \kO_{\PP^1}(1)^{n_2}$
such that $\widetilde{\mm}'
= \tilde\eta^*(F)^{-1} \circ \widetilde{\mm} \circ \tilde\pi^*(f)$.
For any $\widetilde{\mm}$, using $f= \mm_{0}^{-1}$ and $F = \id$, we find an
equivalent matrix $\widetilde{\mm}'$ with $\mm'_{0}=\id$.
In order to reduce $\widetilde{\mm}=\id+\varepsilon \mm_{\varepsilon}$
further, we split $m_{\varepsilon}$ into blocks
$(\mm_{\varepsilon})_{ij}\in\Mat_{n_{i}\times  n_{j}}(\CC)$
and let
$$
F =
\left(
\begin{array}{c|c}
I_{n_{1}} & 0 \\
\hline
F_{21} & I_{n_{2}}
\end{array}
\right)
$$
be the automorphism of $\widetilde\kV$ with
$F_{21} =z_{0} (\mm_{\varepsilon})_{21}\in\Mat_{n_2 \times n_1}(\CC[z_0, z_1]_1)$.
With this choice of $F$ and $f=\id$, a straightforward calculation,
which uses
$$
\tilde\eta^*(F) =
\left(
\begin{array}{c|c}
I_{n_{1}}  & 0 \\
\hline
0 & I_{n_{2}}
\end{array}
\right)
+
\varepsilon
\left(
\begin{array}{c|c}
0  & 0 \\
\hline
(\mm_{\varepsilon})_{21} & 0
\end{array}
\right),
$$
shows that we can reduce
$\widetilde{\mm}=\id+\varepsilon \mm_{\varepsilon}$ further to the form
$$
\widetilde{\mm} =
\left(
\begin{array}{c|c}
I_{n_{1}} & 0 \\
\hline
0 & I_{n_{2}}
\end{array}
\right) +\varepsilon
\left(
\begin{array}{c|c}
M_{11} & M_{12} \\
\hline
0 & M_{22}
\end{array}
\right).
$$
Therefore, triples
$\bigl(\kO_{\PP^1}^{n_1} \oplus \kO_{\PP^1}(1)^{n_2},
\CC_{s}^{n_{1}+n_{2}}, \widetilde{\mm}\bigr)$
with such $\widetilde{\mm}$ form a category which is equivalent to
$\VB^{(0,1)}_{n_{1},n_{2}}(E)$.
This motivates the following definition.
\begin{definition}\label{D:MPcusp}
  Let $E$ be a cuspidal cubic curve and $n_{1}>0$, $n_{2}\ge0$ integers.
  The category  $\MPcusp(n_{1},n_{2})$  is defined as follows.
  \begin{itemize}
  \item Its objects are ``matrices'' with three blocks
    $M_{ij}\in\Mat_{n_{i}\times n_{j}}(\CC)$
    $$
    M =  \left(
      \begin{array}{c|c}
        M_{11} & M_{12} \\
        \hline
        \times   & M_{22}
      \end{array}
    \right),
    $$
    where $\times$ is an ``empty'' or ``non-existing'' block.
  \item Morphisms are given by ``matrices''
    $$
    \Hom_{\MPcusp}(M, M') = \left\{ S \mid S M = M' S\right\}.
    $$
    with obvious composition and such that $S=\left(
      \begin{array}{c|c}
        S_{11} & \times  \\
        \hline
        S_{21} & S_{22}
      \end{array}
    \right)$
    has blocks of the same size as the blocks of $M$ and $M'$.
    The condition $S M = M' S$ means that
    $$
    \begin{array}{rcl}
      S_{11} M_{11} & =  & M'_{11} S_{11} + M'_{12} S_{21} \\
      S_{11} M_{12} & = & M'_{12} S_{22} \\
      S_{21} M_{12} + S_{22} M_{22} & = & M'_{22} S_{22},
    \end{array}
    $$
    in other words: we ignore the lower left block in $SM$ and $M'S$.
  \end{itemize}
As in the nodal case, we denote by $\MPcusp^{s}(n_{1},n_{2})$ the full subcategory
of \emph{simple} objects of $\MPcusp(n_{1},n_{2})$.
\end{definition}

\begin{remark}\label{R:degree-zero-cusp}
  If $n_{2}=0$ the block structure and the non-existing block disappear and we
  end up in a situation similar to the one described in Remark
  \ref{R:degree-zero}. Here we have $\MPcusp(n,0)=\Mat_{n\times n}(\CC)$ and
  $\Hom_{\MPcusp(n,0)}(M, M') = \left\{ S \mid S M = M' S\right\}$.
  The indecomposable objects in this category are precisely
  those which are isomorphic to a Jordan block $J_{n}(\lambda)$,
  $\lambda\in\CC$. As before, this implies $\MPcusp^{s}(n,0)=\emptyset$ if
  $n>1$ and $\MPcusp^{s}(1,0)=\CC$.
\end{remark}

\begin{lemma}\label{L:VBequivMPcusp}
  For any pair of non-negative integers $(n_{1},n_{2})$
  with $n_{1}>0$, the categories $\VB^{(0,1)}_{n_{1},n_{2}}(E)$ and
  $\MPcusp(n_{1},n_{2})$ are equivalent.
  Under this equivalence, simple vector bundles correspond to objects of
  $\MPcusp^s(n_{1},n_{2})$.
\end{lemma}

\begin{proof} Sending $M\in\MPcusp^{s}(n_{1},n_{2})$ to
  $I_{n_{1}+n_{2}}+\varepsilon M\in\GL_{n_{1}+n_{2}}(\kO_{\widetilde{Z}})$,
  with inherited block structure from $M$, gives an equivalence between
  $\MPcusp^{s}(n_{1},n_{2})$ and $\BM_{n_{1},n_{2}}(E)$.
  The proof of the lemma is now completely parallel to
  the case of a nodal cubic curve and is, therefore, left to the reader.
\end{proof}

Next, we wish to find a canonical form for objects
$M\in\MPcusp^{s}(n_{1},n_{2})$. Similar to the nodal case,
in each isomorphism class of $\MPcusp^{s}(n_{1},n_{2})$, we are going to
describe  a unique object with a particularly ``simple'' structure. Again, the
reduction procedure described below is based on an easy lemma.

\begin{lemma}\label{L:invblock-cusp}
  The block $M_{12}$ has  full rank, if $M \in \MPcusp^s(n_1, n_2)$ is simple.
\end{lemma}

\begin{proof}
  Just as in the nodal case (Lemma \ref{L:invblock}), if $M_{12}$ does not
  have full rank, the matrix $M$ can be reduced to the form
  $$
  M =
  \left(
    \begin{array}{c|c||c|c}
      M_1 & M_2 & 0 & 0 \\
      \hline
      0   &  0  & I & 0 \\
      \hline
      \hline
      \times  & \times & M_3 & M_4 \\
      \hline
      \times & \times & 0 & M_5 \\
    \end{array}
  \right)
  $$
  and we obtain a non-scalar endomorphism
  $$
  S =
  \left(
    \begin{array}{c|c||c|c}
      I & 0 & \times & \times \\
      \hline
      0   &  I  & \times & \times \\
      \hline
      \hline
      0 & 0  & I  & 0  \\
      \hline
      W & 0 & 0  & I \\
    \end{array}
  \right),
  $$
  where $W$ is an arbitrary matrix of  appropriate size.
\end{proof}

\begin{example}
  For any $\lambda \in \CC$, the triple
  $\bigl(\kO_{\PP^1}\oplus \kO_{\PP^1}(1), \CC^2_s, \widetilde{\mm}\bigr)$ with
  $$\widetilde{\mm} = \left(
    \begin{array}{c|c}
      1 & 0 \\
      \hline
      0 & 1
    \end{array}
  \right)
  +
  \varepsilon
  \left(
    \begin{array}{c|c}
      \lambda & 1 \\
      \hline
      0 & 0
    \end{array}
  \right)
  $$
  defines a simple  vector bundle of rank
  $2$ and degree $1$ on a cuspidal cubic curve $E$. It corresponds to
  $M_{1,1}(\lambda) :=
  \left(\begin{smallmatrix}\lambda&1\\
      \times&0\end{smallmatrix}\right)
  \in \MPcusp^s(1,1)$.
\end{example}

\begin{theorem}[see \cite{BodnarchukDrozd}]\label{T:detcusp}
  Let $E$ be a cuspidal cubic curve and denote by $\Spl^{(n,d)}(E)$ the set of
  all isomorphism classes of simple vector bundles of rank $n$ and degree $d$
  on $E$.
  If $\gcd(n,d)=1$ the map
  $\det:\Spl^{(n,d)}(E) \rightarrow \Pic^{d}(E)\cong\CC$ is
  bijective. If $\gcd(n,d)>1$ we have $\Spl^{(n,d)}(E) = \emptyset$.
\end{theorem}

\begin{proof}
  The proof very similar  to the proof of Theorem \ref{T:nodecurve}.
  A first difference is that if $n_1 = n_2$, we can transform the matrix
  $M\in\MPcusp^{s}(n,n)$ to the form
  $$
  M = \left(
    \begin{array}{c|c}
      M_{11}   & I \\
      \hline
      \times & 0
    \end{array}
  \right),
  $$
  because the block $M_{12}$ is square and invertible.
  We can further reduce the block $M_{11}$ to its  Jordan canonical form
  keeping the block $M_{12} = I $ unchanged. This implies that $M$ splits into
  a direct sum  of objects of the form
  $$
  \left(
    \begin{array}{c|c}
      J_m(\lambda)  & I_m \\
      \hline
      \times       & 0
    \end{array}
  \right)
  $$
  which are simple in $\MPcusp(m,m)$ if and only if $m = 1$.

  The other difference is that a simple object $M \in \MPcusp(n_1, n_2)$
  can be reduced to
  $$M =
  \left(
    \begin{array}{c|cc}
      0 & I  & 0 \\
      \hline
      \times & M'_{11} & M'_{12} \\
      \times & 0 & M'_{22}
    \end{array}
  \right)
  \text{ if }n_{2}> n_{1},\; \text{ or to }\quad
  M =
  \left(
    \begin{array}{cc|c}
      M'_{11} & M'_{12} & 0 \\
      0 & M'_{22} & I \\
      \hline
      \times & \times & 0
    \end{array}
  \right)\text{ if } n_{1}> n_{2}.
  $$
  A straightforward calculation shows that the matrix
  $$
  M' = \left(
    \begin{array}{c|c}
      M'_{11} & M'_{12} \\
      \hline
      \times        & M'_{22}
    \end{array}
  \right)
  $$
  is an object of $\MPcusp^s(n_1, n_2 - n_1)$ or $\MPcusp^s(n_1 - n_2, n_2)$
  respectively and that
  $$\det(\id+\varepsilon M) = \det(\id+\varepsilon M').$$
  If $\gcd(n_{1},n_{2})=1$, we end up with an equivalence
  $\MPcusp^s(1,1)\rightarrow\MPcusp^s(n_1,n_2)$ just as in the nodal case.
  Using Lemma \ref{L:invblock-cusp} we see that each object in
  $\MPcusp^{s}(1,1)$ is isomorphic to
  $$
  M_{1,1}(\lambda) =
  \left(
    \begin{array}{c|c}
      \lambda & 1 \\
      \hline
      \times  & 0
    \end{array}
  \right),\quad \lambda\in\CC.
  $$
  By definition, $M_{n_1, n_2}(\lambda)\in\MPcusp^{s}(n_{1},n_{2})$ denotes
  the image of $M_{1,1}(\lambda)$ under the equivalence described above.
  Again, $M_{n_1, n_2}(\lambda) \cong M_{n_1, n_2}(\lambda')$ in
  $\MPcusp^{s}(n_{1},n_{2})$ if and only if $\lambda = \lambda'$. From the
  identity $\det\bigl(\id+\varepsilon M_{1,1}(\lambda)\bigr)=1+\varepsilon\lambda$, the
  bijectivity of the determinant map follows.
\end{proof}

\begin{remark}
  Reversing the reduction step in the proof of Theorem \ref{T:detcusp} gives
  us an algorithm similar to Algorithm \ref{A:node} which produces the matrix
  $M_{n_1, n_2}(\lambda)$ starting with $M_{1,1}(\lambda)$. The only non-zero
  diagonal element of $M_{n_1, n_2}(\lambda)$ will be the moduli parameter
  $\lambda \in \CC$, i.e.\ $\lambda=\tr\bigl(M_{n_1, n_2}(\lambda)\bigr)$.
\end{remark}

\begin{lemma}\label{L:mult-diag-cusp}
  Let $\lambda\in\CC$, $A=\mathrm{diag}(\alpha_1,\alpha_2,\dots,\alpha_n)$
  with $\alpha_{i}\in\CC$ and $n=n_{1}+n_{2}$, then
  $
  A+M_{n_{1},n_{2}}(\lambda)\cong
  M_{n_{1},n_{2}}(\lambda + \alpha_1 + \alpha_2 + \dots + \alpha_n)
  $
  as objects of $\MPcusp$.
\end{lemma}

\begin{proof}
  We proceed by induction on $n$, the size of the matrix
  $M_{n_{1},n_{2}}(\lambda)$. The case $n_1 = n_2 = 1$ is an easy calculation.
  Assume the statement is true for all pairs $(n_1, n_2)$ of positive integers
  such that $n_1 + n_2 < n$. We shall deal with the case $n_{1}>n_{2}$, the
  opposite case is similar and left to the reader.

  From the proof of Theorem \ref{T:detcusp} we know that
  $M_{n_{1},n_{2}}(\lambda)$ has the  structure
  $$  M_{n_{1},n_{2}}(\lambda) =
  \left(
    \begin{array}{cc|c}
      M'_{11} & M'_{12} & 0 \\
      0 & M'_{22} & I_{n_{2}} \\
      \hline
      \times & \times & 0
    \end{array}
  \right)\quad\text{ so that }\quad
  M_{n_{1}-n_{2},n_{2}}(\lambda)=
    \left(
    \begin{array}{c|c}
      M'_{11} & M'_{12} \\
      \hline
      \times  & M'_{22} \\
    \end{array}
  \right),
  $$
  with $M'_{11}, M'_{22}$ being square matrices of sizes $n_{1}-n_{2}$ and
  $n_{2}$ respectively.
  If we write
  $$A=\mathrm{diag}(\alpha_1,\alpha_2,\dots,\alpha_{n_{1}-n_{2}} |
  \alpha_{n_{1}-n_{2}+1},\ldots,\alpha_{n_{1}} |
  \alpha_{n_{1}+1},\ldots,\alpha_n)  = \mathrm{diag}(A_{1}|A_{2}|A_{3}),$$ with
  $A_{i}$ being diagonal blocks, we obtain
  $$A+M_{n_{1},n_{2}}(\lambda) =
  \left(
    \begin{array}{cc|c}
      M'_{11}+A_{1} & M'_{12} & 0 \\
      0 & M'_{22}+A_{2} & I_{n_{2}} \\
      \hline
      \times & \times & A_{3}
    \end{array}
  \right).
  $$
  A straightforward calculation shows that the matrix
  $$
  S=
  \left(
    \begin{array}{cc|c}
      I_{n_{1}-n_{2}} & 0      & \times \\
      0          & I_{n_{2}} & \times \\
      \hline
      0          & -A_{3} & I_{n_{2}}
    \end{array}
  \right)
  $$
  defines an isomorphism in the category $\MPcusp(n_{1},n_{2})$ between
  $A+M_{n_{1},n_{2}}(\lambda)$ and
  $$
  \left(
    \begin{array}{cc|c}
      M'_{11}+A_{1} & M'_{12} & 0 \\
      0 & M'_{22}+A_{2}+A_{3} & I_{n_{2}} \\
      \hline
      \times & \times & 0
    \end{array}
  \right).
  $$
  The inductive hypothesis implies
  that $M_{n_{1}-n_{2},n_{2}}(\lambda) + \mathrm{diag}(A_{1}|A_{2}+A_{3})$ is
  a simple object of $\MPcusp(n_{1}-n_{2},n_{2})$ which is isomorphic to
  $M_{n_{1}-n_{2},n_{2}}(\lambda+\alpha_{1}+\ldots+\alpha_{n})$. This implies
  that $A+M_{n_{1},n_{2}}(\lambda) \cong
  M_{n_{1},n_{2}}(\lambda+\alpha_{1}+\ldots+\alpha_{n})$ in
  $\MPcusp(n_{1},n_{2})$.
  It was not possible to give a direct proof like for Lemma
  \ref{L:mult-diag-node} in the nodal case, because it is not obvious at the
  beginning that $A + M_{n_{1},n_{2}}(\lambda)\in\MPcusp$ is again simple.
\end{proof}

\begin{remark}\label{R:eqvcusp}
  This lemma implies that the object $N_{n_1, n_2}(\lambda)\in\MPcusp(n_1,
  n_2)$, obtained from
  $M_{n_1, n_2}(\lambda)$ by replacing each diagonal entry by
  $\frac{\lambda}{n}$, is isomorphic to $M_{n_1, n_2}(\lambda)$  in
  $\MPcusp(n_1, n_2)$. Moreover, this matrix is \emph{compatible with the
    action of the Jacobian} in the sense that for all $\lambda,\beta\in\CC$ we
  have $\beta I_{n}+ N_{n_1, n_2}(\lambda) = N_{n_1, n_2}(n\beta+\lambda)$.
  This is equivalent to
  $(1+\varepsilon\beta)(I_{n}+\varepsilon N_{n_1, n_2}(\lambda))
  = I_{n}+\varepsilon N_{n_1, n_2}(n\beta+\lambda)$.
  The precise meaning of this condition will be clarified in Subsection
  \ref{SS:spectriv}.
\end{remark}

\begin{remark}\label{R:detvbcusp}
  Because simple vector bundles on Weierstra{\ss} curves are stable
  (\cite[Cor.\ 4.5]{BK3}), we have $\Spl^{(n,d)}(E)=M_{E}^{(n,d)}$ and
  Theorem \ref{T:detcusp} provides another proof of the part of
  Theorem \ref{T:ellfibrrepr} which says that two stable vector bundles
  $\kV_1$ and $\kV_2$ of the same rank on a cuspidal Weierstra\ss{} curve are
  isomorphic if and only if $\det(\kV_1) \cong \det(\kV_2)$.

  Moreover, because the group $\Pic^0(E)\cong \CC$ is torsion free and
  divisible, it follows from
  $\det(\kV \otimes \kL) \cong \det(\kV) \otimes \kL^{\otimes n}$
  that the action of the Jacobian $\Pic^0(E)$ on the set $M_E^{(n, d)}$ of
  stable vector bundles of rank $n$ and degree $d$ is simply transitive.
\end{remark}

\begin{example}
The family of vector bundles on a cuspidal cubic curve, described by the
triples
$\bigl(\mathcal{O}_{\PP^{1}}\oplus \mathcal{O}_{\PP^{1}}(1),
\CC^{2}_{s}, \widetilde{\mm}\bigr)$
with
$$\renewcommand{\arraystretch}{1.2}
\widetilde{\mm} = I_{2}+\varepsilon N_{1,1}(\lambda)=
\left(
\begin{array}{c|c}
1 & 0 \\
\hline
0 & 1
\end{array}
\right)
+
\varepsilon
\left(
\begin{array}{c|c}
\frac{\lambda}{2} & 1 \\
\hline
0 & \frac{\lambda}{2}
\end{array}
\right),
$$
defines a universal family of stable vector bundles of rank $2$ and degree
$1$. The family of matrices
$$\renewcommand{\arraystretch}{1.2}
N_{n_{1},n_{2}}(\lambda)=\left(
\begin{array}{c|c}
\frac{\lambda}{2} & 1 \\
\hline
0 & \frac{\lambda}{2}
\end{array}
\right)
$$
is compatible with the action of $\Pic^0(E)$.
\end{example}

\medskip
\subsection{Universal families and their trivializations}\label{SS:spectriv}
The goal of this subsection is an explicit description of a
universal family on the moduli space of stable vector bundles on a singular
Weierstra\ss{} cubic curve $E$.

\medskip
\noindent
For a \emph{reduced} complex space $B$ consider the Cartesian diagram
$$
\begin{CD}
  \widetilde{Z} \times B
  @>{\tilde{\eta}_{B} = \tilde{\eta} \times \id}>>
  \PP^1 \times B  \\
  @V{\tilde{\pi}_{B} = \tilde{\pi} \times \id}VV
  @VV{\pi_{B} = \pi \times \id}V\\
  Z \times B
  @>{\eta_{B} = \eta \times \id}>>
  E\times B
\end{CD}
$$
and abbreviate $\nu_{B}=\pi_{B}\circ \tilde{\eta}_{B} = \eta_{B}\circ
\tilde{\pi}_{B}$. If $B$ is a point we omit the subscript $B$ in the notation
introduced.

Let us fix homogeneous coordinates $(z_0: z_1)$ on $\PP^1$ and   denote
$\kO_{\PP^1 \times B}(l) = \pr_1^*\kO_{\PP^1}(l)$, where
$\pr_1: \PP^1 \times B \to \PP^1$  is the projection map.
Let $p = p(z_0, z_1) \in H^0\bigl(\kO_{\PP^1}(1)\bigr)$ be a section, which
is non-vanishing on $\widetilde Z$. The recipe of Remark
\ref{R:singvssmooth} gives a family of isomorphisms
$\zeta_l:  \tilde\eta^*\bigl(\kO_{\PP^1}(l)\bigr) \lar \kO_{\widetilde Z}$.
Pulling back to $\PP^1 \times B$, we obtain isomorphisms
$\zeta_l:  \tilde\eta_B^*\bigl(\kO_{\PP^1 \times B}(l)\bigr)
\lar \kO_{\widetilde{Z}\times B}$ denoted for the sake of simplicity by the
same letters.

Let $\widetilde\kA = \bigoplus_{l \in \mathbb{Z}}\kO_{\PP^1 \times B}(l)^{n_l}$
be a vector bundle of rank $n$ on $\PP^1 \times B$.
Then there are induced  isomorphisms
$\zeta^{\widetilde \kA}: \tilde\eta_B^* \widetilde\kA \lar
\kO_{\widetilde{Z}\times B}^n$ and
$\zeta^{\widetilde \kA}: \nu_{B\ast}\tilde\eta_B^* \widetilde\kA \lar
\nu_{B\ast}\kO_{\widetilde{Z}\times B}^n$, denoted again by the same letters.
For any point $b \in B$ we have:
$\widetilde{\kA}_b := \widetilde{\kA}|_{\PP^1 \times \{b\}} =
\bigoplus_{l \in \mathbb{Z}}\kO_{\PP^1}(l)^{n_l}$. In a similar way, we obtain
isomorphisms
$\zeta^{\widetilde{\kA}_b}: \tilde\eta^* \widetilde{\kA}_b \lar \kO_{\widetilde Z}^n$ and
$\zeta^{\widetilde{\kA}_b}: \nu_{\ast}\tilde\eta^* \widetilde{\kA}_b \lar
\nu_{\ast}\kO_{\widetilde Z}^n$.
Note that the isomorphisms $\zeta^{\widetilde\kA}$ and
$\zeta^{\widetilde{\kA}_b}$ are related by the canonical base change diagrams.

Now we proceed with our construction of  vector bundles on $E\times B$.
We start with an invertible matrix $M \in \GL_n(\kO_{\widetilde{Z}\times B})$.
If convenient, we may start instead with a holomorphic function
$M:B\lar\Mat_{n\times n}(\CC)$, denoted by the same letter. The corresponding
element in $\GL_n(\kO_{\widetilde{Z}\times B})$ will then be $(M,\id)$ in the
nodal case and $1_{n}+\varepsilon M$ in the cuspidal case (see Lemmas
\ref{L:VBequivMP} and \ref{L:VBequivMPcusp}).
For any point $b \in B$ we denote by $M(b)$ the corresponding matrix
in  $\GL_n(\kO_{\widetilde{Z}})$ or $\Mat_{n\times n}(\CC)$ respectively.
Following the notation of Subsection \ref{SS:catoftriples}, we denote
\[\mm: \eta_{B \ast} \kO_{Z \times B}^n \stackrel{\can}\lar
\nu_{B\ast} \kO_{\widetilde{Z} \times B}^n \xrightarrow{\nu_{B\ast}(M)}
\nu_{B\ast} \kO_{\widetilde{Z} \times B}^n\] and let $\widetilde\mm$ be the
unique map which makes the following diagram of isomorphisms commutative:
\[
\xymatrix{
\kO_{\widetilde{Z} \times B}^n \ar[r]^{\widetilde\mm} \ar@/_{2.5ex}/[rr]_{M}&
\widetilde{\eta}_{B}^{\ast}\widetilde\kA \ar[r]^{\zeta^{\widetilde\kA}}&
\kO_{\widetilde{Z} \times B}^n
}
\]
In a similar way, let
$\mm_b$ and $\widetilde{\mm}_b$ be the morphisms determined by the matrix
$M(b)$. The following theorem is a mild generalization of Theorem
\ref{thm:Drozd-Greuel}.

\medskip

\begin{theorem}\label{thm:DGrel} Let
$\bar\kA = \oplus_{l \in \mathbb{Z}}\kO_{\PP^1}(l)^{n_l}$ be a vector bundle
of rank $n$  on $\PP^1$, $\zeta: \kO(1)_{\PP^1}\big|_{\widetilde Z}
\lar \kO_{\widetilde Z}$ be the isomorphism induced by a section
$p = p_\zeta(z_0, z_1) \in H^0\bigl(\kO_{\PP^1}(1)\bigr)$, and
$\widetilde\kA = \oplus_{l \in \mathbb{Z}}\kO_{\PP^1 \times B}(l)^{n_l}$
be the  pull-back of $\bar\kA$ to $\PP^1 \times B$.

\hspace{1mm}

\noindent
$\bullet$
Consider the  coherent sheaf   $\kA$ on $E\times B$  given by the exact
sequence
\begin{equation}
  \label{eq:flat-family-node}
  0\rightarrow \kA \xrightarrow
  {
  \left(
  \begin{smallmatrix}
  \ii \\
  \pp
  \end{smallmatrix}
  \right)
  }
  \pi_{B\ast}\widetilde{\kA} \oplus \eta_{B\ast}\mathcal{O}_{Z\times B}^{n}
  \xrightarrow{\left(\begin{smallmatrix}{\zeta}^{\widetilde\kA}
        & \mm\end{smallmatrix}\right)}
  \nu_{B\ast}\mathcal{O}_{\widetilde{Z}\times B}^{n}
  \rightarrow 0.
\end{equation}
Then $\kA \in \Coh(E\times B)$ is locally free and for each $b\in B$ we have:
$\mathcal{A}_{b}
=  \kA|_{E \times \{b\}}
\cong \mathbb{G}\bigl(\widetilde\kA\big|_{E \times \{b\}}, \, \kO_Z^n, \,
\widetilde{\mm}(b)\bigr)$, where $\mathbb{G}$ is the functor described in
Theorem \ref{thm:Drozd-Greuel}. In what follows, we shall use the notation
$\kA = \mathbb{G}\bigl(\widetilde\kA, \kO_{Z \times B}^n, \widetilde \mm\bigr)$.

\hspace{1mm}

\noindent
$\bullet$
Let $\widetilde\kB = \oplus_{l \in \mathbb{Z}}\kO_{\PP^1 \times B}(l)^{m_l}
$ be a vector bundle of rank $m$ on $\PP^1 \times B$,
$N \in \GL_m(\kO_{\widetilde{Z} \times B})$,
$\nn:  \eta_{B \ast} \kO_{Z \times B}^m \stackrel{\can}\lar
\nu_{B\ast} \kO_{\widetilde{Z} \times B}^m \xrightarrow{\nu_{B\ast}(N)}
\nu_{B\ast} \kO_{\widetilde{Z} \times B}^m$
and $\kB = \mathbb{G}\bigl(\widetilde\kA, \kO_{Z \times B}^m, \widetilde \nn\bigr)$.
Then
$$
\Hom_{E \times B}(\kA, \, \, \kB) \subseteq
\Hom_{\PP^1 \times B}(\widetilde\kA, \widetilde\kB) \times
\Mat_{m \times n}(\kO_{Z \times B})
$$
consists of those pairs $(F,f)$ for which the following diagram is
commutative:
\[
\xymatrix
{
\kO_{\widetilde{Z}\times B}^n
\ar[d]_{\tilde\pi_B^*(f)}
\ar[rr]_M
\ar@/^3ex/[rrrr]^{\widetilde\mm}
& & \kO_{\widetilde{Z}\times B}^n
& & \tilde\eta_B^* \widetilde\kA
\ar[ll]^{\zeta^{\widetilde\kA}}
\ar[d]^{\tilde\eta_B^*(F)}
\\
\kO_{\widetilde{Z}\times B}^m
\ar[rr]^N
\ar@/_3ex/[rrrr]_{\widetilde\nn}
& & \kO_{\widetilde{Z}\times B}^m
& & \tilde\eta_B^* \widetilde\kB
\ar[ll]_{\zeta^{\widetilde\kB}}.
}
\]

\hspace{1mm}

\noindent
$\bullet$ Let $f: B' \lar B$ be a holomorphic map between reduced analytic
spaces.
Let $M' \in \GL_n(\kO_{\widetilde{Z}\times B'})$ be the image of the matrix
$M$ under the morphism induced by
$H^0(\kO_{\widetilde{Z}\times B}) \rightarrow H^0(\kO_{\widetilde{Z}\times B'})$,
which is given by pull-back under $f$. Let $\kA' = (\id \times f)^* \kA$,
$\widetilde\kA' = (\id \times f)^* \widetilde\kA$ etc. Then we have:
$\kA' \cong
\mathbb{G}\bigl(\widetilde\kA', \kO_{Z \times B'}^n, \widetilde{\mm}'\bigr)$.
In other words, this construction is compatible with base change.

 \hspace{1mm}

 \noindent
 $\bullet$ Let $N \in H^0(\kO^{\ast}_{\widetilde{Z} \times B})$,
 $\widetilde\kL = \kO_{\PP^1 \times B}(c)$ for some $c \in \ZZ$
 and $\kL = \mathbb{G}\bigl(\widetilde\kL, \kO_{Z \times B}, \widetilde{\nn}\bigr)$
 be the corresponding line bundle, i.e.\ we have an exact sequence
 $$
 0\rightarrow \kL \xrightarrow
 {\left(
     \begin{smallmatrix}
       \jj \\
       \qq
     \end{smallmatrix}
   \right)
 }
 \pi_{B\ast}\widetilde{\kL} \oplus \eta_{B\ast}\mathcal{O}_{Z\times B}
 \xrightarrow{\left(\begin{smallmatrix}{\zeta}^{\widetilde\kL} &
       \nn\end{smallmatrix}\right)}
 \nu_{B\ast}\mathcal{O}_{\widetilde{Z}\times B}
 \rightarrow 0.
 $$
 Then the following sequence is exact:
 $$
 0\rightarrow \kA \otimes \kL
 \xrightarrow{
   \left(
     \begin{smallmatrix}
       \ii \boxtimes \jj \\
       \pp \boxtimes \qq
     \end{smallmatrix}
   \right)
 }
 \pi_{B\ast}\bigl(\widetilde{\kA} \otimes \widetilde{\kL}\bigr)
 \oplus \eta_{B\ast}\mathcal{O}_{Z\times B}^{n}
 \xrightarrow{
   \left(
     \begin{smallmatrix}
       {\zeta}^{\widetilde\kA \otimes \widetilde\kL} & \mm \boxtimes \nn
     \end{smallmatrix}
   \right)
 }
 \nu_{B\ast}\mathcal{O}_{\widetilde{Z}\times B}^{n}
 \rightarrow 0,
 $$
 where the morphism $\mm \boxtimes \nn$ is induced by the matrix $M \cdot N$,
 $\ii \boxtimes \jj$ is the morphism
 $\kA \otimes \kL  \xrightarrow{\ii \otimes \jj} \pi_{B\ast} \widetilde\kA
 \otimes \pi_{B\ast} \widetilde\kL
 \stackrel{\can}\lar  \pi_{B\ast}\bigl(\widetilde\kA  \otimes \widetilde\kL\bigr)$
 and $\pp \boxtimes \qq$ is defined is a similar way.
 This means that the tensor product of a vector bundle with a line bundle is
 given by the product of the corresponding matrices. More generally, the
 tensor product of two vector bundles corresponds to the tensor product of the
 defining matrices.

\hspace{1mm}

 \noindent
 $\bullet$ Finally,  we have the following short exact sequence:
 $$
  0\rightarrow \det\kA \lar
  \pi_{B\ast}\det\widetilde{\kA} \oplus \eta_{B\ast}\mathcal{O}_{Z\times B}
  \xrightarrow{\left(\begin{smallmatrix} \zeta^{\det\widetilde\kA} &
  \det(\mm) \end{smallmatrix}\right)}
  \nu_{B\ast}\mathcal{O}_{\widetilde{Z}\times B}
  \rightarrow 0,
$$
where the morphism $\det(\mm)$ corresponds to the determinant
$\det(M) \in H^0(\kO_{\widetilde{Z} \times B})$.
\end{theorem}

\begin{proof}
  To prove the first part of the statement, note that the sheaf   $\nu_{B\ast}\mathcal{O}_{\widetilde{Z}\times B}^{n}$  is
  $B$-flat. Hence,   for any point $b \in B$ the restriction of the sequence
  \eqref{eq:flat-family-node} to $E \times \{b\}$
  $$
  0\rightarrow \kA_b \rightarrow
  \pi_{\ast}\widetilde{\kA}_b  \oplus \eta_{\ast}\mathcal{O}_{Z}^{n}
  \xrightarrow{\left({\zeta}^{\widetilde\kA_b} \, \,  \mm_b\right)}
  \nu_{\ast}\mathcal{O}_{\widetilde{Z}}^{n}
  \rightarrow 0,
  $$
   is exact again.
  By  Theorem \ref{thm:Drozd-Greuel} the coherent sheaf  $\kA_b$ is
  locally free. Since $B$ is reduced, $\kA$ is locally free, too.

  The description of morphisms between $\kA$ and $\kB$ in terms of morphisms
  between $\widetilde\kA$ and $\widetilde\kB$ and matrices $M$ and $N$ follows
  from the universal property of the kernel.

  In order to show the base change property for $\GG$ we use again that
   $\nu_{B\ast}\kO_{\widetilde{Z} \times B}^n$ is flat
  over $B$. Let $\tilde{f} = \id \times f: E \times B' \lar E \times B$ then the
  functor $\tilde{f}^*$ induces a short exact sequence
  $$
  0\rightarrow \tilde{f}^*\kA \xrightarrow
  {\tilde{f}^*
  \left(
  \begin{smallmatrix}
  \ii \\
  \pp
  \end{smallmatrix}
  \right)
  }
  \tilde{f}^*\bigl(\pi_{B\ast}\widetilde{\kA} \oplus
  \eta_{B\ast}\mathcal{O}_{Z\times B}^{n}\bigr)
  \xrightarrow{\tilde{f}^*\left(\begin{smallmatrix}{\zeta}^{\widetilde\kA} &
  \mm
  \end{smallmatrix}\right)}
  \tilde{f}^*\bigl(\nu_{B\ast}\mathcal{O}_{\widetilde{Z}\times B}^{n}\bigr)
  \rightarrow 0.
  $$
  Consider the following commutative diagram:
  $$
  \xymatrix
  {
    \PP^1 \times B' \ar[rr]^{\id \times f} \ar[d]_{\pi_{B'}} & &
    \PP^1 \times B \ar[rr]^{\pr_1} \ar[d]^{\pi_{B}} & & \PP^1 \ar[d]^{\pi}\\
    E \times B' \ar[rr]^{\tilde f} \ar[d]_{\pr_2} & &
    E \times B \ar[rr]^{\pr_1} \ar[d]^{\pr_2} & & E \\
    B' \ar[rr]^{ f} & &  B & &
  }
  $$
  It implies that the base-change morphism
  $\tilde{f}^* \pi_{B\ast} \widetilde\kA \lar
  \pi_{B'\ast} (\id \times f)^*  \widetilde\kA$
  is an isomorphism. Indeed, it can be identified with the composition
  of isomorphisms
  $$
  \tilde{f}^*  \pi_{B\ast} \widetilde\kA
  \cong \tilde{f}^* \pi_{B\ast} \pr_1^* \bar\kA
  \cong \tilde{f}^*   \pr_1^*  \pi_{\ast} \bar\kA
  \cong \pi_{B'\ast} (\id \times f)^* \pr_1^* \bar\kA
  \cong \pi_{B'\ast} (\id \times f)^*  \widetilde\kA,
  $$
  given by the flat base change.  Denote $\widetilde\kA'  =
  (\id \times f)^* \widetilde\kA$. Then it is not difficult to show that the
  following diagrams are commutative:
  $$
  \xymatrix
  {
    \tilde{f}^* \pi_{B\ast} \widetilde\kA
    \ar[rr]^{\tilde{f}^*(\zeta^{\widetilde\kA})} \ar[d]_\can & &
    \tilde{f}^* \nu_{B \ast} \kO_{\widetilde{Z} \times B}^n \ar[d]^\can\\
    \pi_{B'\ast} \widetilde\kA' \ar[rr]^{\zeta^{\widetilde\kA'}} & &
    \nu_{B' \ast} \kO_{\widetilde{Z} \times B'}^n
  }
  \quad
  \xymatrix
  {
    \tilde{f}^* \eta_{B\ast}  \kO_{Z \times B}^n
    \ar[rr]^{\tilde{f}^*(\mm)} \ar[d]_\can & &
    \tilde{f}^* \nu_{B \ast} \kO_{\widetilde{Z} \times B}^n \ar[d]^\can\\
    \eta_{B'\ast} \kO_{Z \times B'}^n \ar[rr]^{\mm'} &
    &
    \nu_{B' \ast} \kO_{\widetilde{Z} \times B'}^n\;,
  }
  $$
  in which the vertical morphisms are induced by the base change. This implies
  that we have the following commutative diagram
  $$
  \xymatrix
  {
  0 \ar[r] &  \tilde
  {f}^*\kA \ar[d]  \ar[rr]^-{\tilde{f}^*
  \left(
  \begin{smallmatrix}
  \ii \\
  \pp
  \end{smallmatrix}
  \right)
  }
   &  &
   \tilde{f}^*\bigl(\pi_{B\ast}\widetilde{\kA} \oplus
   \eta_{B\ast}\mathcal{O}_{Z\times B}^{n}\bigr)
  \ar[d]_\can
  \ar[rr]^-{\tilde{f}^*\left(
      \begin{smallmatrix}
        {\zeta}^{\widetilde\kA} & \mm
      \end{smallmatrix}\right)} & &
  \tilde{f}^*\bigl(\nu_{B\ast}\mathcal{O}_{\widetilde{Z}\times B}^{n}\bigr)
  \ar[d]^\can \ar[r]  &  0 \\
  0 \ar[r] &  \tilde{f}^*\kA \ar[rr]
   & &
  \pi_{B'\ast}\widetilde{\kA}' \oplus \eta_{B'\ast}\mathcal{O}_{Z\times B'}^{n}
  \ar[rr]^-{\left(
      \begin{smallmatrix}{\zeta}^{\widetilde\kA'} & \mm'
      \end{smallmatrix}\right)} & &
  \nu_{B'\ast}\mathcal{O}_{\widetilde{Z}\times B'}^{n}
  \ar[r]  &  0.
  }
  $$
  Finally, the  compatibility of $\mathbb{G}$  with tensor products
  can be proven along similar lines as in the absolute case, see
  \cite[Kapitel 2]{Thesis} for more details.
\end{proof}

\medskip
It turns out that the description of simple vector bundles in terms of objects
of $\MPnode^{s}(n_{1},n_{2})$ and $\MPcusp^{s}(n_{1},n_{2})$ respectively,
allows us to give an explicit description of a universal family of stable
vector bundles of rank $n$ and degree $d$ on nodal and cuspidal Weierstra\ss{}
cubic curves.

\begin{proposition}\label{P:tripuniversalfam}
  Let $E$ be either a nodal or a cuspidal Weierstra\ss{} cubic curve,
  $0\le d<n$ be coprime integers and $G = \CC^*$ if $E$ is nodal and $G = \CC$
  if $E$ is cuspidal.
  If $M:G\lar \Mat_{n\times n}(\CC)$ is a holomorphic function
  such that the image of $M$ contains exactly one representative of
  each isomorphism class in $\MPnode^{s}(n-d,d)$ or
  $\MPcusp^{s}(n-d,d)$ respectively, then
  $$\kP :=
  \mathbb{G}\bigl(\kO_{\PP^1 \times G}^{n-d} \oplus \kO_{\PP^1 \times G}(1)^{d},
  \kO_{Z \times G}^n, \widetilde \mm\bigr)$$ is a universal family of stable
  vector bundles of rank $n$ and degree $d$ on $E$.
\end{proposition}

\begin{proof}
By Lemma \ref{L:VBequivMP} and Lemma \ref{L:VBequivMPcusp} and because each
simple vector bundle is stable (\cite[Cor.\ 4.5]{BK3}), for each stable vector
bundle $\kV$ of rank $n$ and degree $d$ there exists a unique $b\in G$ such
that $\kV\cong\kP_{b}$.

Let  $\kQ \in \VB\bigl(E\times M_E^{(n,d)}\bigr)$ be  a universal family of
stable vector bundles of rank $n$ and degree $d$ on $E$. The universal
property implies that there exists a unique morphism $f: G \to M_E^{(n,d)}$
such that $\kP = (\id \times f)^* \kQ \otimes \pr_{2}^{\ast}\kL$. Restricting on
$E\times\{b\}$ shows that $f(b)=[\kP_{b}]$, the point in
$M^{(n,d)}_E$ which corresponds to the isomorphism class of the vector bundle
$\kP_{b}$.
Since  $\kP_{b_1} \not\cong \kP_{b_2}$ for $b_1 \ne b_2 $, the map  $f$ is
bijective. Because $M_E^{(n,d)}$ is known to be smooth, bijectivity of $f$
implies that $f$ is an isomorphism. Hence, the pair $(G, \kP)$ represents the
moduli functor.
\end{proof}

\begin{corollary}\label{C:constrofuniv}
Let $E$ be a singular Weierstra\ss{} cubic curve,
$G = \CC^*$ if $E$ is nodal and $G = \CC$ if $E$ is cuspidal, $0 \le d < n$ be
coprime integers,
$\widetilde\kP = \kO_{\PP^1 \times G}^{n-d} \oplus \kO_{\PP^1 \times G}(1)^{d}$, and
$M$ be the  canonical form from the proofs of Theorems \ref{T:nodecurve} and
\ref{T:detcusp} respectively. Then the coherent sheaf
$$
\kP = \ker\bigl(\pi_{G\ast}\widetilde{\kP} \oplus
\eta_{G\ast}\mathcal{O}_{Z\times G}^{n}
\xrightarrow{\left({\zeta}^{\widetilde\kP} \, \, \mm\right)}
\nu_{G\ast}\mathcal{O}_{\widetilde{Z}\times G}^{n}\bigr)
$$
is a universal family of stable vector bundles of rank $n$ and degree $d$ on
the curve $E$.
\end{corollary}

To construct a trivialization of a vector bundle $\kA$ given by a matrix
$A\in\GL_{n}\bigl(\kO_{\widetilde{Z}\times G}\bigr)$ via the sequence
\eqref{eq:flat-family-node}, we pick a
holomorphic section $p = p_\xi \in H^0\bigl(\kO_{\PP^1}(1)\bigr)$
(in our applications we shall have $p = z_1$). This section induces a family of
trivializations  $\left.\left\{\xi_l: \kO_{\PP^1}(l)\right|_{\widetilde{U}}
\lar \kO_{\widetilde{U}}\right\}_{l \in \mathbb{Z}}$, compatible with tensor
products, from which we obtain isomorphisms
$\xi^{\widetilde\kA}: \widetilde\kA|_{\widetilde{U} \times G}
\lar \kO_{\widetilde{U} \times G}^n$, because
$\widetilde{\kA}=\bigoplus\kO_{\PP^{1} \times G}(l)^{k_{l}}$.
Here, $\widetilde{U}\subset\PP^{1}$ could be any subset on which
$p_\xi$ does not vanish, but we shall assume
$\widetilde{U}\cap\widetilde{Z}=\emptyset$, which implies that
$\pi_{G}|_{\widetilde{U}\times G}:\widetilde{U}\times G \lar U\times G$ is an
isomorphism.
Restricting the sequence \eqref{eq:flat-family-node} to the open subset
$U\times G\subset E\times G$, we obtain a trivialization $\xi^\kA$ of the
family $\kA$
\begin{equation}\label{E:trivofunifam}
\xi^\kA: \kA|_{U \times G} \stackrel{\ii}\lar \pi_{G\ast}\widetilde{\kA}|_{U \times G}
\xrightarrow{\pi_{G\ast}\left(\xi^{\widetilde\kA}\right)}  \kO_{U \times G}^n.
\end{equation}

\begin{remark}
Note that in the construction of all our families of stable vector bundles on
a singular Weierstra\ss{} cubic curve we have chosen  \emph{two} sections
$p_\zeta, p_\xi \in H^0\bigl(\kO_{\PP^1}(1)\bigr)$.  These choices are
independent of each other!
The section $p_{\zeta}$ is used to \emph{define} a family $\kA$ associated to
a matrix $A\in\GL_{n}\bigl(\kO_{\widetilde{Z}\times G}\bigr)$, whereas
$p_{\xi}$ is  used to \emph{trivialize} it on $U \times G$.
\end{remark}

Even though in our application of Theorem \ref{T:goodtrivsing} we shall
come back to the framework of Section \ref{S:GeomYB}, we consider here a more
general setting: $J$ and $M$ are arbitrary reduced complex spaces. We denote
the canonical projections as before by
$p: E \times J \times M \lar E \times M$,
$q: E \times J \times M \lar E \times J$,
$p_{\widetilde{Z}}: \widetilde{Z}\times J\times M \lar \widetilde{Z}\times M$
and
$q_{\widetilde{Z}}: \widetilde{Z}\times J\times M \lar \widetilde{Z}\times J$.

\begin{theorem}\label{T:goodtrivsing}
  Let $\tau:J\times M \lar M$ be a flat morphism,
  $\tilde\tau = \id_{E} \times \tau$,
  $\tau_{\widetilde{Z}}=\id_{\widetilde{Z}} \times \tau$ and fix
  $P\in\GL_{n}\bigl(\kO_{\widetilde{Z}\times M}\bigr)$ and
  $N\in\GL_{1}\bigl(\kO_{\widetilde{Z}\times M}\bigr)$ such that
  \begin{equation}
    \label{eq:matrixcompat}
    \tau_{\widetilde{Z}}^{\ast}P = q_{\widetilde{Z}}^{\ast}N \cdot
    p_{\widetilde{Z}}^{\ast} P\;.
  \end{equation}
  Denote by $\kP\in M_{E}^{(n,d)}$ and $\kN\in\Pic^{c}(E)$ the bundles defined
  by $P$ and $N$ respectively.
  Let $p_\xi \in H^0\bigl(\kO_{\PP^1}(1)\bigr)$ be a section which gives
  trivializations $\xi^{\kP}$ of $\kP$ on $U \times M$ and $\xi^{\kN}$ of
  $\kN$ on $U\times J$.

  Then there exists an isomorphism $\varphi:q^{\ast}\kN \otimes p^{\ast}\kP
  \lar \tilde\tau^{\ast}\kP$ which is represented by the identity  with
  respect to the trivializations $\xi^{\kP}$ and $\xi^{\kN}$  on
  $U\times J\times M$ (compare with Proposition \ref{P:goodtriv}).
\end{theorem}

\begin{proof}
  The bundles $\kP$ and $\kN$ are defined by the short exact sequences
  \[
  0\lar  \kN
  \xrightarrow{
    \left(
      \begin{smallmatrix}
        \jj  \\
        \qq
      \end{smallmatrix}
    \right)
  }
  \pi_{J\ast}\widetilde\kN \oplus \eta_{J\ast}\mathcal{O}_{Z\times J}
  \xrightarrow{
    \left(
      \begin{smallmatrix}
        {\zeta^{\widetilde\kN}} & \nn
      \end{smallmatrix}
    \right)
  }
  \nu_{J\ast}\mathcal{O}_{\widetilde{Z}\times J}
  \lar  0
  \]
  and
 $$
  0\lar  \kP
  \xrightarrow{
    \left(
      \begin{smallmatrix}
        \ii  \\
        \rr
      \end{smallmatrix}
    \right)
  }
  \pi_{M\ast}\widetilde\kP\oplus \eta_{M\ast}\mathcal{O}_{Z\times M}^n
  \xrightarrow{
    \left(
      \begin{smallmatrix}
        {\zeta^{\widetilde\kP}} &  \pp
      \end{smallmatrix}
    \right)
  }
  \nu_{M\ast}\mathcal{O}_{\widetilde{Z}\times M}^n
  \lar  0,
  $$
  where
  $\widetilde\kP = \bigoplus_{l\in\ZZ}\kO_{\PP^1 \times M}(l)^{k_{l}}$
  and $\widetilde\kN = \kO_{\PP^1 \times J}(c)$.
  Let
  $\widehat\kP = \bigoplus_{l\in\ZZ}\kO_{\PP^1 \times J \times M}(l)^{k_{l}}$
  and
  $\widehat\kN = \kO_{\PP^1 \times J \times M}(c)$.
  By Theorem \ref{thm:DGrel} we have short exact sequences
  \[
  0\lar  \tilde\tau^*\kP
  \lar
  (\pi_{J \times M})_\ast
  \widehat\kP\oplus (\eta_{J \times M})_\ast\mathcal{O}_{Z\times J \times M}^n
  \xrightarrow{\left(
      \begin{smallmatrix}
        {\zeta^{\widehat\kP}} &  \tilde\tau^*(\pp)
      \end{smallmatrix}
    \right)
  }
  (\nu_{J \times M})_\ast \mathcal{O}_{\widetilde{Z}\times J \times M}^n
  \lar  0,
  \]
  \[
  0\lar  p^*\kP
  \lar
  (\pi_{J \times M})_\ast
  \widehat\kP\oplus (\eta_{J \times M})_\ast\mathcal{O}_{Z\times J \times M}^n
  \xrightarrow{\left(
      \begin{smallmatrix}
        {\zeta^{\widehat\kP}} &  p^*(\pp)
      \end{smallmatrix}
    \right)
  }
  (\nu_{J \times M})_\ast \mathcal{O}_{\widetilde{Z}\times J \times M}^n
  \lar  0,
  \]
  and
  \[
  0\lar  q^*\kN
  \lar
  (\pi_{J \times M})_\ast \widehat\kN\oplus (\eta_{J \times M})_\ast\mathcal{O}_{Z\times J \times M}
  \xrightarrow{\left(
      \begin{smallmatrix}
        {\zeta^{\widehat\kN}} &  q^*(\nn)
      \end{smallmatrix}
    \right)
  }
  (\nu_{J \times M})_\ast \mathcal{O}_{\widetilde{Z}\times J \times M}
  \lar  0.
  \]
  Moreover, we also know that the sequence
 $$
  0\lar  q^*\kN \otimes p^*\kP
  \lar
  (\pi_{J \times M})_\ast
  \widehat\kP\oplus (\eta_{J \times M})_\ast\mathcal{O}_{Z\times J \times M}^n
  \stackrel{\beta}{\rightarrow}
  (\nu_{J \times M})_\ast \mathcal{O}_{\widetilde{Z}\times J \times M}^n
  \lar  0
  $$
  is exact, where we abbreviated $\beta=\left(
    \begin{smallmatrix}
      {\zeta^{\widehat\kP}} &  q^*(\nn) \boxtimes p^*(\pp)
    \end{smallmatrix}\right)$.
  Because $\tau$ is flat,
  the morphism $\tilde\tau^*(\pp)$ is induced by the
  matrix $\tau_{\widetilde{Z}}^*(P)$ and the morphism
  $q^*(\nn) \boxtimes p^*(\pp)$ by the matrix
  $q_{\widetilde{Z}}^*(N) \cdot p_{\widetilde{Z}}^*(P)$.
  Hence, using \eqref{eq:matrixcompat}, we obtain a commutative diagram
  $$
  \xymatrix
  {
    (\pi_{J \times M})_\ast \widehat\kP \oplus
    (\eta_{J \times M})_\ast \kO_{Z \times J \times M}^n
    \ar[rrrr]^{\bigl(\zeta^{\widehat\kP} \, \, \tilde\tau^*(\mm)\bigr)}
    & & & &  (\nu_{J \times M})_\ast \kO_{\widetilde{Z} \times J \times M}^n
    \\
    & & & & \\
    (\pi_{J \times M})_\ast
    \widehat\kP  \oplus (\eta_{J \times M})_\ast \kO_{Z \times J \times M}^n
    \ar[rrrr]^{\bigl(\zeta^{\widehat\kP} \, \,
      q^*(\nn) \boxtimes p^*(\pp) \bigr)}
    \ar[uu]^{
      \left(
        \begin{smallmatrix}
          \id & 0 \\
          0 & \id
        \end{smallmatrix}
      \right)
    }
    & & & &
    (\nu_{J \times M})_\ast \kO_{\widetilde{Z} \times J \times M}^n
    \ar[uu]_\id
  }
  $$
  But this implies that we have an induced morphism
  $\varphi: q^*\kN  \otimes p^*\kP \lar \tilde\tau^*\kP $.

  It remains to verify that the isomorphism $\varphi$ is the \emph{identity}
  with respect to the trivializations induced by $\xi^\kP$ and $\xi^\kN$.
  By Theorem \ref{thm:DGrel}, we have a commutative diagram
  $$
  \xymatrix
  {
    \tilde\tau^*\kP \ar[rrr]^{\tilde\tau^*(\ii)} \ar[d]_\id  & &  &
    \tilde\tau^*(\pi_M)_\ast \widetilde\kP \ar[d]^\can \\
    \tilde\tau^*\kP \ar[rrr]  & &  &
    (\pi_{J \times M})_\ast \widehat\kP \\
    q^*\kN \otimes p^*\kP \ar[rrr]^{q^*(\jj) \boxtimes p^*(\ii)}
    \ar[rrrd]_{q^*(\jj) \otimes p^*(\ii)}
    \ar[u]^\varphi
    & & &
    (\pi_{J \times M})_\ast \widehat\kP
    \ar[u]_\id \\
    & & & q^* \pi_{J \ast} \widetilde\kN \otimes p^* \pi_{M\ast} \widetilde\kP
    \ar[u]_\can
  }
  $$
  The key point is now that the trivializations
  $\xi^{\widetilde\kP}$ and $\xi^{\widetilde\kN}$ are the pull backs of
  trivializations $\xi^{\bar\kP}_{\widetilde{U}}:
  \bigoplus \kO_{\PP^{1}}(l)^{k_{l}}|_{\widetilde{U}} \lar \kO^n_{\widetilde{U}}$
  and $\xi^{\bar\kN}_{\widetilde{U}}: \kO_{\PP^{1}}(c)|_{\widetilde{U}}
  \lar \kO_{\widetilde{U}}$, on $\PP^1$.
  This implies that the base-change morphism $\can:
  \tilde\tau^*(\pi_M)_\ast \widetilde\kP \lar (\pi_{J \times M})_\ast \widehat\kP$
  is the \emph{identity} with respect to the trivializations
  $\xi^{\widetilde\kP}$ and $\xi^{\widehat\kP}$.
  It follows that, in these trivializations, the isomorphism $\varphi$ is the
  identity.
\end{proof}

After choosing representing pairs $(J^{d},\kL^{(d)})$, $(J,\kL)$ and $(M,\kP)$
for the functors $\underline{\Pic}^{d}$, $\underline{\Pic}^{0}$ and
$\underline{\Mf}_{E}^{n,d}$, in Section \ref{S:GeomYB} we have constructed a
morphism $\tau:J\times M\lar M$. To make this explicit,
let $G = \CC^*$ if $E$ is nodal and $G = \CC$ if $E$ is cuspidal. For
simplicity, we assume again $0\le d<n$.

We define $J=J^{d}=M=G$ with universal bundles $\kL^{(d)}$ and $\kL$ both
given by $(y, \id)$ in the nodal case and by $1+\varepsilon y$ in the cuspidal
case. We define the universal bundle $\kP$ to be given by
$(M_{n-d,d},\id)$ in the nodal case (proof of Theorem \ref{T:nodecurve})
and $1_{n}+\varepsilon M_{n-d,d}$ in the cuspidal case
(proof of Theorem \ref{T:detcusp}). Universality was shown in Corollary
\ref{C:constrofuniv}.

Using the notation of diagram \eqref{diag:jacob} in Section
\ref{S:GeomYB}, we obtain now $t^{e}=\id_{G}$ and $\det=\pm \id_{G}$, the sign
depending on $(n,d)$ only. Moreover, because
$\bigl((y_{1}),(1)\bigr)\bigl((y_{2}),(1)\bigr)=\bigl((y_{1}y_{2}),(1)\bigr)$,
the group structure on $J=G$ is
multiplication in the nodal case and because
$(1+y_{1}\varepsilon)(1+y_{2}\varepsilon)=1+(y_{1}+y_{2}\varepsilon)$, the
group structure is addition in the cuspidal case.
Therefore, $\tau=\sigma':G\times G\lar G$ has the
description $\tau(a,b)=a^{n}b$ in the nodal case and $\tau(a,b)=na+b$ in the
cuspidal case.

We also consider the vector bundle $\kP'$ of rank $n$ and degree $d$ on
$G\times E$, which is given by $(\widetilde{N}_{n-d,d},\id)$ in the nodal case
(Remark \ref{R:noderank2}) and $1_{n}+\varepsilon \widetilde{N}_{n-d,d}$
in the cuspidal case (Remark \ref{R:eqvcusp}). These were constructed in such
a way that \eqref{eq:matrixcompat} holds with respect to the morphism
$\tau':G\times G\lar G$ given by $\tau'(a,b)=ab$
(respectively $\tau'(a,b)=a+b$).

If $f_{n}:G\lar G$ is given by $f_{n}(t)=t^{n}$ in the nodal case and by
$f_{n}(t)=nt$ in the cuspidal case, we have
$(\id_{E}\times f_{n})\circ \tilde\tau' =
\tilde\tau\circ(\id_{E\times J}\times f_{n})$, because $J=G$ is abelian.

From Remark \ref{R:noderank2} and Remark \ref{R:eqvcusp} respectively, it is
clear that $(\id_{E}\times f_{n})^{\ast} \kP$ and $\kP'$ are isomorphic after
restriction to a fibre. This implies that these two bundles are locally
isomorphic (with respect to the basis $G$, which is reduced). Equivalently, up
to a twist by the pull-back of a line bundle on $G$, $(\id_{E}\times
f_{n})^{\ast} \kP$ and $\kP'$ are isomorphic. As $G$ is a non-compact
Riemann surface, we even get $(\id_{E}\times f_{n})^{\ast} \kP \cong
\kP'$, but we do not need this in the sequel.

\begin{corollary}\label{C:goodtrivsing}
  The morphism $\tau'$ and the bundles $\kP'$ and $\kL$ satisfy the properties
  of Theorem \ref{T:goodtrivsing}. Moreover, each point of $G$ has an open
  neighbourhood $M'\subset G$ such that there exits an isomorphism
  $\varphi:  q^*\kL \otimes p^*\kP|_{E \times J \times M'} \lar
  \tilde\tau^*\kP|_{E \times J \times M'}$,
  which is represented by the identity  with
  respect to the trivializations $\xi^{\kP}$ and $\xi^{\kL}$  on
  $U\times J\times M'$.
\end{corollary}
\begin{proof}
  Because, up to a local isomorphism, $\tau$ is isomorphic to
  addition of complex numbers, it is flat. Now, the first statement is clear
  from the above.
  The prove the second, chose a sufficiently small open neighbourhood
  $M'\subset G$ around a given point on $G$ such that $f_{n}$ restricted to
  $M'$ is an isomorphism and $(\id_{E}\times f_{n})^{\ast} \kP$ is isomorphic
  to $\kP'$ over $E\times M'$. Because $q\circ (\id_{E\times J}\times f_{n}) = q$,
  $p\circ (\id_{E\times J}\times f_{n}) = (\id_{E}\times f_{n})\circ p$ and
  $(\id_{E}\times f_{n})\circ \tilde\tau' =
  \tilde\tau\circ(\id_{E\times J}\times f_{n})$, with the aid of the
  isomorphism $\id_{E\times J} \times f_{n}\big|_{M'}$ the claim follows from
  the first part of the corollary.
\end{proof}

In the cuspidal case, $f_{n}$ is an isomorphism, hence we may choose $M'=G$.
In the nodal case, however, the matrix $\widetilde{N}_{n_{1},n_{2}}$ does not
descent to $M=G$, as it involves taking $n$-th roots; $f_{n}$ is an unramified
$n$-fold cover.
Therefore, the isomorphism $\varphi$ exists only locally on $M$ in this case.

\begin{remark}
  The map $x\mapsto\kO_{E}(x)$ gives a canonical isomorphism
  $\breve{E} \rightarrow \Pic^1_E$. Let $e\in E$  be the point which has
  coordinate $z=1$ in the nodal case (Subsection \ref{SS:node}) and $z=0$ in
  the cuspidal case (Subsection \ref{SS:cusp}).
  This point corresponds to the neutral element of $J$ under the isomorphisms
  $\breve{E}\stackrel{\can}{\lar} J^{1} \stackrel{t^{e}}{\longleftarrow} J$
  (see Section \ref{S:GeomYB}). This isomorphism and our choice of the
  representing pair $(J,\kL)$ with $J=G$ induce coordinates on $\breve{E}$.

  Lemma \ref{L:idenonnode} shows that these coordinates coincide
  with the ones induced by the coordinates on $\PP^1$ and the normalization
  morphism $\PP^1 \setminus \widetilde{Z} \stackrel{\pi}\lar
  \breve{E}$. However, by Lemma \ref{L:idenoncusp} we see that in the case of
  a cuspidal curve these two choices are \emph{different}; they are
  related by the involution of $\CC$ mapping $z$ to $-z$.
\end{remark}

\begin{remark}
  If $(n, d) \in \mathbb{Z}^+ \times \mathbb{Z}$ are coprime integers, we let
  $c$ be the unique integer for which $n_{1}=(1+c)n-d>0$ and
  $n_{2}=d-cn\ge 0$. With the aid of the equivalence between
  $\VB^{(c,c+1)}_{n_{1},n_{2}}(E)$ and $\VB^{(0,1)}_{n_{1},n_{2}}(E)$, given
  by the tensor product with $\kO_{E}(ce)$, it can be shown that all the
  results of this section are also valid for such pairs $(n,d)$.
\end{remark}

\section{Computations of $r$-matrices for singular Weierstra\ss{} curves}
\label{S:singrm}

Let $E$ be a singular Weierstra\ss{} cubic curve, $\Omega_E$ the sheaf of
regular holomorphic  1-forms, $\omega  \in H^0(\Omega_E)$ a no-where vanishing
global section. As usual, for a pair of coprime integers $(n, d) \in
\mathbb{Z}^+ \times \mathbb{Z}$, let $M = M_E^{(n, d)}$ be the moduli space
of stable holomorphic vector bundles of rank $n$ and degree $d$ on $E$,
$\kP = \kP(n, d) \in \VB(E \times M)$ be a universal family and
$\xi^\kP: \kP|_{U \times M'} \lar \kO|_{U \times M'}^n$ its trivialization, as
constructed in the previous section.
Recall that these data define the germ of  a meromorphic function
$$
\tilde{r} = \tilde{r}^\xi : \bigl(M \times M \times E \times E, o\bigr) \lar
\Mat_{n \times n}(\CC) \otimes \Mat_{n \times n}(\CC),
$$ whose value at the  point $(v_1, v_2; y_1, y_2)$,
where $v_1 \ne v_1$ and $y_1 \ne y_2$,  is defined via the commutative diagram
(\ref{E:maindiag}).
 Our next goal is to get explicit formulae to calculate
the morphisms $\res^{\kP^{v_1}, \kP^{v_2}}_{y_1}(\omega)$ and
$\ev^{\kP^{v_1},\kP^{v_2}(y_1) }_{y_2}$ in the case of nodal and cuspidal
 Weierstra\ss{} cubic curves.  To do this, we consider first the case of vector
bundles on a projective line $\PP^1$.

\subsection{Residue and evaluation morphisms on $\PP^1$}\label{SS:resandevonp1}
Let $(z_0: z_1)$ be homogeneous coordinates on $\PP^1$, $0 = (1: 0)$ and
$U = \bigl\{(z_0: z_1)| z_0 \ne 0 \bigr\}$. Let $z = \frac{\displaystyle z_1}{\displaystyle z_0}$ be a local coordinate
on $U$. In what follows, we shall use the identification
$\CC \stackrel{\cong}\lar U$ mapping $x \in \CC$ to $(1:x) \in U$.

Let $\kV = \bigoplus_{j \in \ZZ} \kO_{\PP^1}(j)^{n_j}$ and $\kW = \bigoplus_{i \in \ZZ} \kO_{\PP^1}(i)^{m_i}$
be a pair of vector bundles on $\PP^1$ of ranks $n$ and $m$ respectively. Then a morphism
$F \in \Hom_{\PP^1}(\kV, \kW)$ can be written in  matrix form: $F = (F_{ij})$, where
$F_{ij} \in \Mat_{m_i \times n_j}\bigl(\CC[z_0, z_1]_{i-j}\bigr)$.

Consider the set of trivializations $\bigl\{\xi_l: \kO_{\PP^1}(l)|_{U} \lar \kO_U\bigr\}_{l \in \ZZ}$
mapping a local section $p$ of $\kO_{\PP^1}(l)$ to the holomorphic function
$\frac{\displaystyle p}{\displaystyle z_0^l|_U}$. They induce trivializations
$\xi^\kV: \kV|_U \lar \kO_U^n$ and $\xi^\kW: \kW|_U \lar \kO_U^m$. Let $x, y \in U$ and
$\omega$ be a meromorphic differential form on $\PP^1$ holomorphic at $x$.
Let  $\bar\xi^\kV$ be the morphism $$\kV|_U \otimes \CC_x
 \xrightarrow{\xi^\kV \otimes \mathsf{id}}
 \kO_U^n \otimes \CC_x
\stackrel{\can}\lar \CC_x^n;$$
 the morphism $\bar\xi^\kW$ is defined in a similar way.
Our goal is to get
explicit formulae to calculate the morphisms:
$$
\Hom_{\PP^1}\bigl(\kV, \kW(x)\bigr)
\xrightarrow{\res^{\kV, \kW}_x(\omega)}
\Hom_{\PP^1}\bigl(\kV \otimes \CC_x, \kW \otimes \CC_x\bigr)
\xrightarrow{\conj(\bar\xi^{\kV}, \bar\xi^{\kW})} \Mat_{m \times n}(\CC)
$$
and
$$
\Hom_{\PP^1}\bigl(\kV, \kW(x)\bigr)
\xrightarrow{\ev^{\kV, \kW(x)}_y}
\Hom_{\PP^1}\bigl(\kV \otimes \CC_y, \kW \otimes \CC_y\bigr)
\xrightarrow{\conj(\bar\xi^{\kV}, \bar\xi^{\kW})} \Mat_{m \times n}(\CC).
$$
Let $\sigma = z_1 - xz_0 \in H^0\bigl(\kO_{\PP^1}(1)\bigr)$ be a section such that
$\mathsf{div}(\sigma) = [x]$. Then $\sigma$ defines an isomorphism
$\kO_{\PP^1}(1) \lar \kO_{\PP^1}(x)$ mapping a global section $p = p(z_0, z_1)$ to
the meromorphic function $\frac{\displaystyle p}{\displaystyle \sigma}$. Moreover, it induces
an isomorphism
$$
t_\sigma: \Hom_{\PP^1}\bigl(\kV, \kW(1)\bigr) \lar \Hom_{\PP^1}\bigl(\kV, \kW(x)\bigr).
$$
Let $\kV' = \kV|_U$ and $\kW' = \kW|_U$.  By Proposition \ref{P:functofres} we have a commutative diagram
$$
\xymatrix
{
 \Hom_{\PP^1}\bigl(\kV, \kW(x)\bigr) \ar[rrr]^-{\res_x^{\kV, \kW}(\omega)}
 \ar[d]_{\can}
 & & &
 \Hom_{\PP^1}\bigl(\kV\otimes \CC_x, \kW\otimes \CC_x\bigr) \ar[d]^\can
 \\
 \Hom_{U}\bigl(\kV', \kW'(x)\bigr) \ar[rrr]^-{\res_x^{\kV', \kW'}(\omega)}
 \ar[d]_{\conj\bigl(\xi^\kV, \, \xi^\kW(x)\bigr)}
 & & &
 \Hom_{U}\bigl(\kV'\otimes \CC_x, \kW'\otimes \CC_x\bigr)
 \ar[d]^{\conj\bigl(\bar\xi^\kV, \, \bar\xi^\kW\bigr)}
 \\
 \Hom_{U}\bigl(\kO_U^n, \kO_U^m(x)\bigr) \ar[rrr]^-{\res_x(\omega)} \ar[d]_\can
 & & &
 \Hom_{U}\bigl(\CC_x^n, \CC_x^m\bigr) \ar[d]^\can
 \\
 \Mat_{m \times n}\bigl(O(x)\bigr) \ar[rrr]^-{\res_x(\omega)}
 & &  &
 \Mat_{m \times n}(\CC),
}
$$
where $O(x)$ denotes the vector space of meromorphic functions on $U$ which
have at most a pole of order one at $x$.
Let $T$ be the composition $$\Hom_{\PP^1}\bigl(\kV, \kW(1)\bigr)
\stackrel{t_\sigma}\lar \Hom_{\PP^1}\bigl(\kV, \kW(x)\bigr)
\xrightarrow{\conj\bigl(\xi^\kV, \, \xi^{\kW}(x)\bigr)} \Mat_{m \times n}\bigl(O(x)\bigr).
$$
Then we  have the following result:
if $F \in \Hom_{\PP^1}\bigl(\kV, \kW(1)\bigr)$ then $T(F) =
\frac{\displaystyle F(1, z)}{\displaystyle z - x}$. Let $\omega = g(z) dz$ then  by Lemma \ref{L:resexplit} we have:
$$
\res_x(\omega)\Bigl(\frac{F(1,z)}{z-x}\Bigr) =  g(x) F(1, x).
$$
\begin{corollary}
In the above notation, the morphism $$\overline{\res}_x := \res_x(\omega) \, \circ \, T: \quad \Hom_{\PP^1}\bigl(\kV, \kW(1)\bigr)
\lar \Mat_{m \times n}(\CC)$$ has the following form:
\begin{itemize}
\item if $\omega = \frac{\displaystyle dz}{\displaystyle z}$ then $F(z_0, z_1)$ is mapped to
$\frac{\displaystyle 1}{\displaystyle x} F(1, x)$;
\item if $\omega = dz$ then $F(z_0, z_1)$ is mapped to
$F(1, x)$.
\end{itemize}
\end{corollary}

In a similar way, we compute the morphism $\ev^{\kV, \kW(x)}_y$. Indeed, by Proposition \ref{L:functofev}
we have a commutative diagram
$$
\xymatrix
{
  \Hom_{\PP^1}\bigl(\kV, \kW(x)\bigr) \ar[rrr]^-{\ev_y^{\kV, \kW(x)}}
  \ar[d]_{\can}
  & & &
  \Hom_{\PP^1}\bigl(\kV\otimes \CC_y, \kW\otimes \CC_y\bigr) \ar[d]^\can
  \\
  \Hom_{U}\bigl(\kV', \kW'(x)\bigr) \ar[rrr]^-{\ev_y^{\kV', \kW'(x)}}
  \ar[d]_{\conj\bigl(\zeta^\kV, \,  \zeta^\kW(x)\bigr)}
  & & &
  \Hom_{U}\bigl(\kV'\otimes \CC_y, \kW'\otimes \CC_y\bigr)
  \ar[d]^{\conj\bigl(\bar\zeta^\kV, \, \bar\zeta^\kW\bigr)}
  \\
  \Hom_{U}\bigl(\kO_U^n, \kO_U^m(x)\bigr) \ar[rrr]^-{\ev_y} \ar[d]_\can
  & & &
  \Hom_{U}\bigl(\CC_y^n, \CC_y^m\bigr) \ar[d]^\can
  \\
  \Mat_{m \times n}\bigl(O(x)\bigr) \ar[rrr]^-{\ev_y}
  & &  &
  \Mat_{m \times n}(\CC)
}
$$
and for $A(z) \in \Mat_{m \times n}\bigl(O(x)\bigr)$ we have
$
\ev_y\bigl(A(z)\bigr)  =  A(y).
$
\begin{corollary}
In the above notations, the morphism $$\overline{\ev}_y := \ev_y  \, \circ\,  T: \quad \Hom_{\PP^1}\bigl(\kV, \kW(1)\bigr)
\lar \Mat_{m \times n}(\CC)$$
maps a matrix $F(z_0, z_1) \in \Hom_{\PP^1}\bigl(\kV, \kW(1)\bigr)$ to $\frac{\displaystyle 1}{\displaystyle y-x}F(1, y)$.
\end{corollary}

\begin{remark}
The above morphisms can be included into the following diagram:
$$
\xymatrix
{
 & \Hom_{\PP^1}\bigl(\kV, \kW(1)\bigr) \ar[d]_{t_\sigma}
 \ar@/_10pt/[ddl]_{\overline{\res}_x} \ar@/^10pt/[ddr]^{\overline{\ev}_y} & \\
 & \Hom_{\PP^1}\bigl(\kV, \kW(x)\bigr) \ar[ld]^{\res_x(\omega)} \ar[rd]_{\ev_y}& \\
 \Mat_{m \times n}(\CC) &  & \Mat_{m \times n}(\CC)
}
$$
In particular, the linear maps  $\overline{\res}_x$ and $\overline{\ev}_y$ depend on the choice of
a section $\sigma \in H^0\bigl(\kO_{\PP^1}(1)\bigr)$ vanishing at $x$: such a $\sigma$ is
determined uniquely only up to a non-zero constant. However,
as we shall see below, this choice does not affect the final formula to
compute the triple Massey product $\tilde{r}^{\kP^{v_1}, \kP^{v_1}}_{y_1, y_2}(\omega)$
on a cubic curve.
 \end{remark}

\medskip
 \subsection{Residue and evaluations maps on singular Weierstra\ss{} curves}
Let $E$ be a singular Weierstra\ss{} cubic curve, $s \in E$ its unique singular point, $\pi: \PP^1 \lar E$ normalization of $E$. We choose homogeneous coordinates on $E$ in such a way that
$\pi^{-1}(s) = \{(0:1), (1:0)\}$ if $E$ is nodal and $\pi^{-1}(s) = \{(0:1)\}$ if $E$ is cuspidal.
Let $\breve{E}$ be the regular part of $E$ then the isomorphism
$\pi: \PP^1 \setminus \pi^{-1}(s) \lar \breve{E}$ induces local coordinates on
$\breve{E}$. In what follows, we shall identify a point $y \in \breve{E} = G$ with
its preimage $\tilde{y} = (1:y) \in \PP^1$.
Let $0 < d < n$ be a pair of mutually prime integers and $\kP \in \VB(E \times M)$ a universal
family of stable vector bundles of rank $n$ and degree $d$ on $E$.
Recall that for $v_1 \ne v_2 \in M$  and $y_1 \ne y_2 \in \breve{E}$ we have a commutative diagram
$$
\small{
\xymatrix
{
\Hom_E\bigl(\kP^{v_1}|_{{y}_1}, \kP^{v_2}|_{y_1}\bigr)
\ar[rr]^{\pi^*} & &
\Hom_{\PP^1}\bigl(
\pi^*\kP^{v_1}|_{\tilde{y}_1}, \pi^*\kP^{v_2}|_{\tilde{y}_1}\bigr)
\\
 \Hom_E\bigl(\kP^{v_1}, \kP^{v_2}(y_1)\bigr)
\ar[u]^{\res_{y_1}^{\kP^{v_1}, \kP^{v_2}}(\omega)}
\ar[rr]^{\pi^*}
\ar[d]_{\ev_{y_2}^{\kP^{v_1}, \kP^{v_2}(y_1)}} &  &
\Hom_{\PP^1}\bigl(\pi^*\kP^{v_1}, \pi^*\kP^{v_2}(\tilde{y}_1)\bigr)
\ar[u]_{\res_{\tilde{y}_1}^{\pi^*\kP^{v_1}, \pi^*\kP^{v_2}}(\tilde\omega)}
\ar[d]^{\ev_{\tilde{y}_2}^{\pi^*\kP^{v_1}, \pi^*\kP^{v_2}(\tilde{y}_1)}} \\
\Hom_E\bigl(\kP^{v_1} \otimes \CC_{y_2}, \kP^{v_2} \otimes \CC_{y_2}\bigr)
\ar[rr]^{\pi^*} & &
\Hom_{\PP^1}\bigl(
\pi^*\kP^{v_1}|_{\tilde{y}_2}, \pi^*\kP^{v_2}|_{\tilde{y}_2}\bigr).
}
}
$$
In particular, the computation of the morphisms $\res_{y_1}^{\kP^{v_1}, \kP^{v_2}}(\omega)$ and
$\ev_{y_2}^{\kP^{v_1}, \kP^{v_2}(y_1)}$ can be reduced to an analogous computation
on $\PP^1$.

Let $E'$ and $M'$ be open neighbourhoods of $e \in E$ and $m \in M$,
$\xi: \kP|_{E' \times M'} \lar \kO_{E' \times M'}^n$ be a trivialization of the universal family
$\kP$ which is compatible with the action of the Jacobian. For $v \in M'$ let
$\kP^v = \kP|_{E \times v}$ and
$\xi^v$ be the induced trivialization $\xi^v: \kP^v|_{E'} \lar \kO_{E'}^n$.
Next,   for $y \in E'$ let $\bar\zeta^v$ be the corresponding isomorphism
$\kP^v \otimes \CC_y \lar \CC_y^n$.  Our goal is to compute the value of the geometric
$r$--matrix $r^\xi$ at the point $(v_1, v_2; y_1, y_2)$. This linear morphism
$\tilde{r}^\xi(v_1, v_2; y_1, y_2)$
is defined via the commutative diagram
$$
\xymatrix
{
  \Hom_E(\kP^{v_1} \otimes \CC_{y_1}, \kP^{v_2} \otimes \CC_{y_1})
  \ar[rrrr]^-{\conj\bigl(\bar\xi^{v_1}, \bar\xi^{v_2}\bigr)}
  \ar[d]_{\tilde{r}^{\kP^{v_1}, \kP^{v_2}}_{y_1, y_2}(\omega)}
  & & & &
  \Mat_{n \times n}(\CC)  \ar[d]^{\tilde{r}^\xi(v_1, v_2; y_1, y_2)}
  \\
  \Hom_E(\kP^{v_1} \otimes \CC_{y_2}, \kP^{v_2} \otimes \CC_{y_2})
  \ar[rrrr]^-{\conj\bigl(\bar\xi^{v_1}, \bar\xi^{v_2}\bigr)}
  & & & &
  \Mat_{n \times n}(\CC).
}
$$
Let $\widetilde\kP = \kO_{\PP^1}^{n-d} \oplus \kO_{\PP^1}^{d}(1)$ and
$\widehat\kP$ be the pull-back of $\widetilde\kP$ to $\PP^1 \times M$. Let
$\zeta_1: \kO_{\PP^1}(1)|_{\widetilde{Z}} \lar \kO_{\widetilde{Z}}$ be the isomorphism
used in Section \ref{S:triples} in the construction of the category $\Tri(E)$.

Recall that the universal family $\kP = \kP(n, d)$ was defined using the following short exact sequence
$$
  0\lar  \kP
  \xrightarrow{
\left(
\begin{smallmatrix}
\ii  \\
\pp
\end{smallmatrix}
\right)
}
  \pi_{M\ast}\widehat{\kP} \oplus \eta_{M\ast}\mathcal{O}_{Z\times M}^{n}
  \xrightarrow{\left(\begin{smallmatrix}{\zeta}^{\widehat\kP} \, \, \mm\end{smallmatrix}\right)}
  \nu_{M \ast}\mathcal{O}_{\widetilde{Z}\times M}^{n}
  \lar  0,
$$
In particular, the trivialization $\hat{\xi}^{\widehat\kP}: \widehat\kP|_{U \times M'} \lar
\kO^n_{U \times M'}$ induces the trivialization
$$
\xi^\kP: \kP|_{E' \times M'} \stackrel{\ii}\lar \pi_{M\ast}\widehat{\kP}|_{E' \times M'}
\stackrel{\xi^{\widehat\kP}}\lar \kO_{E' \times M'}^n.
$$
Moreover, by Theorem \ref{thm:Drozd-Greuel} we know that the morphism
$$
\tilde{\ii}: \pi_M^* \kP \xrightarrow{\pi^*(\ii)} \pi_M^* \pi_{M\ast}\widehat{\kP}
\stackrel{\can}\lar \widehat{\kP}
$$
is an isomorphism. Hence, we have the following commutative diagram
$$
\scriptsize{
\xymatrix
{
\Hom_E\bigl(\kP^{v_1} \otimes \CC_{{y}_1}, \kP^{v_2} \otimes \CC_{y_1}\bigr)
\ar[rrd]_{\conj\bigl(\bar\xi^{v_1}, \,  \bar\xi^{v_2}\bigr)}
\ar[rrrr]^{\conj \, \circ \,  \pi^*} \ar[dddd]_{\tilde{r}^{\kP^{v_1}, \kP^{v_2}}_{y_1, y_2}(\omega)}
& & & &
\Hom_{\PP^1}(\widetilde\kP\otimes \CC_{\tilde{y}_1}, \widetilde\kP\otimes \CC_{\tilde{y}_1})
\ar[lld]^{\conj\bigl(\tilde\xi^{v_1}, \,  \tilde\xi^{v_2}\bigr)} \\
& & \Mat_{n \times n}(\CC) \ar[dd]_{\tilde{r}^\xi(v_1, v_2; y_1, y_2)} & & \\
& & & & \Pi^{v_1, v_2}_{y_1} \ar[uu]_{\res_{\tilde{y}_1}} \ar[dd]^{\ev_{\tilde{y}_2}}\\
& & \Mat_{n \times n}(\CC) & & \\
\Hom_E\bigl(\kP^{v_1} \otimes \CC_{{y}_2}, \kP^{v_2} \otimes \CC_{y_2}\bigr)
\ar[rrrr]^{\conj \, \circ \,  \pi^*} \ar[rru]^{\conj\bigl(\bar\xi^{v_1}, \,  \bar\xi^{v_2}\bigr)}
& & & &
\Hom_{\PP^1}(\widetilde\kP\otimes \CC_{\tilde{y}_2}, \widetilde\kP\otimes \CC_{\tilde{y}_2}) \ar[llu]_{\conj\bigl(\tilde\xi^{v_1}, \,  \tilde\xi^{v_2}\bigr)}\\
}
}
$$
where $\Pi^{v_1, v_2}_{y_1} = \mathsf{Im}\Bigl(\Hom_E\bigl(\kP^{v_1}, \kP^{v_2}(y_1)\bigr)
\xrightarrow{\conj\bigl(\tilde{\ii}^{v_1}, \, \tilde{\ii}^{v_2}(y_1)\bigr) \, \circ \, \pi^*}
\Hom_{\PP^1}\bigl(\widetilde\kP, \widetilde\kP(\tilde{y}_1)\bigr)\Bigr)$.

The morphisms
$\res_{\tilde{y}_1}: \Pi^{v_1, v_2}_{y_1} \lar \Hom_{\PP^1}(\widetilde\kP\otimes \CC_{\tilde{y}_1}, \widetilde\kP\otimes \CC_{\tilde{y}_1})$ and
$\ev_{\tilde{y}_2}: \Pi^{v_1, v_2}_{y_1} \lar \Hom_{\PP^1}(\widetilde\kP\otimes \CC_{\tilde{y}_2}, \widetilde\kP\otimes \CC_{\tilde{y}_2})$ are isomorphisms. Hence, in order to compute
the linear map $r^\xi(v_1, v_2; y_1, y_2)$,  it suffices to get  explicit formulae
for  the morphisms
$$
\Hom_{\PP^1}(\widetilde\kP\otimes \CC_{\tilde{y}_1}, \widetilde\kP\otimes \CC_{\tilde{y}_1})
\xrightarrow{\res_{\tilde{y}_1}^{-1}} \Pi^{v_1, v_2}_{y_1}
\xrightarrow{\ev_{\tilde{y}_2}} \Hom_{\PP^1}(\widetilde\kP\otimes \CC_{\tilde{y}_2}, \widetilde\kP\otimes \CC_{\tilde{y}_2}).
$$
Consider the short exact sequence
$$
0 \lar \kO_E(y_1)
\xrightarrow{
\left(
\begin{smallmatrix}
\jj \\
\qq
\end{smallmatrix}
\right)
} \pi_*\bigl(\kO_{\PP^1}(1)\bigr) \oplus \eta_*\kO_Z
\xrightarrow{
\left(
\begin{smallmatrix}
\zeta_1 &  \nn_y
\end{smallmatrix}
\right)
}
\nu_*\kO_{\widetilde{Z}} \lar 0.
$$
By Theorem \ref{thm:Drozd-Greuel} we know that the morphism $$
\tilde{\jj}: \kO_{\PP^1}(\tilde{y}_1) =
\pi^*\bigl(\kO_E(y_1)\bigr) \xrightarrow{\pi^*(\jj)} \pi^* \pi_* \kO_{\PP^1}(1)
\xrightarrow{\can} \kO_{\PP^1}(1)$$
 is an isomorphism. Let $\sigma = \tilde{\jj}(1) \in
H^0\bigl(\kO_{\PP^1}(1)\bigr)$. Then the morphism
$$t_\sigma: \Hom_{\PP^1}\bigl(\widetilde\kP, \widetilde\kP(1)\bigr)
\lar \Hom_{\PP^1}\bigl(\widetilde\kP, \widetilde\kP(y_1)\bigr)$$ introduced in Subsection
\ref{SS:resandevonp1}
is the inverse of the morphism
$$\tilde{\jj}_*: \Hom_{\PP^1}\bigl(\widetilde\kP, \widetilde\kP(y_1)\bigr) \lar \Hom_{\PP^1}\bigl(\widetilde\kP, \widetilde\kP(1)\bigr)$$ induced by $\tilde{\jj}$.
Consider the short exact sequence
$$
0 \lar \kP^{v_2}
\xrightarrow{
\left(
\begin{smallmatrix}
\ii^{v_2} \\
\qq^{v_2}
\end{smallmatrix}
\right)
} \pi_*\widetilde\kP \oplus \eta_{*}\kO_Z^n
\xrightarrow{
\left(
\begin{smallmatrix}
\zeta^{\widetilde\kP} &  \mm(v_2)
\end{smallmatrix}
\right)
}
\nu_*\kO_{\widetilde{Z}}^n  \lar 0.
$$
We know that the sequence
$$
0 \lar \kP^{v_2}(y_1)
\xrightarrow{
\left(
\begin{smallmatrix}
\mathsf{k} \\
\mathsf{l}
\end{smallmatrix}
\right)
} \pi_*\widetilde{\kP}(1) \oplus \eta_{*}\kO_Z^n
\xrightarrow{
\left(
\begin{smallmatrix}
\zeta^{\widetilde{\kP}(1)} &  \mm
\end{smallmatrix}
\right)
}
\nu_*\kO_{\widetilde{Z}}^n  \lar 0
$$
is exact, where $\mathsf{k}$ is the morphism
$
\kP^{v_2} \otimes \kO_E(y_1)
\xrightarrow{\ii^{v_1} \otimes \jj}
\pi_*\widetilde\kP \otimes \pi_* \kO_{\PP^1}(1)
\stackrel{\can}\lar \pi_*\bigl(\widetilde{\kP}(1)\bigr)
$
and $\mm$ corresponds to the tensor product of matrices $\mm(v_2)$ and $\nn$.

\begin{lemma}
In the above notation, the following diagram is commutative:
$$
\xymatrix
{
\widetilde{\kP}(\tilde{y}_1) \ar[rr]^{\tilde \jj } & & \widetilde{\kP}(1) \\
& \pi^*\bigl(\kP^{v_2}(y_1)\bigr) \ar[ru]_{\tilde{\mathsf{k}}} \ar[lu]^{\tilde{\jj}(\tilde{y}_1)}&
}
$$
\end{lemma}

\begin{proof}\label{L:somerel}
This follows from the definition of the morphism $k$ and the fact that the
diagram
$$
\xymatrix
{
  \pi^*\bigl(\kP(y_1)\bigr) \ar[rrr] \ar[d]_\can
  & & &
  \pi^* \pi_* \widetilde\kP \otimes \pi^* \pi_* \kO_{\PP^1}(1)
  \ar[rrr]^-{\pi^*(\can)} \ar[d]_=
  & & &\pi^* \pi_* \widetilde\kP(1) \ar[d]^\can
  \\
  \pi^* \kP \otimes \kO_{\PP^1}(\tilde{y}_1) \ar[rrr]
  & & &
  \pi^* \pi_* \widetilde\kP \otimes \pi^* \pi_* \kO_{\PP^1}(1)
  \ar[rrr]^-{\pi^*(\can) \otimes \pi^*(\can)}
  & & &
  \widetilde\kP \otimes \kO_{\PP^1}(1)
}
$$
is commutative.
\end{proof}

\noindent
Using Lemma \ref{L:somerel},  we obtain  the following result.
\begin{proposition}
The following diagram is commutative:
$$
\xymatrix
{
  \Hom_{\Tri(E)}\bigl(T_{\kP^{v_1}}, T_{\kP^{v_2}(y_1)}\bigr)
  \ar[rrr]^-{\mathsf{For}} \ar[d]_{\mathbb{G}}
  & & &
  \Hom_{\PP^1}\bigl(\widetilde\kP, \widetilde\kP(1)\bigr) \ar[d]^{t_\sigma}
  \\
  \Hom_E\bigl(\kP^{v_1}, \kP^{v_2}(y_1)\bigr)
  \ar[rrr]^-{\conj\bigl(\tilde{\ii}^{v_1}, \,
    \tilde{\ii}^{v_2}(\tilde{y}_1)\bigr) \, \circ\, \pi^*}
  & & & \Hom_{\PP^1}\bigl(\widetilde\kP, \widetilde\kP(\tilde{y}_1)\bigr),
}
$$
where $T_{\kP^v} = \FF(\kP^v)$ for all $v \in M$.
\end{proposition}

\begin{proof}
First note that by Theorem \ref{thm:Drozd-Greuel} the diagram
$$
\xymatrix
{
& \Hom_{\Tri(E)}\bigl(T_{\kP^{v_1}}, T_{\kP^{v_2}(y_1)}\bigr) \ar[ddl]_{\mathbb{G}}
\ar[ddr]^{\mathsf{For}}& \\
& & \\
\Hom_E\bigl(\kP^{v_1}, \kP^{v_2}(y_1)\bigr) \ar[rr]^{\conj\bigl(\tilde{\ii}^{v_1}, \, \tilde{\mathsf{k}} \bigr) \, \circ\, \pi^*}  & & \Hom_{\PP^1}\bigl(\widetilde\kP, \widetilde\kP(1)\bigr)
}
$$
is commutative. By Lemma \ref{L:somerel} the diagram
$$
\xymatrix
{
  &  \Hom_E\bigl(\kP^{v_1}, \kP^{v_2}(y_1)\bigr)
  \ar[ddl]_{\conj\bigl(\tilde{\ii}^{v_1}, \tilde{\ii}^{v_2}(\tilde{y}_1)\bigr) \, \circ \, \pi^*}
\ar[ddr]^{\conj\bigl(\tilde{\ii}^{v_1}, \tilde{\mathsf{k}}\bigr) \, \circ \, \pi^*} & \\
& & \\
\Hom_{\PP^1}\bigl(\widetilde\kP, \widetilde\kP(\tilde{y}_1)\bigr)
 \ar[rr]^{\tilde{\jj}_*}  & & \Hom_{\PP^1}\bigl(\widetilde\kP, \widetilde\kP(1)\bigr)
}
$$
is commutative, too. Patching both diagrams together and using that
$\tilde{\jj}_* = t_\sigma^{-1}$, we get the claim.
\end{proof}

\begin{corollary}
Let $\widetilde\Pi^{v_1, v_2}_{y_1} =
 \mathsf{Im}\Bigl(\Hom_{\Tri(E)}\bigl(T_{\kP^{v_1}}, T_{\kP^{v_2}(y_1)}\bigr)
\xrightarrow{\mathsf{For}}
\Hom_{\PP^1}\bigl(\widetilde\kP, \widetilde\kP(1)\bigr)\Bigr)$. Then the following diagram is commutative:
$$
\xymatrix
{
\widetilde\Pi^{v_1, v_2}_{y_1}  \ar[rr] \ar[d]_{t_\sigma} & & \Hom_{\PP^1}\bigl(\widetilde\kP, \widetilde\kP(1)\bigr) \ar[d]^{t_\sigma} \\
\Pi^{v_1, v_2}_{y_1}  \ar[rr] & & \Hom_{\PP^1}\bigl(\widetilde\kP, \widetilde\kP(\tilde{y}_1)\bigr) \\
}
$$
\end{corollary}

\noindent
Collecting everything together, we get the following algorithm for computing associative
 $r$--matrices coming from a singular Weierestra\ss{} cubic curve $E$.

\begin{algorithm}
Let $\omega \in H^0(\Omega^{1, R}_E) \subset H^0(\Omega^{1, M}_{\PP^1})$ be the global regular
differential one form on $\PP^1$ equal to
$\dfrac{dz}{z}$ if $E$ is nodal and $dz$ if $E$ is cuspidal. The linear morphism
$\tilde{r}^\xi_\omega(v_1, v_2; y_1, y_2)$ can be computed in the following way.
\begin{itemize}
\item Compute the vector space $\widetilde\Pi^{v_1, v_2}_{y_1} \subseteq
\Hom_{\PP^1}\bigl(\widetilde\kP, \widetilde\kP(1)\bigr)$.
\item Consider the morphism
$\overline{\res}_{y_1}: \Hom_{\PP^1}\bigl(\widetilde\kP, \widetilde\kP(1)\bigr) \lar
\Mat_{n \times n}(\CC)$ given by
$$
\overline{\res}_{y_1}(F) =
\left\{
\begin{array}{cl}
\frac{\displaystyle 1}{\displaystyle y_1}F(1, y_1) & \mbox{\textrm{if}} \, \, $E$  \, \, \mbox{\textrm{is nodal}} \\
F(1, y_1) & \mbox{\textrm{if}} \, \, $E$  \, \, \mbox{\textrm{is cuspidal.}}
\end{array}
\right.
$$
\item If $E$ is either  nodal or cuspidal, we set $\ev_{y_2}:
\Hom_{\PP^1}\bigl(\widetilde\kP, \widetilde\kP(1)\bigr) \lar
\Mat_{n \times n}(\CC)$ to be given by the formula
$$
\overline{\ev}_{y_2}(F) = \frac{1}{y_2 - y_1}F(1, y_2).
$$
\item The linear morphism $\tilde{r}^\xi(v_1, v_2; y_1, y_2):
\Mat_{n \times n}(\CC) \lar \Mat_{n \times n}(\CC)$ can be computed  as the composition
$$
\Mat_{n \times n}(\CC) \xrightarrow{\overline{\res}_{y_1}^{-1}} \widetilde\Pi^{v_1, v_2}_{y_1}
\xrightarrow{\overline{\ev}_{y_2}} \Mat_{n \times n}(\CC).
$$
\end{itemize}
\end{algorithm}

\subsection{A trigonometric solution  obtained from a nodal cubic curve}

Let $E$ be a nodal Weierstra\ss{} cubic curve. In this subsection we calculate
the associative $r$-matrices  corresponding to the moduli spaces  of rank two
(semi-)stable vector bundles on a nodal Weierstra\ss{} curve. We use the
notation from Subsection \ref{SS:node}.

\medskip
We start with the case of the moduli space of stable vector
bundles of rank $2$ and degree $1$, $M= M_E^{(2, 1)} = \CC^*$. It is convenient to
use the local homeomorphism $\sigma: \CC^* \to \CC^*$
given by $\sigma(z) = z^2$, because, according to Example \ref{E:rank2} and
Remark \ref{R:noderank2}, the family of stable vector bundles $(1 \times
\sigma)^*\kP(2,1)$ is then given by the triple
$\bigl(\kO_{\PP^1} \oplus \kO_{\PP^1}(1), \CC_s^2, \mm\bigr)$, where
$$
\mm(0) =
\left(
\begin{array}{cc}
0 & \lambda \\
\lambda & 0
\end{array}
\right), \, \, \lambda \in \CC^*
\quad \mbox{\textrm{and}} \quad
\mm(\infty) =
\left(
\begin{array}{cc}
1 & 0  \\
0 & 1
\end{array}
\right).
$$
Our goal is to compute the map
$$
r^{\lambda_1, \lambda_2}_{y_1, y_2}: \Mat_{2 \times 2}(\CC) \lar \Mat_{2 \times 2}(\CC), \quad
\left(
\begin{array}{cc}
a & b \\
c & d
\end{array}
\right)
\mapsto
\left(
\begin{array}{cc}
\varphi & \psi \\
\eta & \xi
\end{array}
\right).
$$

\noindent
\underline{Step 1}.
In order to calculate the entries  $\varphi,  \psi, \eta, \xi$
we first need  to describe the subspace
$
\Pi^{\lambda_1, \lambda_2}_{y_1} \subset
\Hom_{\PP^1}\bigl(\kO_{\PP^1} \oplus \kO_{\PP^1}(1), \kO_{\PP^1}(1)
\oplus \kO_{\PP^1}(2)\bigr).
$
Following the recipe of  Subsection \ref{SS:node}, we take the section
$p_\zeta = z_1 - z_0$ to evaluate a morphism
$$
F =
\left(
\begin{array}{cc}
a' z_0 + a'' z_1 & t \\
b'z_0^2 + b'' z_0 z_1 + b''' z_1^2 & d'z_0 + d'' z_1
\end{array}
\right) \in \Hom_{\PP^1}\bigl(\kO_{\PP^1} \oplus \kO_{\PP^1}(1), \kO_{\PP^1}(1)
\oplus \kO_{\PP^1}(2)\bigr).
$$
This gives the following  evaluation rule:
$$
F(0) =
\left(
\begin{array}{cc}
-a'  & t \\
b' &  -d'
\end{array}
\right),
\qquad
F(\infty) =
\left(
\begin{array}{cc}
a''  & t \\
b''' &  d''
\end{array}
\right).
$$
From the definition of the category of triples we see  that
$F$ belongs to $\Pi^{\lambda_1, \lambda_2}_{y_1}$ if and only if there exists a matrix $\varphi \in
\Mat_{2 \times 2}(\CC)$ making the following diagram commutative:
{\scriptsize
$$
{
\xymatrix
{\CC_0^2 \oplus \CC_\infty^2
\ar[rrrr]^{\begin{pmatrix} F(0) & 0  \\ 0 & F(\infty) \end{pmatrix}}
& & & & \CC_0^2 \oplus \CC_\infty^2 \\
& & & & \\
& & & & \\
& & & & \\
\CC_0^2 \oplus \CC_\infty^2
\ar[uuuu]^{
\begin{pmatrix}
\left(\begin{smallmatrix}
0  & \lambda_1  \\
\lambda_1  & 0
\end{smallmatrix}\right)
& 0 \\  0  &
\left(\begin{smallmatrix}
1  & 0   \\
0  & 1
\end{smallmatrix}\right)
\end{pmatrix}}
\ar[rrrr]^{
\begin{pmatrix}
\varphi & 0 \\
0  & \varphi
\end{pmatrix}}
& & & & \CC_0^2 \oplus \CC_\infty^2  \ar[uuuu]_{
\begin{pmatrix}
\left(\begin{smallmatrix}
0  & \lambda_2 y_1  \\
\lambda_2 y_1  & 0
\end{smallmatrix}\right)
& 0 \\ 0  &
\left(\begin{smallmatrix}
1  & 0   \\ 0  & 1
\end{smallmatrix}\right)
\end{pmatrix}}
 \\
}
}
$$
}
This is equivalent to the equations:
$$
F(\infty)  =   \varphi  \quad \mbox{\textrm{and}} \quad
F(0)
\left(
\begin{array}{cc}
0 & \lambda_1 \\
\lambda_1 & 0
\end{array}
\right)
 =
\left(
\begin{array}{cc}
0 & \lambda_2 y_1 \\
\lambda_2 y_1 & 0
\end{array}
\right)
\varphi.
$$
Taking $a'', b'', b''', d''$
as free variables and solving the above system we get
$$
\left\{
\begin{array}{cccl}
a' & =  &  - \lambda y_1 d'' &\\
d' & =  &  - \lambda y_1 a'' &\\
t & = &  \lambda y_1 b'''&\\
b' & = & (\lambda y_1)^2 b''',&\quad\text{ where }
\lambda = \dfrac{\lambda_2}{\lambda_1}.
\end{array}
\right.
$$

\medskip
\noindent
\underline{Step 2}.
Next, the equation
$
\overline{\res}_{y_1}(F) =
\left(
\begin{array}{cc}
a & b \\
c & d
\end{array}
\right)
$
reads as
$$
\left(
\begin{array}{cc}
a' + a'' y_1 & t \\
b' + b'' y_1 + b''' y_1^2 & d' + d'' y_1
\end{array}
\right)
=
y_1
\left(
\begin{array}{cc}
a & b \\
c & d
\end{array}
\right).
$$
From this we obtain
$$
\left\{
\begin{array}{ccl}
a' & = & - \frac{\displaystyle \lambda y_1}{\displaystyle 1 -
           \lambda^2}(d + \lambda a)\\
a''& = &  \frac{\displaystyle 1}{\displaystyle 1 - \lambda^2}(a + \lambda d)\\
b' & = &  \lambda y_1^2 b \\
b''& = &  c - \frac{\displaystyle \lambda^2 + 1}{\displaystyle \lambda} y_1 b\\
b''' & = & \frac{\displaystyle 1}{\displaystyle \lambda} b  \\
d'& = &  - \frac{\displaystyle \lambda y_1}{\displaystyle 1 -
            \lambda^2}(a + \lambda d)\\
d''& = &  \frac{\displaystyle 1}{\displaystyle 1 - \lambda^2}(d + \lambda a)\\
t & = & y_1 b.\\
\end{array}
\right.
$$

\noindent
\underline{Step 3}.
By the formula for the evaluation map we get:
$$
\overline{\ev}_{y_2}(F) =
\frac{1}{y_2 - y_1}
\left(
\begin{array}{cc}
a' + a'' y_2 & t \\
b' + b'' y_2 + b'''y_2^2 & d' + d'' y_2
\end{array}
\right) =
\frac{1}{y}\left(
\begin{array}{cc}
\varphi & \psi \\
\eta & \xi
\end{array}
\right),
$$
where we denote
$$
\left\{
\begin{array}{ccl}
\varphi & = &
   \frac{\displaystyle y_2 - \lambda^2 y_1}{\displaystyle 1 - \lambda^2}a +
   \frac{\displaystyle \lambda(y_2 - y_1)}{\displaystyle 1 - \lambda^2}d \\
\psi & = & y_1 b \\
\xi & = &
   \frac{\displaystyle \lambda(y_2 - y_1)}{\displaystyle 1 - \lambda^2}a +
   \frac{\displaystyle y_2 - \lambda^2 y_1}{\displaystyle 1 - \lambda^2}d \\
\eta & = &
\frac{\displaystyle (y_2 - y_1)(y_2 - \lambda^2 y_1)}{\displaystyle \lambda}b
     + y_2 c.
\end{array}
\right.
$$
In order to calculate  the corresponding solution
of the associative Yang--Baxter equation
 we use the inverse of the canonical isomorphism
$$
\Mat_{2 \times 2}(\CC) \otimes \Mat_{2 \times 2}(\CC) \lar
\Lin\bigl(\Mat_{2 \times 2}(\CC), \Mat_{2 \times 2}(\CC)\bigr)
$$
given by $X \otimes Y \mapsto \tr(X \circ -) Y$.
It is easy to see that under this inverse
$$
\Lin\bigl(\Mat_{2 \times 2}(\CC), \Mat_{2 \times 2}(\CC)\bigr) \lar
\Mat_{2 \times 2}(\CC) \otimes \Mat_{2 \times 2}(\CC)
$$
a  linear function $e_{ij} \mapsto \alpha_{ij}^{kl} e_{kl}, \alpha_{ij}^{kl}
\in \CC$ corresponds to the tensor
$\alpha_{ij}^{kl}e_{ji} \otimes e_{kl}$. Having this
rule in mind we obtain the
desired associative $r$--matrix:
\begin{align*}
r(\lambda; y_1, y_2) &=
\dfrac{y_2 - \lambda^2 y_1}{(y_2 - y_1)(1 - \lambda^2)}(e_{11}\otimes e_{11} +
e_{22} \otimes e_{22}) + \!
\dfrac{\lambda}{1 - \lambda^2}(e_{11}\otimes e_{22}  +
e_{22} \otimes e_{11}) \\
&+
\dfrac{y_1}{y_2 - y_1} e_{21} \otimes e_{12} +
\dfrac{y_2}{y_2 - y_1} e_{12} \otimes e_{21} +
\dfrac{y_2 - \lambda^2 y_1}{\lambda} e_{21} \otimes e_{21}.
\end{align*}
The gauge transformation
$\varphi(z) = \varphi(\mu; z): (\CC^2, 0) \lar \Aut\bigl(\Mat_n(\CC)\bigr)$
(see  Definition \ref{D:gaugeAYBE})  given by
$$
\left(
\begin{array}{cc}
a & b \\
c & d
\end{array}
\right)
\mapsto
\left(
\begin{array}{cc}
\scriptstyle{\sqrt{z}} & 0  \\
0 & 1
\end{array}
\right)
\left(
\begin{array}{cc}
a & b \\
c & d
\end{array}
\right)
\left(
\begin{array}{cc}
\frac{1}{\sqrt{z}} & 0  \\
0 & 1
\end{array}
\right)
$$
yields    the transformation
$$
\left\{
\begin{array}{l}
e_{ii} \otimes e_{jj} \mapsto e_{ii} \otimes e_{jj}, \quad i,j \in \{1, 2\} \\
e_{21} \otimes e_{12} \mapsto  \sqrt{\frac{y_2}{y_1}} e_{21} \otimes e_{12} \\
e_{12} \otimes e_{21} \mapsto  \sqrt{\frac{y_1}{y_2}}  e_{12} \otimes e_{21} \\
e_{21} \otimes e_{21} \mapsto \frac{1}{\sqrt{y_1 y_2}} e_{21} \otimes e_{21}.
\end{array}
\right.
$$
Thus, we end up with the solution
\begin{align*}
  r(\lambda, y) &=
\dfrac{y-\lambda^2}{(y-1)(1-\lambda^2)}(e_{11} \otimes e_{11} +
e_{22} \otimes e_{22}) +
\dfrac{\lambda}{1 - \lambda^2}(e_{11} \otimes e_{22} + e_{22} \otimes e_{11}) +\\
&+
\dfrac{\sqrt{y}}{y-1}(e_{12} \otimes e_{21} + e_{21} \otimes e_{12}) +
\left(\dfrac{\sqrt{y}}{\lambda} - \frac{\lambda}{\sqrt{y}}\right)
e_{21} \otimes e_{21},
\end{align*}
where $y = \dfrac{y_2}{y_1}$. Using the notation $\mathbbm{1} = e_{11} +
e_{22}$, this can be rewritten as
\begin{align*}
  r(\lambda, y) &= \dfrac{\mathbbm{1} \otimes \mathbbm{1}}{1 - \lambda^2} +
\dfrac{1}{y-1}(e_{11} \otimes e_{11} +
e_{22} \otimes e_{22}) -
\dfrac{1}{\lambda + 1}(e_{11} \otimes e_{22} +
e_{22} \otimes e_{11})
\\
&+
\dfrac{\sqrt{y}}{y-1}(e_{12} \otimes e_{21} + e_{21} \otimes e_{12}) +
\left(\dfrac{\sqrt{y}}{\lambda} - \frac{\lambda}{\sqrt{y}}\right)
e_{21} \otimes e_{21}.
\end{align*}
This is a solution of the associative Yang--Baxter equation of type
(\ref{E:AnsatzAYBE}), and by Theorem \ref{T:QYBEandCYBE} this  tensor also
satisfies the quantum  Yang--Baxter equation.

In order to rewrite $r(\lambda; y)$ in the  additive form, we make
the change of variables $y = \exp(2iz), \lambda = \exp(iv)$. Making a gauge
transformation we can multiply the tensor
$e_{21}\otimes e_{12}$ with an arbitrary scalar without changing the
coefficients of the other tensors. Therefore, we obtain
\begin{align*}
  2r_{\mathrm{trg}}(v, z) &=
\dfrac{\sin(z+v)}{\sin(z)\sin(v)}(e_{11}\otimes e_{11}+e_{22}\otimes e_{22}) +
\dfrac{1}{\sin(v)}(e_{11} \otimes e_{22} + e_{22} \otimes e_{11}) +
\\
&+ \dfrac{1}{\sin(z)}(e_{12}\otimes e_{21} + e_{21}\otimes e_{12}) +
\sin(z+v) e_{21} \otimes e_{21}.
\end{align*}
Up to a scalar, the corresponding solution
$\bar{r}(z)  := \lim_{v \to 0}(\mathrm{pr} \otimes \mathrm{pr})
r(v; z)$
of the classical Yang--Baxter equation is the  trigonometric solution of
Cherednik:
$$
\bar{r}_{\textrm{trg}}(z) = \frac{1}{2} \cot(z)h\otimes h +
\frac{1}{\sin(z)}(e_{12}\otimes e_{21} + e_{21} \otimes e_{12}) +
\sin(z) e_{21} \otimes e_{21}.
$$

\medskip
\subsection{Trigonometric solutions coming from semi-stable vector bundles}\label{SS:solfromSS}
Our next goal is to construct a solution $r(v; y)$ of the associative
Yang--Baxter equation (\ref{E:AYBE2}) having a higher-order pole with respect
to $v$.
The triple $\bigl(\kO_{\PP^1}^2, \CC_s^2, \mm\bigr)$ with
$$
\mm(0) =
\left(
\begin{array}{cc}
\lambda & \lambda \\
0       & \lambda
\end{array}
\right), \lambda \in \CC^*
\quad
\mbox{\textrm{and}} \quad
\mm(\infty) =
\left(
\begin{array}{cc}
1 & 1 \\
0       & 1
\end{array}
\right)
$$
describes a universal family of semi-stable indecomposable vector bundles of
rank two and degree one, having locally free Jordan-H\"older factors.

\medskip
\noindent
\underline{Step 1}. First we compute the subspace
$
\Pi^{\lambda_1, \lambda_2}_{y_1} \subset
\Hom_{\PP^1}\bigl(\kO_{\PP^1} \oplus\kO_{\PP^1},
\kO_{\PP^1}(1) \oplus \kO_{\PP^1}(1)\bigr).
$
Recall that for a morphism
$$
F =
\left(
\begin{array}{cc}
a' z_0 + a'' z_1 & b' z_0 + b'' z_1 \\
c' z_0 + c'' z_1 & d' z_0 + d'' z_1
\end{array}
\right) \in \Hom_{\PP^1}\bigl(\kO_{\PP^1} \oplus\kO_{\PP^1},
\kO_{\PP^1}(1) \oplus \kO_{\PP^1}(1)\bigr)
$$
we take the evaluation rule induced by the section $p_\zeta = z_1 - z_0$:
$$
F(0) = - F' :=
- \left(
\begin{array}{cc}
a' & b' \\
c' & d'
\end{array}
\right)
\quad
F(\infty) = F'' :=
 \left(
\begin{array}{cc}
a'' & b'' \\
c'' & d''
\end{array}
\right).
$$
Thus, $F$ belongs to $\Pi^{\lambda_1, \lambda_2}_{y_1}$  if and only if we have
$$
F(0)
\left(
\begin{array}{cc}
\lambda_1 & \lambda_1 \\
0         & \lambda_1
\end{array}
\right) =
\left(
\begin{array}{cc}
\lambda_2 y_1 & \lambda_2 y_1 \\
0         & \lambda_2 y_1
\end{array}
\right)
F(\infty).
$$
This implies that
$$
F' =  - \lambda y_1 F'' +
\lambda y_1
\left(
\begin{array}{cc}
- c'' & a'' + c'' - d'' \\
0 & c''
\end{array}
\right).
$$

\noindent
\underline{Step 2}. The equation $\overline{\res}_{y_1}(F) = y_1
\left(
\begin{array}{cc}
a & b \\
c & d
\end{array}
\right)
$ reads $F' + y_1 F'' = y_1
\left(
\begin{array}{cc}
a & b \\
c & d
\end{array}
\right)$.

\noindent
Solving this equation we obtain
$$
\left\{
\begin{array}{lll}
a'' & = & \dfrac{1}{1 - \lambda} a +
\dfrac{\lambda}{(1 - \lambda)^2} c
 \\
b'' & =  & \dfrac{\lambda}{1 - \lambda} a +
\dfrac{1}{1 - \lambda} b
- \dfrac{\lambda(\lambda + 1)}{(1 - \lambda)^3} c
+ \dfrac{\lambda}{(1 - \lambda)^2} d \\
c'' & =  & \dfrac{1}{1 - \lambda} c \\
d'' & = & -\dfrac{\lambda}{(1 - \lambda)^2} c +
\dfrac{1}{1 - \lambda} d.
\end{array}
\right.
$$

\noindent
\underline{Step 3}. From the formula
$\overline{\ev}_{y_2}(F) = \dfrac{1}{y_2 - y_1}(F' + y_2 F'')$ we obtain:
$$
\tilde{r}^{\lambda_1, \lambda_2}_{y_1, y_2}
\left(
\begin{array}{cc}
a & b \\
c & d
\end{array}
\right)
=
\left(
\begin{array}{cc}
\varphi & \psi  \\
\eta  & \zeta
\end{array}
\right),
$$
where
$$
\left\{
\begin{array}{l}
\varphi =
\dfrac{y - \lambda}{(y - 1)(1-\lambda)} a +
\dfrac{\lambda}{(1-\lambda)^2} c \\
\eta = \dfrac{y - \lambda}{(y - 1)(1-\lambda)} c \\
\xi = - \dfrac{\lambda}{(1-\lambda)^2} c +
\dfrac{y - \lambda}{(y - 1)(1-\lambda)} d \\
\psi =
- \dfrac{\lambda}{(1-\lambda)^2} a +
\dfrac{y - \lambda}{(y - 1)(1-\lambda)} b
- \dfrac{\lambda(1 + \lambda)}{(1-\lambda)^3} c +
\dfrac{\lambda}{(1-\lambda)^2} d
\end{array}
\right.
$$
and  $y = \dfrac{y_2}{y_1}$,
$\lambda = \dfrac{\lambda_2}{\lambda_1}$.
Hence, we  obtain the associative $r$--matrix
$$
r(\lambda; y) =
\dfrac{y - \lambda}{(y - 1)(1-\lambda)}
\bigl(
e_{11} \otimes e_{11} + e_{22} \otimes e_{22} + e_{21}
\otimes e_{12} + e_{12} \otimes
e_{21}
\bigr) +
$$
$$
+
\frac{\displaystyle \lambda}{(1-\lambda)^2} \bigl(e_{12} \otimes h -
h \otimes e_{12}\bigr) -
\frac{\displaystyle \lambda(1 + \lambda)}{(1-\lambda)^3} e_{12} \otimes e_{12}.
$$
Denoting $y = \exp(2iz), \lambda = \exp(-2iv)$
and making  a gauge transformation
$$e_{11} \mapsto e_{11}, e_{22} \mapsto e_{22},
 e_{12} \mapsto 2 e_{12} \text{ and } e_{21} \mapsto \frac{1}{2}e_{21}$$
we finally end up with an associative $r$--matrix
\begin{align*}
  r(v;z) &=
\dfrac{\sin(z +v)}{2\sin(z)\sin(v)}
\bigl(
e_{11} \otimes e_{11} + e_{22} \otimes e_{22} + e_{21} \otimes e_{12} +
e_{12} \otimes e_{21}
\bigr)
\\
&+
\dfrac{1}{2 \sin^2(v)} \bigl(e_{12} \otimes h - h \otimes e_{12}\bigr) -
\dfrac{\cos(v)}{\sin^3(v)} e_{12} \otimes e_{12}.
\end{align*}
\begin{remark}
Since  $\lim\limits_{v \to 0} \bigl(\pr\otimes \pr (r(v; z))\bigr)$ does
not exist, the family of indecomposable semi-stable vector bundles of
rank two and degree zero on a nodal Weierstra\ss{} curve $E$,  whose
Jordan-H\"older factors are locally free, does not give a solution of the
classical Yang--Baxter equation.
\end{remark}

\subsection{A rational solution obtained from a cuspidal cubic curve}

In this subsection we shall calculate the rational solution
of the classical Yang--Baxter equation, obtained from a universal family
of stable vector bundles  of rank $2$ and degree $1$ on a cuspidal cubic
curve.
In terms of Subsection \ref{SS:cusp} the universal family is described by the
family of triples
$\bigl(\kO_{\PP^1} \oplus \kO_{\PP^1}(1), \CC_s^2, \mm\bigr)$, where
$$
\mm = \mm_0 + \varepsilon \mm_\varepsilon =
\left(
\begin{array}{cc}
1 & 0 \\
0 & 1
\end{array}
\right)
+
\varepsilon
\left(
\begin{array}{cc}
\lambda & 1  \\
0 & \lambda
\end{array}
\right), \, \, \lambda \in \CC.
$$
As in the previous subsection, let
$$
\left(
\begin{array}{cc}
a & b \\
c & d
\end{array}
\right) \in \Mat_{2\times 2}(\CC) \quad \mbox{\textrm{and}} \quad
\left(
\begin{array}{cc}
\varphi & \psi \\
\eta & \xi
\end{array}
\right) =
\tilde{r}^{\lambda_1, \lambda_2}_{y_1, y_2}
\left(
\begin{array}{cc}
a & b \\
c & d
\end{array}
\right).
$$

\noindent
\underline{Step 1}.
Again, we start by  calculating the linear subspace
$
\Pi^{\lambda_1, \lambda_2}_{y_1} \subset
\Hom_{\PP^1}\bigl(\kO_{\PP^1} \oplus \kO_{\PP^1}(1), \kO_{\PP^1}(1)
\oplus \kO_{\PP^1}(2)\bigr).
$
Recall that, in the case of a cuspidal curve,  a morphism
$$
F=
\begin{pmatrix}
a' z_0 + a'' z_1 & t \\
b'z_0^2 + b'' z_0 z_1 + b''' z_1^2 & d'z_0 + d'' z_1
\end{pmatrix}
\in
\Hom_{\PP^1}\bigl(\kO_{\PP^1} \oplus \kO_{\PP^1}(1), \kO_{\PP^1}(1)
\oplus \kO_{\PP^1}(2)\bigr)
$$
is evaluated
on the analytic subspace  $\widetilde{Z}$ using the section $p_\zeta = z_1$. This gives
the following evaluation rule:
$$
F
\mapsto
\left(
\begin{array}{cc}
a'' + a' \varepsilon  & t \\
b''' + b'' \varepsilon  &  d'' + d'\varepsilon
\end{array}
\right).
$$
From the definition of the category $\Tri(E)$ we see  that
$F$ belongs to $\Pi^{\lambda_1, \lambda_2}_{y_1}$
if and only if there exists a matrix $f \in \Mat_{2\times 2}(\CC)$ making
the following diagram commutative
{
$$
\xymatrix
{\fR^2
\ar[rrrr]^{\begin{pmatrix} a'' & t \\ b'''  & d''\end{pmatrix} +
{\displaystyle \varepsilon}
\begin{pmatrix} a' & 0  \\ b''  & d'\end{pmatrix}
}
& & & & \fR^2 \\
 & & & & \\
%& & & & \\
& & & & \\
\fR^2
\ar[uuu]^{\begin{pmatrix} 1 & 0 \\ 0  & 1\end{pmatrix} +
{\displaystyle \varepsilon}
\begin{pmatrix} \lambda_1 &  1  \\ 0  &  \lambda_1\end{pmatrix}
}
\ar[rrrr]^{f}
& & & & \fR^2,
\ar[uuu]_{\begin{pmatrix} 1 & 0 \\ 0  & 1\end{pmatrix} +
{\displaystyle \varepsilon}
\begin{pmatrix} \lambda_2 - y_1 &  1  \\ 0  &
\lambda_2 -y_1\end{pmatrix}
}
}
$$
}
where $\fR = \CC[\varepsilon]/\varepsilon^2$.
 This leads to the equality
$$
\begin{pmatrix}
a' & 0 \\
b'' & d'
\end{pmatrix}
+
\begin{pmatrix}
a'' & t \\
b''' & d''
\end{pmatrix}
\begin{pmatrix}
\lambda_1 & 1 \\
0  & \lambda_1
\end{pmatrix} =
\begin{pmatrix}
\lambda_2 - y_1 & 1 \\
0  & \lambda_2 - y_1
\end{pmatrix}
\begin{pmatrix}
a'' & t \\
b''' & d''
\end{pmatrix}.
$$
Taking $a'', b', b'''$ and $t$ as free variables we obtain
$$
\left\{
\begin{array}{l}
a' = (\lambda - y_1) a'' + b''' \\
b'' = (\lambda - y_1) b''' \\
d' = (\lambda - y_1) a'' - b''' - (\lambda - y_1)^2 t \\
d'' = a'' - (\lambda - y_1) t.
\end{array}
\right.
$$

\noindent
\underline{Step 2}. By the formula for the residue map $\overline{\res}_{y_1}$ we have:
$$
\overline{\res}_{y_1}(F) =
\begin{pmatrix}
a' + a'' y_1   & t \\
b' + b''y_1 +  b'''y_1^2 & d' + d'' y_1
\end{pmatrix} =
\begin{pmatrix}
a & b \\
c & d
\end{pmatrix}
$$
from which we get:
$$
\left\{
\begin{array}{ccl}
t & = & b \\[2mm]
a'' & = & \dfrac{1}{2 \lambda} a +
\dfrac{\lambda - y_1}{2} b +
\dfrac{1}{2 \lambda} d \\[2mm]
b' & = &  - \dfrac{\lambda y_1}{2} a +
\dfrac{\lambda^2y_1(\lambda - y_1)}{2} b + c
+ \dfrac{\lambda y_1}{2} d \\[2mm]
b''' & = & \dfrac{1}{2} a -
\dfrac{\lambda(\lambda - y_1)}{2} b
- \dfrac{1}{2} d
\end{array}
\right.
$$

\noindent
\underline{Step 3}.
Since the formula for the map $\overline{\ev}_{y_2}$ is given by:
$$
\overline{\ev}_{y_2}(F) =
\frac{1}{y_2 - y_1}
\left(
\begin{array}{cc}
a' + a'' y_2   & t \\
b' + b''y_2 +  b'''y_2^2 & d' + d'' y_2
\end{array}
\right) =
\frac{1}{y}
\left(
\begin{array}{cc}
\varphi & \psi \\
\eta & \xi
\end{array}
\right)
$$
we obtain:
$$
\left\{
\begin{array}{ll}
\varphi &= (1 + \dfrac{y_2 - y_1}{2\lambda}) a +
\dfrac{(\lambda  - y_1)(y_2 - y_1)}{2} b +
\dfrac{y_2 - y_1}{2\lambda} d \\[2mm]
\psi &= t \\[2mm]
\eta &=
  \dfrac{(y_2 - y_1)(y_2 - y_1 + \lambda)}{2} a -
  \dfrac{\lambda(\lambda - y_1)(\lambda  +  y_2)(y_2 - y_1)}{2} b +\\
  &\phantom{=}+ c -\dfrac{(y_2 - y_1)(\lambda + y_2)}{2}d \\[2mm]
\xi &= \dfrac{y_2 - y_1}{2\lambda} a -
\dfrac{(y_2 - y_1)(y_1 - \lambda)}{2} b +
\left(1 + \dfrac{y_2 - y_1}{2\lambda}\right) d.
\end{array}
\right.
$$
From this we get the following associative $r$-matrix:
\begin{equation}\label{E:ratsol}
  r(\lambda, y_1, y_2)  =
\dfrac{1}{2 \lambda}
\mathbbm{1}\otimes \mathbbm{1} +
\dfrac{1}{y_2 - y_1}\Bigl(e_{11}
\otimes e_{11} + e_{22}\otimes e_{22} +
e_{12}\otimes e_{21} + e_{21}\otimes e_{12}\Bigr) +
\end{equation}
$$
+ \dfrac{\lambda - y_1}{2} e_{21} \otimes h +
\dfrac{\lambda + y_2}{2}  h \otimes e_{21} -
\dfrac{\lambda(\lambda - y_1)(\lambda + y_2)}{2}
e_{21} \otimes e_{21}.
$$
Projecting this matrix to
$\mathfrak{sl}_2(\CC) \otimes \mathfrak{sl}_2(\CC)$ we obtain the rational
solution of the classical Yang--Baxter equation
\begin{equation}\label{E:Stolinmatr}
\bar{r}(y_1, y_2) = \dfrac{1}{y_2 - y_1}\left(
\dfrac{1}{2} h \otimes h + e_{12}
\otimes e_{21} + e_{21} \otimes e_{12}\right)
+ \dfrac{y_2}{2} h \otimes e_{21} -
\dfrac{y_1}{2} e_{21} \otimes h,
\end{equation}
found for the first time by Stolin in \cite{Stolin}. It is easy to
check that $\bar{r}(y_1, y_2)$ does not have infinitesimal symmetries, hence by
Theorem \ref{T:QYBEandCYBE} the tensor $r(\lambda, y_1, y_2)$ satisfies
 the Quantum Yang--Baxter equation.  This solution was recently found by
Khoroshkin, Stolin and Tolstoy \cite{KST}.

\begin{remark}
By a result of Belavin and Drinfeld \cite{BelavinDrinfeld2} it is known that
the $r$--matrix (\ref{E:Stolinmatr}) is equivalent to a solution depending on the difference of spectral parameters. However, we were not able to find the corresponding
gauge transformation in the literature and the
following form of Stolin's solution seems to be new.
Consider the  gauge transformation
$
\phi: (\CC, 0) \lar \Aut\bigl(\mathfrak{sl}_2(\CC)\bigr)$
given by the formula
$$
\phi(y)  h = h - 2 y^2 e_{21}, \quad  \phi(y)  e_{12} = -\frac{y^2}{2}  h +  \frac{1}{4} e_{12}
- \frac{y^4}{4}  e_{21} \quad   \mbox{\textrm and} \quad
\phi(y)  e_{21} = 4 e_{21}.
$$
Then we have: $\bigl(\phi(y_1) \otimes \phi(y_2)\bigr) r(y_1, y_2) = s(y_2 - y_1)$, where
\begin{equation}\label{E:ratStol}
s(y)  =
\frac{1}{y}\Bigl(\frac{1}{2} h \otimes h + e_{12} \otimes e_{21} +
e_{21} \otimes e_{12}\Bigr) +
y(e_{21} \otimes h + h \otimes e_{21}) - y^3 e_{21} \otimes e_{21}.
\end{equation}
\end{remark}

\medskip
\begin{remark}
For any $t \in \CC^*$ consider the constant gauge transformation
$$
e_{11} \mapsto e_{11}, \quad e_{22} \mapsto e_{22}, \quad e_{21} \mapsto t \, e_{21}, \quad
e_{12} \mapsto \dfrac{1}{t} \, e_{12}.
$$
Then the associative $r$--matrix (\ref{E:ratsol}) transforms into the solution
\begin{align*}
  r_t(\lambda, y_1, y_2)  & =
\dfrac{1}{2 \lambda}
\mathbbm{1}\otimes \mathbbm{1} +
\dfrac{1}{y_2 - y_1}\Bigl(e_{11}
\otimes e_{11} + e_{22}\otimes e_{22} +
e_{12}\otimes e_{21} + e_{21}\otimes e_{12}\Bigr) + \\
& + t\Bigl(\dfrac{\lambda - y_1}{2} e_{21} \otimes h +
\dfrac{\lambda + y_2}{2}  h \otimes e_{21}\Bigr)  +
\dfrac{t^2\lambda(\lambda - y_1)(\lambda + y_2)}{2}
e_{21} \otimes e_{21}.
\end{align*}
Taking the limit $t \to  0$, we get the following solution of the associative Yang--Baxter equation (\ref{E:AYBE3}):
\begin{equation}
r(\lambda, y) = \dfrac{1}{2 \lambda}
\mathbbm{1}\otimes \mathbbm{1} +
\dfrac{1}{y}\bigl(e_{11}
\otimes e_{11} + e_{22}\otimes e_{22} +
e_{12}\otimes e_{21} + e_{21}\otimes e_{12}\bigr).
\end{equation}
Note that the corresponding solution of the classical Yang--Baxter equation (\ref{E:CYBE1})
is the rational solution of Yang.
\end{remark}

\section{Summary}\label{S:summary}
Let us summarise the main analytical results obtained is this article.
We have shown that for any pair of coprime integers $0 < d < n$ and an
irreducible reduced projective curve $E$ with trivial dualizing sheaf
one can \emph{canonically} attach the germ of a tensor-valued function
$$
r^{(n,d)}_E(v; y_1, y_2): (\CC^3, 0) \lar \Mat_{n \times n}(\CC)  \otimes \Mat_{n \times n}(\CC)
$$
satisfying the associative Yang--Baxter equation
$$
r(u; y_1, y_2)^{12} r(u+v; y_2, y_3)^{23}  =
r(u+v;y_1, y_3)^{13} r(-v;y_1, y_2)^{12}
$$
$$
+ r(v; y_2, y_3)^{23} r(u; y_1, y_3)^{13}.
$$
By Proposition \ref{P:Kirillov}, the  tensor $r^{(n,d)}_E(v; y_1, y_2)$
defines a family of commuting first order differential operators.  Moreover,
under certain conditions (which are always fulfilled at least for elliptic
curves and for nodal cubic curves) it also satisfies the quantum Yang--Baxter
equation with respect to the spectral variables
$y_1 $ and $y_2$ and a fixed value of $v \ne 0$.  For $(n, d) = (2,1)$ these tensors $r^{(2,1)}_E $ have the following explicit form

\medskip
\noindent
$\bullet$  For  an elliptic curve $E = E_\tau = \CC/\langle 1,  \tau\rangle$ we get
 \begin{align*}
 r^{(2, 1)}_{\mathrm{ell}}(v; y) =
\frac{\theta'_1(0|\tau)}{\theta_1(y|\tau)}
&\left[
\frac{\theta_1(y + v|\tau)}{\theta_1(v|\tau)}
\mathbbm{1} \otimes \mathbbm{1} +
\frac{\theta_2(y + v|\tau)}{\theta_2(v|\tau)} h  \otimes h+\right.
\\
&\left.+
\frac{\theta_3(y + v|\tau)}{\theta_3(v|\tau)} \sigma  \otimes \sigma
+
\frac{\theta_4(y + v|\tau)}{\theta_4(v|\tau)} \gamma
\otimes \gamma\right],
\end{align*}
where $\mathbbm{1} =  e_{11} + e_{22}$, $h = e_{11} - e_{22}, \sigma =
i(e_{21} - e_{12})$ and  $\gamma = e_{21} + e_{12}$. This solution
is a quantization of the elliptic solution of the classical Yang--Baxter equation
$$
\bar{r}^{(2, 1)}_{\mathrm{ell}}(y) =
\frac{1}{2}\left(\frac{\cn(y)}{\sn(y)} h\otimes h +
\frac{1}{\sn(y)} \gamma  \otimes \gamma  +
 \frac{\dn(y)}{\sn(y)} \sigma \otimes \sigma\right).
$$
studied by Baxter, Belavin and Sklyanin.

\noindent
$\bullet$
For the  plane nodal cubic curve $E = V(zy^2 - x^3 - zx^2) \subset \PP^2$ we get a trigonometric
solution
\begin{align*}
  r^{(2,1)}_{\mathrm{trg}}(v, y) &= \frac{\sin(y+v)}{\sin(y)\sin(v)}(e_{11}\otimes e_{11} +
e_{22} \otimes e_{22}) + \frac{1}{\sin(v)}(e_{11} \otimes e_{22} + e_{22} \otimes
e_{11}) +
\\
&+ \frac{1}{\sin(y)}(e_{12}\otimes e_{21} + e_{21}\otimes e_{12}) + \sin(y+v)
e_{21} \otimes e_{21}.
\end{align*}

\noindent
This solution is a quantization of   the  trigonometric solution of
Cherednik:
$$
\bar{r}^{(2,1)}_{\textrm{trg}}(y) = \frac{1}{2} \cot(z)h\otimes h +
\frac{1}{\sin(y)}(e_{12}\otimes e_{21} + e_{21} \otimes e_{12}) +
\sin(y) e_{21} \otimes e_{21}.
$$

\noindent
$\bullet$
For the cuspidal nodal cubic curve $E = V(zy^2 - x^3) \subset \PP^2$ we get a rational
solution
\begin{align*}
  r^{(2,1)}_{\mathrm{rat}}(v, y_1, y_2) &=
\frac{1}{v} \mathbbm{1} \otimes \mathbbm{1} +
\frac{2}{y_2 - y_1}(e_{11} \otimes e_{11} + e_{22}\otimes e_{22} +
e_{12} \otimes e_{21} + e_{21}\otimes e_{12}) +
\\
& + (v - y_1)  e_{21} \otimes h +
(v + y_2)  h  \otimes e_{21} -
v(v - y_1)(v + y_2)  e_{21}\otimes e_{21}.
\end{align*}

\noindent
This solution of the quantum Yang--Baxter equation is a quantization of the rational solution
of Stolin
$$
\bar{r}(y_1, y_2) = \dfrac{1}{y_2 - y_1}\left(
\dfrac{1}{2} h \otimes h + e_{12}
\otimes e_{21} + e_{21} \otimes e_{12}\right)
+ \dfrac{y_2}{2} h \otimes e_{21} -
\dfrac{y_1}{2} e_{21} \otimes h.
$$
Note that this solution is gauge equivalent to the solution
$$
s(y)  =
\frac{1}{y}\Bigl(\frac{1}{2} h \otimes h + e_{12} \otimes e_{21} +
e_{21} \otimes e_{12}\Bigr) +
y(e_{21} \otimes h + h \otimes e_{21}) - y^3 e_{21} \otimes e_{21}
$$
depending on the difference of spectral parameters.

\noindent
$\bullet$
Next, the following solution $r^{(2,1)}_{\mathrm{rat-deg}}(v; y)$ of the associative Yang--Baxter equation is a degeneration
of the rational solution $ r^{(2,1)}_{\mathrm{rat}}(v, y_1, y_2)$:
$$ r^{(2,1)}_{\mathrm{rat-deg}}(v, y) = \dfrac{1}{2 v}
\mathbbm{1}\otimes \mathbbm{1} +
\dfrac{1}{y}\bigl(e_{11}
\otimes e_{11} + e_{22}\otimes e_{22} +
e_{12}\otimes e_{21} + e_{21}\otimes e_{12}\bigr).
 $$
 The corresponding solution of the classical Yang--Baxter equation is the rational solution of Yang:
 $$
r^{(2,1)}_{\mathrm{rat-deg}}(y)  =
\frac{1}{y}\Bigl(\frac{1}{2} h \otimes h + e_{12} \otimes e_{21} +
e_{21} \otimes e_{12}\Bigr).
$$

\noindent
$\bullet$
In the case of a nodal cubic curve $E$, the universal family of \emph{indecomposable
semi-stable} vector
bundles of rank $2$ and degree $0$ having \emph{locally free} Jordan-H\"older factors  gives
the following solution of the associative Yang--Baxter equation:
\begin{align*}
  r^{(2,0)}_{\mathrm{trg}}(v;y) &=
\dfrac{\sin(y +v)}{2\sin(y)\sin(v)}
\bigl(
e_{11} \otimes e_{11} + e_{22} \otimes e_{22} + e_{21} \otimes e_{12} +
e_{12} \otimes e_{21}
\bigr)
\\
&+
\dfrac{1}{2 \sin^2(v)} \bigl(e_{12} \otimes h - h \otimes e_{12}\bigr) -
\dfrac{\cos(v)}{\sin^3(v)} e_{12} \otimes e_{12}.
\end{align*}
This solution has higher order poles with respect to the spectral variable $v$ and does not project to a solution of the classical Yang--Baxter equation. However, it still yields a family of commuting
Dunkl operators.

\medskip
The second analytic application of methods developed in this article, is the following.
Consider the  Weierstra\ss{} family of plane cubic curves
$zy^2 = 4 x^3 - g_2 xz^2 - g_3 z^3$, where $g_2, g_3 \in \CC$.
Recall the following classical result.

\begin{proposition}[see Section II.4 in \cite{HurwitzCourant}]
Let $\tau \in \CC \setminus \mathbb{R}$ and $\Lambda_\tau = \mathbb{Z} + \tau
\mathbb{Z} \subset \CC^2$ be the corresponding lattice.
 Then the complex torus $\CC/\Lambda_\tau$ is isomorphic
to the projective cubic curve
$
zy^2 = 4 x^3 - g_2 xz^2 - g_3 z^3,
$
where
\begin{equation}\label{E:periods}
g_2 = 60 \sum\limits_{(m',m'') \in \mathbb{Z}^2\setminus \{(0,0)\}}
\frac{1}{(m' +  m'' \tau)^4}, \quad
g_3 = 140  \sum\limits_{(m',m'') \in \mathbb{Z}^2\setminus \{(0,0)\}}
\frac{1}{(m' +  m'' \tau)^6}.
\end{equation}
Conversely, for any pair $(g_2, g_3) \in \CC^2$ such that
$\Delta(g_2, g_3) = g_2^3 - 27 g_3^2  \ne 0$
there exists a unique $\tau$ from the domain $D$ given below,
such that $(g_2, g_3) = \bigl(g_2(\tau), g_3(\tau)\bigr)$.
$$
D= \left\{\tau \in \CC\left|\; |\mathrm{Re}(\tau)| \le \frac{1}{2},\quad
|\tau| \ge 1 \,\,  \text{ if } \,\,
\mathrm{Re}(\tau) \le 0, \quad |\tau| > 1 \,\, \text{ if } \,\,
  \mathrm{Re}(\tau) > 0 \right.\right\}
$$
\end{proposition}

Let  $(n, d) \in \mathbb{N} \times \mathbb{Z}$ be a pair of coprime integers
and $M = M_{E/T}^{(n,d)} \cong \breve{E}$ be the moduli space of relatively
stable vector bundles on $E$ of rank $n$ and degree $d$ with universal family
$\kP(n,d)$. Let $t = (g_2, g_3)$  and $o = \bigl((0:1:0), t\bigr) \in E$,
$m \in M$ be the point corresponding to $o$ and $\xi$ be some trivialization
of $\kP(n,d)$ in a neighbourhood of $(o, m) \in E \times_T M$ and
$\omega \in H^0(\omega_{E/T})$ be a nowhere vanishing regular  one-form.
Then we get the  germ of a meromorphic function
$$
r^\xi := \bigl(r^{(n,d)}_{E/T}(\omega)\bigr)^\xi:
\bigl(M\times_T  M\times_T\breve{E}\times_T\breve{E}, \hat{o}\bigr)
\lar \Mat_{n \times n}(\CC) \times \Mat_{n \times n}(\CC)
$$
which satisfies  the associative Yang--Baxter equation
$$
r^\xi(t; v_1, v_2; y_1, y_2)^{12} r^\xi(t; v_1, v_3; y_2, y_3)^{23} =
r^\xi(t; v_1, v_3; y_1, y_3)^{13} r^\xi(t; v_3, v_2; y_1, y_2)^{12} +
$$
$$
+ r^\xi(t; v_2, v_3; y_2, y_3)^{23} r^\xi(t; v_1, v_2; y_1, y_3)^{13}
$$
and its ``dual''
$$
r^\xi(t; v_2, v_3; y_1, y_2)^{23} r^\xi(t; v_1, v_3; y_1, y_2)^{12} =
r^\xi(t; v_1, v_2; y_1, y_2)^{12} r^\xi(t; v_2, v_3; y_1, y_3)^{13} +
$$
$$
+ r^\xi(t; v_1, v_3; y_1, y_3)^{13} r^\xi(t; v_2, v_1; y_2, y_3)^{23}.
$$
Moreover, it fulfills  the  unitarity condition
$$
r^\xi(t; v_1, v_2; y_1, y_2) = - \tau\bigl(r^\xi(t; v_2, v_1; y_2, y_1) \bigr),
$$
where $\tau(a \otimes b) = b \otimes a$.
The function $r^\xi(t; v_1, v_2; y_1, y_2)$ depends analytically on the parameter $t \in T$ and
its poles lie on the hypersurfaces $v_1 = v_2$ and $y_1 = y_2$.

Next, different choices of trivializations of the universal family $\kP$
lead to equivalent solutions: if $\zeta$ is another trivialization of $\kP$ and
$\phi = \zeta \, \circ \,  \xi^{-1}: (M \times_T E, o) \lar \GL_n(\CC)$ is the corresponding holomorphic function,
then we have:
$$
r^\zeta =
\bigl(\phi(t, v_1, y_1) \otimes \phi(t, v_2, y_2)\bigr)
 r^\xi\bigl(\phi(t, v_2, y_1)^{-1}
\otimes \phi(t, v_1, y_2)^{-1}\bigr).
$$
We have shown that $r^{(2,1)}_{\mathrm{ell}}(v; y)$ is equivalent
to the solution $r^{(2,1)}_{E/T}(t; v; y)$ for all $t = (g_2, g_3)$ such that
 $\Delta(t) \ne  0$. This equivalence relation is
  generated by the gauge transformations and  coordinate changes.

The trigonometric solution $r^{(2,1)}_{\mathrm{trg}}(v; y)$ is  equivalent
to the solution $r^{(2,1)}_{E/T}(t; v; y)$ for $ t = (g_2, g_3)  \ne (0, 0)$ but such that $\Delta(t) = 0$.
Finally, the rational solution  $r^{(2,1)}_{\mathrm{rat}}(v; y_1, y_2)$ is  equivalent
to the solution $r^{(2,1)}_{E/T}\bigl((0,0); v; y\bigr)$.

In other words, we get the following result, which seems to be difficult to show by direct computations:
the rational solution $r^{(2,1)}_{\mathrm{rat}}(v; y_1, y_2)$ is
\emph{equivalent}  to a solution, which is a \emph{degeneration} of the
trigonometric solution $r^{(2,1)}_{\mathrm{trg}}(v; y)$ and of the
elliptic solution $r^{(2,1)}_{\mathrm{ell}}(v; y)$.

\end{document}